\documentclass[12pt,a4paper]{article}

\usepackage{amsfonts,amssymb,amsthm,amsmath,latexsym}
\usepackage{eqsection,indent}
\usepackage{subeqnarray, eepicemu, url,cite,bm}
\usepackage{multirow}  
\usepackage{tensor}  
\usepackage{graphicx,graphbox}
\usepackage{float}  
\usepackage{mycaption}
\usepackage{tikz}

\date{July 9, 2018 \\[1mm] revised May 17, 2020 \\[2.5mm]
      To appear in {\em Memoirs of the American Mathematical Society}\/}   

\begin{document}

\title{\vspace*{-2cm}
       Lattice paths and branched continued fractions:  \\[5mm]
       \hspace*{-4mm}An infinite sequence of generalizations \\
       \hspace*{-8mm}of the Stieltjes--Rogers and Thron--Rogers polynomials,
       \hspace*{-3mm}with coefficientwise Hankel-total positivity
      }

\author{ \\
      \hspace*{-2cm}
      {\large Mathias P\'etr\'eolle${}^1$, Alan D.~Sokal${}^{1,2}$,
                Bao-Xuan Zhu${}^{3,1}$}
   \\[5mm]
     \hspace*{-2.2cm}
      \normalsize
           ${}^1$Department of Mathematics, University College London,
                    London WC1E 6BT, UK   \\[1mm]
     \hspace*{-3.8cm}
      \normalsize
           ${}^2$Department of Physics, New York University,
                    New York, NY 10003, USA     \\
     \hspace*{-1.5cm}
      \normalsize
           ${}^3$School of Mathematical Sciences, Jiangsu Normal University,
                    Xuzhou 221116, CHINA
       \\[5mm]
     \hspace*{-1cm}
     {\tt mathias.petreolle@gmail.com}, {\tt sokal@nyu.edu},
     {\tt bxzhu@jsnu.edu.cn} \\[1cm]
}

\maketitle
\thispagestyle{empty}   

\begin{abstract}
We define an infinite sequence of generalizations,
parametrized by an integer $m \ge 1$,
of the Stieltjes--Rogers and Thron--Rogers polynomials;
they arise as the power-series expansions
of some branched continued fractions,
and as the generating polynomials for $m$-Dyck and $m$-Schr\"oder paths
with height-dependent weights.
We prove that all of these sequences of polynomials are
coefficientwise Hankel-totally positive,
jointly in all the (infinitely many) indeterminates.
We then apply this theory to prove the
coefficientwise Hankel-total positivity
for combinatorially interesting sequences of polynomials.
Enumeration of unlabeled ordered trees and forests
gives rise to multivariate Fuss--Narayana polynomials
and Fuss--Narayana symmetric functions.
Enumeration of increasing (labeled) ordered trees and forests
gives rise to multivariate Eulerian polynomials
and Eulerian symmetric functions,
which include the univariate $m$th-order Eulerian polynomials
as specializations.
We also find branched continued fractions
for ratios of contiguous hypergeometric series ${}_r \! F_s$
for arbitrary $r$ and $s$,
which generalize Gauss' continued fraction
for ratios of contiguous ${}_2 \! F_1$;
and for $s=0$ we prove the coefficientwise Hankel-total positivity.
Finally, we extend the branched continued fractions
to ratios of contiguous basic hypergeometric series ${}_r \! \phi_s$.
\end{abstract}

\bigskip
\noindent
{\bf Key Words:}  Dyck path, $m$-Dyck path,
Schr\"oder path, $m$-Schr\"oder path,
Motzkin path, \L{}ukasiewicz path,
Catalan numbers, Fuss--Catalan numbers, Schr\"oder numbers,
continued fraction, S-fraction, T-fraction, J-fraction,
branched continued fraction,
Stieltjes--Rogers polynomials, Thron--Rogers polynomials,
Jacobi--Rogers polynomials,
production matrix,
totally positive matrix, total positivity, Hankel matrix,
Lindstr\"om--Gessel--Viennot lemma, Stieltjes moment problem,
Fuss--Narayana polynomial, Fuss--Narayana symmetric function,
Eulerian polynomial, Eulerian symmetric function, Stirling permutation,
hypergeometric series, basic hypergeometric series, contiguous relation.

\bigskip
\bigskip
\noindent
{\bf Mathematics Subject Classification (MSC 2010) codes:}
05A15 (Primary);
05A19, 05A20, 05C30, 05E05, 15B05, 15B48, 30B70, 30E05,
33C05, 33C20, 33D05, 44A60
(Secondary).

\vspace*{1cm}

\newtheorem{theorem}{Theorem}[section]
\newtheorem{proposition}[theorem]{Proposition}
\newtheorem{lemma}[theorem]{Lemma}
\newtheorem{corollary}[theorem]{Corollary}
\newtheorem{definition}[theorem]{Definition}
\newtheorem{conjecture}[theorem]{Conjecture}
\newtheorem{question}[theorem]{Question}
\newtheorem{problem}[theorem]{Problem}
\newtheorem{openproblem}[theorem]{Open Problem}
\newtheorem{example}[theorem]{Example}

\renewcommand{\theenumi}{\alph{enumi}}
\renewcommand{\labelenumi}{(\theenumi)}
\def\eop{\hbox{\kern1pt\vrule height6pt width4pt
depth1pt\kern1pt}\medskip}
\def\prf{\par\noindent{\bf Proof.\enspace}\rm}
\def\rmk{\par\medskip\noindent{\bf Remark\enspace}\rm}

\newcommand{\textbfit}[1]{\textbf{\textit{#1}}}

\newcommand{\bigdash}{%
\smallskip\begin{center} \rule{5cm}{0.1mm} \end{center}\smallskip}

\newcommand{\safepar}{ {\protect\hfill\protect\break\hspace*{5mm}} }

\newcommand{\be}{\begin{equation}}
\newcommand{\ee}{\end{equation}}
\newcommand{\<}{\langle}
\renewcommand{\>}{\rangle}
\newcommand{\widebar}{\overline}
\def\reff#1{(\protect\ref{#1})}
\def\spose#1{\hbox to 0pt{#1\hss}}
\def\ltapprox{\mathrel{\spose{\lower 3pt\hbox{$\mathchar"218$}}
    \raise 2.0pt\hbox{$\mathchar"13C$}}}
\def\gtapprox{\mathrel{\spose{\lower 3pt\hbox{$\mathchar"218$}}
    \raise 2.0pt\hbox{$\mathchar"13E$}}}
\def\textprime{${}^\prime$}
\def\proof{\par\medskip\noindent{\sc Proof.\ }}
\def\firstproof{\par\medskip\noindent{\sc First Proof.\ }}
\def\secondproof{\par\medskip\noindent{\sc Second Proof.\ }}
\def\alternateproof{\par\medskip\noindent{\sc Alternate Proof.\ }}
\def\algebraicproof{\par\medskip\noindent{\sc Algebraic Proof.\ }}
\def\combinatorialproof{\par\medskip\noindent{\sc Combinatorial Proof.\ }}
\def\proofof#1{\bigskip\noindent{\sc Proof of #1.\ }}
\def\firstproofof#1{\bigskip\noindent{\sc First Proof of #1.\ }}
\def\secondproofof#1{\bigskip\noindent{\sc Second Proof of #1.\ }}
\def\thirdproofof#1{\bigskip\noindent{\sc Third Proof of #1.\ }}
\def\algebraicproofof#1{\bigskip\noindent{\sc Algebraic Proof of #1.\ }}
\def\combinatorialproofof#1{\bigskip\noindent{\sc Combinatorial Proof of #1.\ }}
\def\sketchofproof{\par\medskip\noindent{\sc Sketch of proof.\ }}
\renewcommand{\qed}{ $\square$ \bigskip}
\newcommand{\myendremark}{ $\blacksquare$ \bigskip}
\def\half{ {1 \over 2} }
\def\third{ {1 \over 3} }
\def\twothird{ {2 \over 3} }
\def\smfrac#1#2{{\textstyle{#1\over #2}}}
\def\smhalf{ {\smfrac{1}{2}} }
\newcommand{\real}{\mathop{\rm Re}\nolimits}
\renewcommand{\Re}{\mathop{\rm Re}\nolimits}
\newcommand{\imag}{\mathop{\rm Im}\nolimits}
\renewcommand{\Im}{\mathop{\rm Im}\nolimits}
\newcommand{\sgn}{\mathop{\rm sgn}\nolimits}
\newcommand{\tr}{\mathop{\rm tr}\nolimits}
\newcommand{\supp}{\mathop{\rm supp}\nolimits}
\newcommand{\disc}{\mathop{\rm disc}\nolimits}
\newcommand{\diag}{\mathop{\rm diag}\nolimits}
\newcommand{\tridiag}{\mathop{\rm tridiag}\nolimits}
\newcommand{\AZ}{\mathop{\rm AZ}\nolimits}
\newcommand{\NC}{\mathop{\rm NC}\nolimits}
\newcommand{\PF}{{\rm PF}}
\newcommand{\rk}{\mathop{\rm rk}\nolimits}
\newcommand{\perm}{\mathop{\rm perm}\nolimits}
\def\hboxscript#1{ {\hbox{\scriptsize\em #1}} }
\renewcommand{\emptyset}{\varnothing}
\newcommand{\eqdef}{\stackrel{\rm def}{=}}

\newcommand{\restrict}{\upharpoonright}

\newcommand{\compinv}{{\langle -1 \rangle}}   

\newcommand{\scra}{{\mathcal{A}}}
\newcommand{\scrb}{{\mathcal{B}}}
\newcommand{\scrc}{{\mathcal{C}}}
\newcommand{\scrd}{{\mathcal{D}}}
\newcommand{\scre}{{\mathcal{E}}}
\newcommand{\scrf}{{\mathcal{F}}}
\newcommand{\scrg}{{\mathcal{G}}}
\newcommand{\scrh}{{\mathcal{H}}}
\newcommand{\scri}{{\mathcal{I}}}
\newcommand{\scrj}{{\mathcal{J}}}
\newcommand{\scrk}{{\mathcal{K}}}
\newcommand{\scrl}{{\mathcal{L}}}
\newcommand{\scrm}{{\mathcal{M}}}
\newcommand{\scrn}{{\mathcal{N}}}
\newcommand{\scro}{{\mathcal{O}}}
\newcommand\scroo{
  \mathchoice
    {{\scriptstyle\mathcal{O}}}
    {{\scriptstyle\mathcal{O}}}
    {{\scriptscriptstyle\mathcal{O}}}
    {\scalebox{0.6}{$\scriptscriptstyle\mathcal{O}$}}
  }
\newcommand{\scrp}{{\mathcal{P}}}
\newcommand{\scrq}{{\mathcal{Q}}}
\newcommand{\scrr}{{\mathcal{R}}}
\newcommand{\scrs}{{\mathcal{S}}}
\newcommand{\scrt}{{\mathcal{T}}}
\newcommand{\scrv}{{\mathcal{V}}}
\newcommand{\scrw}{{\mathcal{W}}}
\newcommand{\scrz}{{\mathcal{Z}}}

\newcommand{\bfa}{{\mathbf{a}}}
\newcommand{\bfb}{{\mathbf{b}}}
\newcommand{\bfc}{{\mathbf{c}}}
\newcommand{\bfd}{{\mathbf{d}}}
\newcommand{\bfe}{{\mathbf{e}}}
\newcommand{\bfh}{{\mathbf{h}}}
\newcommand{\bfj}{{\mathbf{j}}}
\newcommand{\bfi}{{\mathbf{i}}}
\newcommand{\bfk}{{\mathbf{k}}}
\newcommand{\bfl}{{\mathbf{l}}}
\newcommand{\bfL}{{\mathbf{L}}}
\newcommand{\bfm}{{\mathbf{m}}}
\newcommand{\bfn}{{\mathbf{n}}}
\newcommand{\bfp}{{\mathbf{p}}}
\newcommand{\bfr}{{\mathbf{r}}}
\newcommand{\bfu}{{\mathbf{u}}}
\newcommand{\bfv}{{\mathbf{v}}}
\newcommand{\bfw}{{\mathbf{w}}}
\newcommand{\bfx}{{\mathbf{x}}}
\newcommand{\bfy}{{\mathbf{y}}}
\newcommand{\bfz}{{\mathbf{z}}}
\renewcommand{\k}{{\mathbf{k}}}
\newcommand{\n}{{\mathbf{n}}}
\newcommand{\vv}{{\mathbf{v}}}
\newcommand{\bv}{{\mathbf{v}}}
\newcommand{\w}{{\mathbf{w}}}
\newcommand{\x}{{\mathbf{x}}}
\newcommand{\y}{{\mathbf{y}}}
\newcommand{\cc}{{\mathbf{c}}}
\newcommand{\zero}{{\mathbf{0}}}
\newcommand{\one}{{\mathbf{1}}}
\newcommand{\bmm}{{\mathbf{m}}}

\newcommand{\ahat}{{\widehat{a}}}
\newcommand{\Zhat}{{\widehat{Z}}}

\newcommand{\C}{{\mathbb C}}
\newcommand{\D}{{\mathbb D}}
\newcommand{\Z}{{\mathbb Z}}
\newcommand{\N}{{\mathbb N}}
\newcommand{\Q}{{\mathbb Q}}
\newcommand{\PP}{{\mathbb P}}
\newcommand{\R}{{\mathbb R}}
\newcommand{\RR}{{\mathbb R}}
\newcommand{\E}{{\mathbb E}}

\newcommand{\Sym}{{\mathfrak{S}}}
\newcommand{\SymB}{{\mathfrak{B}}}
\newcommand{\Alt}{{\mathrm{Alt}}}

\newcommand{\germanA}{{\mathfrak{A}}}
\newcommand{\germanB}{{\mathfrak{B}}}
\newcommand{\germanQ}{{\mathfrak{Q}}}
\newcommand{\germanh}{{\mathfrak{h}}}

\newcommand{\myle}{\preceq}
\newcommand{\myge}{\succeq}
\newcommand{\mygt}{\succ}

\newcommand{\B}{{\sf B}}
\newcommand{\OB}{B^{\rm ord}}
\newcommand{\OS}{{\sf OS}}
\newcommand{\OO}{{\sf O}}
\newcommand{\SP}{{\sf SP}}
\newcommand{\OSP}{{\sf OSP}}
\newcommand{\Eu}{{\sf Eu}}
\newcommand{\ERR}{{\sf ERR}}
\newcommand{\sfB}{{\sf B}}
\newcommand{\sfD}{{\sf D}}
\newcommand{\sfE}{{\sf E}}
\newcommand{\sfG}{{\sf G}}
\newcommand{\sfJ}{{\sf J}}
\newcommand{\sfP}{{\sf P}}
\newcommand{\sfQ}{{\sf Q}}
\newcommand{\sfS}{{\sf S}}
\newcommand{\sfT}{{\sf T}}
\newcommand{\sfW}{{\sf W}}
\newcommand{\sfMV}{{\sf MV}}
\newcommand{\AMV}{{\sf AMV}}
\newcommand{\BM}{{\sf BM}}
\newcommand{\emIB}{B^{\rm irr}}
\newcommand{\emIP}{P^{\rm irr}}
\newcommand{\emOB}{B^{\rm ord}}
\newcommand{\emCB}{B^{\rm cyc}}
\newcommand{\emSC}{P^{\rm cyc}}

\newcommand{\lev}{{\rm lev}}
\newcommand{\stat}{{\rm stat}}
\newcommand{\cyc}{{\rm cyc}}
\newcommand{\mysteryone}{{\rm mys1}}
\newcommand{\mysterytwo}{{\rm mys2}}
\newcommand{\Asc}{{\rm Asc}}
\newcommand{\asc}{{\rm asc}}
\newcommand{\Des}{{\rm Des}}
\newcommand{\des}{{\rm des}}
\newcommand{\Exc}{{\rm Exc}}
\newcommand{\exc}{{\rm exc}}
\newcommand{\Wex}{{\rm Wex}}
\newcommand{\wex}{{\rm wex}}
\newcommand{\Fix}{{\rm Fix}}
\newcommand{\fix}{{\rm fix}}
\newcommand{\lrmax}{{\rm lrmax}}
\newcommand{\rlmax}{{\rm rlmax}}
\newcommand{\Rec}{{\rm Rec}}
\newcommand{\rec}{{\rm rec}}
\newcommand{\Arec}{{\rm Arec}}
\newcommand{\arec}{{\rm arec}}
\newcommand{\ERec}{{\rm ERec}}
\newcommand{\erec}{{\rm erec}}
\newcommand{\EArec}{{\rm EArec}}
\newcommand{\earec}{{\rm earec}}
\newcommand{\recarec}{{\rm recarec}}
\newcommand{\nonrec}{{\rm nonrec}}
\newcommand{\Cpeak}{{\rm Cpeak}}
\newcommand{\cpeak}{{\rm cpeak}}
\newcommand{\Cval}{{\rm Cval}}
\newcommand{\cval}{{\rm cval}}
\newcommand{\Cdasc}{{\rm Cdasc}}
\newcommand{\cdasc}{{\rm cdasc}}
\newcommand{\Cddes}{{\rm Cddes}}
\newcommand{\cddes}{{\rm cddes}}
\newcommand{\cdrise}{{\rm cdrise}}
\newcommand{\cdfall}{{\rm cdfall}}
\newcommand{\Peak}{{\rm Peak}}
\newcommand{\peak}{{\rm peak}}
\newcommand{\Val}{{\rm Val}}
\newcommand{\val}{{\rm val}}
\newcommand{\Dasc}{{\rm Dasc}}
\newcommand{\dasc}{{\rm dasc}}
\newcommand{\Ddes}{{\rm Ddes}}
\newcommand{\ddes}{{\rm ddes}}
\newcommand{\inv}{{\rm inv}}
\newcommand{\maj}{{\rm maj}}
\newcommand{\rs}{{\rm rs}}
\newcommand{\cross}{{\rm cr}}
\newcommand{\crosshat}{{\widehat{\rm cr}}}
\newcommand{\nest}{{\rm ne}}
\newcommand{\rodd}{{\rm rodd}}
\newcommand{\reven}{{\rm reven}}
\newcommand{\lodd}{{\rm lodd}}
\newcommand{\leven}{{\rm leven}}
\newcommand{\sg}{{\rm sg}}
\newcommand{\bl}{{\rm bl}}
\newcommand{\tran}{{\rm tr}}
\newcommand{\area}{{\rm area}}
\newcommand{\ret}{{\rm ret}}
\newcommand{\peaks}{{\rm peaks}}
\newcommand{\hl}{{\rm hl}}
\newcommand{\sll}{{\rm sl}}
\newcommand{\negg}{{\rm neg}}
\newcommand{\imp}{{\rm imp}}
\newcommand{\osg}{{\rm osg}}
\newcommand{\ons}{{\rm ons}}
\newcommand{\isg}{{\rm isg}}
\newcommand{\ins}{{\rm ins}}
\newcommand{\LL}{{\rm LL}}
\newcommand{\height}{{\rm ht}}
\newcommand{\as}{{\rm as}}

\newcommand{\ba}{{\bm{a}}}
\newcommand{\bahat}{{\widehat{\bm{a}}}}
\newcommand{\sfa}{{{\sf a}}}
\newcommand{\bb}{{\bm{b}}}
\newcommand{\bc}{{\bm{c}}}
\newcommand{\bchat}{{\widehat{\bm{c}}}}
\newcommand{\bd}{{\bm{d}}}
\newcommand{\bee}{{\bm{e}}}
\newcommand{\beh}{{\bm{eh}}}
\newcommand{\bff}{{\bm{f}}}
\newcommand{\bg}{{\bm{g}}}
\newcommand{\bh}{{\bm{h}}}
\newcommand{\bll}{{\bm{\ell}}}
\newcommand{\bp}{{\bm{p}}}
\newcommand{\br}{{\bm{r}}}
\newcommand{\bs}{{\bm{s}}}
\newcommand{\bu}{{\bm{u}}}
\newcommand{\bw}{{\bm{w}}}
\newcommand{\bx}{{\bm{x}}}
\newcommand{\by}{{\bm{y}}}
\newcommand{\bz}{{\bm{z}}}
\newcommand{\bA}{{\bm{A}}}
\newcommand{\bB}{{\bm{B}}}
\newcommand{\bC}{{\bm{C}}}
\newcommand{\bE}{{\bm{E}}}
\newcommand{\bF}{{\bm{F}}}
\newcommand{\bG}{{\bm{G}}}
\newcommand{\bH}{{\bm{H}}}
\newcommand{\bI}{{\bm{I}}}
\newcommand{\bJ}{{\bm{J}}}
\newcommand{\bM}{{\bm{M}}}
\newcommand{\bN}{{\bm{N}}}
\newcommand{\bP}{{\bm{P}}}
\newcommand{\bQ}{{\bm{Q}}}
\newcommand{\bR}{{\bm{R}}}
\newcommand{\bS}{{\bm{S}}}
\newcommand{\bT}{{\bm{T}}}
\newcommand{\bW}{{\bm{W}}}
\newcommand{\bX}{{\bm{X}}}
\newcommand{\bY}{{\bm{Y}}}
\newcommand{\bIB}{{\bm{B}^{\rm irr}}}
\newcommand{\bOB}{{\bm{B}^{\rm ord}}}
\newcommand{\bOS}{{\bm{OS}}}
\newcommand{\bERR}{{\bm{ERR}}}
\newcommand{\bSP}{{\bm{SP}}}
\newcommand{\bMV}{{\bm{MV}}}
\newcommand{\bBM}{{\bm{BM}}}
\newcommand{\balpha}{{\bm{\alpha}}}
\newcommand{\balphapre}{{\bm{\alpha}^{\rm pre}}}
\newcommand{\bbeta}{{\bm{\beta}}}
\newcommand{\bgamma}{{\bm{\gamma}}}
\newcommand{\bdelta}{{\bm{\delta}}}
\newcommand{\bkappa}{{\bm{\kappa}}}
\newcommand{\bmu}{{\bm{\mu}}}
\newcommand{\bomega}{{\bm{\omega}}}
\newcommand{\bsigma}{{\bm{\sigma}}}
\newcommand{\btau}{{\bm{\tau}}}
\newcommand{\bphi}{{\bm{\phi}}}
\newcommand{\bpsi}{{\bm{\psi}}}
\newcommand{\bzeta}{{\bm{\zeta}}}
\newcommand{\bone}{{\bm{1}}}
\newcommand{\bzero}{{\bm{0}}}

\newcommand{\Cbar}{{\overline{C}}}
\newcommand{\Dbar}{{\overline{D}}}
\newcommand{\dbar}{{\overline{d}}}
\def\Ctilde{{\widetilde{C}}}
\def\Ftilde{{\widetilde{F}}}
\def\Gtilde{{\widetilde{G}}}
\def\Htilde{{\widetilde{H}}}
\def\Ptilde{{\widetilde{P}}}
\def\Chat{{\widehat{C}}}
\def\ctilde{{\widetilde{c}}}
\def\zbar{{\overline{Z}}}
\def\pitilde{{\widetilde{\pi}}}

\newcommand{\sech}{{\rm sech}}

%
%
\newcommand{\sn}{{\rm sn}}
\newcommand{\cn}{{\rm cn}}
\newcommand{\dn}{{\rm dn}}
\newcommand{\sm}{{\rm sm}}
\newcommand{\cm}{{\rm cm}}

%
%
\newcommand{\zfz}{ {{}_0 \! F_0} }
\newcommand{\zfo}{ {{}_0  F_1} }
\newcommand{\ofz}{ {{}_1 \! F_0} }
\newcommand{\ofo}{ {{}_1 \! F_1} }
\newcommand{\oft}{ {{}_1 \! F_2} }

%
%
\newcommand{\FHyper}[2]{ {\tensor[_{#1 \!}]{F}{_{#2}}\!} }
\newcommand{\FHYPER}[5]{ {\FHyper{#1}{#2} \!\biggl(
   \!\!\begin{array}{c} #3 \\[1mm] #4 \end{array}\! \bigg|\, #5 \! \biggr)} }
\newcommand{\tfo}{ {\FHyper{2}{1}} }
\newcommand{\tfz}{ {\FHyper{2}{0}} }
\newcommand{\threefz}{ {\FHyper{3}{0}} }
\newcommand{\FHYPERbottomzero}[3]{ {\FHyper{#1}{0} \hspace*{-0mm}\biggl(
   \!\!\begin{array}{c} #2 \\[1mm] \hbox{---} \end{array}\! \bigg|\, #3 \! \biggr)} }
\newcommand{\FHYPERtopzero}[3]{ {\FHyper{0}{#1} \hspace*{-0mm}\biggl(
   \!\!\begin{array}{c} \hbox{---} \\[1mm] #2 \end{array}\! \bigg|\, #3 \! \biggr)} }

\newcommand{\phiHyper}[2]{ {\tensor[_{#1}]{\phi}{_{#2}}} }
\newcommand{\psiHyper}[2]{ {\tensor[_{#1}]{\psi}{_{#2}}} }
\newcommand{\PhiHyper}[2]{ {\tensor[_{#1}]{\Phi}{_{#2}}} }
\newcommand{\PsiHyper}[2]{ {\tensor[_{#1}]{\Psi}{_{#2}}} }
\newcommand{\phiHYPER}[6]{ {\phiHyper{#1}{#2} \!\left(
   \!\!\begin{array}{c} #3 \\ #4 \end{array}\! ;\, #5, \, #6 \! \right)\!} }
\newcommand{\psiHYPER}[6]{ {\psiHyper{#1}{#2} \!\left(
   \!\!\begin{array}{c} #3 \\ #4 \end{array}\! ;\, #5, \, #6 \! \right)} }
\newcommand{\PhiHYPER}[5]{ {\PhiHyper{#1}{#2} \!\left(
   \!\!\begin{array}{c} #3 \\ #4 \end{array}\! ;\, #5 \! \right)\!} }
\newcommand{\PsiHYPER}[5]{ {\PsiHyper{#1}{#2} \!\left(
   \!\!\begin{array}{c} #3 \\ #4 \end{array}\! ;\, #5 \! \right)\!} }
\newcommand{\zerophizero}{ {\phiHyper{0}{0}} }
\newcommand{\ophizero}{ {\phiHyper{1}{0}} }
\newcommand{\zphio}{ {\phiHyper{0}{1}} }
\newcommand{\ophio}{ {\phiHyper{1}{1}} }
\newcommand{\tphio}{ {\phiHyper{2}{1}} }
\newcommand{\tphiz}{ {\phiHyper{2}{0}} }
\newcommand{\tPhio}{ {\PhiHyper{2}{1}} }
\newcommand{\opsio}{ {\psiHyper{1}{1}} }

%
%
\newcommand{\stirlingsubset}[2]{\genfrac{\{}{\}}{0pt}{}{#1}{#2}}
\newcommand{\stirlingcycle}[2]{\genfrac{[}{]}{0pt}{}{#1}{#2}}
\newcommand{\assocstirlingsubset}[3]{{\genfrac{\{}{\}}{0pt}{}{#1}{#2}}_{\! \ge #3}}
\newcommand{\genstirlingsubset}[4]{{\genfrac{\{}{\}}{0pt}{}{#1}{#2}}_{\! #3,#4}}
\newcommand{\irredstirlingsubset}[2]{{\genfrac{\{}{\}}{0pt}{}{#1}{#2}}^{\!\rm irr}}
\newcommand{\euler}[2]{\genfrac{\langle}{\rangle}{0pt}{}{#1}{#2}}
\newcommand{\eulergen}[3]{{\genfrac{\langle}{\rangle}{0pt}{}{#1}{#2}}_{\! #3}}
\newcommand{\eulersecond}[2]{\left\langle\!\! \euler{#1}{#2} \!\!\right\rangle}
\newcommand{\eulersecondgen}[3]{{\left\langle\!\! \euler{#1}{#2} \!\!\right\rangle}_{\! #3}}
\newcommand{\binomvert}[2]{\genfrac{\vert}{\vert}{0pt}{}{#1}{#2}}
\newcommand{\binomsquare}[2]{\genfrac{[}{]}{0pt}{}{#1}{#2}}
\newcommand{\doublebinom}[2]{\left(\!\! \binom{#1}{#2} \!\!\right)}


\newenvironment{sarray}{
             \textfont0=\scriptfont0
             \scriptfont0=\scriptscriptfont0
             \textfont1=\scriptfont1
             \scriptfont1=\scriptscriptfont1
             \textfont2=\scriptfont2
             \scriptfont2=\scriptscriptfont2
             \textfont3=\scriptfont3
             \scriptfont3=\scriptscriptfont3
           \renewcommand{\arraystretch}{0.7}
           \begin{array}{l}}{\end{array}}

\newenvironment{scarray}{
             \textfont0=\scriptfont0
             \scriptfont0=\scriptscriptfont0
             \textfont1=\scriptfont1
             \scriptfont1=\scriptscriptfont1
             \textfont2=\scriptfont2
             \scriptfont2=\scriptscriptfont2
             \textfont3=\scriptfont3
             \scriptfont3=\scriptscriptfont3
           \renewcommand{\arraystretch}{0.7}
           \begin{array}{c}}{\end{array}}


\newcommand*\circled[1]{\tikz[baseline=(char.base)]{
  \node[shape=circle,draw,inner sep=1pt] (char) {#1};}}
\newcommand{\ostar}{{\circledast}}
\newcommand{\ostarN}{{\,\circledast_{\vphantom{\dot{N}}N}\,}}
\newcommand{\ostarPsi}{{\,\circledast_{\vphantom{\dot{\Psi}}\Psi}\,}}
\newcommand{\starN}{{\,\ast_{\vphantom{\dot{N}}N}\,}}
\newcommand{\starpsi}{{\,\ast_{\vphantom{\dot{\bpsi}}\!\bpsi}\,}}
\newcommand{\starone}{{\,\ast_{\vphantom{\dot{1}}1}\,}}
\newcommand{\startwo}{{\,\ast_{\vphantom{\dot{2}}2}\,}}
\newcommand{\starinfty}{{\,\ast_{\vphantom{\dot{\infty}}\infty}\,}}
\newcommand{\starT}{{\,\ast_{\vphantom{\dot{T}}T}\,}}

\newcommand*{\Scale}[2][4]{\scalebox{#1}{$#2$}}

\newcommand*{\Scaletext}[2][4]{\scalebox{#1}{#2}} 

\tableofcontents

\clearpage

\section{Introduction}

In a seminal 1980 paper, Flajolet \cite{Flajolet_80}
showed that the coefficients in the Taylor expansion
of the generic Stieltjes-type (resp.\ Jacobi-type) continued fraction
--- which he called the {\em Stieltjes--Rogers}\/
 (resp.\ {\em Jacobi--Rogers}\/) {\em polynomials}\/ ---
can be interpreted as the generating polynomials
for Dyck (resp.\ Motzkin) paths with specified height-dependent weights.
Very recently it was independently discovered by several authors
\cite{Fusy_15,Oste_15,Josuat-Verges_18,Sokal_totalpos}
that Thron-type continued fractions also have an interpretation of this kind:
namely, their Taylor coefficients
--- which we shall call, by analogy, the {\em Thron--Rogers polynomials}\/ ---
can be interpreted as the generating polynomials 
for Schr\"oder paths with specified height-dependent weights.
(All this will be explained in more detail in Section~\ref{sec2} below.)

The purpose of the present paper is to present an infinite sequence
of generalizations of the Stieltjes--Rogers and Thron--Rogers polynomials:
these generalizations are parametrized by an integer $m \ge 1$
and reduce to the classical Stieltjes--Rogers and Thron--Rogers polynomials
when $m=1$;
they are the generating polynomials of $m$-Dyck and $m$-Schr\"oder paths,
respectively, with height-dependent weights.
One fundamental feature of these generalizations
is that they all possess the property of
coefficientwise Hankel-total positivity \cite{Sokal_flajolet,Sokal_totalpos},
jointly in all the (infinitely many) indeterminates.
We will give two proofs of this fact:
a combinatorial proof based on the Lindstr\"om--Gessel--Viennot lemma,
and an algebraic proof based on the theory of production matrices.\footnote{
   More precisely, the combinatorial proof will apply to both the
   Stieltjes and Thron cases,
   while the algebraic proof will apply to the Stieltjes case
   and to a {\em subclass}\/ of the Thron case.
 \label{footnote.1}
}
These facts were known when $m = 1$ \cite{Sokal_flajolet,Sokal_totalpos}
but are new when $m > 1$.
By specializing the indeterminates we can give many examples
of Hankel-totally positive sequences
whose generating functions do not possess nice classical continued fractions.
(The concept of Hankel-total positivity \cite{Sokal_flajolet,Sokal_totalpos}
 will be explained in more detail later in this Introduction.)

In a similar way we will generalize the Jacobi--Rogers polynomials
to $m$-Jacobi--Rogers polynomials.
Just as when $m=1$ \cite{Sokal_flajolet,Sokal_totalpos},
these polynomials are {\em not}\/ in general Hankel-totally positive
--- in contrast to the Stieltjes--Rogers and Thron--Rogers cases ---
but they do provide a useful framework.
In particular, the standard formulae
for the contraction of an S-fraction to a J-fraction
\cite[p.~21]{Wall_48} \cite[p.~V-31]{Viennot_83}
have simple generalizations in which the $m$-Stieltjes--Rogers polynomials
(and also some specializations of the $m$-Thron--Rogers polynomials)
can be rewritten as $m$-Jacobi--Rogers polynomials
with suitably ``contracted'' weights.

We will furthermore introduce triangular arrays of generalized
$m$-Stieltjes--Rogers, $m$-Thron--Rogers and $m$-Jacobi--Rogers polynomials,
whose first columns are given by the ordinary
$m$-Stieltjes--Rogers, $m$-Thron--Rogers and $m$-Jacobi--Rogers polynomials,
respectively.
In the Stieltjes and Thron cases,
these triangular arrays turn out also to be coefficientwise totally positive,
jointly in all the indeterminates.
We will again give two proofs of this fact:
one combinatorial and one algebraic.\footnote{
   With the same clarification as in footnote~\ref{footnote.1}.
}

Just as the classical
Stieltjes--Rogers, Thron--Rogers and Jacobi--Rogers polynomials
can be interpreted as the Taylor coefficients
of Stieltjes-type, Thron-type or Jacobi-type continued fractions, respectively,
so the
$m$-Stieltjes--Rogers, $m$-Thron--Rogers and $m$-Jacobi--Rogers polynomials
can be interpreted as the Taylor coefficients
of certain branched continued fractions.
But these branched continued fractions seem very difficult to work with,
and we will not make much use of them.
(It would, however, be an interesting project for someone interested in
 the analytic theory of continued fractions
 \cite{Perron,Wall_48,Jones_80,Lorentzen_92,Cuyt_08}
 to study the convergence properties of these branched continued fractions
 when the variable $t$ lies in the complex plane.)

Let us now explain in more detail the context of this work
\cite{Sokal_flajolet,Sokal_totalpos},
which is basically combinatorial and linear-algebraic
but is strongly motivated by a classical branch of analysis:
the Stieltjes moment problem
\cite{Stieltjes_1894,Shohat_43,Akhiezer_65,Simon_98,Schmudgen_17}.
Recall first that a (finite or infinite) matrix of real numbers is called
{\em totally positive}\/ (TP) if all its minors are nonnegative,
and {\em strictly totally positive}\/ (STP)
if all its minors are strictly positive.\footnote{
   {\bf Warning:}  Many authors
   (e.g.\ \cite{Gantmakher_37,Gantmacher_02,Fomin_00,Fallat_11})
   use the terms ``totally nonnegative'' and ``totally positive''
   for what we have termed ``totally positive'' and
   ``strictly totally positive'', respectively.
   So it is very important, when seeing any claim about
   ``totally positive'' matrices, to ascertain which sense of
   ``totally positive'' is being used.
   (This is especially important because many theorems in this subject
    require {\em strict}\/ total positivity for their validity.)
}
Background information on totally positive matrices can be found
in \cite{Karlin_68,Ando_87,Gantmacher_02,Pinkus_10,Fallat_11};
they have applications to combinatorics
\cite{Brenti_89,Brenti_95,Brenti_96,Fomin_00,Skandera_03},
stochastic processes \cite{Karlin_59,Karlin_68},
statistics \cite{Karlin_68},
the mechanics of oscillatory systems \cite{Gantmakher_37,Gantmacher_02},
the zeros of polynomials and entire functions
\cite{Karlin_68,Asner_70,Kemperman_82,Holtz_03,Pinkus_10,Dyachenko_14},
spline interpolation \cite{Schoenberg_53,Karlin_68,Gasca_96},
Lie theory \cite{Lusztig_94,Lusztig_98,Fomin_99,Lusztig_08}
and cluster algebras \cite{Fomin_10,Fomin_forthcoming},
the representation theory of the infinite symmetric group
\cite{Thoma_64,Borodin_17},
the theory of immanants \cite{Stembridge_91},
planar discrete potential theory \cite{Curtis_98,Fomin_01}
and the planar Ising model \cite{Lis_17},
and several other areas of pure and applied mathematics \cite{Gasca_96}.

Now let $\ba = (a_n)_{n \ge 0}$ be a sequence of real numbers,
and let $H_\infty(\ba) = (a_{i+j})_{i,j \ge 0}$ be its
associated infinite Hankel matrix.
More generally, for each $m \ge 0$,
let $H_\infty^{(m)}(\ba) = (a_{i+j+m})_{i,j \ge 0}$
be the $m$-shifted infinite Hankel matrix.
The fundamental fact about total positivity of Hankel matrices
is the following:

\begin{theorem}[Stieltjes--Gantmakher--Krein]
   \label{thm.stieltjes}
For a sequence $\ba = (a_n)_{n \ge 0}$ of real numbers,
the following are equivalent:
\begin{itemize}
   \item[(a)]  $H^{(0)}_\infty(\ba)$ is totally positive.
      [That is, all the minors of $H^{(0)}_\infty(\ba)$ are nonnegative.]
   \item[(b)]  Both $H^{(0)}_\infty(\ba)$ and $H^{(1)}_\infty(\ba)$
      are positive-semidefinite.
      [That is, all the {\em principal} minors of $H^{(0)}_\infty(\ba)$
       and $H^{(1)}_\infty(\ba)$ are nonnegative.]
   \item[(c)]  There exists a positive measure $\mu$ on $[0,\infty)$
      such that $a_n = \int x^n \, d\mu(x)$ for all $n \ge 0$.
      [That is, $\ba$ is a Stieltjes moment sequence.]
   \item[(d)]  There exist numbers $\alpha_0,\alpha_1,\ldots \ge 0$
      such that
\be
   \sum_{n=0}^{\infty} a_n t^n
   \;=\;
   \cfrac{\alpha_0}{1 - \cfrac{\alpha_1 t}{1 - \cfrac{\alpha_2 t}{1 - \cdots}}}
   \label{eq.thm.hankelreal.infty.Stype}
\ee
in the sense of formal power series.
    [That is, the ordinary generating function
     $f(t) = \sum\limits_{n=0}^{\infty} a_n t^n$
    can be represented
    as a Stieltjes-type continued fraction with nonnegative coefficients.]
\end{itemize}
\end{theorem}

\noindent
Here (a)$\implies$(b) is trivial;
(b)$\iff$(c) is the standard necessary and sufficient condition
for a Stieltjes moment sequence
\cite{Stieltjes_1894,Shohat_43,Akhiezer_65,Simon_98,Schmudgen_17};
(c)$\iff$(d) is due to Stieltjes \cite{Stieltjes_1894} in 1894;
and (c)$\implies$(a) is due to Gantmakher and Krein \cite{Gantmakher_37}
in 1937.

So every Stieltjes moment sequence has a continued-fraction expansion
\reff{eq.thm.hankelreal.infty.Stype} with nonnegative coefficients $\balpha$.
For some Stieltjes moment sequences these coefficients are ``nice'',
but for others they look horrendous.
For instance, $n! = \int_0^\infty x^n \, e^{-x} \, dx$
is a Stieltjes moment sequence,
and Euler \cite{Euler_1760} showed in 1746 that\footnote{
   The paper \cite{Euler_1760},
   which is E247 in Enestr\"om's \cite{Enestrom_13} catalogue,
   was probably written circa 1746;
   it~was presented to the St.~Petersburg Academy in 1753,
   and published in 1760.
 \label{footnote_Euler_1760}
}
\be
   \sum_{n=0}^\infty n! \: t^n
   \;=\;
   \cfrac{1}{1 - \cfrac{1t}{1 - \cfrac{1t}{1 - \cfrac{2t}{1- \cfrac{2t}{1- \cdots}}}}}
   \;\:,
 \label{eq.nfact.contfrac}
\ee
i.e.\ it has a continued fraction \reff{eq.thm.hankelreal.infty.Stype}
with coefficients $\alpha_0 = 1$ and $\alpha_{2k-1} = \alpha_{2k} = k$
for $k \ge 1$.
Now the entrywise product of two Stieltjes moment sequences
is easily seen to be a Stieltjes moment sequence;
in particular, any positive-integer entrywise power
of a Stieltjes moment sequence is a Stieltjes moment sequence.
So $( (n!)^2 )_{n \ge 0}$ is a Stieltjes moment sequence;
therefore its ordinary generating function
has a continued fraction \reff{eq.thm.hankelreal.infty.Stype}
with nonnegative coefficients.  Straightforward computation gives
\be
   \alpha_1, \alpha_2, \ldots
   \;=\;
   1,\, 3,\, \smfrac{20}{3},\, \smfrac{164}{15},\, \smfrac{3537}{205},\,
     \smfrac{127845}{5371},\, \smfrac{4065232}{124057},\,
     \smfrac{244181904}{5868559},\, \smfrac{38418582575}{721944303},\,
     \ldots
   \;.
\ee
Perhaps the reader can see a pattern in these coefficients,
but we, alas, cannot!\footnote{
   And neither, at least as far as we checked, can {\tt gfun} \cite{Salvy_94}.
}
Likewise, any arithmetic-progression subsequence of a
Stieltjes moment sequence is a Stieltjes moment sequence;
thus, for instance, $( (2n)! )_{n \ge 0}$ is a Stieltjes moment sequence.
Its continued-fraction expansion has
\be
   \alpha_1, \alpha_2, \ldots
   \;=\;
   2,\, 10,\, \smfrac{108}{5},\, \smfrac{596}{15},\, \smfrac{27250}{447},\,
     \smfrac{1448550}{16241},\, \smfrac{351090488}{2923925},\,
     \smfrac{25537748648}{161245075},\, \smfrac{35020650343950}{176040604499},\,
     \ldots
   \;,
\ee
which looks, if anything, a bit worse.

Here we will show that although $(n!)^2$ and $(2n)!$
have ugly classical continued fractions,
they have nice $m$-branched continued fractions with $m=2$:
the coefficients are, respectively,
\be
   \balpha
   \;=\;
   1,1,2,4,4,6,9,9,12,\ldots
 \label{eq.mCF.nfactsq}
\ee
(corresponding to products of successive pairs of the ``pre-alphas''
 $\balphapre = 1,1,1,2,2,2,$ {}$3,3,3,\ldots$)
and
\be
   \balpha
   \;=\;
   2,4,6,12,16,20,30,36,42,\ldots
 \label{eq.mCF.2nfact}
\ee
(corresponding to products of successive pairs of
 $\balphapre = 1,2,2,3,4,4,5,6,6,\ldots$).
More generally, for any integer $m \ge 1$,
the sequences $(n!)^m$ and $(mn)!$ have nice $m$-branched continued fractions;
and yet more generally, for any integers $a,b \ge 1$,
the sequences $(an)!^b$, $(2an-1)!!^b$, $(3an-1)!!!^b$, $(3an-2)!!!^b$, etc.\ 
have nice $ab$-branched continued fractions
(see Section~\ref{sec.hyper.rF0}).

But this is only the beginning of the story,
because we are here principally concerned,
not with sequences and matrices of real numbers,
but with sequences and matrices of {\em polynomials}\/
(with integer or real coefficients) in one or more indeterminates $\bfx$:
in applications they will typically be generating polynomials that enumerate
some combinatorial objects with respect to one or more statistics.
We equip the polynomial ring $\R[\bfx]$ with the coefficientwise
partial order:  that is, we say that $P$ is nonnegative
(and write $P \myge 0$)
in case $P$ is a polynomial with nonnegative coefficients.
We then say that a matrix with entries in $\R[\bfx]$ is
{\em coefficientwise totally positive}\/
if all its minors are polynomials with nonnegative coefficients;
and we say that a sequence $\ba = (a_n)_{n \ge 0}$ with entries in $\R[\bfx]$
is {\em coefficientwise Hankel-totally positive}\/
if its associated infinite Hankel matrix is coefficientwise totally positive.
Similar definitions apply to the formal-power-series ring $\R[[\bfx]]$.
Most generally, we can consider sequences and matrices
with entries in an arbitrary partially ordered commutative ring;
total positivity and Hankel-total positivity
are then defined in the obvious way.

Now, in a general partially ordered commutative ring
--- and even in the univariate polynomial ring $\R[x]$
with the coefficientwise order ---
the analogue of Theorem~\ref{thm.stieltjes} fails to hold.
Indeed, there are many sequences $\ba = (a_n)_{n \ge 0}$ in $\R[x]$
that are (provably or conjecturally) coefficientwise Hankel-totally positive,
but for which there does not exist any continued fraction
\reff{eq.thm.hankelreal.infty.Stype} with coefficients $\balpha$ in $\R[x]$
\cite{Sokal_flajolet,Sokal_totalpos}:
in general the coefficients $\balpha$ are rational functions of $x$,
not polynomials.
Nevertheless, the implication (d)$\implies$(a) does hold
in an arbitrary partially ordered commutative ring
\cite{Sokal_flajolet,Sokal_totalpos}:
the existence of a continued fraction \reff{eq.thm.hankelreal.infty.Stype}
with nonnegative coefficients is a {\em sufficient}\/
(though far from necessary) condition for Hankel-total positivity.

In this paper we will prove an analogous result for the
$m$-branched continued fractions of Stieltjes and Thron type, for all $m \ge 1$.
This will allow us to prove the coefficientwise Hankel-total positivity
of some sequences of polynomials that do not have a classical
continued fraction \reff{eq.thm.hankelreal.infty.Stype}
in the polynomial ring $\R[\bfx]$,
but which do have an $m$-branched continued fraction (for suitably chosen $m$)
in which the $\balpha$ are polynomials in $\bfx$ with nonnegative coefficients.
Among these will be some multivariate polynomials
that refine the sequences $(n!)^m$ and $(mn)!$
(see Section~\ref{sec.hyper.rF0}).
In particular, for any integers $r,s \ge 0$
we will find a branched continued fraction
for ratios of contiguous hypergeometric series $\FHyper{r}{s}$,
which generalizes Gauss' \cite{Gauss_1813}
continued fraction for ratios of contiguous $\tfo$;
and for $s=0$ we will prove the coefficientwise Hankel-total positivity
of the corresponding sequence of polynomials
(see Sections~\ref{sec.hyper.rF0} and \ref{sec.hyper.rFs}).
We will also generalize these results to
basic hypergeometric series $\phiHyper{r}{s}$
(see Section~\ref{sec.hyper.rphis}).

Let us conclude this introduction by mentioning some precursors of our work.
There has been some study of branched continued fractions
in the analysis literature\footnote{
   See e.g.\ \cite[pp.~274--280]{Lorentzen_92}
   and the work cited on \cite[p.~285]{Lorentzen_92}
   and \cite[p.~28]{Cuyt_08}.
},
but this work seems rather distant from our concerns here.
More relevant are a few papers in which combinatorialists
have studied classes of lattice paths leading to branched continued fractions.
Already a quarter-century ago, Viennot \cite[section~V.6]{Viennot_83}
briefly considered the branched continued fractions
({\em fractions multicontinu\'ees}\/)
generated by \L{}ukasiewicz paths;
these correspond to our $m$-Jacobi--Rogers polynomials with $m=\infty$.
This work was then carried forward in the Ph.D.~theses of
Roblet \cite{Roblet_94} and Varvak \cite{Varvak_04}.
Gouyou-Beauchamps \cite{Gouyou-Beauchamps_98}
used the $m$-Jacobi--Rogers polynomials with $m=2,3,4$
to enumerate certain classes of convex polyominoes,
while Arqu\`es and Fran\c{c}on \cite{Arques_84}
employed a different type of branched continued fraction
to enumerate well-labeled trees.
Finally --- and closest to our own work ---
Albenque and Bouttier \cite{Albenque_12}
introduced the branched continued fractions generated by $m$-Dyck paths,
which correspond to our $m$-Stieltjes--Rogers polynomials,
and proved many interesting results about them
(but did not consider total positivity).
We will comment further on these papers in the subsequent sections.

The plan of this paper is as follows:
In Section~\ref{sec.m-SR_and_m-TR} we define the
$m$-Stieltjes--Rogers and $m$-Thron--Rogers polynomials
and derive the fundamental recurrences that they satisfy.
In Section~\ref{sec.different_m} we analyze the
relation between the $m$-Stieltjes--Rogers polynomials
for different values of $m$,
and similarly for the $m$-Thron--Rogers polynomials.
In Section~\ref{sec.m-JR} we define the $m$-Jacobi--Rogers polynomials.
In Section~\ref{sec.generalized}
we introduce the triangular arrays of generalized
$m$-Stieltjes--Rogers, $m$-Thron--Rogers and $m$-Jacobi--Rogers polynomials.
In Section~\ref{sec.mJR.trees} we interpret the generalized
$m$-Jacobi--Rogers polynomials as generating polynomials
for ordered forests of ordered trees.
In Section~\ref{sec.contraction}
we show how the $m$-Stieltjes--Rogers polynomials
(and also some special cases of the $m$-Thron--Rogers polynomials)
can be rewritten as $m$-Jacobi--Rogers polynomials
with suitably ``contracted'' weights.
In Section~\ref{sec.production} we review the theory of production matrices,
and then exhibit the production matrices for the triangular arrays of
generalized
$m$-Stieltjes--Rogers, $m$-Thron--Rogers and $m$-Jacobi--Rogers polynomials.
In Section~\ref{sec.totalpos}
--- which is the theoretical heart of the paper ---
we prove the coefficientwise total positivity of the Hankel matrices
associated to the $m$-Stieltjes--Rogers and $m$-Thron--Rogers polynomials,
and of the lower-triangular matrices of generalized
$m$-Stieltjes--Rogers and $m$-Thron--Rogers polynomials.

The remainder of the paper is devoted to illustrating the theory
and applying it to cases of combinatorial interest.
In Section~\ref{sec.periodic} we compute explicit formulae
for the $m$-Stieltjes--Rogers polynomials when the weights $\balpha$
are periodic of period~$m+1$ or $m$;
these give rise to multivariate Fuss--Narayana polynomials
and Fuss--Narayana symmetric functions.
In Section~\ref{sec.eventually_periodic} we generalize this
to treat some cases in which the weights $\balpha$ are eventually periodic
(i.e.\ a finite sequence followed by a periodic sequence);
these include the multivariate Aval polynomials.
In Section~\ref{sec.quasi-affine} we take one step up in complexity,
and treat some cases in which the weights $\balpha$
are ``quasi-affine'' of period~$m+1$ or $m$,
i.e.\ $\alpha_i$ is affine in~$i$
within each residue class of $i \bmod m+1$ or $m$;
these give rise to multivariate Eulerian polynomials
and Eulerian symmetric functions,
which include the classical (univariate) $m$th-order Eulerian polynomials
as specializations.\footnote{
   The Fuss--Narayana symmetric functions were known \cite{Stanley_97}
   (under the name of parking-function symmetric functions),
   though we here give them a simple (and possibly new)
   combinatorial interpretation in terms of ordered trees.
   The Eulerian symmetric functions, by contrast,
   appear not to have been introduced previously.
}
In Section~\ref{sec.hyper.rF0} we take another step up,
and prove an $m$-branched continued fraction
for ratios of contiguous hypergeometric series $\FHyper{m+1}{0}$;
here the weights $\balpha$ are not quasi-affine,
but they are products of $m$ successive ``pre-alphas'' that are quasi-affine.
We show that this $m$-branched continued fraction
implies the coefficientwise Hankel-total positivity
of a sequence of multivariate polynomials
that refines the sequences $(n!)^m$, $(mn)!$, $(2n-1)!!^m$ and many others.
In Section~\ref{sec.hyper.rFs} we extend this
to find a branched continued fraction
for ratios of contiguous hypergeometric series $\FHyper{r}{s}$
for arbitrary $r,s$,
which generalizes Gauss' \cite{Gauss_1813}
continued fraction for ratios of contiguous $\tfo$;
and we discuss the implications for Hankel-total positivity.
In Section~\ref{sec.hyper.rphis} we extend the branched continued fractions
to the basic hypergeometric series $\phiHyper{r}{s}$.
In Section~\ref{sec.final} we make some final remarks,
including a conjecture on higher-order Genocchi numbers.

We hope that this paper will be read (or at least readable)
by mathematicians working in a variety of subfields:
combinatorics and linear algebra, to be sure,
but also analysis and special functions.
We have therefore tried hard to make the paper self-contained
and to make our arguments understandable to non-specialists.
We apologize in advance to experts for boring them from time to time
with overly detailed explanations of elementary facts.

We suggest \cite{Stanley_86,Stanley_99}
as general references on enumerative combinatorics;
\cite{Gessel_16} for Lagrange inversion;
\cite[Chapter~1]{Macdonald_95} and \cite[Chapter~7]{Stanley_99}
for symmetric functions;
\cite{Karlin_68,Gantmacher_02,Pinkus_10,Fallat_11} for total positivity;
\cite{Rainville_60,Slater_66,Andrews_99} for hypergeometric series;
and \cite{Gasper_04} for basic hypergeometric series.

\section[The $\bm{m}$-Stieltjes--Rogers and $\bm{m}$-Thron--Rogers polynomials]{The $\bm{m}$-Stieltjes--Rogers and $\bm{m}$-Thron--Rogers \break polynomials}
   \label{sec2}
   \label{sec.m-SR_and_m-TR}

We begin by reviewing briefly some well-known facts
concerning Dyck paths, Catalan numbers, Stieltjes--Rogers polynomials
and Stieltjes-type continued fractions;
then we review some slightly less well-known facts
concerning Schr\"oder paths, large Schr\"oder numbers,
Thron--Rogers polynomials and Thron-type continued fractions;
and finally, we present our generalizations,
which are parametrized by an integer $m \ge 1$
and which reduce to the classical Stieltjes--Rogers and Thron--Rogers
polynomials when $m=1$.

\subsection{Dyck paths}

A {\em Dyck path}\/ is a path in the upper half-plane $\Z \times \N$,
starting and ending on the horizontal axis,
using steps $(1,1)$ [``rise'' or ``up step'']
and $(1,-1)$ [``fall'' or ``down step''].
More generally, a {\em Dyck path at level $k$}\/
is a path in $\Z \times \N_{\ge k}$,
starting and ending at height $k$,
using steps $(1,1)$ and $(1,-1)$.
Clearly a Dyck path must be of even length;
we denote by $\scrd_{2n}$ the set of Dyck paths from $(0,0)$ to $(2n,0)$.
The ordinary generating function
$D(t) = \sum_{n=0}^\infty |\scrd_{2n}| \, t^n$
satisfies the functional equation
\be
   D(t)  \;=\;  1 \,+\, t D(t)^2
   \;,
 \label{eq.dyck.1}
\ee
as can be seen by splitting the Dyck path
(if it is of nonzero length)
at its next-to-last return to the horizontal axis
(i.e.\ its last visit to the horizontal axis excepting the final point);
then the path is of the form $\scrp_0 U \scrp_1 D$
where $\scrp_i$ is an arbitrary Dyck path at level $i$,
$U$ is an up step and $D$ is a down step.
The functional equation \reff{eq.dyck.1} can equivalently be written in the form
\be
   D(t)  \;=\;  {1 \over 1 \,-\, t D(t)}
   \;,
 \label{eq.dyck.2}
\ee
which also has a nice combinatorial interpretation:
a Dyck path can be uniquely decomposed as the concatenation
of zero or more {\em irreducible}\/ Dyck paths
(that is, Dyck paths of nonzero length that touch the horizontal axis
 only at the starting and ending points);
and an irreducible Dyck path is of the form $U \scrp_1 D$
where $\scrp_1$ is a Dyck path at level 1.
It follows easily from \reff{eq.dyck.1}/\reff{eq.dyck.2} that
\be
   D(t)  \;=\;  {1 \,-\, \sqrt{1-4t}  \over 2t}
 \label{eq.dyck.ogf}
\ee
and hence (by binomial expansion) that $|\scrd_{2n}|$
equals the {\em Catalan number}\/ \cite{Roman_15,Stanley_15}
\be
   C_n  \;=\;  {1 \over n+1} \binom{2n}{n}
 \label{def.catalan}
\ee
\!\!\cite[A000108]{OEIS}.
All this is extremely well-known and is presented
in numerous textbooks on enumerative combinatorics.

Now let $\balpha = (\alpha_i)_{i \ge 1}$ be an infinite set of indeterminates,
and let $S_n(\balpha)$ be the generating polynomial
for Dyck paths of length~$2n$ in which each rise gets weight~1
and each fall from height~$i$ gets weight $\alpha_i$.
Clearly $S_n(\balpha)$ is a homogeneous polynomial
of degree~$n$ with nonnegative integer coefficients;
following Flajolet \cite{Flajolet_80},
we call it the \textbfit{Stieltjes--Rogers polynomial} of order~$n$.
The first few are
\begin{subeqnarray}
   S_0  & = &   1          \\
   S_1  & = &   \alpha_1   \\
   S_2  & = &   \alpha_1^2 + \alpha_1 \alpha_2 \\
   S_3  & = &   \alpha_1^3 + 2\alpha_1^2 \alpha_2 + \alpha_1 \alpha_2^2 +
                   \alpha_1 \alpha_2 \alpha_3 \\
   S_4  & = &   \alpha_1^4 + 3 \alpha_1^3 \alpha_2 + 3 \alpha_1^2 \alpha_2^2 +
                   \alpha_1 \alpha_2^3 + 2 \alpha_1^2 \alpha_2 \alpha_3 +
                   2 \alpha_1 \alpha_2^2 \alpha_3 +
                   \alpha_1 \alpha_2 \alpha_3^2 +
                   \alpha_1 \alpha_2 \alpha_3 \alpha_4
   \nonumber \\
\end{subeqnarray}

Let $f_0(t) = \sum_{n=0}^\infty S_n(\balpha) \, t^n$
be the ordinary generating function for Dyck paths with these weights
(considered as a formal power series in $t$);
and more generally, let $f_k(t)$ be the ordinary generating function
for Dyck paths at level $k$ with these same weights.
(Obviously $f_k$ is just $f_0$ with each $\alpha_i$ replaced by $\alpha_{i+k}$;
 but we shall not explicitly use this fact.)
The same combinatorial arguments used earlier give the functional equations
\be
   f_k(t)  \;=\;  1 \:+\: \alpha_{k+1} t \, f_k(t) \, f_{k+1}(t)
 \label{eq.SRfk.1}
\ee
and
\be
   f_k(t)  \;=\;  {1 \over 1 \:-\: \alpha_{k+1} t \, f_{k+1}(t)}
   \;,
 \label{eq.SRfk.2}
\ee
which generalize \reff{eq.dyck.1}/\reff{eq.dyck.2}
and reduce to them when $\balpha = \bone$.
Iterating \reff{eq.SRfk.2}, we see immediately that $f_k$
is given by the continued fraction
\be
   f_k(t)
   \;=\;
   \cfrac{1}{1 - \cfrac{\alpha_{k+1} t}{1 - \cfrac{\alpha_{k+2} t}{1- \cfrac{\alpha_{k+3} t}{1- \cdots}}}}
 \label{eq.fk.Sfrac}
\ee
and in particular that $f_0$ is given by
\be
   f_0(t)
   \;=\;
   \cfrac{1}{1 - \cfrac{\alpha_{1} t}{1 - \cfrac{\alpha_{2} t}{1- \cfrac{\alpha_{3} t}{1- \cdots}}}}
   \;\,.
 \label{eq.f0.Sfrac}
\ee
The right-hand sides of \reff{eq.fk.Sfrac}/\reff{eq.f0.Sfrac}
are called \textbfit{Stieltjes-type continued fractions},
or \textbfit{S-fractions} for short.
This combinatorial interpretation of S-fractions
in terms of weighted Dyck paths
is due to Flajolet \cite{Flajolet_80}.


Let us also consider the combinatorial meaning of
a product of generating functions $f_0 f_1 \cdots f_\ell$ ($\ell \ge 0$).
We use the term {\em partial Dyck path}\/
to denote a path in the upper half-plane $\Z \times \N$,
using steps $(1,1)$ and $(1,-1)$,
that starts on the horizontal axis
but is allowed to end anywhere in the upper half-plane.\footnote{
    Here we follow the terminology of Stanley \cite[p.~146]{Stanley_15}.
    Some authors \cite{Banderier_02} \cite[p.~77]{Flajolet_09}
    call this a {\em meander}\/.
}
We then claim that the coefficient of $t^n$ in $f_0 f_1 \cdots f_\ell$
is the generating polynomial for partial Dyck paths from
$(0,0)$ to $(2n+\ell,\ell)$
[with the usual weights 1 for a rise and $\alpha_i$ for a fall],
which we will denote by $S_{n|\ell}(\balpha)$.
To see this, we split the partial Dyck path at its last return to level 0;
then it takes an up step;
we split the remaining path at its last return to level 1, and so forth.
So the path is of the form $\scrp_0 U \scrp_1 U \cdots \scrp_\ell$
where each $\scrp_i$ is an arbitrary Dyck path at level $i$.\footnote{
   See also \cite[pp.~II-7--II-8]{Viennot_83} \cite[pp.~295--296]{Goulden_83}
   for a similar argument applied to Motzkin paths.
}

And finally, let us mention an important method for proving
continued fractions of the form \reff{eq.fk.Sfrac}/\reff{eq.f0.Sfrac},
which was employed implicitly by Euler \cite[section~21]{Euler_1760}
for proving \reff{eq.nfact.contfrac}
and explicitly by Gauss \cite[sections~12--14]{Gauss_1813}
for proving his continued fraction for ratios of contiguous $\tfo$
(see also \cite{Sokal_alg_contfrac} for further discussion).
Namely, let $(g_k(t))_{k \ge -1}$ be a sequence of formal power series
(with coefficients in some commutative ring $R$)
with constant term 1,
and suppose that this sequence satisfies
a linear three-term recurrence of the form
\be
   g_k(t) - g_{k-1}(t)  \;=\; \alpha_{k+1} t \, g_{k+1}(t)
   \qquad\hbox{for } k \ge 0
 \label{eq.recurrence.gk1.0}
\ee
for some coefficients $\balpha = (\alpha_i)_{i \ge 1}$ in $R$.
If we define $f_k(t) = g_k(t)/g_{k-1}(t)$ for $k \ge 0$,
then \reff{eq.recurrence.gk1.0} can be rewritten as
\be
   f_k(t)  \;=\;  1 \:+\: \alpha_{k+1} t \, f_k(t) \, f_{k+1}(t)
   \;,
 \label{eq.SRfk.1.bis.0}
\ee
which is precisely the recurrence \reff{eq.SRfk.1}.
It follows that $f_k(t)$ is the ordinary generating function
for Dyck paths at level $k$ with weights $\balpha$
and is therefore given by the S-fraction \reff{eq.fk.Sfrac};
in particular,
$f_0(t) \eqdef g_0(t)/g_{-1}(t) = \sum_{n=0}^\infty S_n(\balpha) \, t^n$
where $S_n(\balpha)$ is the Stieltjes--Rogers polynomial
evaluated at the specified values $\balpha$;
and more generally $f_0(t) \cdots f_k(t) \eqdef g_k(t)/g_{-1}(t)
 = \sum_{n=0}^\infty S_{n|k}(\balpha) \, t^n$
as defined in the preceding paragraph.
This supplies a combinatorial interpretation of $g_k$,
at least when $g_{-1} = 1$ (which occurs in some but not all applications).
And conversely, if $(f_k)_{k \ge 0}$ satisfy \reff{eq.SRfk.1}/\reff{eq.SRfk.2},
then for any $g_{-1}(t)$ with constant term 1,
the series $g_k \eqdef g_{-1} f_0 f_1 \cdots f_k$
satisfy \reff{eq.recurrence.gk1.0}.

\subsection{Schr\"oder paths}

A {\em Schr\"oder path}\/ is a path in the upper half-plane $\Z \times \N$,
starting and ending on the horizontal axis,
using steps $(1,1)$ [``rise'' or ``up step''],
$(1,-1)$ [``fall'' or ``down step''],
and $(2,0)$ [``long level step''].
More generally, a {\em Schr\"oder path at level $k$}\/
is a path in $\Z \times \N_{\ge k}$,
starting and ending at height $k$,
using steps $(1,1)$, $(1,-1)$ and $(2,0)$.
We define the {\em length}\/ of a Schr\"oder path
to be the number of rises plus the number of falls
plus twice the number of long level steps.
Clearly a Schr\"oder path must be of even length;
we denote by $\scrs_{2n}$ the set of Schr\"oder paths from $(0,0)$ to $(2n,0)$.
The ordinary generating function
$S(t) = \sum_{n=0}^\infty |\scrs_{2n}| \, t^n$
satisfies the functional equation
\be
   S(t)  \;=\;  1 \,+\, t S(t) \,+\, t S(t)^2
   \;,
 \label{eq.schroder.1}
\ee
as can be seen by splitting the Schr\"oder path
(if it is of nonzero length)
at its next-to-last return to the horizontal axis:
then the path is either of the form $\scrp_0 L$ or $\scrp_0 U \scrp_1 D$,
where $\scrp_i$ is an arbitrary Schr\"oder path at level $i$,
$L$ is a long level step,
$U$ is an up step and $D$ is a down step.
The functional equation \reff{eq.schroder.1}
can equivalently be written in the form
\be
   S(t)  \;=\;  {1 \over 1 \,-\, t \,-\, t S(t)}
   \;,
 \label{eq.schroder.2}
\ee
which can be interpreted as saying that
a Schr\"oder path can be uniquely decomposed as the concatenation
of zero or more irreducible Schr\"oder paths;
and an irreducible Schr\"oder path is either a long level step
or else is of the form $U \scrp_1 D$.
It follows easily from \reff{eq.schroder.1}/\reff{eq.schroder.2} that
\be
   S(t)
   \;=\;
   {1 - t - \sqrt{1-6t+t^2} \over 2t}
   \;.
 \label{eq.schroder.ogf}
\ee
The coefficient $[t^n] \, S(t) = |\scrs_{2n}|$
is called the {\em large Schr\"oder number}\/ $r_n$
\cite[A006318]{OEIS};
it has the explicit expression (among others)
\be
   r_n  \;=\; \sum_{k=0}^n \binom{2n-k}{k} C_{n-k}
   \;,
\ee
since a Schr\"oder path of length $2n$ with $k$ long level steps
can be obtained from a Dyck path of length $2n-2k$
by inserting the long level steps arbitrarily among the Dyck steps.

Now let $\balpha = (\alpha_i)_{i \ge 1}$ and $\bdelta = (\delta_i)_{i \ge 1}$
be infinite sets of indeterminates,
and let $T_n(\balpha,\bdelta)$ be the generating polynomial
for Schr\"oder paths of length~$2n$ in which each rise gets weight~1,
each fall from height~$i$ gets weight $\alpha_i$,
and each long level step at height~$i$ gets weight $\delta_{i+1}$.
Clearly $T_n(\balpha,\bdelta)$ is a homogeneous polynomial
of degree~$n$ with nonnegative integer coefficients;
we call it the \textbfit{Thron--Rogers polynomial} of order~$n$.
The first few are
\begin{subeqnarray}
   T_0  & = &   1          \\
   T_1  & = &   \alpha_1 + \delta_1   \\
   T_2  & = &   (\alpha_1 + \delta_1)^2 + \alpha_1 (\alpha_2 + \delta_2)  \\
   T_3  & = &   (\alpha_1+\delta_1)^3 +
                   2\alpha_1 (\alpha_1 + \delta_1) (\alpha_2 + \delta_2) +
                   \alpha_1 (\alpha_2 + \delta_2)^2 +
                   \alpha_1 \alpha_2 (\alpha_3 + \delta_3)
          \nonumber \\
\end{subeqnarray}
When $\bdelta = \bzero$,
the Thron--Rogers polynomials reduce to the Stieltjes--Rogers polynomials.

Let $f_0(t) = \sum_{n=0}^\infty T_n(\balpha,\bdelta) \, t^n$
be the ordinary generating function for Schr\"oder paths with these weights;
and more generally, let $f_k(t)$ be the ordinary generating function
for Schr\"oder paths at level $k$ with these same weights.
(Obviously $f_k$ is just $f_0$ with each $\alpha_i$ replaced by $\alpha_{i+k}$
 and each $\delta_i$ replaced by $\delta_{i+k}$;
 but we shall not explicitly use this fact.)
The same combinatorial arguments used earlier give the functional equations
\be
   f_k(t)
   \;=\;
   1 \:+\: \delta_{k+1} t \, f_k(t) \:+\: \alpha_{k+1} t \, f_k(t) \, f_{k+1}(t)
 \label{eq.TRfk.1}
\ee
and
\be
   f_k(t)  \;=\;  {1 \over 1 \:-\: \delta_{k+1} t 
                             \:-\: \alpha_{k+1} t \, f_{k+1}(t)}
   \;,
 \label{eq.TRfk.2}
\ee
which generalize \reff{eq.schroder.1}/\reff{eq.schroder.2}
and reduce to them when $\balpha = \bdelta = \bone$;
they also generalize \reff{eq.SRfk.1}/\reff{eq.SRfk.2}
and reduce to them when  $\bdelta = \bzero$.
Iterating \reff{eq.TRfk.2}, we see immediately that $f_k$
is given by the continued fraction
\be
   f_k(t)
   \;=\;
   \cfrac{1}{1 - \delta_{k+1} t - \cfrac{\alpha_{k+1} t}{1 - \delta_{k+2} t - \cfrac{\alpha_{k+2} t}{1- \delta_{k+3} t - \cfrac{\alpha_{k+3} t}{1- \cdots}}}}
 \label{eq.fk.Tfrac}
\ee
and in particular that $f_0$ is given by
\be
   f_0(t)
   \;=\;
   \cfrac{1}{1 - \delta_{1} t - \cfrac{\alpha_{1} t}{1 - \delta_{2} t - \cfrac{\alpha_{2} t}{1- \delta_{3} t - \cfrac{\alpha_{3} t}{1- \cdots}}}}
   \;\,.
 \label{eq.f0.Tfrac}
\ee
The right-hand sides of \reff{eq.fk.Tfrac}/\reff{eq.f0.Tfrac}
are called \textbfit{Thron-type continued fractions},
or \textbfit{T-fractions} for short.
This combinatorial interpretation of T-fractions
in terms of weighted Schr\"oder paths
was independently discovered by several authors
over the past few years
\cite{Fusy_15,Oste_15,Josuat-Verges_18,Sokal_totalpos}.
(Two decades earlier, Roblet and Viennot \cite{Roblet_96}
 gave an alternate interpretation of T-fractions
 in terms of weighted Dyck paths: see Remark~2 at the end of this section.)

As for Dyck paths, we can consider the combinatorial meaning of
a product $f_0 f_1 \cdots f_\ell$.
We use the term {\em partial Schr\"oder path}\/
to denote a path in the upper half-plane $\Z \times \N$,
using steps $(1,1)$, $(1,-1)$ and $(2,0)$,
that starts on the horizontal axis 
but is allowed to end anywhere in the upper half-plane.
We then claim that the coefficient of $t^n$ in $f_0 f_1 \cdots f_\ell$
is the generating polynomial for partial Schr\"oder paths from
$(0,0)$ to $(2n+\ell,\ell)$ [with the usual weights],
which we denote by $T_{n|\ell}(\balpha,\bdelta)$.
The argument is the same as for Dyck paths, with one small addition:
we observe that the step after the subwalk $\scrp_i$
cannot be a down step or a long level step,
because that would either take us below the horizontal axis
or else violate the fact that we have already seen
the last returns to levels $\le i$.

And finally, the Euler--Gauss method for proving S-fractions,
discussed in \reff{eq.recurrence.gk1.0}--\reff{eq.SRfk.1.bis.0} above,
has a straightforward generalization to T-fractions.
Namely, let $(g_k(t))_{k \ge -1}$ be as before,
and suppose that this sequence satisfies
a linear three-term recurrence of the form
\be
   g_k(t) - g_{k-1}(t)
   \;=\;
   \delta_{k+1} t g_k(t) \:+\: \alpha_{k+1} t \, g_{k+1}(t)
   \qquad\hbox{for } k \ge 0
 \label{eq.recurrence.gk1.T}
\ee
for some coefficients $\balpha = (\alpha_i)_{i \ge 1}$
and $\bdelta = (\delta_i)_{i \ge 1}$ in $R$.
Defining $f_k(t) = g_k(t)/g_{k-1}(t)$ as before,
\reff{eq.recurrence.gk1.T} can be rewritten as \reff{eq.TRfk.1}.
So $f_k(t)$ is the ordinary generating function
for Schr\"oder paths at level $k$ with weights $\balpha$ and $\bdelta$
and is therefore given by the T-fraction \reff{eq.fk.Tfrac};
in particular, $f_0(t) \eqdef g_0(t)/g_{-1}(t)
            = \sum_{n=0}^\infty T_n(\balpha,\bdelta) \, t^n$
where $T_n(\balpha,\bdelta)$ is the Stieltjes--Rogers polynomial
evaluated at the specified values $\balpha$ and~$\bdelta$;
and more generally $f_0(t) \cdots f_k(t) \eqdef g_k(t)/g_{-1}(t)
 = \sum_{n=0}^\infty T_{n|k}(\balpha,\bdelta) \, t^n$.
This method was used by Euler
to derive a T-fraction for ratios of contiguous $\tfo$
\cite[p.~98]{Andrews_99} \cite{Berndt_85,Askey_85,Ramanathan_87}.

\subsection[$m$-Dyck paths]{$\bm{m}$-Dyck paths}

\begin{figure}[!ht]
\begin{center}
\includegraphics[scale=1.5]{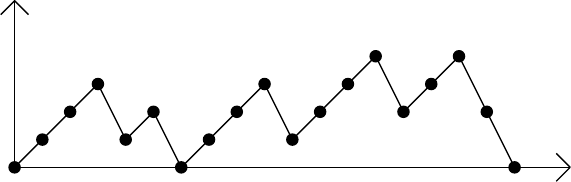}
\caption{A 2-Dyck path of length 18.}
\label{fig.2dyck}
\end{center}
\end{figure}

We recall the following definition \cite{Aval_08,Cameron_16,Prodinger_16}:

\begin{definition}
   \label{definition.mdyck}
\rm
Fix an integer $m \ge 1$.
An \textbfit{$\bm{m}$-Dyck path}
is a path in the upper half-plane $\Z \times \N$,
starting and ending on the horizontal axis,
using steps $(1,1)$ [``rise'' or ``up step'']
and $(1,-m)$ [``$m$-fall'' or ``down step'']:
see Figure~\ref{fig.2dyck} for an example.
More generally, an \textbfit{$\bm{m}$-Dyck path at level $\bm{k}$}
is a path in $\Z \times \N_{\ge k}$,
starting and ending at height $k$,
using steps $(1,1)$ and $(1,-m)$.
\end{definition}

\begin{figure}[!ht]
\begin{center}
\includegraphics[scale=1.5]{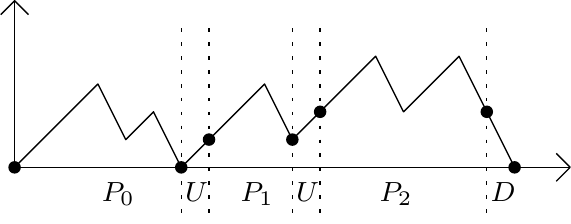}
\caption{\label{fig.2dyckdecomposition}
         The decomposition of a 2-Dyck path.
         Only the important vertices for this decomposition
         are shown explicitly.
}
\end{center}
\end{figure}

Since the number of up steps must equal $m$ times the number of down steps,
the length of an $m$-Dyck path must be a multiple of $m+1$;
we denote by $\scrd^{(m)}_{(m+1)n}$
the set of $m$-Dyck paths from $(0,0)$ to $((m+1)n,0)$.
The ordinary generating function
$D_m(t) = \sum_{n=0}^\infty |\scrd^{(m)}_{(m+1)n}| \, t^n$
satisfies the functional equation
\be
   D_m(t)  \;=\;  1 \,+\, t D_m(t)^{m+1}
   \;,
 \label{eq.mdyck.1}
\ee
as can be seen by splitting the $m$-Dyck path
(if it is of nonzero length)
at its next-to-last return to the horizontal axis;
then we further split the remaining part of the path
at its last return to height 1,
then its last return to height 2, \ldots,
and finally its last return to height $m$;
so the complete path is of the form
$\scrp_0 U \scrp_1 U \scrp_2 U \cdots \scrp_m D$
where $\scrp_i$ is an arbitrary $m$-Dyck path at level $i$
(see Figure~\ref{fig.2dyckdecomposition} for an example
 of this decomposition).

\noindent
The functional equation \reff{eq.mdyck.1}
can equivalently be written in the form
\be
   D_m(t)  \;=\;  {1 \over 1 \,-\, t D_m(t)^m}
   \;,
 \label{eq.mdyck.2}
\ee
which also has a nice combinatorial interpretation:
an $m$-Dyck path can be uniquely decomposed as the concatenation
of zero or more irreducible $m$-Dyck paths;
and an irreducible $m$-Dyck path is of the form
$U \scrp_1 U \scrp_2 \cdots U \scrp_m D$
where $\scrp_1,\ldots,\scrp_m$ are as above.
By applying Lagrange inversion to \reff{eq.mdyck.1},
one finds that $|\scrd^{(m)}_{(m+1)n}|$ equals
the {\em Fuss--Catalan number of order $p=m+1$}\/:
\be
   C_n^{(p)}
  \;=\;  {1 \over (p-1)n+1} \binom{pn}{n}
  \;=\;  {1 \over pn+1} \binom{pn+1}{n}
 \label{def.Cpn}
\ee
\cite[A000108, A001764, A002293, A002294, A002295, A002296
      for $m=1,2,3,4,5,6$]{OEIS}.

We now come to the central object of this paper:

\begin{definition}
   \label{definition.mSR}
\rm
Fix an integer $m \ge 1$,
and let $\balpha = (\alpha_i)_{i \ge m}$ be an infinite set of indeterminates.
The \textbfit{$\bm{m}$-Stieltjes--Rogers polynomial} of order~$n$,
denoted $S^{(m)}_n(\balpha)$, is the generating polynomial
for $m$-Dyck paths of length~$(m+1)n$ in which each rise gets weight~1
and each $m$-fall from height~$i$ gets weight $\alpha_i$.
\end{definition}

Clearly $S_n^{(m)}(\balpha)$ is a homogeneous polynomial
of degree~$n$ with nonnegative integer coefficients.
For instance, for $m=2$ the first few are
\begin{subeqnarray}
   S^{(2)}_0  & = &   1          \\
   S^{(2)}_1  & = &   \alpha_2   \\
   S^{(2)}_2  & = &   \alpha_2^2 + \alpha_2 \alpha_3 + \alpha_2 \alpha_4  \\
   S^{(2)}_3  & = &   \alpha_2^3 + 2 \alpha_2^2 \alpha_3 + \alpha_2 \alpha_3^2 +
 2 \alpha_2^2 \alpha_4 + 2 \alpha_2 \alpha_3 \alpha_4 + 
 \alpha_2 \alpha_4^2 
         \nonumber \\
              & & \quad\;  +\, \alpha_2 \alpha_3 \alpha_5 + 
 \alpha_2 \alpha_4 \alpha_5 + \alpha_2 \alpha_4 \alpha_6
 \label{eq.mSR.m=2.examples}
\end{subeqnarray}

Let $f_0(t) = \sum_{n=0}^\infty S^{(m)}_n(\balpha) \, t^n$
be the ordinary generating function for $m$-Dyck paths with these weights;
and more generally, let $f_k(t)$ be the ordinary generating function
for $m$-Dyck paths at level $k$ with these same weights.
(Obviously $f_k$ is just $f_0$ with each $\alpha_i$ replaced by $\alpha_{i+k}$;
 but we shall not explicitly use this fact.)
The same combinatorial arguments used earlier give the functional equations
\be
   f_k(t)  \;=\;  1 \:+\: \alpha_{k+m} t \, f_k(t) \, f_{k+1}(t) \,\cdots\, f_{k+m}(t)
 \label{eq.mSRfk.1}
\ee
and
\be
   f_k(t)  \;=\;  {1 \over 1 \:-\: \alpha_{k+m} t \, f_{k+1}(t) \,\cdots\, f_{k+m}(t)}
   \;,
 \label{eq.mSRfk.2}
\ee
which generalize \reff{eq.mdyck.1}/\reff{eq.mdyck.2}
and reduce to them when $\balpha = \bone$;
they also generalize \reff{eq.SRfk.1}/\reff{eq.SRfk.2}
and reduce to them when $m=1$.
Iterating \reff{eq.mSRfk.2}, we see immediately that $f_k$
is given by the branched continued fraction
\begin{subeqnarray}
   f_k(t)
   & = &
   \cfrac{1}
         {1 \,-\, \alpha_{k+m} t
            \prod\limits_{i_1=1}^{m}
                 \cfrac{1}
            {1 \,-\, \alpha_{k+m+i_1} t
               \prod\limits_{i_2=1}^{m}
               \cfrac{1}
            {1 \,-\, \alpha_{k+m+i_1+i_2} t
               \prod\limits_{i_3=1}^{m}
               \cfrac{1}{1 - \cdots}
            }
           }
         }
%
      \slabel{eq.fk.mSfrac.a} \\[2mm]
   & = &
\Scale[0.6]{
   \cfrac{1}{1 - \cfrac{\alpha_{k+m} t}{
     \Biggl( 1 - \cfrac{\alpha_{k+m+1} t}{
        \Bigl( 1  - \cfrac{\alpha_{k+m+2} t}{(\cdots) \,\cdots\, (\cdots)} \Bigr)
        \,\cdots\,
        \Bigl( 1  - \cfrac{\alpha_{k+2m+1} t}{(\cdots) \,\cdots\, (\cdots)} \Bigr)
       }
     \Biggr)
     \,\cdots\,
     \Biggl( 1 - \cfrac{\alpha_{k+2m} t}{
        \Bigl( 1  - \cfrac{\alpha_{k+2m+1} t}{(\cdots) \,\cdots\, (\cdots)} \Bigr)
        \,\cdots\,
        \Bigl( 1  - \cfrac{\alpha_{k+3m} t}{(\cdots) \,\cdots\, (\cdots)} \Bigr)
       }
     \Biggr)
    }
   }
}
     \nonumber \\
 \slabel{eq.fk.mSfrac.b}
 \label{eq.fk.mSfrac}
\end{subeqnarray}
and in particular that $f_0$ is given by
the specialization of \reff{eq.fk.mSfrac} to $k=0$.
We shall call the right-hand side of \reff{eq.fk.mSfrac}
an \textbfit{$\bm{m}$-branched Stieltjes-type continued fraction},
or \textbfit{$\bm{m}$-S-fraction} for short.

As for Dyck and Schr\"oder paths,
we can consider a product $f_0 f_1 \cdots f_\ell$.
We use the term {\em partial $m$-Dyck path}\/
to denote a path in the upper half-plane $\Z \times \N$,
using steps $(1,1)$ and $(1,-m)$,
that starts on the horizontal axis 
but is allowed to end anywhere in the upper half-plane.
We claim that the coefficient of $t^n$ in $f_0 f_1 \cdots f_\ell$
is the generating polynomial for partial $m$-Dyck paths from
$(0,0)$ to $((m+1)n+\ell,\ell)$ [with the usual weights],
which we denote by $S^{(m)}_{n|\ell}(\balpha)$.
The argument is the same as for Dyck and Schr\"oder paths.

\begin{figure}[!ht]
\begin{center}
\includegraphics[scale=1.5]{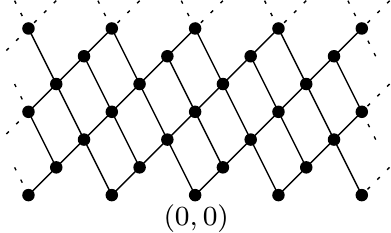}
\caption{\label{fig.Gm}
   The directed graph $G_m$ for $m=2$.
   All edges point towards the right.
}
\end{center}
\end{figure}

All the (partial) $m$-Dyck paths to be considered in this paper
will live in the directed graph $G_m = (V_m, E_m)$ with vertex set
\be
   V_m  \;=\;  \{ (x,y) \in \Z \times \N \colon\:
                  x=y \bmod m+1
               \}
 \label{def.Vm}
\ee
and edge set
\be
   E_m  \;=\;  \bigl\{ \bigl( (x_1,y_1),\, (x_2,y_2) \bigr) \in V_m \times V_m
                       \colon\:
                       x_2 - x_1 = 1 \hbox{ and } y_2 - y_1 \in \{1,-m\}
               \bigr\}
   \;.
 \label{def.Em}
\ee
This is depicted in Figure~\ref{fig.Gm} for $m=2$.
\medskip

\noindent
The key fact about this graph is that it is {\em planar}\/:
this will allow us to apply the Lindstr\"om--Gessel--Viennot lemma
to prove the total positivity of various matrices
associated to the $m$-Stieltjes--Rogers polynomials.

And finally, the Euler--Gauss method for proving S-fractions,
discussed in \reff{eq.recurrence.gk1.0}--\reff{eq.SRfk.1.bis.0} above,
has a straightforward generalization to $m$-S-fractions:

\begin{proposition}[Euler--Gauss recurrence method for $m$-S-fractions]
   \label{prop.euler-gauss.mSR}
Fix an integer $m \ge 1$,
and let $(g_k(t))_{k \ge -1}$ be a sequence of formal power series
(with coefficients in some commutative ring $R$) with constant term 1.
Suppose that this sequence satisfies
a linear three-term recurrence
\be
   g_k(t) - g_{k-1}(t)  \;=\; \alpha_{k+m} t \, g_{k+m}(t)
   \qquad\hbox{for } k \ge 0
 \label{eq.recurrence.gkm.0}
\ee
for some coefficients $\balpha = (\alpha_i)_{i \ge m}$ in $R$.
Defining $f_k(t) = g_k(t)/g_{k-1}(t)$, we have:
\begin{itemize}
   \item[(a)] $f_k(t)$ is the ordinary generating function
for $m$-Dyck paths at level $k$ with weights $\balpha$
and hence is given by the $m$-S-fraction \reff{eq.fk.mSfrac}.
   \item[(b)] In particular,
$f_0(t) = g_0(t)/g_{-1}(t)
            = \sum_{n=0}^\infty S^{(m)}_n(\balpha) \, t^n$
where $S^{(m)}_n(\balpha)$ is the $m$-Stieltjes--Rogers polynomial
evaluated at the specified values $\balpha$.
   \item[(c)] More generally, $f_0(t) \cdots f_k(t) = g_k(t)/g_{-1}(t)
 = \sum_{n=0}^\infty S^{(m)}_{n|k}(\balpha) \, t^n$.
\end{itemize}
\end{proposition}

\noindent
Please note that in the recurrence \reff{eq.recurrence.gkm.0}
we have $k+m$ on the right-hand side,
in place of the $k+1$ that occurred in \reff{eq.recurrence.gk1.0}.

\proof
Using the definition $f_k = g_k/g_{k-1}$,
the recurrence \reff{eq.recurrence.gkm.0}
can be rewritten as \reff{eq.mSRfk.1}.
\qed

We will use this method in Section~\ref{sec.quasi-affine}
to handle some cases with quasi-affine coefficients $\balpha$,
and again in Sections~\ref{sec.hyper.rF0}--\ref{sec.hyper.rphis}
to handle ratios of contiguous hypergeometric series
$\FHyper{r}{s}$ and $\phiHyper{r}{s}$.
The main difficulty in using this method
is that one needs to guess not only the $\balpha$ but also the $(g_k(t))$.

\subsection[$m$-Schr\"oder paths]{$\bm{m}$-Schr\"oder paths}

\begin{figure}[!ht]
\begin{center}
\includegraphics[scale=1.5]{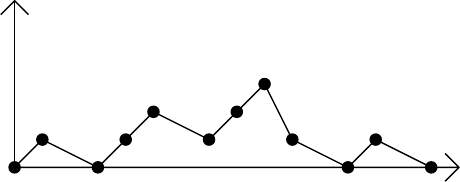}
\caption{\label{fig.schroder} A 2-Schr\"oder path of length 15.}
\end{center}
\end{figure}

In the same way as Dyck paths can be generalized to Schr\"oder paths,
so $m$-Dyck paths can be generalized to $m$-Schr\"oder paths:

\begin{definition}
   \label{definition.mschroder}
\rm
Fix an integer $m \ge 1$.
An \textbfit{$\bm{m}$-Schr\"oder path}
is a path in the upper half-plane $\Z \times \N$,
starting and ending on the horizontal axis,
using steps $(1,1)$ [``rise'' or ``up step''],
$(1,-m)$ [``$m$-fall'' or ``down step'']
and $(2,-(m-1))$ [``$m$-long step'']:
see Figure~\ref{fig.schroder} for an example.
We define the \textbfit{length} of an $m$-Schr\"oder path
to be the number of rises plus the number of $m$-falls
plus twice the number of $m$-long steps.
More generally, an \textbfit{$\bm{m}$-Schr\"oder path at level $\bm{k}$}
is a path in $\Z \times \N_{\ge k}$,
starting and ending at height $k$,
using steps $(1,1)$, $(1,-m)$ and $(2,-(m-1))$.
\end{definition}

It is not difficult to see that
the length of an $m$-Schr\"oder path must be a multiple of $m+1$;
we denote by $\scrs^{(m)}_{(m+1)n}$
the set of $m$-Schr\"oder paths from $(0,0)$ to $((m+1)n,0)$.
The ordinary generating function
$S_m(t) = \sum_{n=0}^\infty |\scrs^{(m)}_{(m+1)n}| \, t^n$
satisfies the functional equation
\be
   S_m(t)  \;=\;  1 \,+\, t S_m(t)^m \,+\, t S_m(t)^{m+1}
   \;,
 \label{eq.mschroder.1}
\ee
as can be seen by splitting the $m$-Schr\"oder path
(if it is of nonzero length)
in the same way as was done for $m$-Dyck paths;
but now, in addition to the path being
$\scrp_0 U \scrp_1 U \scrp_2 U \cdots \scrp_m D$
where $\scrp_i$ is an arbitrary $m$-Schr\"oder path at level $i$,
there is also the possibility that the subwalk $U \scrp_m D$
is replaced by an $m$-long step.
The functional equation \reff{eq.mschroder.1}
can equivalently be written in the form
\be
   S_m(t)  \;=\;  {1 \over 1 \,-\, t S_m(t)^{m-1} \,-\, t S_m(t)^m}
   \;,
 \label{eq.mschroder.2}
\ee
which also has a nice combinatorial interpretation:
an $m$-Schr\"oder path can be uniquely decomposed as the concatenation
of zero or more irreducible $m$-Schr\"oder paths;
and an irreducible $m$-Schr\"oder path
is either $U \scrp_1 U \scrp_2 \cdots U \scrp_m D$
or
\linebreak
$U \scrp_1 U \scrp_2 \cdots U \scrp_{m-1} L$.
By applying Lagrange inversion to \reff{eq.mschroder.1},
one finds that
\be
   |\scrs^{(m)}_{(m+1)n}|
   \;=\;
   {1 \over n}
   \sum_{k=0}^{n-1}
   \binom{mn}{n-1-k}  \binom{n}{k} \, 2^{n-k}
   \quad\hbox{for $n \ge 1$,}
 \label{def.mschroder.n.lagrange}
\ee
or by combinatorial arguments \cite{Song_05} one finds the alternative formula
\be
   |\scrs^{(m)}_{(m+1)n}|
   \;=\;
   \sum_{\ell=0}^n
   {1 \over mn - \ell + 1} \binom{(m+1)n - \ell}{n}  \binom{n}{\ell}
 \label{def.mschroder.n}
\ee
(here $\ell$ corresponds to the number of $m$-long steps).
For $m=1,2,3,4$ these sums are \cite[A006318, A027307, A144097, A260332]{OEIS};
for $m=1,2$ the triangular array \reff{def.mschroder.n}
is \cite[A060693/A088617, A108426]{OEIS}.

\begin{definition}
   \label{definition.mTR}
\rm
Fix an integer $m \ge 1$,
let $\balpha = (\alpha_i)_{i \ge m}$ and $\bdelta = (\delta_i)_{i \ge m}$
be infinite sets of indeterminates.
The \textbfit{$\bm{m}$-Thron--Rogers polynomial} of order~$n$,
denoted $T^{(m)}_n(\balpha,\bdelta)$, is the generating polynomial
for $m$-Schr\"oder paths of length~$(m+1)n$ in which each rise gets weight~1,
each $m$-fall from height~$i$ gets weight $\alpha_i$,
and each $m$-long step from height~$i$ gets weight $\delta_{i+1}$.
\end{definition}

Clearly $T_n^{(m)}(\balpha,\bdelta)$ is a homogeneous polynomial
of degree~$n$ with nonnegative integer coefficients.
For instance, for $m=2$ the first few are
\begin{subeqnarray}
   T^{(2)}_0  & = &   1          \\
   T^{(2)}_1  & = &   \alpha_2 + \delta_2   \\
   T^{(2)}_2  & = &   (\alpha_2 + \delta_2)^2
      + \alpha_2 (\alpha_3 + \alpha_4 + \delta_3 + \delta_4)
      + \delta_2 (\alpha_3 + \delta_3)
\end{subeqnarray}
(and already $T^{(2)}_3$ is quite complicated, containing 41 monomials).

Let $f_0(t) = \sum_{n=0}^\infty T^{(m)}_n(\balpha,\bdelta) \, t^n$
be the ordinary generating function for $m$-Schr\"oder paths
with these weights;
and more generally, let $f_k(t)$ be the ordinary generating function
for $m$-Schr\"oder paths at level $k$ with these same weights.
(Obviously $f_k$ is just $f_0$ with each $\alpha_i$ replaced by $\alpha_{i+k}$
 and each $\delta_i$ replaced by $\delta_{i+k}$;
 but we shall not explicitly use this fact.)
The same combinatorial arguments used earlier give the functional equations
\be
   f_k(t)
   \;=\;
   1 \:+\: \delta_{k+m} t \, f_k(t) \, f_{k+1}(t) \,\cdots\, f_{k+m-1}(t)
     \:+\: \alpha_{k+m} t \, f_k(t) \, f_{k+1}(t) \,\cdots\, f_{k+m}(t)
   \quad
 \label{eq.mTRfk.1}
\ee
and
\be
   f_k(t)
   \;=\;
   {1 \over
    1 \:-\: \delta_{k+m} t \, f_{k+1}(t) \,\cdots\, f_{k+m-1}(t)
      \:-\: \alpha_{k+m} t \, f_{k+1}(t) \,\cdots\, f_{k+m}(t)
   }
   \;,
 \label{eq.mTRfk.2}
\ee
which generalize all the preceding functional equations:
they reduce to \reff{eq.mschroder.1}/\reff{eq.mschroder.2}
when $\balpha = \bdelta = \bone$,
to \reff{eq.mSRfk.1}/\reff{eq.mSRfk.2} when $\bdelta = \bzero$,
and to \reff{eq.TRfk.1}/\reff{eq.TRfk.2} when $m=1$.
Iterating \reff{eq.mTRfk.2}, we see immediately that $f_k$
is given by the branched continued fraction
\be
   f_k(t)
   \;=\;
\Scale[0.88]{
   \cfrac{1}{1 -
 \cfrac{\delta_{k+m} t}{
     \Biggl( 1 - 
\cfrac{\delta_{k+m+1} t}{
        (\cdots) \,\cdots\,
     }
\,-\,
\cfrac{\alpha_{k+m+1} t}{
        (\cdots) \,\cdots\,
     }
     \Biggr)
     \cdots
}
\,-\,
 \cfrac{\alpha_{k+m} t}{
     \Biggl( 1 - 
\cfrac{\delta_{k+m+1} t}{
        (\cdots) \,\cdots\,
     }
-
\cfrac{\alpha_{k+m+1} t}{
        (\cdots) \,\cdots\,
     }
     \Biggr)
     \,\cdots\,
     }
    }
}
 \label{eq.fk.mTfrac}
\ee
and in particular that $f_0$ is given by
the specialization of \reff{eq.fk.mTfrac} to $k=0$.
We shall call the right-hand side of \reff{eq.fk.mTfrac}
an \textbfit{$\bm{m}$-branched Thron-type continued fraction},
or \textbfit{$\bm{m}$-T-fraction} for short.

Let us also consider a product $f_0 f_1 \cdots f_\ell$.
We use the term {\em partial $m$-Schr\"oder path}\/
to denote a path in the upper half-plane $\Z \times \N$,
using steps $(1,1)$, $(1,-m)$ and $(2,-(m-1))$,
that starts on the horizontal axis
but is allowed to end anywhere in the upper half-plane.
We claim that the coefficient of $t^n$ in $f_0 f_1 \cdots f_\ell$
is the generating polynomial for partial $m$-Schr\"oder paths from
$(0,0)$ to $((m+1)n+\ell,\ell)$ [with the usual weights],
which we denote by $T^{(m)}_{n|\ell}(\balpha,\bdelta)$.
The argument is the same as for Dyck and Schr\"oder paths.

\begin{figure}[!ht]
\begin{center}
\includegraphics[scale=1.5]{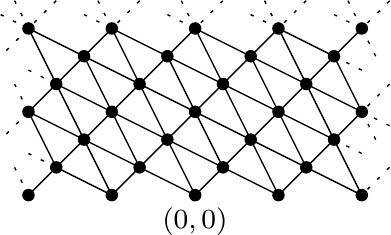}
\caption{\label{fig.Gtildem}
   The directed graph $\widetilde{G}_m$ for $m=2$.
   All edges point towards the right.
}
\end{center}
\end{figure}

All the (partial) $m$-Schr\"oder paths to be considered in this paper
will live in the directed graph $\widetilde{G}_m = (V_m, \widetilde{E}_m)$
with vertex set \reff{def.Vm}
and edge set
\be
   \widetilde{E}_m 
   \;=\;
   E_m  \;\bigcup\;
   \bigl\{ \bigl( (x_1,y_1),\, (x_2,y_2) \bigr) \in V_m \times V_m
                  \colon\:
                  x_2 - x_1 = 2 \hbox{ and } y_2 - y_1 = -(m-1)
   \bigr\}
 \label{def.Etildem}
\ee
where $E_m$ was defined in \reff{def.Em}.
This is depicted in Figure~\ref{fig.Gtildem} for $m=2$.

\noindent
The graph $\widetilde{G}_m$, which augments the graph $G_m$ defined in
\reff{def.Vm}/\reff{def.Em}
by including additional edges associated to the $m$-long steps,
is again planar and acyclic:
this will allow us, once again, to apply the Lindstr\"om--Gessel--Viennot lemma
to prove the total positivity of various matrices
associated to the $m$-Thron--Rogers polynomials.

And finally, the Euler--Gauss method generalized to $m$-S-fractions
as in Proposition~\ref{prop.euler-gauss.mSR} above,
has a straightforward generalization to $m$-T-fractions.
Let $(g_k(t))_{k \ge -1}$ be as before, satisfying a recurrence
\be
   g_k(t) - g_{k-1}(t)
   \;=\;
   \delta_{k+m} t \, g_{k+m-1}(t)  \:+\: \alpha_{k+m} t \, g_{k+m}(t)
   \qquad\hbox{for } k \ge 0
 \label{eq.recurrence.gkm.T}
\ee
for some coefficients $\balpha = (\alpha_i)_{i \ge m}$
and $\bdelta = (\delta_i)_{i \ge m}$.
Defining $f_k(t) = g_k(t)/g_{k-1}(t)$,
\reff{eq.recurrence.gkm.T} can be rewritten as \reff{eq.mTRfk.1}.
So $f_k(t)$ is the ordinary generating function
for $m$-Schr\"oder paths at level $k$ with weights $\balpha$ and $\bdelta$,
and is given by the $m$-T-fraction \reff{eq.fk.mTfrac};
and $f_0(t) \eqdef g_0(t)/g_{-1}(t)
            = \sum_{n=0}^\infty T^{(m)}_n(\balpha,\bdelta) \, t^n$;
and more generally $f_0(t) \cdots f_k(t)$ $\eqdef g_k(t)/g_{-1}(t)
 = \sum_{n=0}^\infty T^{(m)}_{n|k}(\balpha,\bdelta) \, t^n$.

\bigskip

{\bf Remark.}
When $m=1$, a T-fraction with $\balpha = \bzero$ is trivial:
$f_0(t) = 1/(1 - \delta_1 t)$.
Combinatorially this is because a Schr\"oder path without falls
cannot have rises:  it can only be a sequence of zero or more long level steps.
But when $m > 1$, the situation is different:
we have in fact $T_n^{(m)}(\bzero,\bdelta) = S_n^{(m-1)}(\bdelta)$.
This is immediate from the functional equation \reff{eq.mTRfk.1}
for $m$-T-fractions:
when $\balpha = 0$ it reduces to the functional equation \reff{eq.mSRfk.1}
for $(m-1)$-S-fractions with $\balpha$ replaced by $\bdelta$.
It can also be seen combinatorially:
$m$-Schr\"oder paths of length $(m+1)n$ without $m$-falls
are in bijection with $(m-1)$-Dyck paths of length $mn$,
by replacing each $m$-long step by an $(m-1)$-fall.
So the $m$-Thron--Rogers polynomials $T_n^{(m)}(\balpha,\bdelta)$
contain both $T_n^{(m)}(\balpha,\bzero) = S_n^{(m)}(\balpha)$
and $T_n^{(m)}(\bzero,\bdelta) = S_n^{(m-1)}(\bdelta)$ as specializations.
\myendremark

\subsection{Some further remarks}

1.  Many authors, in presenting Dyck and Schr\"oder paths,
rotate our picture by $+45^\circ$:
then the roles of up, down and long level steps are played, respectively,
by $(0,1)$ [``north''], $(1,0)$ [``east'']
and $(1,1)$ [``diagonal'' or ``northeast''],
and the path runs from $(0,0)$ to $(n,n)$
while staying in the region $y \ge x$.
Alternatively, the roles of north and east can be reversed:
then the path stays in the region $y \le x$.
Following Duchon \cite{Duchon_00},
we can call this the ``northeastbound'' representation
of Dyck and Schr\"oder paths.

Likewise, in presenting $m$-Dyck and $m$-Schr\"oder paths,
some authors \cite{Song_05,Bergeron_12}
use an analogously turned picture
(although for $m \ne 1$ the required linear transformation is no longer
angle-preserving):
the roles of up, down and long steps are again played
by $(0,1)$, $(1,0)$ and $(1,1)$,
and now the path runs from $(0,0)$ to $(n,mn)$
while staying in the region $y \ge mx$.
Alternatively, the roles of north and east can be reversed:
then the path runs from $(0,0)$ to $(mn,n)$
while staying in the region $y \le x/m$.
We again call this the ``northeastbound'' representation.

However, we prefer the ``horizontal'' (or ``eastbound'')
presentation employed here,
which makes clear that the paths are {\em directed}\/
along the horizontal (``time'') axis.
See \cite{Duchon_00} for a comparison of the
northeastbound and eastbound representations;
and see \cite{Banderier_02} for a general discussion
of enumeration of directed lattice paths.

2. Roblet and Viennot \cite{Roblet_96} gave an alternate interpretation
of the Thron--Rogers polynomials $T_n(\balpha,\bdelta)$
in terms of Dyck paths.
A {\em peak}\/ of a Dyck path is a vertex
that is arrived at by a rise and departed from by a fall;
thus, a {\em fall from a peak}\/ is a fall preceded by a rise,
while a {\em fall from a non-peak}\/ is a fall preceded by a fall.
(Note that the {\em first}\/ step of a Dyck path {\em cannot}\/ be a fall!)
Roblet and Viennot proved \cite[Proposition~1]{Roblet_96}
that $T_n(\balpha,\bdelta)$
is the generating polynomial for Dyck paths of length $2n$ in which
each rise gets weight 1,
each fall from a non-peak starting at height~$i$ gets weight $\alpha_i$,
and each fall from a peak starting at height~$i$
gets weight $\alpha_i + \delta_i$.
This is in fact an easy consequence of the Schr\"oder-path interpretation
that we have already discussed:
it suffices to apply the surjection of Schr\"oder paths onto Dyck paths
defined by replacing each long level step by a pair $UD$
(thus creating a new peak).
When $\bdelta = \bzero$,
falls from peaks and non-peaks get the same weight,
and we recover $T_n(\balpha,\bzero) = S_n(\balpha)$.

A similar mapping shows that the $m$-Thron--Rogers polynomial
$T_n^{(m)}(\balpha,\bdelta)$ 
is the generating polynomial for $m$-Dyck paths of length $(m+1)n$ in which
each rise gets weight 1,
each $m$-fall from a non-peak starting at height~$i$ gets weight $\alpha_i$,
and each $m$-fall from a peak starting at height~$i$ 
gets weight $\alpha_i + \delta_i$.

3. For weighted $m$-Dyck paths with $m=2$,
equations \reff{eq.mSRfk.1}--\reff{eq.fk.mSfrac}
were presented briefly by Varvak \cite[pp.~3--4]{Varvak_04};
but she did not dwell on them, because her main goal was
to motivate a more general class of paths ---
namely, \L{}ukasiewicz paths, in which $m$-falls of
{\em all}\/ $m \ge 0$ are simultaneously allowed ---
and the corresponding branched continued fractions.
We will discuss this construction in Section~\ref{sec.m-JR}.
In Section~\ref{sec.totalpos} we will show that the $m$-Dyck paths possess
a desirable property not possessed by the more general \L{}ukasiewicz paths,
namely, Hankel-total positivity.

4.  Stieltjes-type continued fractions \reff{eq.f0.Sfrac}
are an extraordinarily useful tool in combinatorics:
one reason for this is that ``generically''
(i.e.\ barring certain degenerate cases)
there is a one-to-one correspondence between
formal power series $f_0(t) = \sum_{n=0}^\infty a_n t^n$ with $a_0 = 1$
and S-fractions \reff{eq.f0.Sfrac}.
Thron-type continued fractions \reff{eq.f0.Tfrac},
by contrast, are more difficult to work with,
because there is a high degree of nonuniqueness:
the number of free parameters at each level is twice as many as needed.
The $m$-S-fractions and $m$-T-fractions
suffer from an even greater nonuniqueness:
indeed, already from the definition in terms of $m$-Dyck paths
one can see that the $m$-Stieltjes--Rogers polynomial $S_n^{(m)}(\balpha)$
depends on the coefficients $\alpha_i$ for $m \le i \le mn$
[cf.\ \reff{eq.mSR.m=2.examples} for $m=2$ and $n \le 3$],
so that the number of free parameters at each level
is roughly $m$ times as many as needed;
and the $m$-Thron--Rogers polynomials are twice as bad as this.
So it is highly nontrivial to know whether a given sequence
$\ba = (a_n)_{n \ge 0}$ has a nice $m$-S-fraction for some given $m \ge 2$;
this has to be done largely by trial and error.

In addition, the branched continued fractions
\reff{eq.fk.mSfrac} and \reff{eq.fk.mTfrac} are irredeemably ugly,
and seem almost impossible to work with (or even to write!).
We shall therefore hardly ever write out a branched continued fraction;
we shall instead work directly with the $m$-Dyck and $m$-Schr\"oder paths
and/or with the recurrences
\reff{eq.mSRfk.1}/\reff{eq.mSRfk.2} and \reff{eq.mTRfk.1}/\reff{eq.mTRfk.2}
that their generating functions satisfy.

5.  In a very interesting paper, Albenque and Bouttier \cite{Albenque_12}
introduced some polynomials that are essentially equivalent to our
$m$-Stieltjes--Rogers polynomials, along with the corresponding
branched continued fraction.  Their expressions are slightly different
from ours, because they chose to weight the rises rather than the $m$-falls.
(They also interchanged rises and falls compared to our notation,
 but that is a trivial change, which is equivalent to reversing the path.)
They noted, as we just did, that the $m$-S-fractions have roughly
$m$ times as many unknowns as equations;
and they showed how to restore the uniqueness by specifying
not only the generating polynomials of $m$-Dyck paths
but also the generating polynomials of a suitable set of
partial $m$-Dyck paths
\cite[Theorem~3 and Corollary~4]{Albenque_12}.
Finally, they applied these branched continued fractions
to the enumeration of certain classes of planar maps.

6.  Our terminology for continued fractions
follows the general practice in the combinatorial literature,
starting with Flajolet \cite{Flajolet_80}.
The classical analytic literature on continued fractions
\cite{Perron,Wall_48,Jones_80,Lorentzen_92,Cuyt_08}
generally uses a different terminology.
For instance, Jones and Thron \cite[pp.~128--129, 386--389]{Jones_80}
use the term ``regular C-fraction''
for (a minor variant of) what we~have called an S-fraction;
they call it an ``S-fraction'' if all $\alpha_n < 0$.
\myendremark

\section{Relation between different values of $\bm{m}$}
   \label{sec.different_m}

From the definitions given in the preceding section,
it might seem that there is no relation between the
$m$-Stieltjes--Rogers polynomials for different values of $m \ge 1$:
each value of $m$ seems to provide a distinct generalization
of the classical Stieltjes--Rogers polynomials.
But appearances are deceptive:
the actual situation
is that the 1-Stieltjes--Rogers polynomials are contained in the
2-Stieltjes--Rogers polynomials, which are in turn contained in the
3-Stieltjes--Rogers polynomials, and so forth;
thus, the $m$-Stieltjes--Rogers polynomials become increasingly
general as $m$ grows.
To see this, let us show explicitly how
the $m$-Stieltjes--Rogers polynomials are equal to
a suitable specialization of the $m'$-Stieltjes--Rogers polynomials
for any $m' > m$.

Given $1 \le m < m'$, let us say that a subset $I \subseteq [m',\infty)$
is \textbfit{$\bm{(m,m')}$-good} if
\begin{itemize}
   \item[(a)]  $m' \notin I$;  and
   \item[(b)]  For each $j \ge m'$, either $j \in I$ or
      exactly $m'-m$ of the elements $j+1,\ldots,j+m'$ belong to $I$.
\end{itemize}
We then have:

\begin{proposition}
   \label{prop.reduction}
Fix $1 \le m < m'$,
and let $I$ be an $(m,m')$-good subset of $[m',\infty)$.
Let $\balpha = (\alpha_i)_{i \ge m}$ be an infinite set of indeterminates,
and define $\balpha' = (\alpha'_i)_{i \ge m'}$ by
setting $\alpha'_i = 0$ for $i \in I$,
and then setting $(\alpha'_i)_{i \in [m',\infty) \setminus I}$
to be $(\alpha_i)_{i \ge m}$ in increasing order.
Then $S_n^{(m)}(\balpha) = S_n^{(m')}(\balpha')$.
\end{proposition}

We will give two proofs of Proposition~\ref{prop.reduction}:
a first proof based on manipulating the branched continued fraction
\reff{eq.fk.mSfrac},
and a second proof based on the functional equation \reff{eq.mSRfk.2}.
For the special case in which the set $I$ is periodic of period~${m+1}$,
we will also give (in Section~\ref{subsec.production.SR.TR.JR}) a third proof,
based on consideration of the corresponding production matrices.

\smallskip

\firstproof
Consider the $m'$-branched continued fraction \reff{eq.fk.mSfrac.b}
with $k=0$ (and $m$ replaced by $m'$, and $\balpha$ by $\balpha'$).
Since $I$ is an $(m,m')$-good subset of $[m',\infty)$,
condition (a) guarantees that $\alpha'_{m'} = \alpha_m$.
Then condition (b) guarantees that exactly $m'-m$ of the subsequent branches
are killed (because the index $i$ corresponding to the root of that branch
belongs to the set $I$);
and then, for the $m$ branches that are not killed,
exactly $m'-m$ of their sub-branches will be killed;
and so on at all levels.
The result is an $m$-branched continued fraction \reff{eq.fk.mSfrac.b}
with the weights $\balpha$.
\qed

\secondproof
For each $k \ge 0$, let $f'_k(t)$ be the ordinary generating function
for $m'$-Dyck paths at level $k$ with the weights $\balpha'$.
[The prime here is {\em not}\/ intended as a derivative!]
Then $(f'_k)_{k \ge 0}$ satisfies the functional equation \reff{eq.mSRfk.2}
with $m$ replaced by $m'$ and $\alpha_{k+m}$ by $\alpha'_{k+m'}$.
Since $\alpha'_i = 0$ for $i \in I$,
it follows that $f'_k = 1$ for $k \in I-m'$.
[Note that $0 \notin I-m'$ by hypothesis~(a).]
By hypothesis~(b), for each $k \in \N \setminus (I-m')$,
exactly $m'-m$ of the indices $i \in \{k+1,\ldots,k+m'\}$ belong to $I-m'$
and hence have $f'_i = 1$.
So if we define $(f_k)_{k \ge 0} = (f'_k)_{k \in \N \setminus (I-m')}$
where both sequences are taken in increasing order,
then the functional equation \reff{eq.mSRfk.2} for $(f'_k)_{k \ge 0}$
[with $m$ replaced by $m'$ and $\alpha_{k+m}$ by $\alpha'_{k+m'}$]
is equivalent to the functional equation \reff{eq.mSRfk.2} for $(f_k)_{k \ge 0}$
[with $m$ and $\balpha$ as is].
\qed

\bigskip

{\bf Remarks.}
1.  The first proof of Proposition~\ref{prop.reduction}
is the {\em only}\/ place in this paper where we will find it convenient
to use the branched continued fraction \reff{eq.fk.mSfrac}.

2.  The second proof can be reinterpreted combinatorially
as a bijection from $m$-Dyck paths to the subclass of $m'$-Dyck paths
that have no falls from heights in the set $I$.
We leave the details to the reader.

3.  For the case $(m,m') = (1,2)$,
it is easy to give a ``constructive'' explanation
of what it means for a set $I$ to be $(m,m')$-good:
namely, $I \subseteq [2,\infty)$ is $(1,2)$-good
if and only if the sequence $2,3,\ldots$ starts with an element of $I^c$
(that is, $2 \in I^c$)
and consists of blocks of either one or two elements of $I^c$
alternating with single elements of $I$.

It follows that, given any classical S-fraction with $\alpha_i \ne 0$,
there are {\em uncountably many}\/ 2-S-fractions that are equivalent to it!
This illustrates the extreme nonuniqueness of $m$-S-fractions with $m \ge 2$.
(Later we will see also many examples of this nonuniqueness
 that do not merely involve inserting zeroes into an existing sequence.)
\myendremark

\bigskip

For the $m$-Thron--Rogers polynomials we have a more restrictive result,
in which the set $I$ must be taken to be periodic of period $m'$:

\begin{proposition}
   \label{prop.reduction.thron}
Fix $1 \le m < m'$,
and let $J \subseteq \{1,\ldots,m'-1\}$ be a set of cardinality $m'-m$.
Then let $I = \{ i \ge m' \colon\: \hbox{there exists } j \in J \hbox{ with }
 i \equiv j \bmod m'\}$.
Let $\balpha = (\alpha_i)_{i \ge m}$ and $\bdelta = (\delta_i)_{i \ge m}$
be infinite sets of indeterminates,
and define $\balpha' = (\alpha'_i)_{i \ge m'}$ by
setting $\alpha'_i = 0$ for $i \in I$,
and then setting $(\alpha'_i)_{i \in [m',\infty) \setminus I}$
to be $(\alpha_i)_{i \ge m}$ in increasing order,
and likewise for $\bdelta'$.
Then $T_n^{(m)}(\balpha,\bdelta) = T_n^{(m')}(\balpha',\bdelta')$.
\end{proposition}

\noindent
Note that $I$ does not contain any multiple of $m'$.

\proof
We imitate the second proof of Proposition~\ref{prop.reduction}.
For each $k \ge 0$, let $f'_k(t)$ be the ordinary generating function
for $m'$-Schr\"oder paths at level $k$ with the weights $\balpha'$
and $\bdelta'$.
Then $(f'_k)_{k \ge 0}$ satisfies the functional equation \reff{eq.mTRfk.2}
with $m$ replaced by $m'$, $\alpha_{k+m}$ by $\alpha'_{k+m'}$,
and  $\delta_{k+m}$ by $\delta'_{k+m'}$.
Since $\alpha'_i = \delta'_i = 0$ for $i \in I$,
it follows that $f'_k = 1$ for $k \in I-m'$.
[Note that $0 \notin I-m'$.]
Then, for each $k \in \N \setminus (I-m')$,
exactly $m'-m$ of the indices $i \in \{k+1,\ldots,k+m'-1\}$ belong to $I-m'$,
and moreover $k+m' \notin I-m'$.
So, whenever $k \in \N \setminus (I-m')$,
{\em both}\/ of the terms
$\delta'_{k+m'} t f'_{k+1} \cdots f'_{k+m'-1}$
and $\alpha'_{k+m'} t f'_{k+1} \cdots f'_{k+m'}$
will reduce as desired.
So if we define $(f_k)_{k \ge 0} = (f'_k)_{k \in \N \setminus (I-m')}$
where both sequences are taken in increasing order,
then the functional equation \reff{eq.mTRfk.2} for $(f'_k)_{k \ge 0}$
[with $m$ replaced by $m'$, etc.]\ 
is equivalent to the functional equation
\reff{eq.mTRfk.2} for $(f_k)_{k \ge 0}$.
\qed

\section{The $\bm{m}$-Jacobi--Rogers polynomials}
   \label{sec.m-JR}

\begin{figure}[!ht]
\begin{center}
\includegraphics[scale=1.5]{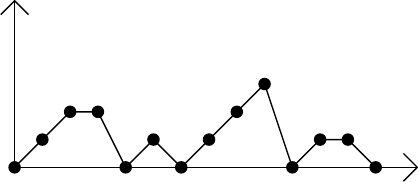}
\caption{\label{fig.lukasiewicz} A \L{}ukasiewicz path of length 13.}
\end{center}
\end{figure}

The following class of paths is well known:

\begin{definition}
   \label{definition.lukasiewicz}
\rm
A \textbfit{\L{}ukasiewicz path}
is a path in the upper half-plane $\Z \times \N$,
starting and ending on the horizontal axis,
using steps $(1,r)$ with $r \le 1$:
the allowed steps are thus $r=1$ (``rise''),
$r=0$ (``level step''), and $r = -\ell$ for any $\ell > 0$ (``$\ell$-fall).
See Figure~\ref{fig.lukasiewicz} for an example.
More generally, a \textbfit{\L{}ukasiewicz  path at level~$\bm{k}$}
is a path in $\Z \times \N_{\ge k}$, starting and ending at height~$k$,
using these same steps.
A \L{}ukasiewicz path is called an \textbfit{$\bm{m}$-\L{}ukasiewicz path}
if it uses only the steps $(1,r)$ with $-m \le r \le 1$.
A 1-\L{}ukasiewicz path is called a \textbfit{Motzkin path}.
We also use the term ``$\infty$-\L{}ukasiewicz path''
as a synonym of ``\L{}ukasiewicz path''.
For $1 \le m \le \infty$,
we denote by $\scrl_n^{(m)}$ the set of $m$-\L{}ukasiewicz paths
from $(0,0)$ to $(n,0)$.\footnote{
   {\bf Warning:}  Many authors define \L{}ukasiewicz paths
   by interchanging rises and falls compared to our definition:
   thus, for them a ``\L{}ukasiewicz path''
   (resp.\ ``$m$-\L{}ukasiewicz path'')
   uses steps $(1,r)$ with $r \ge -1$ (resp.\ $-1 \le r \le m$).
   We shall call this a
   {\em reversed \L{}ukasiewicz}\/ (resp.\ {\em reversed $m$-\L{}ukasiewicz}\/)
   {\em path}\/.
   Obviously, $\omega = (\omega_0,\ldots,\omega_n)$
   is a \L{}ukasiewicz (resp.\ $m$-\L{}ukasiewicz) path
   if and only if its reversal $\omega^{\rm rev} = (\omega_n,\ldots,\omega_0)$
   is a reversed \L{}ukasiewicz (resp.\ reversed $m$-\L{}ukasiewicz) path.
}
\end{definition}

\begin{definition}
   \label{definition.mJR}
\rm
Fix $1 \le m \le \infty$,
and let $\bbeta = (\beta_i^{(\ell)})_{0 \le \ell \le m, \, i \ge \ell}$
be indeterminates.
The \textbfit{$\bm{m}$-Jacobi--Rogers polynomial} of order~$n$,
denoted $J^{(m)}_{n}(\bbeta)$, is the generating polynomial
for $m$-\L{}ukasiewicz paths from $(0,0)$ to $(n,0)$
in which each rise gets weight~1,
each level step at height~$i$ gets weight $\beta_i^{(0)}$,
and each $\ell$-fall from height~$i$ gets weight $\beta_i^{(\ell)}$.
\end{definition}

For $m=1$ the $m$-Jacobi--Rogers polynomials
reduce to the standard Jacobi--Rogers polynomials
as introduced by Flajolet \cite{Flajolet_80}.
When $m=1$ we will also use the notation
 $\beta_i = \beta_i^{(1)}$ and $\gamma_i = \beta_i^{(0)}$,
 and write the Jacobi--Rogers polynomials as
 $J_n(\bbeta,\bgamma)$ where
 $\bbeta = (\beta_i)_{i \ge 1}$ and $\bgamma = (\gamma_i)_{i \ge 0}$.

Let $f_0(t) = \sum_{n=0}^\infty J^{(m)}_{n}(\bbeta) \, t^n$
be the ordinary generating function for $m$-\L{}ukasiewicz paths
with these weights;
and more generally, let $f_k(t)$ be the ordinary generating function
for $m$-\L{}ukasiewicz paths at level $k$ with these same weights.
(Obviously $f_k$ is just $f_0$ with each $\beta_i^{(\ell)}$
replaced by $\beta_{i+k}^{(\ell)}$;
 but we shall not explicitly use this fact.)
We can derive a recurrence for the $f_k$ as follows:
Split an $m$-\L{}ukasiewicz path $\scrp$ at level $k$
(if it is of nonzero length) at its next-to-last return to level $k$:
it is of the form $\scrp = \scrp_k \scrp'$
where $\scrp_k$ is an arbitrary $m$-\L{}ukasiewicz path at level $k$
and $\scrp'$ is an irreducible $m$-\L{}ukasiewicz path at level $k$.
Now $\scrp'$ is either a single level step,
or else it starts with a rise and ends with an $\ell$-fall
for some $\ell$ ($1 \le \ell \le m$).
In the latter case, split $\scrp'$ at its {\em last}\/ return to level $k+1$,
then its last return to level $k+2$, and so on,
until its last return to level $k+\ell$; then it is of the form
$\scrp' = U \scrp_{k+1} U \scrp_{k+2} \cdots U \scrp_{k+\ell} D$
where $\scrp_i$ is an arbitrary $m$-\L{}ukasiewicz path at level $i$
and $D$ is an $\ell$-fall.
It follows that \cite[pp.~V-38--V-39]{Viennot_83}
\cite[pp.~22, 143]{Roblet_94} \cite[p.~5]{Varvak_04}
\be
   f_k(t)  \;=\;  1 \:+\: \sum_{\ell=0}^m \beta_{k+\ell}^{(\ell)} \, t^{\ell+1}
                   \, f_k(t) \, f_{k+1}(t) \,\cdots\, f_{k+\ell}(t)
 \label{eq.mJRfk.1}
\ee
and hence
\be
   f_k(t)
   \;=\;
   {1 \over 1 \:-\: \sum\limits_{\ell=0}^m \beta_{k+\ell}^{(\ell)} \, t^{\ell+1}  
                   \, f_{k+1}(t) \,\cdots\, f_{k+\ell}(t)}
   \;.
 \label{eq.mJRfk.2}
\ee
When $m=1$ these reduce to the usual functional equations
for weighted Motzkin paths \cite{Flajolet_80};
and when $\beta_i^{(m)} = \alpha_i$
and $\beta_i^{(\ell)} = 0$ for $\ell \ne m$,
they reduce to the functional equations \reff{eq.mSRfk.1}/\reff{eq.mSRfk.2}
for weighted $m$-Dyck paths after a change of variables $t \to t^{m+1}$.

We can iterate \reff{eq.mJRfk.2} to obtain a branched continued fraction
\cite[pp.~{V-39}--\hbox{V-40}]{Viennot_83}
\cite[pp.~22--24, 143]{Roblet_94} \cite[p.~5]{Varvak_04}:
\be
   f_k(t)
   \;=\;
\Scale[0.93]{
   \cfrac{1}
         {1 \,-\, \beta_{k}^{(0)} t
            \,-\, \sum\limits_{\ell_1=1}^m \beta_{k+\ell_1}^{(\ell_1)} t^{\ell_1+1}
            \prod\limits_{i_1=1}^{\ell_1}
                 \cfrac{1}
            {1 \,-\, \beta_{k+i_1}^{(0)} t
               \,-\, \sum\limits_{\ell_2=1}^m \beta_{k+i_1+\ell_2}^{(\ell_2)} t^{\ell_2+1}
               \prod\limits_{i_2=1}^{\ell_2}
               \cfrac{1}{1 - \cdots}
            }
         }
}
   \;.
 \label{eq.fk.mJfrac}
\ee
We shall call the right-hand side of \reff{eq.fk.mJfrac}
an \textbfit{$\bm{m}$-branched Jacobi-type continued fraction},
or \textbfit{$\bm{m}$-J-fraction} for short.
Some authors \cite{Viennot_83,Roblet_94}
call it a {\em \L{}ukasiewicz}\/ (or {\em $m$-\L{}ukasiewicz}\/)
{\em continued fraction}\/,
or {\em L-fraction}\/ for short.
When $m=1$ it reduces to the classical J-fraction
studied by Flajolet \cite{Flajolet_80}.

\section{The generalized $\bm{m}$-Stieltjes--Rogers, $\bm{m}$-Thron--Rogers
   and $\bm{m}$-Jacobi--Rogers polynomials}   \label{sec.generalized}

In this section we introduce triangular arrays of generalized
$m$-Stieltjes--Rogers, $m$-Thron--Rogers and $m$-Jacobi--Rogers polynomials,
whose first columns are given by the ordinary
$m$-Stieltjes--Rogers, $m$-Thron--Rogers and $m$-Jacobi--Rogers polynomials,
respectively.
These generalized polynomials will play an important role
when we discuss production matrices
and their connection with total positivity
(Sections~\ref{sec.production} and \ref{sec.totalpos}).

\bigskip

{\bf Generalized $\bm{m}$-Stieltjes--Rogers polynomials.}
Fix an integer $m \ge 1$.
Recall that a {\em partial $m$-Dyck path}\/
is a path in the upper half-plane $\Z \times \N$,
starting on the horizontal axis but ending anywhere,
using steps $(1,1)$ [``rise'']
and $(1,-m)$ [``$m$-fall''].
A partial $m$-Dyck path starting at $(0,0)$
must stay always within the set $V_m$ defined in \reff{def.Vm}:
that is, every point $(x,y)$ of the path
--- and in particular the final point ---
must satisfy $x=y \bmod m+1$.
We denote by $\scrd^{(m)}_{(m+1)n,(m+1)k}$
the set of partial $m$-Dyck paths from $(0,0)$ to $((m+1)n,(m+1)k)$.

Now let $\balpha = (\alpha_i)_{i \ge m}$ be an infinite set of indeterminates,
and let $S^{(m)}_{n,k}(\balpha)$ be the generating polynomial
for partial $m$-Dyck paths from $(0,0)$ to ${((m+1)n,(m+1)k)}$
in~which each rise gets weight~1
and each $m$-fall from height~$i$ gets weight $\alpha_i$.
We call the $S^{(m)}_{n,k}$ the
\textbfit{generalized $\bm{m}$-Stieltjes--Rogers polynomials}.
Obviously $S^{(m)}_{n,k}$ is nonvanishing only for $0 \le k \le n$,
and $S^{(m)}_{n,n} = 1$.
We therefore have an infinite unit-lower-triangular array
$\sfS^{(m)} = \big( S^{(m)}_{n,k}(\balpha) \big)_{n,k \ge 0}$
in which the first ($k=0$) column displays
the ordinary $m$-Stieltjes--Rogers polynomials $S^{(m)}_{n,0} = S^{(m)}_n$.

Later (in Section~\ref{subsec.periodic.row-generating})
we will want to study the \textbfit{row-generating polynomials}
associated to the lower-triangular matrix $\sfS^{(m)}$, namely
\be
   S^{(m)}_n(\balpha;\xi)
   \;\eqdef\;
   \sum_{k=0}^n S^{(m)}_{n,k}(\balpha) \: \xi^k
 \label{def.Sn.rowgen}
\ee
where $\xi$ is an indeterminate.
More generally, it will prove useful to define \cite{Chang_16,Mu_17a,Zhu_17a}
the \textbfit{row-generating matrix}
$\sfS^{(m)}(\xi) = (S^{(m)}_{n,k}(\balpha;\xi))_{n,k \ge 0}$ by
\be
   S^{(m)}_{n,k}(\balpha;\xi)
   \;\eqdef\;
   \sum_{\ell=k}^n S^{(m)}_{n,\ell}(\balpha) \: \xi^{\ell-k}
   \;,
 \label{def.Snk.rowgen}
\ee
which reduces to the row-generating polynomials when $k=0$.
Please note that the definition \reff{def.Snk.rowgen}
can be written in matrix form as
\be
   \sfS^{(m)}(\xi)  \;=\;  \sfS^{(m)} \, T_\xi
\ee
where $T_\xi$ is the lower-triangular Toeplitz matrix of powers of $\xi$:
\be
   (T_\xi)_{\ell k}  \;=\;  \begin{cases}
                           \xi^{\ell-k}  & \textrm{if $\ell \ge k$}  \\
                           0           & \textrm{if $\ell < k$}
                      \end{cases}
 \label{def.Txi}
\ee

\bigskip

{\bf Generalized $\bm{m}$-Thron--Rogers polynomials.}
Fix an integer $m \ge 1$.
Recall that a {\em partial $m$-Schr\"oder path}\/
is a path in the upper half-plane $\Z \times \N$,
starting on the horizontal axis but ending anywhere,
using steps $(1,1)$ [``rise''],
$(1,-m)$ [``$m$-fall'']
and $(2,-(m-1))$ [``$m$-long step''].
A partial $m$-Schr\"oder path starting at $(0,0)$
must stay always within the set $V_m$ defined in \reff{def.Vm}.
We denote by $\scrs^{(m)}_{(m+1)n,(m+1)k}$
the set of partial $m$-Schr\"oder paths from $(0,0)$ to $((m+1)n,(m+1)k)$.

Now let $\balpha = (\alpha_i)_{i \ge m}$ and $\bdelta = (\delta_i)_{i \ge m}$
be infinite sets of indeterminates,
and let $T^{(m)}_{n,k}(\balpha,\bdelta)$ be the generating polynomial
for partial $m$-Schr\"oder paths from $(0,0)$ to $((m+1)n,(m+1)k)$
in which each rise gets weight~1,
each $m$-fall from height~$i$ gets weight $\alpha_i$,
and each $m$-long step at height~$i$ gets weight $\delta_{i+1}$.
We call the $T^{(m)}_{n,k}$ the
\textbfit{generalized $\bm{m}$-Thron--Rogers polynomials}.
Once again we have an infinite unit-lower-triangular array
$\sfT^{(m)} = \big( T^{(m)}_{n,k}(\balpha,\bdelta) \big)_{n,k \ge 0}$
in which the first ($k=0$) column displays
the ordinary $m$-Thron--Rogers polynomials $T^{(m)}_{n,0} = T^{(m)}_n$.

\bigskip

{\bf Generalized $\bm{m}$-Jacobi--Rogers polynomials.}
Fix $1 \le m \le \infty$.
A {\em partial $m$-\L{}ukasiewicz path}\/
is a path in the upper half-plane $\Z \times \N$,
starting on the horizontal axis but ending anywhere,
using steps $(1,r)$ with $-m \le r \le 1$.
We denote by $\scrl_{n,k}^{(m)}$ the set of partial $m$-\L{}ukasiewicz paths
from $(0,0)$ to $(n,k)$.

Now let $\bbeta = (\beta_i^{(\ell)})_{0 \le \ell \le m, \, i \ge \ell}$
be indeterminates,
and let $J^{(m)}_{n,k}(\bbeta)$ be the generating polynomial
for partial $m$-\L{}ukasiewicz paths from $(0,0)$ to $(n,k)$
in which each rise gets weight~1,
each level step at height~$i$ gets weight $\beta_i^{(0)}$,
and each $\ell$-fall from height~$i$ gets weight $\beta_i^{(\ell)}$.
We call the $J^{(m)}_{n,k}$ the
\textbfit{generalized $\bm{m}$-Jacobi--Rogers polynomials}.
Once again we have an infinite unit-lower-triangular array
$\sfJ^{(m)} = \big( J^{(m)}_{n,k}(\bbeta) \big)_{n,k \ge 0}$
in which the first ($k=0$) column displays
the ordinary $m$-Jacobi--Rogers polynomials $J^{(m)}_{n,0} = J^{(m)}_n$.

\bigskip

{\bf Remark.}
When $m=1$, the generalized Stieltjes--Rogers and Jacobi--Rogers
polynomials play a role that does {\em not}\/ generalize to $m > 1$:
namely, they provide an $LDL^{\rm T}$ factorization
of the Hankel matrices associated to the sequences
$\bS = (S_n(\balpha))_{n \ge 0}$ and $\bJ = (J_n(\bbeta,\bgamma))_{n \ge 0}$
of ordinary Stieltjes--Rogers and Jacobi--Rogers polynomials:
\begin{eqnarray}
   H_\infty(\bS)  & = &  \sfS D \sfS^{\rm T}
     \label{eq.LDLT.Stype}  \\[2mm]
   H_\infty(\bJ)  & = &  \sfJ D' \sfJ^{\rm T}
     \label{eq.LDLT.Jtype}
\end{eqnarray}
for suitable diagonal matrices $D$ and $D'$
(the entries of which are partial products of coefficients
 $\balpha$ or $\bbeta$, respectively).
This factorization is due to Stieltjes \cite{Stieltjes_1889,Stieltjes_1894};
see \cite{Sokal_totalpos} for further discussion.
The combinatorial proof of \reff{eq.LDLT.Stype}
[resp.\ \reff{eq.LDLT.Jtype}]
is based on splitting a Dyck path of length $2(n+n')$
[resp.\ a Motzkin path of length $n+n'$]
into its first $2n$ (resp.\ $n$) steps
and its last $2n'$ (resp.\ $n'$) steps
and then imagining the second part run backwards;
this works since the reverse of a Dyck (resp.\ Motzkin) path
is again a Dyck (resp.\ Motzkin) path,
with rises and falls interchanged
(this interchange of rises and falls gives rise to the factors
 $\balpha$ or $\bbeta$ occurring in the matrices $D$ or $D'$).
By contrast, for $m > 1$ the reverse of an $m$-Dyck
(resp.\ $m$-\L{}ukasiewicz) path
is {\em not}\/ an $m$-Dyck (resp.\ $m$-\L{}ukasiewicz) path,
and this method fails.
In fact, it {\em must}\/ fail,
since it can be shown \cite{Sokal_totalpos}
(based on ideas from \cite{Peart_00,Woan_01,Aigner_01a})
that for a unit-lower-triangular matrix $L$,
there exists a nontrivial diagonal matrix $D$
such that $LDL^{\rm T}$ is Hankel
{\em if and only if}\/ $L$ is the specialization
of the generalized Jacobi--Rogers polynomials
$\sfJ = \big( J_{n,k}(\bbeta,\bgamma) \big)_{n,k \ge 0}$
[that is, the traditional ones with $m=1$]\ 
to some values $\bbeta$ and $\bgamma$;
and for a sequence $\ba = (a_n)_{n \ge 0}$,
the Hankel matrix $H_\infty(\ba)$ has an $LDL^{\rm T}$ factorization
with $L$ unit-lower-triangular and $D$ nontrivial diagonal
{\em if and only if}\/ $\ba$ is the specialization
of the Jacobi--Rogers polynomials
$\bJ = \big( J_{n}(\bbeta,\bgamma) \big)_{n \ge 0}$
[again the traditional ones with $m=1$]\
to some values $\bbeta$ and $\bgamma$.
(Here ``nontrivial'' means that no diagonal entry is zero
 or a divisor of zero.)
Since for $m>1$ the matrices of generalized $m$-Jacobi--Rogers polynomials
$\sfJ^{(m)} = \big( J^{(m)}_{n,k}(\bbeta) \big)_{n,k \ge 0}$
and generalized $m$-Stieltjes--Rogers polynomials
$\sfS^{(m)} = \big( S^{(m)}_{n,k}(\balpha) \big)_{n,k \ge 0}$
are {\em not}\/ of this form
--- their production matrices are $(m,1)$-banded rather than tridiagonal,
as we shall see in Section~\ref{subsec.production.SR.TR.JR} ---
such an $LDL^{\rm T}$ factorization cannot exist.
\myendremark

\section[The generalized $\bm{m}$-Jacobi--Rogers polynomials in terms of ordered trees and forests]{Generalized $\bm{m}$-Jacobi--Rogers polynomials in \\ terms of ordered trees and forests}
   \label{sec.mJR.trees}

We would now like to give a combinatorial interpretation
of the $m$-Jacobi--Rogers polynomials $J_n^{(m)}(\bbeta)$
as the generating polynomials for certain classes of ordered trees,
and of the generalized $m$-Jacobi--Rogers polynomials $J_{n,k}^{(m)}(\bbeta)$
as the generating polynomials for certain classes of
ordered forests of ordered trees.

Recall first \cite[pp.~294--295]{Stanley_86}
that an {\em ordered tree}\/ (also called {\em plane tree}\/)
is an (unlabeled) rooted tree in which
the children of each vertex are linearly ordered.
An {\em ordered forest of ordered trees}\/ (also called {\em plane forest}\/)
is a linearly ordered collection of ordered trees.

There is a well-known bijection from the set of
ordered forests of ordered trees
with $n+1$ total vertices and $k+1$ components
onto the set of partial \L{}ukasiewicz paths from $(0,0)$ to $(n,k)$,
which can be described as follows:\footnote{
   See \cite[pp.~30--36]{Stanley_99},
   \cite[Chapter~11]{Lothaire_97},
   \cite[section~6.2]{Pitman_06}
   and \cite[proof of Proposition~3.1]{Hackl_18}.
}

1) Given an ordered forest $F$ of ordered trees with $n+1$ total vertices,
label the vertices with integers $1,\ldots,n+1$
in {\em depth-first-search order}\/
(more precisely, {\em preorder traversal}\/,
 i.e.\ parent first, then children from left to right,
 carried out recursively starting at the root):
for each tree this means to walk counterclockwise around the tree,
starting at the root, and to label the vertices in the order
in which they are first seen;
then do this successively for the trees of the forest
\cite[pp.~33--34]{Stanley_99}.
Note that all the children of a vertex~$i$ have labels $> i$.
See Figure~\ref{fig.depthfirstforest}(a) for an example.

2) Now define a lattice path starting at height $h_0 = k$
and taking steps $s_1,\ldots,s_{n+1}$ with $s_i = \deg(i) - 1$,
where $\deg(i)$ is the number of children of vertex $i$.
For a vertex $j$ ($1 \le j \le n+1$)
belonging to the $r$th tree ($1 \le r \le k+1$),
the height $h_j = k+ \sum_{i=1}^j s_i$
is the number of children of the vertices $\{1,\ldots,j\}$
whose labels are $> j$, plus $k-r$.\footnote{
   {\sc Proof:}
   By induction on $j$.
   For the base case $j=1$, the claim is clear.
   For $j>1$, vertex $j$ is either the child of another node of the
   $r$th tree, or else is the root of the $r$th tree.

   (i) In the first case, by the inductive hypothesis,
   $h_{j-1}$ is the number of children of the vertices $\{1,\ldots,j-1\}$
   whose labels are $\ge j$, plus $k-r$;
   so $h_{j-1} - 1$ is the number of children of the vertices $\{1,\ldots,j-1\}$
   whose labels are $> j$, plus $k-r$.
   Vertex $j$ has $\deg(j)$ children, all of which have labels $> j$;
   so $h_j = h_{j-1} + s_j = h_{j-1} - 1 + \deg(j)$
   is the number of children of the vertices $\{1,\ldots,j\}$
   whose labels are $> j$, plus $k-r$.

   (ii) In the second case, vertex $j-1$ is the last vertex
   of the $(r-1)$st tree,
   and none of the vertices $\{1,\ldots,j-1\}$
   have a child with label $\ge j$ (by definition of the labeling);
   so $h_{j-1} = k-(r-1)$.
   On the other hand, all of the children of vertex $j$
   have labels $> j$,
   and $h_j = h_{j-1} +  \deg(j) - 1 = k-r + \deg(j)$.
 \label{footnote_DFS_lukasiewicz}
}
In particular, $h_j > k-r$ for every vertex of the $r$th tree
other than the last,
and $h_j = k-r$ for the last vertex of the $r$th tree.
So $h_1,\ldots,h_{n-1} \ge 0$, $h_n = 0$ and $h_{n+1} = -1$.
Let us drop the last step:  then the walk $(h_0,\ldots,h_n)$
is a reversed partial \L{}ukasiewicz path from $(0,k)$ to $(n,0)$.
See Figure~\ref{fig.depthfirstforest}(b) for an example.


%
\begin{figure}[t]
\begin{center}
\begin{tikzpicture}[scale=1]
\filldraw[black] (2,2) circle (3pt);
\node[above] at (2,2.1) {\large 1};
\filldraw[black] (1,1) circle (3pt);
\node[left] at (0.9,1) {\large 2};
\filldraw[black] (0.2,0) circle (3pt);
\node[below] at (0.2,-0.1) {\large 3};
\filldraw[black] (1,0) circle (3pt);
\node[below] at (1,-0.1) {\large 4};
\filldraw[black] (1.8,0) circle (3pt);
\node[below] at (1.8,-0.1) {\large 5};
\filldraw[black] (3,1) circle (3pt);
\node[right] at (3.1,1) {\large 6};
\filldraw[black] (4,0) circle (3pt);
\node[below] at (4,-0.1) {\large 7};
\filldraw[black] (7,2) circle (3pt);
\node[above] at (7,2.1) {\large 8};
\filldraw[black] (7,1) circle (3pt);
\node[right] at (7.1,1) {\large 9};
\filldraw[black] (6.5,0) circle (3pt);
\node[below] at (6.5,-0.1) {\large 10};
\filldraw[black] (7.5,0) circle (3pt);
\node[below] at (7.5,-0.1) {\large 11};
\draw[black, thin] (2,2) -- (1,1);
\draw[black, thin] (1,1) -- (0.2,0);
\draw[black, thin] (1,1) -- (1,0);
\draw[black, thin] (1,1) -- (1.8,0);
\draw[black, thin] (2,2) -- (3,1);
\draw[black, thin] (3,1) -- (4,0);
\draw[black, thin] (7,2) -- (7,1);
\draw[black, thin] (7,1) -- (6.5,0);
\draw[black, thin] (7,1) -- (7.5,0);
\node[below] at (2,-1) {\large $\scrt_1$};
\node[below] at (7,-1) {\large $\scrt_2$};
\end{tikzpicture}
\end{center}
\vspace*{0mm}
\begin{center}
\begin{tikzpicture}[scale=0.9]
\filldraw[black] (0,1) circle (3pt);
\filldraw[black] (1,2) circle (3pt);
\filldraw[black] (2,4) circle (3pt);
\filldraw[black] (3,3) circle (3pt);
\filldraw[black] (4,2) circle (3pt);
\filldraw[black] (5,1) circle (3pt);
\filldraw[black] (6,1) circle (3pt);
\filldraw[black] (7,0) circle (3pt);
\filldraw[black] (8,0) circle (3pt);
\filldraw[black] (9,1) circle (3pt);
\filldraw[black] (10,0) circle (3pt);
\draw[black, thin] (0,1) -- (1,2);
\draw[black, thin] (1,2) -- (2,4);
\draw[black, thin] (2,4) -- (3,3);
\draw[black, thin] (3,3) -- (4,2);
\draw[black, thin] (4,2) -- (5,1);
\draw[black, thin] (5,1) -- (6,1);
\draw[black, thin] (6,1) -- (7,0);
\draw[black, thin] (7,0) -- (8,0);
\draw[black, thin] (8,0) -- (9,1);
\draw[black, thin] (9,1) -- (10,0);
\draw[->, gray, thin] (0,0) -- (11,0);
\draw[->, gray, thin] (0,0) -- (0,5);
\node[below] at (0,-0.1) {\large 0};
\node[below] at (1,-0.1) {\large 1};
\node[below] at (2,-0.1) {\large 2};
\node[below] at (3,-0.1) {\large 3};
\node[below] at (4,-0.1) {\large 4};
\node[below] at (5,-0.1) {\large 5};
\node[below] at (6,-0.1) {\large 6};
\node[below] at (7,-0.1) {\large 7};
\node[below] at (8,-0.1) {\large 8};
\node[below] at (9,-0.1) {\large 9};
\node[below] at (10,-0.1) {\large 10};
\node[left] at (-0.1,0) {\large 0};
\node[left] at (-0.1,1) {\large 1};
\node[left] at (-0.1,2) {\large 2};
\node[left] at (-0.1,3) {\large 3};
\node[left] at (-0.1,4) {\large 4};
\draw[gray, thin] (-0.1,2) -- (0.1,2);
\draw[gray, thin] (-0.1,3) -- (0.1,3);
\draw[gray, thin] (-0.1,4) -- (0.1,4);
\draw[gray, thin] (1,0.1) -- (1,-0.1);
\draw[gray, thin] (2,0.1) -- (2,-0.1);
\draw[gray, thin] (3,0.1) -- (3,-0.1);
\draw[gray, thin] (4,0.1) -- (4,-0.1);
\draw[gray, thin] (5,0.1) -- (5,-0.1);
\draw[gray, thin] (6,0.1) -- (6,-0.1);
\draw[gray, thin] (9,0.1) -- (9,-0.1);
\end{tikzpicture}
\end{center}
\caption{
    (a) An ordered forest consisting of 2 ordered trees
        with 11 total vertices, and its depth-first-search labeling.
    (b) The reversed partial \L{}ukasiewicz path from $(0,1)$ to $(10,0)$
        onto which it maps.
}
\label{fig.depthfirstforest}
\end{figure}

3) We now reverse the path by defining $\widehat{h}_i = h_{n-i}$
and translating its starting abscissa back to the origin:
then $(\widehat{h}_0,\ldots,\widehat{h}_n) = (h_n,\ldots,h_0)$
is a partial \L{}ukasiewicz path from $(0,0)$ to $(n,k)$.

4)  Now define the {\em level}\/ of a vertex $j \in [n]$
belonging to the $r$th tree
to be the number of children of the vertices $1,\ldots,j-1$ that are $> j$,
plus $k+1-r$;
it is precisely the height $h_{j-1}$
in the reversed partial \L{}ukasiewicz path,
as explained in footnote~\ref{footnote_DFS_lukasiewicz},
or equivalently the height $\widehat{h}_{n-j+1}$ in the
partial \L{}ukasiewicz path.
Now, the step in the partial \L{}ukasiewicz path
from $\widehat{h}_{n-j} = h_j$ to $\widehat{h}_{n-j+1} = h_{j-1}$
is an $\ell$-fall with $\ell = h_j - h_{j-1} = s_j = \deg(j) - 1$.
Therefore, a vertex $j \in [n]$ at level $L = h_{j-1}$
with $\deg(j)$ children corresponds to an $\ell$-fall
from height $L+\ell$ to height $L$,
where $\ell = \deg(j) - 1$.
(If $\ell = -1$ this is of course a rise,
 and if $\ell = 0$ it is a level step.)
Finally, vertex $n+1$ is always a leaf.
We have therefore proven:

\begin{proposition}
   \label{prop.Jnk.forests}
Let $\bbeta = (\beta_i^{(\ell)})_{i \ge \ell \ge 0}$ be indeterminates.
Then the generalized $\infty$-Jacobi--Rogers polynomial
$J^{(\infty)}_{n,k}(\bbeta)$ is the generating polynomial
for ordered forests of ordered trees
with $n+1$ total vertices and $k+1$ components
in which each vertex at level $L$ with $c$ children
gets weight 1 if it is a leaf ($c=0$)
and weight $\beta_{L+c-1}^{(c-1)}$ otherwise.
\end{proposition}

\noindent
The $m$-Jacobi--Rogers polynomials with $m < \infty$
correspond to the forests in which each vertex has at most $m+1$ children.


\section{Contraction formulae for $\bm{m}$-branched continued fractions}
   \label{sec.contraction}

In this section we define a natural mapping
from $m$-Dyck paths onto $m$-\L{}ukasiewicz paths,
based on grouping steps in nonoverlapping groups of $m+1$ steps;
for $m=1$ this reduces to Viennot's \cite[pp.~V-30--V-32]{Viennot_83}
mapping of Dyck paths onto Motzkin paths by grouping pairs of steps.
We will then use this mapping to derive a formula
for rewriting the $m$-Stieltjes--Rogers polynomials as
$m$-Jacobi--Rogers polynomials with suitably ``contracted'' weights;
for $m=1$ this reduces to the classical formula
for the ``even contraction'' of an S-fraction to a J-fraction
\cite[p.~21]{Wall_48} \cite[p.~V-31]{Viennot_83}.
We will also extend this construction to a subclass
of $m$-Thron--Rogers polynomials,
and to the generalized
$m$-Stieltjes--Rogers and $m$-Thron--Rogers polynomials.

Finally, we will generalize to $m>1$ the classical formula
for the ``odd contraction'' of an S-fraction to a J-fraction
\cite[p.~V-33]{Viennot_83};
and likewise for the other polynomials.

\subsection{Generalization of even contraction}
   \label{subsec.contraction.even}


The combinatorial essence of even contraction
is contained in the following easy lemma:

\begin{lemma}
   \label{prop.evencontraction.combinatorial}
Fix integers $m \ge 1$ and $n \ge 0$,
and consider the mapping of
the path $\omega = (\omega_0,\ldots,\omega_{(m+1)n})$
to the path $\omega' = (\omega'_0,\ldots,\omega'_n)$ defined by
\be
   \omega'_i  \;=\;  {\omega_{(m+1)i}  \over m+1}
 \label{def.omegaprime}
\ee
[here $\omega_i \in \N$ denotes the height of the path at ``time'' $i$,
 and likewise for $\omega'_i$].
This mapping (which we denote $\Pi_n^{(m)}$) is a surjection
from the set $\scrd_{(m+1)n}^{(m)}$
of $m$-Dyck paths from $(0,0)$ to $((m+1)n,0)$
onto the set $\scrl_n^{(m)}$ of $m$-\L{}ukasiewicz paths
from $(0,0)$ to $(n,0)$.
\end{lemma}

\proof
Here $\omega'_i$ is an integer because the path $\omega$
must stay always within the set $V_m$ defined in \reff{def.Vm};
so when the time is a multiple of $m+1$,
the height must also be a multiple of $m+1$.
Moreover, from the hypothesis $\omega_i - \omega_{i-1} \in \{-m,+1\}$
it follows that $-m \le \omega'_i - \omega'_{i-1} \le 1$
and hence that
$\omega'_i - \omega'_{i-1} \in \{-m,-(m-1),\ldots,0,+1\}$;
so $\omega'$ is an $m$-\L{}ukasiewicz path.
Finally, it is easy to see that every value
$\omega'_i - \omega'_{i-1} \in \{-m,-(m-1),\ldots,0,+1\}$ is achievable;
so the mapping is a surjection.
\qed

To derive the formula for contraction of an $m$-Stieltjes--Rogers polynomial
to an $m$-Jacobi--Rogers polynomial,
we now need only compute, for each path $\omega' \in \scrl_n^{(m)}$,
the total weight of paths $\omega \in (\Pi_n^{(m)})^{-1}(\omega')$
when each $m$-fall from height~$i$ gets weight $\alpha_i$.
It is convenient to state this contraction formula in matrix form.
We define some special matrices $M = (m_{ij})_{i,j \ge 0}$, as follows:
\begin{itemize}
   \item $L(s_1,s_2,\ldots)$ is the lower-bidiagonal matrix
       with 1 on the diagonal and $s_1,s_2,\ldots$ on the subdiagonal:
\be
   L(s_1,s_2,\ldots)
   \;=\;
   \begin{bmatrix}
      1  &     &     &     &    \\
      s_1 & 1  &     &     &    \\
          & s_2 & 1  &     &    \\
          &     & s_3 & 1  &    \\
          &     &     & \ddots & \ddots
   \end{bmatrix}
   \;.
 \label{def.L}
\ee
   \item $U^\star(s_1,s_2,\ldots)$ is the upper-bidiagonal matrix
       with 1 on the superdiagonal and $s_1,s_2,\ldots$ on the diagonal:
\be
   U^\star(s_1,s_2,\ldots)
   \;=\;
   \begin{bmatrix}
      s_1 & 1   &     &     &     &    \\
          & s_2 & 1   &     &     &    \\
          &     & s_3 & 1   &     &    \\
          &     &     & s_4 & 1   &    \\
          &     &     &     & \ddots & \ddots
   \end{bmatrix}
   \;.
 \label{def.Ustar}
\ee
   \item $P^{(m)}(\bbeta)$ is the $(m,1)$-banded lower-Hessenberg
       matrix with 1 on the superdiagonal,
       $\beta_0^{(0)},\beta_1^{(0)},\ldots$ on the diagonal,
       $\beta_1^{(1)},\beta_2^{(1)},\ldots$ on the first subdiagonal,
       $\beta_2^{(2)},\beta_3^{(2)},\ldots$ on the second subdiagonal, etc.,
       and in general
\be
   P^{(m)}(\bbeta)_{ij}
   \;=\;
   \begin{cases}
       1                 & \textrm{if $j=i+1$}  \\[1mm]
       \beta_i^{(i-j)}   & \textrm{if $i-m \le j \le i$}  \\[1mm]
       0                 & \textrm{if $j < i-m$}
   \end{cases}
 \label{def.Pm}
\ee
or in other words
\be
   P^{(m)}(\bbeta)
   \;=\;
   \begin{bmatrix}
      \beta_0^{(0)}   & 1               &                 &     &     &    \\
      \beta_1^{(1)}   & \beta_1^{(0)}   & 1               &     &     &    \\
      \beta_2^{(2)}   & \beta_2^{(1)}   & \beta_2^{(0)}   & 1   &     &    \\
      \beta_3^{(3)}   & \beta_3^{(2)}   & \beta_3^{(1)}   & \beta_3^{(0)} & 1   & \\
      \vdots   &\vdots    & \vdots   & \vdots   & \ddots & \ddots
   \end{bmatrix}
   \;.
 \label{def.Pm.0}
\ee
Of course, $P^{(\infty)}(\bbeta)$ is simply the generic
unit-lower-Hessenberg matrix.
\end{itemize}
The contraction formula is then the following:

\begin{proposition}[Contraction formula for $m$-Stieltjes--Rogers polynomials]
   \label{prop.contraction}
\hfill\break
Fix an integer $m \ge 1$,
and let $\balpha = (\alpha_i)_{i \ge m}$ be indeterminates.
Then the $m$-Stieltjes--Rogers polynomial $S^{(m)}_n(\balpha)$
equals the $m$-Jacobi--Rogers polynomial $J^{(m)}_{n}(\bbeta)$
where $\bbeta = (\beta_i^{(\ell)})_{0 \le \ell \le m, \, i \ge \ell}$
are defined by the matrix equality
$P^{(m)}(\bbeta) = P^{(m)\mathrm{S}}(\balpha)$ with
\begin{eqnarray}
   P^{(m)\mathrm{S}}(\balpha)
   & \eqdef &
   L(\alpha_{m+1}, \alpha_{2m+2}, \alpha_{3m+3}, \ldots)
   \:
   L(\alpha_{m+2}, \alpha_{2m+3}, \alpha_{3m+4}, \ldots)
   \:\cdots\:
   \hspace*{1cm}
       \nonumber \\
   & & \qquad
   L(\alpha_{2m}, \alpha_{3m+1}, \alpha_{4m+2}, \ldots)
   \:
   U^\star(\alpha_m, \alpha_{2m+1}, \alpha_{3m+2}, \ldots)
   \;,
   \hspace*{1cm}
 \label{eq.prop.contraction}
\end{eqnarray}
that is, the product of $m$ factors $L$ and one factor $U^\star$.
\end{proposition}

\noindent
For instance, for $m=1$ we have
\be
\Scale[0.95]{
   \begin{bmatrix}
      \alpha_1          & 1                   &     &   &      \\
      \alpha_1 \alpha_2 & \alpha_2 + \alpha_3 & 1   &   &      \\
                        & \alpha_3 \alpha_4   & \alpha_4 + \alpha_5 & 1   &  \\
                        &                     & \ddots & \ddots & \ddots
   \end{bmatrix}
   \;=\;
   \begin{bmatrix}
      1                 &                     &     &   &      \\
      \alpha_2          & 1                   &     &   &      \\
                        & \alpha_4            & 1   &   &  \\
                        &                     & \ddots & \ddots &
   \end{bmatrix}
   \begin{bmatrix}
      \alpha_1          & 1                   &     &   &      \\
                        & \alpha_3            & 1   &   &      \\
                        &                     & \alpha_5 & 1   &  \\
                        &                     &        & \ddots & \ddots
   \end{bmatrix}
}
  \,,
 \label{def.Snl.production.factorization}
\ee
while for $m=2$ we have
\begin{eqnarray}
   & &
   \hspace*{-13mm}
\Scale[0.95]{
   \begin{bmatrix}
      \alpha_2          & 1                   &     &   &      \\
      \alpha_2 \alpha_3 + \alpha_2 \alpha_4 &
                          \alpha_3 + \alpha_4 + \alpha_5 & 1   &   &      \\
      \alpha_2 \alpha_4 \alpha_6           & \alpha_4 \alpha_6 + \alpha_5 \alpha_6 + \alpha_5 \alpha_7 &
                          \alpha_6 + \alpha_7 + \alpha_8 & 1   &  \\
                        &  \ddots             & \ddots & \ddots & \ddots
   \end{bmatrix}
}
               \nonumber \\[2mm]
   & &
   =\;
\Scale[0.95]{
   \begin{bmatrix}
      1                 &                     &     &   &      \\
      \alpha_3          & 1                   &     &   &      \\
                        & \alpha_6            & 1   &   &  \\
                        &                     & \ddots & \ddots &
   \end{bmatrix}
   \begin{bmatrix}
      1                 &                     &     &   &      \\
      \alpha_4          & 1                   &     &   &      \\
                        & \alpha_7            & 1   &   &  \\
                        &                     & \ddots & \ddots &
   \end{bmatrix}
   \begin{bmatrix}
      \alpha_2          & 1                   &     &   &      \\
                        & \alpha_5            & 1   &   &      \\
                        &                     & \alpha_8 & 1   &  \\
                        &                     &        & \ddots & \ddots
   \end{bmatrix}
}
  \,. \qquad
 \label{def.contraction.m=2}
\end{eqnarray}

\proofof{Proposition~\ref{prop.contraction}}
We need to enumerate the $(m+1)$-step walks in $\N$
going from height $(m+1)i$ to height $(m+1)j$, using steps $(1,1)$ and $(1,-m)$,
with weights $\balpha$;
we will show that the right-hand side of \reff{eq.prop.contraction}
does the job.
The first step starts at height $(m+1)i$
and uses the matrix $L(\alpha_{m+1}, \alpha_{2m+2}, \ldots)$:
it is either a rise, using the diagonal ($i \to i$) with weight~1,
or an $m$-fall, using the subdiagonal ($i \to i-1$)
with weight $\alpha_{(m+1)i}$.
We are now at height $(m+1)i + 1$
[where $i$ is the ``new'' value, i.e.\ either the original $i$
or the original $i$ minus 1],
and the next step uses the matrix $L(\alpha_{m+2}, \alpha_{2m+3}, \ldots)$
in an analogous way; and so forth.
Finally, at the $(m+1)$st step, we are at height $(m+1)i + m$
and we use the matrix $U^\star(\alpha_m, \alpha_{2m+1}, \ldots)$:
this step is either a rise, using the superdiagonal ($i \to i+1$)
with weight~1,
or else an $m$-fall, using the diagonal ($i \to i$)
with weight $\alpha_{(m+1)i + m}$.
The resulting value of $j$ is now interpreted as height $(m+1)j$.
\qed

An identical construction works for the generalized polynomials.
For any pair $0 \le k \le n$, the same formula \reff{def.omegaprime}
defines a surjection
$\Pi_{n,k}^{(m)} \colon\: \scrd_{(m+1)n,(m+1)k}^{(m)} \to \scrl_{n,k}^{(m)}$,
and the computation of the weights is identical to the case $k=0$.
We have therefore proven:

\begin{proposition}[Contraction formula for generalized $m$-Stieltjes--Rogers \hbox{polynomials}]
   \label{prop.contraction.generalized}
\hfill\\[-5mm]
Fix an integer $m \ge 1$,
and let $\balpha = (\alpha_i)_{i \ge m}$ be indeterminates.
Then the generalized $m$-Stieltjes--Rogers polynomial $S^{(m)}_{n,k}(\balpha)$
equals the generalized $m$-Jacobi--Rogers polynomial
$J^{(m)}_{n,k}(\bbeta)$
where $\bbeta$ are defined by $P^{(m)}(\bbeta) = P^{(m)\mathrm{S}}(\balpha)$
and \reff{eq.prop.contraction}.
\end{proposition}

Let us now discuss the extension of this construction
to the $m$-Thron--Rogers polynomials.
We would like to use the same formula \reff{def.omegaprime}
also when $\omega$ is an $m$-Schr\"oder path.
There is, however, a danger:  if the $m$-Schr\"oder path
takes an $m$-long step at a time $(m+1)i-1$,
then the height $\omega_{(m+1)i}$ is not defined.
Since an $m$-Schr\"oder path
must stay always within the set $V_m$ defined in \reff{def.Vm},
the time (i.e.\ abscissa) equals $m \bmod m+1$
if and only if the height equals $m \bmod m+1$.
We therefore define $\scrs^{(m)\prime}_{(m+1)n}$
to be the set of $m$-Schr\"oder paths from $(0,0)$ to $((m+1)n,0)$
in which there are no $m$-long steps starting at a height
equal to $m \bmod m+1$.
Then \reff{def.omegaprime} defines a surjection
$\widetilde{\Pi}_n^{(m)} \colon\: \scrs_{(m+1)n}^{(m)\prime} \to \scrl_n^{(m)}$,
which extends $\Pi_n^{(m)}$.

At the level of the $m$-Thron--Rogers polynomials,
forbidding $m$-long steps starting at a height $m \bmod m+1$
corresponds to setting $\delta_i = 0$ whenever $i$ is a multiple of $m+1$.
For the remaining indeterminates we write
$\bdelta' = (\delta_i)_{i \ge m ,\, i \ne 0 \bmod m+1}$,
and we write $T_n^{(m)\prime}(\balpha,\bdelta')$
for $T_n^{(m)}(\balpha,\bdelta)$ specialized to
$\delta_i = 0$ whenever $i$ is a multiple of $m+1$.
We call $T_n^{(m)\prime}(\balpha,\bdelta')$
the \textbfit{restricted $\bm{m}$-Thron--Rogers polynomial}.

To derive the formula for contraction of a
restricted $m$-Thron--Rogers polynomial to an $m$-Jacobi--Rogers polynomial,
we now need to compute, for each path $\omega' \in \scrl_n^{(m)}$,
the total weight of paths
$\omega \in (\widetilde{\Pi}_n^{(m)})^{-1}(\omega')$
when each $m$-fall from height~$i$ gets weight $\alpha_i$
and each $m$-long step starting at height~$i$
($\ne m \bmod m+1$) gets weight $\delta_{i+1}$.
Unfortunately, we have been unable to find a compact matrix formulation
analogous to \reff{eq.prop.contraction}.
We therefore state the result combinatorially:

\begin{proposition}[Contraction formula for restricted $m$-Thron--Rogers polynomials]
   \label{prop.contraction.thron}
\hfill\break
Fix an integer $m \ge 1$,
and let $\balpha = (\alpha_i)_{i \ge m}$
and $\bdelta' = (\delta_i)_{i \ge m ,\, i \ne 0 \bmod m+1}$
be indeterminates.
Define the matrix $P^{(m)\mathrm{RT}}(\balpha,\bdelta')$
by specifying that the matrix element
$P^{(m)\mathrm{RT}}(\balpha,\bdelta')_{ij}$
is the total weight for walks from $(0,(m+1)i)$ to $(m+1,(m+1)j)$
in the graph $\widetilde{G}_m$
--- or equivalently in the graph $\widetilde{G}_m \cap ([0,m+1] \times \N)$ ---
in which each rise gets weight 1,
each $m$-fall from height~$i$ gets weight $\alpha_i$
and each $m$-long step starting at height~$i$ gets weight $\delta_{i+1}$.
[Note that no such walk can have an $m$-long step
 starting at a height equal to $m \bmod m+1$.]
Then the restricted $m$-Thron--Rogers polynomial
$T_n^{(m)\prime}(\balpha,\bdelta')$
equals the $m$-Jacobi--Rogers polynomial $J^{(m)}_{n}(\bbeta)$
where $\bbeta = (\beta_i^{(\ell)})_{0 \le \ell \le m, \, i \ge \ell}$
are defined by the matrix equality
$P^{(m)}(\bbeta) = P^{(m)\mathrm{RT}}(\balpha,\bdelta')$.
\end{proposition}

An identical construction works for the generalized polynomials,
and we have:

\begin{proposition}[Contraction formula for restricted generalized $m$-Thron--Rogers \hbox{polynomials}]
   \label{prop.contraction.thron.generalized}
Fix an integer $m \ge 1$,
and let $\balpha = (\alpha_i)_{i \ge m}$
and $\bdelta' = (\delta_i)_{i \ge m ,\, i \ne 0 \bmod m+1}$
be indeterminates.
Define the matrix $P^{(m)\mathrm{RT}}(\balpha,\bdelta')$
and the coefficients $\bbeta$
as in Proposition~\ref{prop.contraction.thron}.
Then the restricted generalized $m$-Thron--Rogers polynomial
$T_{n,k}^{(m)\prime}(\balpha,\bdelta')$
equals the generalized $m$-Jacobi--Rogers polynomial
$J^{(m)}_{n,k}(\bbeta)$.
\end{proposition}

\subsection{Generalization of odd contraction}

We shall also need a generalization to $m > 1$ of the
odd contraction \cite[p.~V-33]{Viennot_83}.
Let us begin with the Stieltjes case.
An $m$-Dyck path of length $(m+1)n$ with $n \ge 1$
necessarily begins with $m$ rises and ends with an $m$-fall.
We can therefore remove these steps and define
\be
   \omega''_i  \;=\;  {\omega_{(m+1)i+m} - m  \over m+1}
   \;,
 \label{def.omegadoubleprime}
\ee
which are nonnegative integers satisfying
$\omega''_0 = \omega''_{n-1} = 0$
and $-m \le \omega''_i - \omega''_{i-1} \le 1$.
The path $\omega'' = (\omega''_0,\ldots,\omega''_{n-1})$
is therefore an $m$-\L{}ukasiewicz path of length $n-1$.

\begin{proposition}[Odd contraction formula for $m$-Stieltjes--Rogers polynomials]
   \label{prop.contraction_odd}
\hfill\break
Fix an integer $m \ge 1$,
and let $\balpha = (\alpha_i)_{i \ge m}$ be indeterminates.
Then for $n \ge 1$ we have
\be
   S^{(m)}_n(\balpha)  \;=\;  \alpha_m \, J^{(m)}_{n-1}(\bbeta')
\ee
where $\bbeta' = (\beta_i^{\prime(\ell)})_{0 \le \ell \le m, \, i \ge \ell}$
are defined by $P^{(m)}(\bbeta') = P^{(m)\mathrm{S}\prime}(\balpha)$ and
\begin{eqnarray}
   P^{(m)\mathrm{S}\prime}(\balpha)
   & \eqdef &
   U^\star(\alpha_m, \alpha_{2m+1}, \alpha_{3m+2}, \ldots)
   \:
   L(\alpha_{m+1}, \alpha_{2m+2}, \alpha_{3m+3}, \ldots)
   \:\cdots\:
   \hspace*{1cm}
       \nonumber \\
   & & \qquad
   L(\alpha_{m+2}, \alpha_{2m+3}, \alpha_{3m+4}, \ldots)
   \:
   L(\alpha_{2m}, \alpha_{3m+1}, \alpha_{4m+2}, \ldots)
   \;.
   \hspace*{1cm}
 \label{eq.prop.contraction_odd}
\end{eqnarray}
\end{proposition}

\noindent
Here \reff{eq.prop.contraction_odd} is the same as \reff{eq.prop.contraction}
except that the factor $U^\star$ is in first position rather than last.
Thus, for $m=1$ we have
\be
\Scale[0.95]{
   \begin{bmatrix}
      \alpha_1 + \alpha_2 & 1                   &     &   &      \\
      \alpha_2 \alpha_3   & \alpha_3 + \alpha_4 & 1   &   &      \\
                        & \alpha_4 \alpha_5   & \alpha_5 + \alpha_6 & 1   &  \\
                        &                     & \ddots & \ddots & \ddots
   \end{bmatrix}
   \;=\;
   \begin{bmatrix}
      \alpha_1          & 1                   &     &   &      \\
                        & \alpha_3            & 1   &   &      \\
                        &                     & \alpha_5 & 1   &  \\
                        &                     &        & \ddots & \ddots
   \end{bmatrix}
   \begin{bmatrix}
      1                 &                     &     &   &      \\
      \alpha_2          & 1                   &     &   &      \\
                        & \alpha_4            & 1   &   &  \\
                        &                     & \ddots & \ddots &
   \end{bmatrix}
}
  \,,
 \label{def.Snlprime.production.factorization}
\ee
while for $m=2$ we have
\begin{eqnarray}
   & &
   \hspace*{-13mm}
\Scale[0.95]{
   \begin{bmatrix}
      \alpha_2 + \alpha_3 + \alpha_4 & 1                   &     &   &      \\
      \alpha_3 \alpha_5 + \alpha_4 \alpha_5 + \alpha_4 \alpha_6 &
                          \alpha_5 + \alpha_6 + \alpha_7 & 1   &   &      \\
      \alpha_4 \alpha_6 \alpha_8           & \alpha_6 \alpha_8 + \alpha_7 \alpha_8 + \alpha_7 \alpha_9 &
                          \alpha_8 + \alpha_9 + \alpha_{10} & 1   &  \\
                        &  \ddots             & \ddots & \ddots & \ddots
   \end{bmatrix}
}
               \nonumber \\[2mm]
   & &
   =\;
\Scale[0.95]{
   \begin{bmatrix}
      \alpha_2          & 1                   &     &   &      \\
                        & \alpha_5            & 1   &   &      \\
                        &                     & \alpha_8 & 1   &  \\
                        &                     &        & \ddots & \ddots
   \end{bmatrix}
   \begin{bmatrix}
      1                 &                     &     &   &      \\
      \alpha_3          & 1                   &     &   &      \\
                        & \alpha_6            & 1   &   &  \\
                        &                     & \ddots & \ddots &
   \end{bmatrix}
   \begin{bmatrix}
      1                 &                     &     &   &      \\
      \alpha_4          & 1                   &     &   &      \\
                        & \alpha_7            & 1   &   &  \\
                        &                     & \ddots & \ddots &
   \end{bmatrix}
}
  \,. \qquad
 \label{def.contraction_odd.m=2}
\end{eqnarray}

\proofof{Proposition~\ref{prop.contraction_odd}}
The proof is analogous to that of Proposition~\ref{prop.contraction},
but with a few changes.
We need to enumerate the $(m+1)$-step walks in $\N$
going from height $(m+1)i+m$ to height $(m+1)j+m$,
using steps $(1,1)$ and $(1,-m)$,
with weights $\balpha$;
we will show that the right-hand side of \reff{eq.prop.contraction_odd}
does the job.
The first step starts at height $(m+1)i+m$
and uses the matrix $U^\star(\alpha_{m}, \alpha_{2m+1}, \ldots)$:
it is either a rise, using the superdiagonal ($i \to i+1$) with weight~1,
or an $m$-fall, using the diagonal ($i \to i$)
with weight $\alpha_{(m+1)i+m}$.
We are now at height $(m+1)i$
[where $i$ is the ``new'' value, i.e.\ either the original $i$ plus 1
or the original $i$],
and the next step uses the matrix $L(\alpha_{m+1}, \alpha_{2m+2}, \ldots)$
exactly as in the proof of Proposition~\ref{prop.contraction};
and similarly for the remaining steps.
After the $(m+1)$st step, we are at height $(m+1)j + m$.
\qed

\medskip

{\bf Remarks.}
1.  For $m > 1$ the terminology ``odd contraction'' is of course a misnomer;
it should be called ``$m \bmod m+1$ contraction''.
But since this is a mouthful, we prefer the shorter term ``odd contraction''.

2.  There does not seem to exist any analogous contraction formula for
values of $i \bmod m+1$ other than $i=0$ and $i=m$.
For suppose that we extract from an $m$-Dyck path of length $(m+1)n$
the first $i$ steps and the last $m+1-i$ steps, where $1 \le i \le m-1$.
The first $i$ steps must indeed be rises,
and the last step must be an $m$-fall;
but the $m-i$ penultimate steps could be either rises or $m$-falls.
As a result, the walk with the starting and ending steps deleted
must start at height $i$ but could end at various heights
equal to $i \bmod m+1$, not only $i$.

3.  The odd contraction formula (Proposition~\ref{prop.contraction_odd})
will play an important role in our proof of
Theorem~\ref{thm.factorized} below.
\myendremark

\bigskip

There is also a version of odd contraction for $m$-Thron--Rogers polynomials
restricted to have $\delta_i = 0$ whenever $i = m \bmod m+1$
(so as to forbid $m$-long steps starting at a height equal to $m-1 \bmod m+1$).
We leave the details to the reader.

\section{Production matrices}   \label{sec.production}

The method of production matrices \cite{Deutsch_05,Deutsch_09}
has become in recent years an important tool in enumerative combinatorics.
In the special case of a tridiagonal production matrix,
this construction goes back to Stieltjes' \cite{Stieltjes_1889,Stieltjes_1894}
work on continued fractions:
the production matrix of a classical S-fraction or J-fraction is tridiagonal.
In~the present paper, by contrast,
we shall need production matrices that are lower-Hessenberg
(i.e.\ vanish above the first superdiagonal)
but are not in general tridiagonal.
We therefore begin by reviewing briefly
the basic theory of production matrices.
We then exhibit the production matrices associated to
the lower-triangular arrays of generalized $m$-Stieltjes--Rogers,
$m$-Thron--Rogers and $m$-Jacobi--Rogers polynomials.
The important connection of production matrices with total positivity
will be treated in the next section.

\subsection{General theory of production matrices}

Let $P = (p_{ij})_{i,j \ge 0}$ be an infinite matrix
with entries in a commutative ring $R$.
In~order that powers of $P$ be well-defined,
we shall assume that $P$ is either row-finite
(i.e.\ has only finitely many nonzero entries in each row)
or column-finite.

Let us now define an infinite matrix $A = (a_{nk})_{n,k \ge 0}$ by
\be
   a_{nk}  \;=\;  (P^n)_{0k}
 \label{def.iteration}
\ee
(in particular, $a_{0k} = \delta_{0k}$).
Writing out the matrix multiplications explicitly, we have
\be
   a_{nk}
   \;=\;
   \sum_{i_1,\ldots,i_{n-1}}
      p_{0 i_1} \, p_{i_1 i_2} \, p_{i_2 i_3} \,\cdots\,
        p_{i_{n-2} i_{n-1}} \, p_{i_{n-1} k}
   \;,
 \label{def.iteration.walk}
\ee
so that $a_{nk}$ is the total weight for all $n$-step walks in $\N$
from $i_0 = 0$ to $i_n = k$, in~which the weight of a walk is the
product of the weights of its steps, and a step from $i$ to $j$
gets a weight $p_{ij}$.
Yet another equivalent formulation is to define the entries $a_{nk}$
by the recurrence
\be
   a_{nk}  \;=\;  \sum_{i=0}^\infty a_{n-1,i} \, p_{ik}
   \qquad\hbox{for $n \ge 1$}
 \label{def.iteration.bis}
\ee
with the initial condition $a_{0k} = \delta_{0k}$.

We call $P$ the \textbfit{production matrix}
and $A$ the \textbfit{output matrix},
and we write $A = \scro(P)$.
Note that if $P$ is row-finite, then so is $\scro(P)$;
if $P$ is lower-Hessenberg, then $\scro(P)$ is lower-triangular;
if $P$ is lower-Hessenberg with invertible superdiagonal entries,
then $\scro(P)$ is lower-triangular with invertible diagonal entries;
and if $P$ is unit-lower-Hessenberg
(i.e.\ lower-Hessenberg with entries 1 on the superdiagonal),
then $\scro(P)$ is unit-lower-triangular.
In all the applications in this paper, $P$ will be unit-lower-Hessenberg.

The matrix $P$ can also be interpreted as the adjacency matrix
for a weighted directed graph on the vertex set $\N$
(where the edge $ij$ is omitted whenever $p_{ij}  = 0$).
Then $P$ is row-finite (resp.\ column-finite)
if and only if every vertex has finite out-degree (resp.\ finite in-degree).

This iteration process can be given a compact matrix formulation.
Let $\Delta = (\delta_{i+1,j})_{i,j \ge 0}$
be the matrix with 1 on the superdiagonal and 0 elsewhere.
Then for any matrix $M$ with rows indexed by $\N$,
the product $\Delta M$ is simply $M$ with its zeroth row removed
and all other rows shifted upwards.
(Some authors use the notation $\overline{M} \eqdef \Delta M$.)
The recurrence \reff{def.iteration.bis} can then be written as
\be
   \Delta \, \scro(P)  \;=\;  \scro(P) \, P
   \;.
 \label{def.iteration.bis.matrixform}
\ee
It follows that if $A$ is a row-finite matrix
that has a row-finite inverse $A^{-1}$
and has first row $a_{0k} = \delta_{0k}$,
then $P = A^{-1} \Delta A$ is the unique matrix such that $A = \scro(P)$.
This holds, in particular, if $A$ is lower-triangular with
invertible diagonal entries and $a_{00} = 1$;
then $A^{-1}$ is lower-triangular
and $P = A^{-1} \Delta A$ is lower-Hessenberg.
And if $A$ is unit-lower-triangular,
then $P = A^{-1} \Delta A$ is unit-lower-Hessenberg.

Let us record, for future use, the following easy fact:

\begin{lemma}[Production matrix of a product]
   \label{lemma.production.AB}
Let $A$ and $B$ be infinite lower-triangular matrices
(with entries in a commutative ring $R$)
with invertible diagonal entries.
Then $A$ has production matrix $P_A = A^{-1} \Delta A$,
and $AB$ has production matrix
$P_{AB} = (AB)^{-1} \Delta (AB) = B^{-1} P_A B$.
\end{lemma}

\subsection[Production matrices for generalized $m$-Stieltjes--Rogers, $m$-Thron--Rogers and $m$-Jacobi--Rogers polynomials]{Production matrices for generalized $\bm{m}$-Stieltjes--Rogers, $\bm{m}$-Thron--Rogers and $\bm{m}$-Jacobi--Rogers polynomials}
   \label{subsec.production.SR.TR.JR}

It is now almost trivial to state the production matrices associated to
the lower-triangular arrays of generalized $m$-Stieltjes--Rogers,
$m$-Thron--Rogers and $m$-Jacobi--Rogers polynomials:

\begin{proposition}[Production matrices for generalized $m$-J--R, $m$-S--R and $m$-T--R]
   \label{prop.prod.Sm.Jm.Tm}
\hfill\\[-6mm]
\begin{itemize}
   \item[(a)]  For the unit-lower-triangular matrix
       $\sfJ^{(m)} = \big( J^{(m)}_{n,k}(\bbeta) \big)_{n,k \ge 0}$
       of generalized $m$-Jacobi--Rogers polynomials,
       the production matrix is $P^{(m)}(\bbeta)$
       defined in \reff{def.Pm.0}.
   \item[(b)]  For the unit-lower-triangular matrix
       $\sfS^{(m)} = \big( S^{(m)}_{n,k}(\balpha) \big)_{n,k \ge 0}$
       of generalized $m$-Stieltjes--Rogers polynomials,
       the production matrix is $P^{(m)\mathrm{S}}(\balpha)$
       defined in \reff{eq.prop.contraction}.
   \item[(c)]  For the unit-lower-triangular matrix
       $\sfT^{(m)\prime} = \big( T^{(m)\prime}_{n,k}(\balpha,\bdelta') \big)_{n,k \ge 0}$
       of restricted generalized $m$-Thron--Rogers polynomials,
       the production matrix is $P^{(m)\mathrm{RT}}(\balpha,\bdelta')$
       defined in Proposition~\ref{prop.contraction.thron}.
\end{itemize}
\end{proposition}

\proof
(a) is immediate from the definitions.
(b) follows from (a) together with
Proposition~\ref{prop.contraction.generalized},
while (c) follows from (a) together with
Proposition~\ref{prop.contraction.thron.generalized}.
\qed

We can now give the promised algebraic proof of
Proposition~\ref{prop.reduction}
for the special case in which the set $I$ is periodic of period~$m+1$:

\algebraicproofof{Proposition~\ref{prop.reduction} (special case)}
Consider the production matrix \reff{eq.prop.contraction}
with $m$ replaced by $m'$ and $\balpha$ by $\balpha'$.
When $I$ is periodic of period~$m+1$,
exactly $m'-m$ of the factors $L(\cdots)$ become the identity matrix,
so that we get the production matrix for an $m$-branched S-fraction
with weights $\balpha$.
(Note that since $m' \notin I$, the factor $U^\star(\cdots)$
 is unaltered.)
\qed

\section{Total positivity}
   \label{sec.totalpos}

In this section --- which is the theoretical heart of the paper ---
we prove the coefficientwise total positivity of the Hankel matrices
associated to the $m$-Stieltjes--Rogers and $m$-Thron--Rogers polynomials,
and of the lower-triangular matrices of generalized
$m$-Stieltjes--Rogers and $m$-Thron--Rogers polynomials.
For each Stieltjes result we will give two proofs:
a combinatorial proof using the Lindstr\"om--Gessel--Viennot lemma,
and algebraic proof using the theory of production matrices.
For the Thron results we will have only a combinatorial proof,
although for the restricted Thron--Rogers polynomials
defined in Section~\ref{subsec.contraction.even}
we will have an almost-algebraic proof.
For the case $m=1$, these proofs can be found already
in \cite{Sokal_totalpos}.

We begin (Section~\ref{subsec.totalpos.prelim})
by recalling some general facts
about partially ordered commutative rings and total positivity.
Next (Section~\ref{subsec.totalpos.prodmat})
we present the fundamental theorems
concerning production matrices and total positivity \cite{Sokal_totalpos};
these will be a key tool in the remainder of the paper.
Then we state our main results (Section~\ref{subsec.totalpos.statement})
and give the promised combinatorial and algebraic proofs
(Sections~\ref{subsec.totalpos.combinatorial} and
\ref{subsec.totalpos.algebraic}).

\subsection{Preliminaries concerning partially ordered commutative rings and total positivity}
   \label{subsec.totalpos.prelim}

In this paper all rings will be assumed to have an identity element 1
and to be nontrivial ($1 \ne 0$).

A \textbfit{partially ordered commutative ring} is a pair $(R,\scrp)$ where
$R$ is a commutative ring and $\scrp$ is a subset of $R$ satisfying
\begin{itemize}
   \item[(a)]  $0,1 \in \scrp$.
   \item[(b)]  If $a,b \in \scrp$, then $a+b \in \scrp$ and $ab \in \scrp$.
   \item[(c)]  $\scrp \cap (-\scrp) = \{0\}$.
\end{itemize}
We call $\scrp$ the {\em nonnegative elements}\/ of $R$,
and we define a partial order on $R$ (compatible with the ring structure)
by writing $a \le b$ as a synonym for $b-a \in \scrp$.
Please note that, unlike the practice in real algebraic geometry
\cite{Brumfiel_79,Lam_84,Prestel_01,Marshall_08},
we do {\em not}\/ assume here that squares are nonnegative;
indeed, this property fails completely for our prototypical example,
the ring of polynomials with the coefficientwise order,
since $(1-x)^2 = 1-2x+x^2 \not\myge 0$.

Now let $(R,\scrp)$ be a partially ordered commutative ring
and let $\bfx = \{x_i\}_{i \in I}$ be a collection of indeterminates.
In the polynomial ring $R[\bfx]$ and the formal-power-series ring $R[[\bfx]]$,
let $\scrp[\bfx]$ and $\scrp[[\bfx]]$ be the subsets
consisting of polynomials (resp.\ series) with nonnegative coefficients.
Then $(R[\bfx],\scrp[\bfx])$ and $(R[[\bfx]],\scrp[[\bfx]])$
are partially ordered commutative rings;
we refer to this as the {\em coefficientwise order}\/
on $R[\bfx]$ and $R[[\bfx]]$.

A (finite or infinite) matrix with entries in a
partially ordered commutative ring
is called \textbfit{totally positive} (TP) if all its minors are nonnegative;
it is called \textbfit{totally positive of order~$\bm{r}$} (TP${}_r$)
if all its minors of size $\le r$ are nonnegative.
It follows immediately from the Cauchy--Binet formula that
the product of two TP (resp.\ TP${}_r$) matrices is TP (resp.\ TP${}_r$).

We say that a sequence $\ba = (a_n)_{n \ge 0}$
with entries in a partially ordered commutative ring
is \textbfit{Hankel-totally positive} 
(resp.\ \textbfit{Hankel-totally positive of order~$\bm{r}$})
if its associated infinite Hankel matrix
$H_\infty(\ba) = (a_{i+j})_{i,j \ge 0}$
is TP (resp.\ TP${}_r$).

We will need a few easy facts about the total positivity of special matrices.

\begin{lemma}[Bidiagonal matrices]
  \label{lemma.bidiagonal}
Let $A$ be a matrix with entries in a partially ordered commutative ring,
with the property that all its nonzero entries belong to two consecutive
diagonals.
Then $A$ is totally positive if and only if all its entries are nonnegative.
\end{lemma}

\proof
The nonnegativity of the entries (i.e.\ TP${}_1$)
is obviously a necessary condition for TP.
Conversely, for a matrix of this type it is easy to see that
every nonzero minor is simply a product of some entries.
\qed

\begin{lemma}[Toeplitz matrix of powers]
   \label{lemma.toeplitz.power}
Let $R$ be a partially ordered commutative ring, let $\xi \in R$,
and let $T_\xi$ be the lower-triangular Toeplitz matrix of powers of $\xi$
[cf.\ \reff{def.Txi}].
Then every minor of $T_\xi$ is either zero or else a power of $\xi$.
Hence $T_\xi$ is TP $\iff$ $T_\xi$ is TP${}_1$ $\iff$ $\xi \ge 0$.

In particular, if $\xi$ is an indeterminate, then $T_\xi$
is totally positive in the ring $\Z[\xi]$ equipped with the
coefficientwise order.
\end{lemma}

\proof
Consider a submatrix $A = (T_\xi)_{IJ}$
with rows $I = \{i_1 < \ldots < i_k \}$
and columns $J = \{j_1 < \ldots < j_k \}$.
We will prove the claim by induction on $k$.
It is trivial if $k=0$ or 1.
If $A_{12} = A_{22} = 0$,
then $A_{1s} = A_{2s} = 0$ for all $s \ge 2$ by definition of $T_\xi$,
and $\det A = 0$.
If $A_{12}$ and $A_{22}$ are both nonzero,
then the first column of $A$ is $\xi^{j_2 - j_1}$ times the second column,
and again $\det A = 0$.
Finally, if $A_{12} = 0$ and $A_{22} \ne 0$
(by definition of $T_\xi$ this is the only other possibility),
then $A_{1s} = 0$ for all $s \ge 2$;
we then replace the first column of $A$ by
the first column minus $\xi^{j_2 - j_1}$ times the second column,
so that the new first column has $\xi^{i_1-j_1}$ in its first entry
(or zero if $i_1 < j_1$) and zeroes elsewhere.
Then $\det A$ equals $\xi^{i_1-j_1}$ times
the determinant of its last $k-1$ rows and columns (or zero if $i_1 < j_1$),
so the claim follows from the inductive hypothesis.
\qed

\subsection{Production matrices and total positivity}
   \label{subsec.totalpos.prodmat}

Let $P = (p_{ij})_{i,j \ge 0}$ be a matrix with entries in a
partially ordered commutative ring $R$.
We will use $P$ as a production matrix;
let $A = \scro(P)$ be the corresponding output matrix.
As before, we assume that $P$ is either row-finite or column-finite.

When $P$ is totally positive, it turns out \cite{Sokal_totalpos}
that the output matrix $\scro(P)$ has {\em two}\/ total-positivity properties:
firstly, it is totally positive;
and secondly, its zeroth column is Hankel-totally positive.
Since \cite{Sokal_totalpos} is not yet publicly available,
we shall present briefly here (with proof) the main results
that will be needed in the sequel.

The fundamental fact that drives the whole theory is the following:

\begin{proposition}[Minors of the output matrix]
   \label{prop.iteration.homo}
Every $k \times k$ minor of the output matrix $A = \scro(P)$
can be written as a sum of products of minors of size $\le k$
of the production matrix $P$.
\end{proposition}

In this proposition the matrix elements $\bfp = \{p_{ij}\}_{i,j \ge 0}$
should be interpreted in the first instance as indeterminates:
for instance, we can fix a row-finite or column-finite set
$S \subseteq \N \times \N$
and define the matrix $P^S = (p^S_{ij})_{i,j \in \N}$ with entries
\be
   p^S_{ij}
   \;=\;
   \begin{cases}
       p_{ij}  & \textrm{if $(i,j) \in S$} \\[1mm]
       0       & \textrm{if $(i,j) \notin S$}
   \end{cases}
\ee
Then the entries (and hence also the minors) of both $P$ and $A$
belong to the polynomial ring $\Z[\bfp]$,
and the assertion of Proposition~\ref{prop.iteration.homo} makes sense.
Of course, we can subsequently specialize the indeterminates $\bfp$
to values in any commutative ring $R$.

\proofof{Proposition~\ref{prop.iteration.homo}}
Consider any minor of $A$ involving only the rows 0 through $N$.
We will prove the assertion of the Proposition by induction on $N$.
The statement is obvious for $N=0$.
For $N \ge 1$, let $A_N$ be the matrix consisting of rows
0 through $N-1$ of $A$, and let $A'_N$ be the matrix consisting of rows
1 through $N$ of $A$.  Then we have
\be
   A'_N  \;=\;  A_N P
   \;.
 \label{eq.proof.prop.iteration}
\ee
If the minor in question does not involve row 0,
then obviously it involves only rows 1 through $N$.
If the minor in question does involve row 0,
then it is either zero (in case it does not involve column~0)
or else equal to a minor of $A$ (of one size smaller)
that involves only rows 1 through $N$
(since $a_{0k} = \delta_{0k}$).
Either way it is a minor of $A'_N$;
but by \reff{eq.proof.prop.iteration} and the Cauchy--Binet formula,
every minor of $A'_N$ is a sum of products of minors (of the same size)
of $A_N$ and $P$.
This completes the inductive step.
\qed

If we now specialize the indeterminates $\bfp$
to values in some partially ordered commutative ring $R$,
we can immediately conclude:

\begin{theorem}[Total positivity of the output matrix]
   \label{thm.iteration.homo}
Let $P$ be an infinite matrix that is either row-finite or column-finite,
with entries in a partially ordered commutative ring $R$.
If $P$ is totally positive of order~$r$, then so is $A = \scro(P)$.
\end{theorem}

\medskip

{\bf Remarks.}
1.  In the case $R = \R$, Theorem~\ref{thm.iteration.homo}
is due to Karlin \cite[pp.~132--134]{Karlin_68};
see also \cite[Theorem~1.11]{Pinkus_10}.
Karlin's proof is different from ours.

2.  Our quick inductive proof of Proposition~\ref{prop.iteration.homo}
follows an idea of one of us (Zhu) \cite[proof of Theorem~2.1]{Zhu_13},
which was in turn inspired in part by Aigner \cite[pp.~45--46]{Aigner_99}.
The same idea recurs in recent work of several authors
\cite[Theorem~2.1]{Zhu_14}
\cite[Theorem~2.1(i)]{Chen_15a}
\cite[Theorem~2.3(i)]{Chen_15b}
\cite[Theorem~2.1]{Liang_16}.
However, all of these results concerned only special cases:
\cite{Aigner_99,Zhu_13,Chen_15b,Liang_16}
treated the case in which the production matrix $P$ is tridiagonal;
\cite{Zhu_14} treated a (special) case in which $P$ is upper bidiagonal;
\cite{Chen_15a} treated the case in which
$P$ is the production matrix of a Riordan array.
But the argument is in fact completely general, as we have just seen;
there is no need to assume any special form for the matrix $P$.
\myendremark

\bigskip

Now define 
$\scroo_0(P)$ to be the zeroth-column sequence of $\scro(P)$, i.e.
\be
   \scroo_0(P)_n  \;\eqdef\;  \scro(P)_{n0}  \;\eqdef\;  (P^n)_{00}
   \;.
 \label{def.scroo0}
\ee
Then the Hankel matrix of $\scroo_0(P)$ has matrix elements
\begin{eqnarray}
   & &
   \!\!\!\!\!\!\!
   H_\infty(\scroo_0(P))_{nn'}
   \;=\;
   \scroo_0(P)_{n+n'}
   \;=\;
   (P^{n+n'})_{00}
   \;=\;
   \sum_{k=0}^\infty (P^n)_{0k} \, (P^{n'})_{k0}
   \;=\;
          \nonumber \\
   & &
   \sum_{k=0}^\infty (P^n)_{0k} \, ((P^{\rm T})^{n'})_{0k}
   \;=\;
   \sum_{k=0}^\infty \scro(P)_{nk} \, \scro(P^{\rm T})_{n'k}
   \;=\;
   \big[ \scro(P) \, {\scro(P^{\rm T})}^{\rm T} \big]_{nn'}
   \;.
   \qquad
\end{eqnarray}
(Note that the sum over $k$ has only finitely many nonzero terms:
 if $P$ is row-finite, then there are finitely many nonzero $(P^n)_{0k}$,
 while if $P$ is column-finite,
 there are finitely many nonzero $(P^{n'})_{k0}$.)
We have therefore proven:

\begin{lemma}[Identity for Hankel matrix of the zeroth column]
   \label{lemma.hankel.karlin}
Let $P$ be a row-finite or column-finite matrix
with entries in a commutative ring $R$.
Then
\be
   H_\infty(\scroo_0(P))
   \;=\;
   \scro(P) \, {\scro(P^{\rm T})}^{\rm T}
   \;.
\ee
\end{lemma}

{\bf Remark.}
If $P$ is row-finite, then $\scro(P)$ is row-finite;
$\scro(P^{\rm T})$ need not be row- or column-finite,
but the product $\scro(P) \, {\scro(P^{\rm T})}^{\rm T}$
is anyway well-defined.
\myendremark

\medskip

Combining Proposition~\ref{prop.iteration.homo}
with Lemma~\ref{lemma.hankel.karlin} and the Cauchy--Binet formula,
we obtain:

\begin{corollary}[Hankel minors of the zeroth column]
   \label{cor.iteration2}
Every $k \times k$ minor of the infinite Hankel matrix
$H_\infty(\scroo_0(P)) = ((P^{n+n'})_{00})_{n,n' \ge 0}$
can be written as a sum of products
of the minors of size $\le k$ of the production matrix $P$.
\end{corollary}

And specializing the indeterminates $\bfp$
to nonnegative elements in a partially ordered commutative ring,
in such a way that $P$ is row-finite or column-finite,
we deduce:

\begin{theorem}[Hankel-total positivity of the zeroth column]
   \label{thm.iteration2bis}
Let $P = (p_{ij})_{i,j \ge 0}$ be an infinite row-finite or column-finite
matrix with entries in a partially ordered commutative ring $R$,
and define the infinite Hankel matrix
$H_\infty(\scroo_0(P)) = ((P^{n+n'})_{00})_{n,n' \ge 0}$.
If $P$ is totally positive of order~$r$, then so is $H_\infty(\scroo_0(P))$.
\end{theorem}

\subsection[Statement of results for $m$-Stieltjes--Rogers, $m$-Thron--Rogers and $m$-Jacobi--Rogers polynomials]{Statement of results for $\bm{m}$-Stieltjes--Rogers, $\bm{m}$-Thron--Rogers and $\bm{m}$-Jacobi--Rogers polynomials}
   \label{subsec.totalpos.statement}

Our first main result concerns the Hankel-total positivity of the
$m$-Stieltjes--Rogers and $m$-Thron--Rogers polynomials:

\begin{theorem}[Hankel-total positivity for $m$-Stieltjes--Rogers polynomials]
   \label{thm.Stype.minors}
For each integer $m \ge 1$,
the sequence $\bS^{(m)} = ( S^{(m)}_n(\balpha) )_{n \ge 0}$
of $m$-Stieltjes--Rogers polynomials
is a Hankel-totally positive sequence in the polynomial ring $\Z[\balpha]$
equipped with the coefficientwise partial order.
\end{theorem}

\begin{theorem}[Hankel-total positivity for $m$-Thron--Rogers polynomials]
   \label{thm.Ttype.minors}
For each integer $m \ge 1$,
the sequence $\bT^{(m)} = ( T^{(m)}_n(\balpha,\bdelta) )_{n \ge 0}$
of $m$-Thron--Rogers polynomials
is a Hankel-totally positive sequence in the polynomial ring
$\Z[\balpha,\bdelta]$ equipped with the coefficientwise partial order.
\end{theorem}

Of course, it suffices to prove Theorem~\ref{thm.Ttype.minors},
since Theorem~\ref{thm.Stype.minors} is simply the specialization
to $\bdelta = \bzero$.
We will prove Theorem~\ref{thm.Ttype.minors}
in Section~\ref{subsec.totalpos.combinatorial}
by a combinatorial method using the Lindstr\"om--Gessel--Viennot lemma.
We will also give a second proof of Theorem~\ref{thm.Stype.minors}
in Section~\ref{subsec.totalpos.algebraic},
based on the theory of production matrices.

Note also that,
as a corollary of Theorems~\ref{thm.Stype.minors} and \ref{thm.Ttype.minors},
we can substitute for $\balpha$ and $\bdelta$
nonnegative elements of any partially ordered commutative ring $R$,
and the resulting sequence $\bS^{(m)}$ or $\bT^{(m)}$
will be Hankel-totally positive in $R$.
This will, in fact, be our principal method for constructing examples
of Hankel-totally positive sequences.

Our second main result concerns the total positivity of the
lower-triangular matrices of generalized
$m$-Stieltjes--Rogers and $m$-Thron--Rogers polynomials:

\begin{theorem}[Total positivity for generalized $m$-Stieltjes--Rogers polynomials]
   \label{thm.Stype.generalized.minors}
\hfill\break
The lower-triangular matrix
$\sfS^{(m)} = \big( S^{(m)}_{n,k}(\balpha) \big)_{n,k \ge 0}$
of generalized $m$-Stieltjes--Rogers polynomials
is totally positive in the polynomial ring $\Z[\balpha]$
equipped with the coefficientwise partial order.
\end{theorem}

\begin{theorem}[Total positivity for generalized $m$-Thron--Rogers polynomials]
   \label{thm.Ttype.generalized.minors}
\hfill\break
The lower-triangular matrix
$\sfT^{(m)} = \big( T^{(m)}_{n,k}(\balpha,\bdelta) \big)_{n,k \ge 0}$
of generalized $m$-Thron--Rogers polynomials
is totally positive in the polynomial ring $\Z[\balpha,\bdelta]$
equipped with the coefficientwise partial order.
\end{theorem}

Once again, it suffices to prove Theorem~\ref{thm.Ttype.generalized.minors},
since Theorem~\ref{thm.Stype.generalized.minors} is the specialization
to $\bdelta = \bzero$.
We will prove Theorem~\ref{thm.Ttype.generalized.minors}
in Section~\ref{subsec.totalpos.combinatorial}
using the Lindstr\"om--Gessel--Viennot lemma,
and will give a second proof of Theorem~\ref{thm.Stype.generalized.minors}
in Section~\ref{subsec.totalpos.algebraic}
using the theory of production matrices.

The results on Hankel-total positivity in
Theorems~\ref{thm.Stype.minors} and \ref{thm.Ttype.minors}
can be extended as follows.
Recall from Section~\ref{sec2} that
the coefficient of $t^n$ in the product $f_0 f_1 \cdots f_\ell$
is the generating polynomial for partial $m$-Dyck or $m$-Schr\"oder paths
from $(0,0)$ to $((m+1)n+\ell,\ell)$ [with the usual weights],
which we denote by
$S_{n|\ell}^{(m)}(\balpha)$ or $T_{n|\ell}^{(m)}(\balpha,\bdelta)$,
respectively.
We then have Hankel-total positivity
in the Stieltjes case for $0 \le \ell \le m$
(but {\em not}\/ for $\ell > m$,
 see \reff{eq.ellgtm.counterexample.1}/\reff{eq.ellgtm.counterexample.2}
 below),
and in the Thron case for $0 \le \ell \le m-1$:

\begin{theorem}[Hankel-total positivity for partial $m$-Stieltjes--Rogers polynomials]
   \label{thm.Stype.minors.extended}
For each pair of integers $m \ge 1$ and $0 \le \ell \le m$,
the sequence
$\bS^{(m)}_{|\ell} = ( S^{(m)}_{n|\ell}(\balpha) )_{n \ge 0}$
of partial $m$-Stieltjes--Rogers polynomials
is a Hankel-totally positive sequence in the polynomial ring
$\Z[\balpha]$ equipped with the coefficientwise partial order.
\end{theorem}

\begin{theorem}[Hankel-total positivity for partial $m$-Thron--Rogers polynomials]
   \label{thm.Ttype.minors.extended}
For each pair of integers $m \ge 1$ and $0 \le \ell \le m-1$,
the sequence
$\bT^{(m)}_{|\ell} = ( T^{(m)}_{n|\ell}(\balpha,\bdelta) )_{n \ge 0}$
of partial $m$-Thron--Rogers polynomials
is a Hankel-totally positive sequence in the polynomial ring
$\Z[\balpha,\bdelta]$ equipped with the coefficientwise partial order.
\end{theorem}

We remark that for the Stieltjes case with $\ell = m$,
the recurrence \reff{eq.mSRfk.1} [or a simple direct combinatorial argument]
implies $S_{n|m}^{(m)}(\balpha) = S_{n+1}^{(m)}(\balpha)/\alpha_m$,
so this case of Theorem~\ref{thm.Stype.minors.extended}
is an immediate consequence of Theorem~\ref{thm.Stype.minors}.
But Theorem~\ref{thm.Stype.minors.extended} with $1 \le \ell \le m-1$,
and Theorem~\ref{thm.Ttype.minors.extended} with $1 \le \ell \le m$,
are nontrivial extensions of Theorems~\ref{thm.Stype.minors} and
\ref{thm.Ttype.minors}.

We will prove Theorem~\ref{thm.Ttype.minors.extended}
(and hence also Theorem~\ref{thm.Stype.minors.extended})
in Section~\ref{subsec.totalpos.combinatorial}
using the Lindstr\"om--Gessel--Viennot lemma.

\bigskip

Let us consider, finally, the $m$-Jacobi--Rogers polynomials.
Unlike the Stieltjes and Thron cases,
these polynomials are {\em not}\/ Hankel-totally positive
coefficientwise in $\bbeta$, or even when the coefficients
$\bbeta$ are given arbitrary nonnegative real values;
indeed, they are not even log-convex.
For instance, for the classical Jacobi--Rogers polynomials ($m=1$) we have
(using the notation $\beta_i = \beta_i^{(1)}$ and $\gamma_i = \beta_i^{(0)}$)
\be
   J_1 J_3 - J_2^2  \;=\;  \beta_1 \gamma_0 \gamma_1 - \beta_1^2
   \;,
\ee
which is negative whenever $\beta_1 > \gamma_0 \gamma_1 \ge 0$.
Rather, the Hankel-total positivity of the $m$-Jacobi--Rogers polynomials
--- and also the total positivity of the triangular array of
generalized $m$-Jacobi--Rogers polynomials ---
holds subject to suitable {\em inequalities}\/
on the coefficients $\bbeta$
(assumed to take values in some partially ordered commutative ring).
Indeed, it follows immediately from Theorems~\ref{thm.iteration.homo}
and \ref{thm.iteration2bis} that a {\em sufficient}\/
(though perhaps not necessary) condition for these total positivities to hold
is the total positivity of the production matrix $P^{(m)}(\bbeta)$
defined in \reff{def.Pm}/\reff{def.Pm.0}:

\begin{theorem}[Hankel-total positivity for $m$-Jacobi--Rogers polynomials]
   \label{thm.Jtype.minors}
Fix an integer $m \ge 1$,
and let $\bbeta = (\beta_i^{(\ell)})_{0 \le \ell \le m, \, i \ge \ell}$
be elements in a partially ordered commutative ring.
If the production matrix $P^{(m)}(\bbeta)$
is totally positive of order~$r$,
then the sequence $\bJ^{(m)} = ( J^{(m)}_n(\bbeta) )_{n \ge 0}$
of $m$-Jacobi--Rogers polynomials is Hankel-totally positive of order~$r$.
\end{theorem}

\begin{theorem}[Total positivity for generalized $m$-Jacobi--Rogers polynomials]
   \label{thm.Jtype.generalized.minors}
\hfill\break
Fix an integer $m \ge 1$,
and let $\bbeta = (\beta_i^{(\ell)})_{0 \le \ell \le m, \, i \ge \ell}$
be elements in a partially ordered commutative ring.
If the production matrix $P^{(m)}(\bbeta)$
is totally positive of order~$r$,
then the lower-triangular matrix
$\sfJ^{(m)} = \big( J^{(m)}_{n,k}(\bbeta) \big)_{n,k \ge 0}$
of generalized $m$-Jacobi--Rogers polynomials
is also totally positive of order~$r$.
\end{theorem}

\subsection{Combinatorial proofs using the Lindstr\"om--Gessel--Viennot lemma}
   \label{subsec.totalpos.combinatorial}

We begin by recalling the Lindstr\"om--Gessel--Viennot lemma.
The framework is as follows:
Let $\Gamma = (V, \vec{E})$ be a directed graph;
here the vertex set $V$ and the edge set $\vec{E}$ need not be finite,
but we shall assume for notational simplicity that for each
ordered pair $(i,j) \in V \times V$ there is at most one edge from $i$ to $j$.
Now let ${\bf w} = (w_{ij})_{(i,j) \in \vec{E}}$
be a set of commuting indeterminates associated to the edges of $\Gamma$
(we refer to them as the \emph{edge weights}).
%
For $i,j \in V$, a \emph{walk} from $i$ to $j$ (of length $n \ge 0$)
is a sequence $\gamma = (\gamma_0,\ldots,\gamma_n)$ in $V$
such that $\gamma_0 = i$, $\gamma_n = j$,
and $(\gamma_{k-1},\gamma_k) \in \vec{E}$ for $1 \le k \le n$.
The \emph{weight} of a walk is the product of the weights of its edges:
\be
   W(\gamma)  \;=\;  \prod_{k=1}^n w_{\gamma_{k-1} \gamma_k}
   \;,
\ee
where the empty product is defined as usual to be 1.
Note that for each $i,j \in V$, the total weight of walks from $i$ to $j$,
namely
\be
   b_{ij}  \;=\;  \sum_{\gamma \colon\, i \to j}  W(\gamma)
   \;,
\ee
is a well-defined element of the formal-power-series ring $\Z[[{\bf w}]]$,
since each monomial in the indeterminates ${\bf w}$
corresponds to at most finitely many walks.
We therefore define the \emph{walk matrix} $B = (b_{ij})_{i,j \in V}$
with entries in $\Z[[{\bf w}]]$.

Now let us consider an $r \times r$ minor of the walk matrix,
corresponding to rows $i_1,\ldots,i_r \in V$ (the \emph{source vertices})
and columns $j_1,\ldots,j_r \in V$ (the \emph{sink vertices});
here $i_1,\ldots,i_r$ are all distinct, and $j_1,\ldots,j_r$ are all distinct,
but the two collections are allowed to overlap.
The Lindstr\"om--Gessel--Viennot lemma \cite{Lindstrom_73,Gessel-Viennot_89}
expresses such a minor as a sum over {\em vertex-disjoint}\/ systems of walks
from the source vertices to the sink vertices,
under the condition that the directed graph $\Gamma$
is {\em acyclic}\/ (i.e.\ has no directed cycles):\footnote{
   Lindstr\"om \cite{Lindstrom_73} inadvertently omitted the
   key condition that $\Gamma$ should be acyclic.
   This condition was made explicit by
   Gessel and Viennot \cite{Gessel-Viennot_89}.
}

\begin{lemma}[Lindstr\"om--Gessel--Viennot lemma]
   \label{lemma.gessel-viennot}
Suppose that $\Gamma$ is acyclic.  Then
\be
   \det B\biggl(\!\! \begin{array}{c}
                        i_1,\ldots,i_r \\
                        j_1,\ldots,j_r
                     \end{array}
           \!\!\biggr)
   \;=\;
   \sum_{\sigma \in \Sym_r}  \sgn(\sigma)
   \!\!\!
   \sum_{\begin{scarray}
             \gamma_1 \colon\, i_1 \to j_{\sigma(1)}  \\[-1mm]
                   \vdots \\[-1mm]
             \gamma_r \colon\, i_r \to j_{\sigma(r)}  \\[1mm]
             \gamma_k \cap \gamma_\ell = \emptyset
               \hboxscript{ for } k \ne \ell
         \end{scarray}
        }
   \!\!\!
   \prod_{i=1}^r W(\gamma_i)
 \label{eq.lemma.gessel-viennot}
\ee
where the sum runs over all $r$-tuples of vertex-disjoint walks
$\gamma_1,\ldots,\gamma_r$ connecting the specified vertices;
here ``vertex-disjoint'' means that for $1 \le k < \ell \le r$,
the walks $\gamma_k$ and $\gamma_\ell$ have no vertices in common
(not even the endpoints).
\end{lemma}

\noindent
Please note that the right-hand side of \reff{eq.lemma.gessel-viennot}
is a well-defined element of $\Z[[{\bf w}]]$,
since each monomial corresponds to at most finitely many
systems of walks $\gamma_1,\ldots,\gamma_r$.
See e.g. \cite{Gessel-Viennot_89,Aigner_01b,Aigner_14} for proofs
of the Lindstr\"om--Gessel--Viennot lemma.

Now let $\bfi = (i_1,\ldots,i_r)$
be an ordered $r$-tuple of distinct vertices of $\Gamma$,
and let $\bfj = (j_1,\ldots,j_r)$ be another such ordered $r$-tuple;
we say that the pair $(\bfi, \bfj)$ is \textbfit{nonpermutable}
if the set of vertex-disjoint walk systems
$\bgamma = (\gamma_1,\ldots,\gamma_r)$
satisfying $\gamma_k \colon\, i_k \to j_{\sigma(k)}$
is empty whenever $\sigma$ is not the identity permutation.
In this situation we avoid the sum over permutations
in \reff{eq.lemma.gessel-viennot};
most importantly, we avoid the possibility of terms with $\sgn(\sigma) = -1$.
The nonpermutable case arises frequently in applications:
for instance, if $\Gamma$ is planar, then the nonpermutability
often follows from topological arguments.

Now let $I$ and $J$ be subsets (not necessarily finite) of $V$,
equipped with total orders $<_I$ and $<_J$, respectively.
We say that the pair $\bigl((I,<_I),(J,<_J)\bigr)$
is \textbfit{fully nonpermutable}
if for each $r \ge 1$ and each pair of increasing $r$-tuples
$\bfi = (i_1,\ldots,i_r)$ in $(I,<_I)$
and $\bfj = (j_1,\ldots,j_r)$ in $(J,<_J)$,
the pair $(\bfi, \bfj)$ is nonpermutable.
The Lindstr\"om--Gessel--Viennot lemma then has the immediate consequence:

\begin{corollary}[Lindstr\"om--Gessel--Viennot lemma and total positivity]
   \label{cor.gessel-viennot.totpos}
Suppose that $\Gamma$ is acyclic.
Let $(I,<_I)$ and $(J,<_J)$ be totally ordered subsets of $V$
such that the pair $\bigl((I,<_I),(J,<_J)\bigr)$
is fully nonpermutable.
Then the submatrix $B_{IJ}$
(with rows and columns ordered according to $<_I$ and $<_J$)
is totally positive with respect to
the coefficientwise order on $\Z[[{\bf w}]]$.
\end{corollary}

\noindent
See \cite[pp.~179--180]{Brenti_95} for an example of a nonpermutable pair
$(\bfi, \bfj)$ that is not fully nonpermutable.

We will use repeatedly the following ``topologically obvious'' fact:

\begin{lemma}
   \label{lemma.nonpermutable}
If $\Gamma$ is embedded in the plane
and the vertices of $I \cup J$ lie on the boundary of $\Gamma$
in the order ``first $I$ in reverse order, then $J$ in order'',
then the pair $(I,J)$ is fully nonpermutable.
\end{lemma}

\begin{figure}[t]
\begin{center}
\includegraphics[scale=1.8]{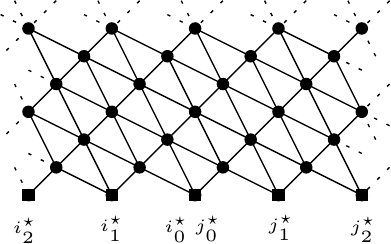}
\caption{\label{fig.Tplanarity}The subsets $I^\star$ and $J^\star$  for $m=2$.
            Vertices of $I^\star \cup J^\star$ are indicated by squares.}
\end{center}
\end{figure}

With these preliminaries in hand, we can now give the promised
combinatorial proofs of
Theorems~\ref{thm.Stype.minors}--\ref{thm.Ttype.generalized.minors}.
In fact, these proofs are very easy, given what we have already developed.
We begin with Theorem~\ref{thm.Ttype.minors} on Hankel minors:

\proofof{Theorem~\ref{thm.Ttype.minors}}
Form the infinite Hankel matrix corresponding to
the sequence $\bT^{(m)} = ( T^{(m)}_n(\balpha,\bdelta) )_{n \ge 0}$
of $m$-Thron--Rogers polynomials:
\be
   H_\infty(\bT^{(m)})
   \;=\;
   \bigl( T^{(m)}_{i+j}(\balpha,\bdelta) \bigr)_{i,j \ge 0}
   \;.
\ee
And consider any $k \times k$ minor of $H_\infty(\bT^{(m)})$:
that is, we choose sets of integers
$I = \{i_1,i_2,\ldots,i_k\}$ with $0 \le i_1 < i_2 < \ldots < i_k$
and $J = \{j_1,j_2,\ldots,j_k\}$ with $0 \le j_1 < j_2 < \ldots < j_k$,
and we consider the $k \times k$ submatrix
\be
   H_{IJ}(\bT^{(m)})
   \;=\;
   \bigl( T^{(m)}_{i_r+j_s}(\balpha,\bdelta) \bigr)_{1 \le r,s \le k}
\ee
and the corresponding minor
\be
   \Delta_{IJ}(\bT^{(m)})  \;=\;  \det H_{IJ}(\bT^{(m)})  \;.
\ee
We can write the elements of the submatrix $H_{IJ}(\bT^{(m)})$
as sums over walks in
the directed graph $\widetilde{G}_m = (V_m, \widetilde{E}_m)$
with vertex set \reff{def.Vm} and edge set \reff{def.Etildem}:
it is easy to see that $T^{(m)}_{i_r+j_s}(\balpha,\bdelta)$
is the sum over walks from $(-(m+1)i_r,0)$ to $((m+1)j_s,0)$,
with a weight~1 on each rising directed edge,
weight~$\alpha_i$ on each $m$-falling directed edge starting at height~$i$,
and weight $\delta_{i+1}$ on each $m$-long directed edge
starting at height~$i$.
It is also easy to see that $\widetilde{G}_m$ is planar and acyclic.
Finally --- and most importantly --- the totally ordered subsets of vertices
$I^\star = \{ (0,0)=i^\star_0 < (-(m+1),0)=i^\star_1 < (-2(m+1),0)=i^\star_2 < \ldots \}$
and $J^\star = \{ (0,0)=j^\star_0 < (m+1,0)=j^\star_1 < (2(m+1),0)=j^\star_2 < \ldots \}$
form a fully nonpermutable pair:
this follows from Lemma~\ref{lemma.nonpermutable}
because the sets $I^\star$ and $J^\star$
lie on the boundary of $\widetilde{G}_m$
and are appropriately ordered
(see Figure~\ref{fig.Tplanarity} for an example with $m=2$).

\noindent
Applying the Lindstr\"om--Gessel--Viennot lemma, we immediately obtain:

\begin{proposition}
   \label{prop.Ttype.minors}
The minor $\Delta_{IJ}(\bT^{(m)})$ is the generating polynomial
for families of vertex-disjoint $m$-Schr\"oder paths $P_1,\ldots,P_k$
where path $P_r$ starts at $(-(m+1)i_r,0)$ and ends at $((m+1)j_r,0)$,
in which each rise gets weight 1,
each $m$-fall from height~$i$ gets weight $\alpha_i$,
and each $m$-long step from height~$i$ gets weight $\delta_{i+1}$.
\end{proposition}

And since every such minor is manifestly a polynomial in $\balpha$
and $\bdelta$ with nonnegative integer coefficients,
we have proven Theorem~\ref{thm.Ttype.minors}.
\qed

\bigskip

{\bf Remark.}
For Dyck paths (with $m=1$), the analogue of
Proposition~\ref{prop.Ttype.minors}
was found already a quarter-century ago
by Viennot \cite[pp.~IV-13--IV-15]{Viennot_83}.
He could easily have deduced from it the Hankel-total positivity
of the classical Stieltjes--Rogers polynomials,
had he thought to pose that question.
\myendremark

\medskip

\begin{figure}[!ht]
\begin{center}
\includegraphics[scale=1.8]{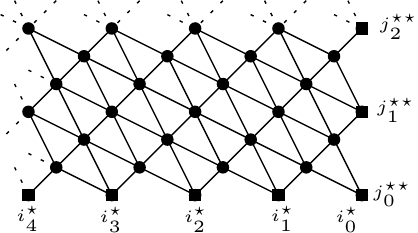}
\caption{\label{fig.thm.Ttype.generalized.minors}The subsets $I^\star$ and $J^{\star\star}$  for $m=2$.}
\end{center}
\end{figure}

We next prove Theorem~\ref{thm.Ttype.generalized.minors}
on the minors of the lower-triangular matrix:

\proofof{Theorem~\ref{thm.Ttype.generalized.minors}}
The proof follows closely that of Theorem~\ref{thm.Ttype.minors};
the only difference is that the sink vertices are now
$J^{\star\star} = \{ (0,0) = j^{\star\star}_0 < (0,m+1) = j^{\star\star}_1 < (0,2(m+1)) = j^{\star\star}_2 < \ldots \}$,
so that the matrix element $T^{(m)}_{i_r,j_s}(\balpha,\bdelta)$
is the sum over walks from $(-(m+1)i_r,0)$ to $(0,(m+1)j_s)$.
Now any path from any $i_r^\star$ to any $j_s^{\star\star}$
in the graph $\widetilde{G}_m$
must actually lie in $\widetilde{G}_m \cap (-\N \times \N)$.
And by Lemma~\ref{lemma.nonpermutable}
the pair $(I^\star, J^{\star\star})$ is fully nonpermutable,
because $I^\star$ and $J^{\star\star}$ lie on the boundary
of $\widetilde{G}_m \cap (-\N \times \N)$
and are appropriately ordered
(see Figure~\ref{fig.thm.Ttype.generalized.minors}).

\noindent
This proves Theorem~\ref{thm.Ttype.generalized.minors}.
\qed

\begin{figure}[t]
\begin{center}
\includegraphics[scale=1.8]{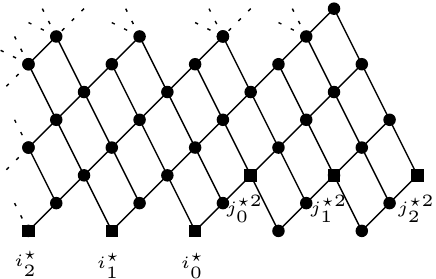}
\caption{\label{fig.thm.Stype.extended.minors}The subsets $I^\star$ and $J^{\star \ell}$  for $m=2$ and $\ell=2$.}
\end{center}
\end{figure}

Let us now prove Theorems~\ref{thm.Stype.minors.extended} and
\ref{thm.Ttype.minors.extended}, which extend
Theorems~\ref{thm.Stype.minors} and \ref{thm.Ttype.minors}
to the polynomials generated by products $f_0 f_1 \cdots f_\ell$
with $0 \le \ell \le m$ (Stieltjes) and $0 \le \ell \le m-1$ (Thron).
We begin with the Stieltjes case:

\proofof{Theorem~\ref{thm.Stype.minors.extended}}
The source vertices are as before, while the sink vertices are now
$J^{\star\ell} = \{ (\ell,\ell) < (m+1+\ell,\ell) < (2(m+1)+\ell,\ell) < \ldots \}$.
When $\ell \le m$, the sink vertices $J^{\star\ell}$
lie on the boundary of the graph $G_m$,
so the pair $(I^\star, J^{\star\ell})$ is fully nonpermutable
by Lemma~\ref{lemma.nonpermutable}.
See Figure~\ref{fig.thm.Stype.extended.minors}.
\qed

\begin{figure}[!ht]
\begin{center}
\includegraphics[scale=1.8]{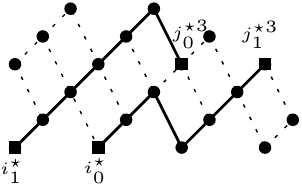}
\caption{\label{fig.thm.Stype.extended.minors.counter}
   A pair of nonintersecting paths from $i_0$ to $j_1^{\star 3}$
   and from $i_1$ to $j_0^{\star 3}$
   in the case of $m=2$ and $\ell=3$,
   showing that the pair $(I^\star, J^{\star 3})$
   is not fully nonpermutable.
}
\end{center}
\end{figure}

Please note that the pair $(I^\star, J^{\star\ell})$
is {\em not}\/ fully nonpermutable if $\ell > m$:
the sink vertices no longer lie on the boundary of the graph,
and there {\em does}\/ exist a vertex-disjoint pair of walks
$i_0^\star \to j_1^{\star\ell}$ and $i_1^\star \to j_0^{\star\ell}$.
See Figure~\ref{fig.thm.Stype.extended.minors.counter} for an example
in the case $(m,\ell) = (2,3)$, which generalizes naturally
to all cases of $1 \le m < \ell$.

\noindent
Indeed, for $\ell > m$ the conclusion of
Theorem~\ref{thm.Stype.minors.extended} is {\em false}\/: we have
\begin{subeqnarray}
   S^{(m)}_{0|\ell}
   & = &
   1
        \\[2mm]
   S^{(m)}_{1|\ell}
   & = &
   \sum_{i=m}^{m+\ell} \alpha_i
        \\[2mm]
   S^{(m)}_{2|\ell}
   & = &
   \sum_{i=m}^{2m+\ell} \sum_{j=\max(i-m,m)}^{m+\ell} \alpha_i \alpha_j
 \label{eq.ellgtm.counterexample.1}
\end{subeqnarray}
and hence for $\ell > m$
\be
   S^{(m)}_{0|\ell} \, S^{(m)}_{2|\ell} \:-\: (S^{(m)}_{1|\ell})^2
   \;=\;
   - \sum_{i=2m+1}^{m+\ell} \sum_{j=m}^{i-m-1} \alpha_i \alpha_j
   \:+\:
   \sum_{i=m+\ell+1}^{2m+\ell} \sum_{j=i-m}^{m+\ell} \alpha_i \alpha_j
   \;,
 \label{eq.ellgtm.counterexample.2}
\ee
which is $\not\myge 0$ since the coefficient of
$\alpha_m \alpha_{2m+1}$ is $-1$.
%
   
\begin{figure}[!ht]
\begin{center}
\includegraphics[scale=1.8]{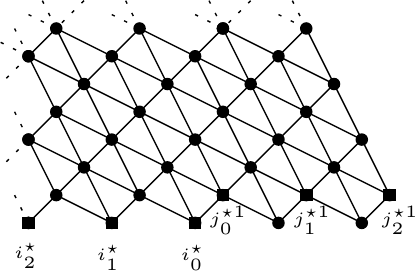}
\caption{\label{fig.thm.Ttype.extended.minors}The subsets $I^\star$ and $J^{\star \ell}$  for $m=2$ and $\ell=1$.}
\end{center}
\end{figure}

\proofof{Theorem~\ref{thm.Ttype.minors.extended}}
The argument is the same as that for Theorem~\ref{thm.Stype.minors.extended},
except that in the Thron case we need
$\ell \le m-1$ for the sink vertices $J^{\star\ell}$
to lie on the boundary of the graph $\widetilde{G}_m$,
because of the extra edges $(2,-(m-1))$ for the $m$-long steps.
See Figure~\ref{fig.thm.Ttype.extended.minors}.

\qed

\begin{figure}[!ht]
\begin{center}
\includegraphics[scale=1.8]{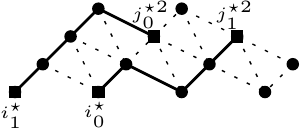}
\caption{\label{fig.thm.Ttype.extended.minors.counter}
   A pair of nonintersecting paths from $i_0$ to $j_1^{\star 2}$
   and from $i_1$ to $j_0^{\star 2}$
   in the case of $m=2$ and $\ell=2$,
   showing that the pair $(I^\star, J^{\star 2})$
   is not fully nonpermutable.
}
\end{center}
\end{figure}

Once again, the pair $(I^\star, J^{\star\ell})$
fails to be fully nonpermutable if $\ell > m-1$:
see Figure~\ref{fig.thm.Ttype.extended.minors.counter} for an example
in the case $(m,\ell) = (2,2)$, which generalizes naturally
to all cases of $1 \le m \le \ell$.

\subsection{Algebraic proofs using production matrices}
   \label{subsec.totalpos.algebraic}

The algebraic proofs of Theorems~\ref{thm.Stype.minors}
and \ref{thm.Stype.generalized.minors}
concerning the $m$-Stieltjes--Rogers polynomials are now easy.
By Proposition~\ref{prop.prod.Sm.Jm.Tm}(b),
the production matrix for the generalized $m$-Stieltjes--Rogers polynomials
is the matrix $P^{(m)\mathrm{S}}(\balpha)$
defined by \reff{eq.prop.contraction};
and by Lemma~\ref{lemma.bidiagonal},
this matrix is totally positive in the polynomial ring $\Z[\balpha]$
equipped with the coefficientwise partial order.
Then Theorem~\ref{thm.iteration.homo} implies that the
lower-triangular matrix of generalized $m$-Stieltjes--Rogers polynomials
is coefficientwise totally positive
(Theorem~\ref{thm.Stype.generalized.minors}),
while Theorem~\ref{thm.iteration2bis} implies that
the sequence of $m$-Stieltjes--Rogers polynomials
is coefficientwise Hankel-totally positive
(Theorem~\ref{thm.Stype.minors}).

As for Theorems~\ref{thm.Ttype.minors} and \ref{thm.Ttype.generalized.minors}
concerning the $m$-Thron--Rogers polynomials,
we can handle by algebraic methods the restricted case
defined in Section~\ref{subsec.contraction.even}
in which we impose $\delta_i = 0$ whenever $i$ is a multiple of $m+1$.
By Proposition~\ref{prop.prod.Sm.Jm.Tm}(c),
the production matrix for the restricted generalized
$m$-Thron--Rogers polynomials
is the matrix $P^{(m)\mathrm{RT}}(\balpha,\bdelta')$
defined in Proposition~\ref{prop.contraction.thron}.
Since we have only a combinatorial definition of this matrix,
not an algebraic one, we will have to prove its total positivity
combinatorially:

\begin{lemma}
   \label{lemma.TP.prodmat.restrictedTR}
Fix an integer $m \ge 1$.
Then the matrix $P^{(m)\mathrm{RT}}(\balpha,\bdelta')$
defined in Proposition~\ref{prop.contraction.thron}
is totally positive, coefficientwise in $\balpha$ and $\bdelta'$.
\end{lemma}

\begin{figure}[!ht]
\begin{center}
\includegraphics[scale=1.8]{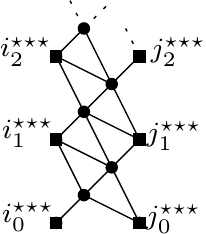}
\caption{
   The subsets $I^{\star\star\star}$ and $J^{\star\star\star}$ for $m=2$.
}
\label{fig.lemma.TP.prodmat.restrictedTR}
\end{center}
\end{figure}

\proof
We apply the Lindstr\"om--Gessel--Viennot lemma
to the graph $\widetilde{G}_m \cap {([0,m+1] \times \N)}$
with
$I^{\star\star\star} = \{ (0,0)=i^{\star\star\star}_0 < (0,m+1)=i^{\star\star\star}_1 < (0,2(m+1))=i^{\star\star\star}_2 < \ldots \}$
and
$J^{\star\star\star} = \{ (m+1,0)=j^{\star\star\star}_0 < (m+1,m+1)=j^{\star\star\star}_1 < (m+1,2(m+1))=j^{\star\star\star}_2 < \ldots \}$.
By Lemma~\ref{lemma.nonpermutable}
the pair $(I^{\star\star\star}, J^{\star\star\star})$ is fully nonpermutable,
because $I^{\star\star\star}$ and $J^{\star\star\star}$ lie on the boundary
of $\widetilde{G}_m \cap ([0,m+1] \times \N)$
and are appropriately ordered
(see Figure~\ref{fig.lemma.TP.prodmat.restrictedTR}).

\noindent
\qed

With this lemma in hand, the desired results then follow from
Theorems~\ref{thm.iteration.homo} and \ref{thm.iteration2bis}
just as in the Stieltjes case.

\section{Weights periodic of period $\bm{m\!+\!1}$ or $\bm{m}$}
   \label{sec.periodic}

Consider the $m$-Stieltjes--Rogers polynomials $S_n^{(m)}(\balpha)$
when the weights $\balpha = (\alpha_i)_{i \ge m}$
are taken to be periodic of period~$p \ge 1$:
\be
   \alpha_{m+j+p k}  \;=\;  x_j
 \label{eq.alpha.periodic}
\ee
where $\bfx = (x_0,\ldots,x_{p-1})$ will be treated as indeterminates.
Let $P_n^{(m,p)}(\bfx)$
be the polynomials obtained by specializing
$S_n^{(m)}(\balpha)$ to these weights.
In this section we will obtain explicit formulae
for the polynomials $P_n^{(m,p)}(\bfx)$
when the period~$p$ is either $m+1$ or $m$.
The results will be {\em multivariate Fuss--Narayana polynomials}\/
\cite{Lenczewski_13,Gessel_16},
which reduce to the ordinary Narayana polynomials when $m=1$
and to the Fuss--Catalan numbers when $\bfx = \bone$.
Furthermore, after obtaining these formulae
it will be natural to pass to the limit $m \to\infty$;
then everything will be expressed in terms of the
{\em Fuss--Narayana symmetric functions}\/ \cite{Stanley_97},
which contain as specializations the multivariate Fuss--Narayana polynomials
for all $m$.


We will perform this computation in two ways:
first, by solving the recurrences \reff{eq.mSRfk.1} and \reff{eq.mSRfk.2}
using Lagrange inversion
(Sections~\ref{subsec.periodic.m+1.recurrence} and
 \ref{subsec.periodic.m.recurrence});
and second, by working with the production matrices
(Sections~\ref{subsec.periodic.m+1.prodmat} and
 \ref{subsec.periodic.row-generating}).
The second method is more complicated but provides additional insight,
and also allows the computation of the generalized
$m$-Stieltjes--Rogers polynomials $S_{n,k}^{(m)}(\balpha)$
and their row-generating polynomials.
However, at present we are able to carry out the production-matrix
calculation only for period $m+1$, not for period $m$.

Along the way (Section~\ref{subsec.periodic.combinatorial})
we will also give combinatorial interpretations
for our multivariate Fuss--Narayana polynomials
and Fuss--Narayana symmetric functions
in terms of ordered trees (of various types),
noncrossing partitions, and parking functions.

We will use Lagrange inversion in the following form \cite{Gessel_16}:
If $\phi(u)$ is a formal power series
with coefficients in a commutative ring $R$ containing the rationals,
then there exists a unique formal power series $f(t)$
with zero constant term satisfying
\be
   f(t)  \;=\;  t \, \phi(f(t))
   \;,
\ee
and it is given by
\be
   [t^n] \, f(t)  \;=\;  {1 \over n} \, [u^{n-1}] \, \phi(u)^n
     \quad\hbox{for $n \ge 1$}
   \;;
\ee
and more generally, if $H(u)$ is any formal power series, then
\be
   [t^n] \, H(f(t))  \;=\;  {1 \over n} \, [u^{n-1}] \, H'(u) \, \phi(u)^n
     \quad\hbox{for $n \ge 1$}
   \;.
 \label{eq.lagrange.H}
\ee

\subsection[Period $m\!+\!1$: Solving the recurrence]{Period $\bm{m\!+\!1}$: Solving the recurrence}
   \label{subsec.periodic.m+1.recurrence}


Let $\balpha$ be given by \reff{eq.alpha.periodic} with period~$p = m+1$,
or in other words $\alpha_i = x_{i+1 \bmod m+1}$
with $\bfx = (x_0,\ldots,x_m)$;
this means that we get a weight $x_i$
for each $m$-fall {\em to}\/ a height equal to $i \bmod m+1$.
We write $P_n^{(m)}$ as a shorthand for $P_n^{(m,m+1)}$.
We use the recurrence \reff{eq.mSRfk.1}.
The periodicity of the $\balpha$ implies the periodicity of the $f_k$;
therefore the product $f_k \cdots f_{k+m}$ is independent of $k$,
let us call it $F$.
Multiplying \reff{eq.mSRfk.1} for $0 \le k \le m$
leads to the functional equation
\be
   F  \;=\;  \prod_{i=0}^m (1 + x_i t F)
   \;.
 \label{eq.functeqn.F}
\ee
Equivalently, defining $f = tF$, the functional equation is
\be
   f  \;=\;  t \, \prod_{i=0}^m (1 + x_i f)
   \;.
 \label{eq.functeqn.f}
\ee
Let us write the solution of \reff{eq.functeqn.F} as
\be
   F(t)  \;=\;  \sum_{n=0}^\infty Q^{(m)}_n(\bfx) \: t^n
 \label{def.F.Qn}
\ee
and more generally
\be
   F(t)^{k+1}  \;=\;  \sum_{n=k}^\infty Q^{(m)}_{n,k}(\bfx) \: t^{n-k}
   \quad\hbox{for integers $k \ge 0$}  \;.
 \label{def.F.Qnk}
\ee
Solving \reff{eq.functeqn.f} for $f$ by Lagrange inversion gives
\cite[eq.~(3.4.3)]{Gessel_16}
\be
   Q^{(m)}_{n,k}(\bfx)
   \;=\;
   [t^{n-k}] \, F(t)^{k+1}
   \;=\;
   [t^{n+1}] \, f(t)^{k+1}
   \;=\;
   {k+1 \over n+1}
   \!\!
   \sum_{\begin{scarray}
            j_0,\ldots,j_m \ge 0 \\
            \sum j_i = n-k
         \end{scarray}}
   \!\!
   \prod_{i=0}^m \binom{n+1}{j_i} \, x_i^{j_i}
   \;.
 \label{def.Qn.x0xm}
\ee
Here $Q^{(m)}_{n,k}(\bfx)$ is a homogeneous polynomial of degree~$n-k$,
which is symmetric in $x_0,\ldots,x_m$.
When $\bfx = \bone$, Chu--Vandermonde gives
$Q^{(m)}_{n,k}(\bone) = \displaystyle {k+1 \over n+1} \binom{(m+1)(n+1)}{n-k}$;
in particular, $Q^{(m)}_{n}(\bone)$
equals the Fuss--Catalan number $C_{n+1}^{(m+1)}$.

Inserting \reff{def.F.Qn}/\reff{def.Qn.x0xm} into \reff{eq.mSRfk.1} gives
\be
   f_0(t)  \;=\;  1 \,+\, x_0 t F(t)
           \;=\;  \sum_{n=0}^\infty P^{(m)}_n(\bfx) \: t^n
 \label{def.f0.Pn}
\ee
where $P^{(m)}_0 = 1$ and
\begin{subeqnarray}
   P^{(m)}_n(\bfx)
   & = &
   x_0 \, Q^{(m)}_{n-1}(\bfx)
     \slabel{def.Pn.x0xm.a} \\[2mm]
   & = &
   {1 \over n}
   \!\!
   \sum_{\begin{scarray}
            j_0,\ldots,j_m \ge 0 \\
            \sum j_i = n
         \end{scarray}}
   \!\!\!
   \binom{n}{j_0 -1} \, x_0^{j_0} \:
   \prod_{i=1}^m \binom{n}{j_i} \, x_i^{j_i}
   \qquad
 \label{def.Pn.x0xm}
\end{subeqnarray}
for $n \ge 1$.
Here $P^{(m)}_n(\bfx)$ is a homogeneous polynomial of degree~$n$:
it~is symmetric in $x_1,\ldots,x_m$, but $x_0$ plays a distinguished role;
however, this distinguished role amounts only to a prefactor $x_0$
for $n \ge 1$, and the remaining polynomial $Q_{n-1}^{(m)}(\bfx)$
is symmetric in $x_0,\ldots,x_m$.
The invariance of $P^{(m)}_n(\bfx)$ under permutations of $x_1,\ldots,x_m$
illustrates once again the nonuniqueness of $m$-S-fractions with $m \ge 2$.
When $\bfx = \bone$, $P^{(m)}_n$ equals the Fuss--Catalan number $C_n^{(m+1)}$.
Since
\be
    Q^{(m)}_n(\bfx)  \;=\;  P^{(m)}_{n+1}(\bfx) / x_0
   \;,
 \label{eq.Qn.Pn+1}
\ee
the sequence $(Q^{(m)}_n)_{n \ge 0}$ is, up to a factor $x_0$,
a subsequence of $(P^{(m)}_n)_{n \ge 0}$.

When $m=1$ these polynomials reduce to the homogenized Narayana polynomials
\begin{eqnarray}
   P^{(1)}_n(x,y)  & = &  \sum_{j=0}^n N(n,j)     \,  x^j y^{n-j}  \\[2mm]
   Q^{(1)}_n(x,y)  & = &  \sum_{j=0}^n N(n+1,j+1) \,  x^j y^{n-j}
\end{eqnarray}
where the Narayana numbers \cite[A001263]{OEIS}
are defined by $N(0,j) = \delta_{j0}$ and\footnote{
   {\bf Warning:}  Many authors (e.g.\ \cite{Petersen_15})
   define the Narayana numbers to be the reversal of ours, i.e.
$$
   N^{\rm theirs}(n,j)  \;=\;
   N(n,n-j)  \;=\;  {1 \over n} \, \binom{n}{j+1} \binom{n}{j}
   \qquad\hbox{for $n \ge 1$}
   \;.
$$
   For $n \ge 1$, our Narayana numbers are nonvanishing for $1 \le j \le n$,
   while in the other convention they are nonvanishing for $0 \le j \le n-1$.
}
\be
   N(n,j)  \;=\;  {1 \over n} \, \binom{n}{j-1} \binom{n}{j}
   \qquad\hbox{for $n \ge 1$}
   \;.
 \label{def.narayana}
\ee
We therefore call $P^{(m)}_n(\bfx)$, $Q^{(m)}_n(\bfx)$
and $Q^{(m)}_{n,k}(\bfx)$
the \textbfit{multivariate Fuss--Narayana polynomials}
\cite{Lenczewski_13,Gessel_16}.

Since the weights \reff{eq.alpha.periodic} are manifestly nonnegative
with respect to the coefficientwise order on $\Z[\bfx]$,
the coefficientwise Hankel-total positivity of the $P^{(m)}_n$
is an immediate consequence of Theorem~\ref{thm.Stype.minors}.
As for $Q^{(m)}_n$, its Hankel-total positivity can either be deduced
from that of $P^{(m)}_n$ by using \reff{eq.Qn.Pn+1},
or else obtained directly as a consequence of
Theorem~\ref{thm.Stype.minors.extended} with $\ell = m$.
We conclude:

\begin{theorem}[Hankel-total positivity of multivariate Fuss--Narayana polynomials]
   \label{thm.QnPn.hankelTP}
For each integer $m \ge 1$,
the sequences $(P^{(m)}_n)_{n \ge 0}$ and $(Q^{(m)}_n)_{n \ge 0}$
are coefficientwise Hankel-totally positive,
jointly in the indeterminates $x_0,\ldots,x_m$.
\end{theorem}

Note now that whenever $m < m'$ we have
$Q^{(m)}_{n,k}(x_0,\ldots,x_m) = Q^{(m')}_{n,k}(x_0,\ldots,x_m,0,\ldots,0)$
and
$P^{(m)}_n(x_0,\ldots,x_m) = P^{(m')}_n(x_0,\ldots,x_m,0,\ldots,0)$:
this follows immediately either from the explicit formulae
\reff{def.Qn.x0xm}/\reff{def.Pn.x0xm}
or from the branched continued fraction
using Proposition~\ref{prop.reduction}.
We can therefore pass to the limit $m \to\infty$ and define
\begin{eqnarray}
   Q^{(\infty)}_{n,k}(\bfx)
   & = &
   {k+1 \over n+1}
   \sum_{\begin{scarray}
            j_0,j_1,\ldots \ge 0 \\
            \sum j_i = n-k
         \end{scarray}}
   \prod_{i=0}^\infty \binom{n+1}{j_i} \, x_i^{j_i}
   \hspace*{1.8cm}
      \label{def.Qinfty}
\end{eqnarray}
\begin{subeqnarray}
   P^{(\infty)}_n(\bfx)
   & = &
   x_0 \, Q^{(\infty)}_{n-1}(\bfx)
          \\[2mm]
   & = &
   {1 \over n}
   \sum_{\begin{scarray}
            j_0,j_1,\ldots \ge 0 \\
            \sum j_i = n
         \end{scarray}}
   \!\!
   \binom{n}{j_0 -1} \, x_0^{j_0} \:
   \prod_{i=1}^\infty \binom{n}{j_i} \, x_i^{j_i}
   \qquad
      \label{def.Pinfty}
\end{subeqnarray}
where $\bfx = (x_i)_{i \ge 0}$ is an infinite set of indeterminates.
Then $Q^{(m)}_{n,k}$ and $P^{(m)}_n$
can be obtained from $Q^{(\infty)}_{n,k}$ and $P^{(\infty)}_n$
by specializing $x_i = 0$ for $i > m$.
Since the sums \reff{def.Qinfty}/\reff{def.Pinfty}
have infinitely many terms,
they have to be interpreted as formal power series
(rather than polynomials) in the indeterminates $\bfx$.
In fact, $Q^{(\infty)}_{n,k}$ belongs to the
ring of symmetric functions in $\bfx$,
while $P^{(\infty)}_n$ is a symmetric function of $\bfx$
up to a prefactor~$x_0$.\footnote{
   Our notation for symmetric functions follows
   Macdonald \cite[Chapter~1]{Macdonald_95}
   and Stanley \cite[Chapter~7]{Stanley_99}.
}
Moreover, the functional equation \reff{eq.functeqn.F} with $m=\infty$
can be rewritten in the language of symmetric functions as
\be
   F(t)  \;=\;  E(tF(t))
\ee
where
\be
   E(t)
   \;=\;
   \prod_{i=0}^\infty (1 + x_i t)
   \;=\;
   \sum_{n=0}^\infty e_n(\bfx) \, t^n
 \label{def.E}
\ee
is the generating function for the elementary symmetric functions.
Then Lagrange inversion gives
\begin{subeqnarray}
   Q^{(\infty)}_{n,k}(\bfx)
   & = &
   [t^{n-k}] \, F(t)^{k+1}   \;=\; {k+1 \over n+1} \, [t^{n-k}] \, E(t)^{n+1}
        \\[2mm]
   &  &
   \!\!=\;
   \!\!\!
   \sum_{\begin{scarray}
           r_1, r_2, \ldots \ge 0 \\
           \sum i r_i = n-k
         \end{scarray}}
 {k+1 \over n+1} \, \binom{n+1}{n+1 - \sum r_i,\, r_1,\, r_2,\, \ldots}
      \, \prod_{j=0}^\infty e_j(\bfx)^{r_j}
    \qquad\qquad
         \slabel{eq.Qinfty.symfun_e.a} \\[2mm]
   &  &
   \!\!=\;
   \sum_{\lambda \vdash n-k} {k+1 \over n+1} \,
       \binom{n+1}{n+1-\ell(\lambda),\, m_1(\lambda),\, m_2(\lambda),\, \ldots}
       \; e_\lambda(\bfx)
   \;,
   \qquad\qquad
         \slabel{eq.Qinfty.symfun_e}
  \label{eq.Qinfty.symfun.all}
\end{subeqnarray}
where the second sum runs over partitions $\lambda$ of the integer $n-k$;
here $\ell(\lambda)$ denotes the number of parts of $\lambda$,
and $m_i(\lambda)$ ($= r_i$) denotes
the number of parts of $\lambda$ equal to~$i$.
Alternatively, we can use
the generating function for the complete homogeneous symmetric functions,
\be
   H(t)  
   \;=\;  
   \prod_{i=0}^\infty {1 \over 1 - x_i t}
   \;=\;  
   \sum_{n=0}^\infty h_n(\bfx) \, t^n
   \;=\;
   {1 \over E(-t)}
   \;,
 \label{def.H}
\ee
to write
\begin{subeqnarray}
   Q^{(\infty)}_{n,k}(\bfx)
   & = &
   {k+1 \over n+1} \, [t^{n-k}] \, H(-t)^{-(n+1)}
        \\[2mm]
   & = &
   \!\!\!
   \sum_{\begin{scarray}
           r_1, r_2, \ldots \ge 0 \\
           \sum i r_i = n-k
         \end{scarray}}
 {k+1 \over n+1} \, \binom{-(n+1)}{-(n+1+\sum r_i),\, r_1,\, r_2,\, \ldots}
      \, \prod_{j=0}^\infty \big[ (-1)^j h_j(\bfx) \bigr]^{r_j}
  \nonumber \\[-6mm] \slabel{eq.Qinfty.symfun_h.a} \\[2mm]
   & = &
   \!\!\!
   \sum_{\begin{scarray}
           r_1, r_2, \ldots \ge 0 \\
           \sum i r_i = n-k
         \end{scarray}}
 {k+1 \over n+1} \, (-1)^{n - \sum r_i}  \,
      \binom{n+ \sum r_i}{n,\, r_1,\, r_2,\, \ldots}
      \, \prod_{j=0}^\infty h_j(\bfx)^{r_j}
  \slabel{eq.Qinfty.symfun_h.b} \\[2mm]
   & = &
   \sum_{\lambda \vdash n-k} {k+1 \over n+1} \,
       (-1)^{n - \ell(\lambda)} \,
       \binom{n+\ell(\lambda)}{n,\, m_1(\lambda),\, m_2(\lambda),\, \ldots}
       \; h_\lambda(\bfx)
   \;.
  \slabel{eq.Qinfty.symfun_h}
  \label{eq.Qinfty.symfun.h.all}
\end{subeqnarray}
Also, the definition \reff{def.Qinfty}
can trivially be rewritten
in terms of the monomial symmetric functions $m_\lambda$ as
\be
   Q^{(\infty)}_{n,k}(\bfx)
   \;=\;
   {k+1 \over n+1}
   \sum_{\lambda \vdash n-k}
   \left[ \prod_i \binom{n+1}{\lambda_i} \right] \, m_\lambda(\bfx)
   \;.
  \label{eq.Qinfty.symfun_m}
\ee
We call $Q^{(\infty)}_{n,k}(\bfx)$
the \textbfit{Fuss--Narayana symmetric functions}.
For $k=0$ they appeared two decades ago, in a different guise,
in the work of Stanley \cite{Stanley_97}
on parking-function symmetric functions:
see Section~\ref{subsubsec.periodic.parking} below.
Stanley also gave (in dual form)
the expansions of $Q^{(\infty)}_n(\bfx)$
in terms of the power-sum symmetric functions $p_\lambda$
and the Schur functions $s_\lambda$ \cite[eqns.~(2) and~(3)]{Stanley_97}.

\bigskip

Since the coefficientwise nonnegativity of Hankel minors
passes trivially to the limit $m \to\infty$,
Theorem~\ref{thm.QnPn.hankelTP} immediately implies:

\begin{theorem}[Hankel-total positivity of $\infty$-variate Fuss--Narayana series]
   \label{thm.QinftyPinfty.hankelTP}
\hfill\break
The sequences $(P^{(\infty)}_n)_{n \ge 0}$ and $(Q^{(\infty)}_n)_{n \ge 0}$
are Hankel-totally positive with respect to the coefficientwise order
on the formal-power-series ring $\Z[[\bfx]]$.
\end{theorem}

\smallskip

Return now to the case of finite $m$,
and suppose that the multiset $\{x_0,\ldots,x_m\}$
consists of elements $\bfy = (y_1,\ldots,y_l)$
with multiplicities $\bfp = (p_1,\ldots,p_l)$
(where of course $\sum p_i = m+1$).
Then the functional equation \reff{eq.functeqn.F} becomes
\be
   F  \;=\;  \prod_{i=1}^l (1 + y_i t F)^{p_i}
   \;,
 \label{eq.functeqn.F.bis}
\ee
and solving it by Lagrange inversion gives \cite[eq.~(3.4.3)]{Gessel_16}
\be
   F(t)^{k+1}
   \;=\;
   \sum_{n=k}^\infty \widetilde{Q}^{(m)}_{n,k}(\bfy;\bfp) \: t^{n-k}
\ee
where
\be
   \widetilde{Q}^{(m)}_{n,k}(\bfy;\bfp)
   \;=\;
   {k+1 \over n+1}
   \!\!
   \sum_{\begin{scarray}
            j_1,\ldots,j_l \ge 0 \\
            j_1+\ldots+j_l = n-k
         \end{scarray}}
   \!\!
   \prod_{i=1}^l \binom{p_i (n+1)}{j_i} \, y_i^{j_i}
   \;.
 \label{def.Qtilden.y1yl}
\ee
Then, if we make the convention that $x_0 = y_1$, we have
\be
   f_0(t)
   \;=\;
   1 \,+\, y_1 t F(t)
   \;=\;
   \sum_{n=0}^\infty \widetilde{P}^{(m)}_n(\bfy;\bfp) \: t^n
\ee
where $\widetilde{P}^{(m)}_0 = 1$ and
\begin{subeqnarray}
   \widetilde{P}^{(m)}_n(\bfy;\bfp)
   & = &
   y_1 \, \widetilde{Q}^{(m)}_{n-1,0}(\bfy;\bfp)
            \\[2mm]
   & = &
   {1 \over n}
   \!\!
   \sum_{\begin{scarray}
            j_1,\ldots,j_l \ge 0 \\
            j_1+\ldots+j_l = n
         \end{scarray}}
   \!\!
   \binom{p_1 n}{j_1 -1} \, y_1^{j_1} \:
   \prod_{i=2}^l \binom{p_i n}{j_i} \, y_i^{j_i}
 \slabel{def.Ptilden.y1yl}
\end{subeqnarray}
for $n \ge 1$.

\begin{example}[Interpolating between $p$-Fuss--Catalan and $p'$-Fuss--Catalan]
   \label{exam.pp'}
\rm
Let ${m' \ge 1}$ and $0 \le m \le m'$, and set $p = m+1$, $p' = m'+1$.
Let us consider an $m'$-branched continued fraction
with $x_0,\ldots,x_m = x$ and $x_{m+1},\ldots,x_{m'} = y$,
hence $y_1 = x$, $y_2 = y$, $p_1 = p$, $p_2 = p'-p$.
Then \reff{def.Ptilden.y1yl} becomes
\be
   \widetilde{P}^{(p,p')}_n(x,y)
   \;=\;
   \sum_{k=0}^n T^{(p,p')}_{n,j} \, x^k y^{n-j}
\ee
where $T^{(p,p')}_{0,j} = \delta_{j0}$ and
\be
   T^{(p,p')}_{n,j}
   \;=\;
   {1 \over n} \binom{pn}{j-1} \binom{(p'-p)n}{n-j}
   \qquad\hbox{for $n \ge 1$}
 \label{eq.Tpp.nk}
\ee
are generalized Narayana numbers.
For $(p,p') = (1,2)$ they reduce to the Narayana numbers \reff{def.narayana};
for $(p,p') = (1,3)$ they are \cite[A120986]{OEIS};
for $(p,p') = (2,3)$, $(3,4)$, $(4,5)$
they are \cite[A108767, A173020, A173621]{OEIS}.
Please note that
the triangular array $T^{(p,p')}$ has the $p$-Fuss--Catalan numbers
on the diagonal ($j=n$) and the $p'$-Fuss--Catalan numbers as the row sums;
moreover, if $p=p'$ we have $T^{(p,p)}_{n,j}  \;=\;  C_n^{(p)} \, \delta_{nj}$.
These facts are easily derived either from \reff{eq.Tpp.nk}
or from the branched continued fraction
(using Proposition~\ref{prop.reduction} to handle $y=0$).
It follows from Theorem~\ref{thm.QnPn.hankelTP}
that for each pair $1 \le p \le p'$,
the sequence of polynomials
$\big( \widetilde{P}^{(p,p')}_n(x,y) \big)_{n \ge 0}$
is coefficientwise Hankel-totally positive.

Alternatively, we can make the substitution $x = 1+z$, $y=z$:
simple algebra then yields
\be
   \widetilde{P}^{(p,p')\star}_n(z)
   \;\eqdef\;
   \widetilde{P}^{(p,p')}_n(1+z,z)
   \;=\;
   \sum_{j=0}^n T^{(p,p')\star}_{n,j} \, z^j
\ee
where $T^{(p,p')\star}_{0,j} = \delta_{j0}$ and
\be
   T^{(p,p')\star}_{n,j}
   \;=\;
   {1 \over n}
   \sum_{i=0}^j \binom{pn}{n-i-1} \binom{(p'-p)n}{i} \binom{n-i}{n-j}
   \qquad\hbox{for $n \ge 1$}
   \;.
 \label{eq.Tppstar.nk}
\ee
The triangular array $T^{(p,p')\star}$ has the $p$-Fuss--Catalan numbers
in the first ($j=0$) column and the $p'$-Fuss--Catalan numbers
on the diagonal;
moreover, if $p=p'$ we have $T^{(p,p)}_{n,j}  \;=\;  C_n^{(p)} \, \binom{n}{j}$.
These formulae are easily derived either from \reff{eq.Tppstar.nk}
or from the branched continued fraction.
Also, when $p' = p+1$, it is not hard to show from \reff{eq.Tppstar.nk} that
\be
   T^{(p,p+1)\star}_{n,j}
   \;=\;
   {1 \over (p-1)n+j+1}  \binom{pn+j}{n} \binom{n}{j}
   \;.
 \label{eq.Tppstar.nk.p+1}
\ee
For $(p,p') = (1,2)$ this is \cite[A088617/A060693]{OEIS};
for $(p,p') = (2,3)$ it is \cite[A104978]{OEIS}.
Writing $j = n-\ell$ in \reff{eq.Tppstar.nk.p+1},
we see that it is equivalent to \reff{def.mschroder.n};
therefore, the $(m+1)$-branched S-fraction of period $m+2$
with $x_0,\ldots,x_m = 1+z$ and $x_{m+1} = z$
is also equivalent to an $(m+1)$-branched T-fraction
with $\alpha_i = z$ and $\delta_i = 1$ for all $i$.
For $m=0$ this is a well-known fact about Narayana polynomials.
%
\myendremark
\end{example}

\subsection[Period $m$: Solving the recurrence]{Period $\bm{m}$: Solving the recurrence}
   \label{subsec.periodic.m.recurrence}

Now let $\balpha$ be given by \reff{eq.alpha.periodic} with period~$p = m$,
or in other words $\alpha_i = x_{i \bmod m}$
with $\bfx = (x_0,\ldots,x_{m-1})$;
this means that we get a weight $x_i$
for each $m$-fall {\em to}\/ a height equal to $i \bmod m$.
We write $P_n^{(m)-}$ as a shorthand for $P_n^{(m,m)}$.
This time we use the recurrence \reff{eq.mSRfk.2}.
The periodicity of the $\balpha$ implies the periodicity of the $f_k$;
therefore the product $f_{k+1} \cdots f_{k+m}$ is independent of $k$,
let us call it $G$.
Multiplying \reff{eq.mSRfk.2} for $0 \le k \le m-1$
leads to the functional equation
\be
   G  \;=\;  \prod_{i=0}^{m-1} (1 - x_i t G)^{-1}
   \;,
 \label{eq.functeqn.G}
\ee
or equivalently, defining $g = tG$,
\be
   g  \;=\;  t \prod_{i=0}^{r-1} (1 - x_i g)^{-1}
   \;.
 \label{eq.functeqn.g}
\ee
Writing
\begin{eqnarray}
   G(t)^{k+1}  & = &  \sum_{n=k}^\infty Q^{(m)-}_{n,k}(\bfx) \: t^{n-k}
     \label{def.G.Qn-k}
\end{eqnarray}
(with $Q_n^{(m)-} \eqdef Q^{(m)-}_{n,0}$)
and solving \reff{eq.functeqn.g} by Lagrange inversion gives
\cite[eq.~(3.4.5)]{Gessel_16}
where
\be
   Q_{n,k}^{(m)-}(\bfx)
   \;=\;
   {k+1 \over n+1}
   \!\!
   \sum_{\begin{scarray}
            j_0,\ldots,j_{m-1} \ge 0 \\
            \sum j_i = n-k
         \end{scarray}}
   \!\!
   \prod_{i=0}^{m-1} \binom{n+j_i}{j_i} \, x_i^{j_i}
   \;.
 \label{def.Qn-.x0xm-1}
\ee
When $\bfx = \bone$ we can use $\binom{n+j}{j} = (-1)^j \binom{-n-1}{j}$
together with Chu--Vandermonde to obtain
$Q^{(m)-}_{n,k}(\bone) =
 \displaystyle {k+1 \over n+1} \binom{(m+1)(n+1)-k-2}{n-k}$.

A further application of Lagrange inversion yields after some algebra
\be
   f_0(t)  \;=\;  {1 \over 1 \,-\, x_0 t G(t)}
           \;=\;  \sum_{n=0}^\infty P_n^{(m)-}(\bfx) \: t^n
 \label{def.f0.Pn-}
\ee
where $P_0^{(m)-} = 1$ and
\begin{eqnarray}
   P_n^{(m)-}(\bfx)
   & = &
   {1 \over n}
   \!\!
   \sum_{\begin{scarray}
            j_0,\ldots,j_{m-1} \ge 0 \\
            \sum j_i = n
         \end{scarray}}
   \!\!
   \binom{n+j_0}{j_0 -1} \, x_0^{j_0}
   \prod_{i=1}^{m-1} \binom{n+j_i-1}{j_i} \, x_i^{j_i}
   \qquad\;
 \label{def.Pn-.x0xm-1}
\end{eqnarray}
for $n \ge 1$.
When $\bfx = \bone$ this equals $C_n^{(m+1)}$.
A combinatorial proof of \reff{def.Pn-.x0xm-1}
was given by M\l{}otkowski \cite[Theorem~1.3]{Mlotkowski_12}.
We call $P^{(m)-}_n(\bfx)$, $Q^{(m)-}_n(\bfx)$ and $Q^{(m)-}_{n,k}(\bfx)$
the \textbfit{multivariate Fuss--Narayana polynomials of negative type}.
Please note that, in contrast to \reff{def.Pn.x0xm.a},
the relation between $P_n^{(m)-}$ and $Q_n^{(m)-}$ is not trivial.

When $m=1$, we have $Q_n^{(1)-}(x_0) = P_n^{(1)-}(x_0) = C_n \, x_0^n$.
When $m=2$, we have
\begin{eqnarray}
   Q_n^{(2)-}(x_0,x_1)
   & = &
   {1 \over n+1}
     \sum_{j=0}^n \binom{n+j}{j} \binom{2n-j}{n-j} \, x_0^j x_1^{n-j}
          \\[2mm]
   P_n^{(2)-}(x_0,x_1)
   & = &
   {1 \over n}
     \sum_{j=0}^n \binom{n+j}{j-1} \binom{2n-j-1}{n-j} \, x_0^j x_1^{n-j}
\end{eqnarray}
To our surprise, these coefficient arrays are not currently in \cite{OEIS}.

Once again we have coefficientwise Hankel-total positivity:
for $P^{(m)-}_n$ this follows from Theorem~\ref{thm.Stype.minors},
while for $Q^{(m)-}_n$ it follows from Theorem~\ref{thm.Stype.minors.extended}
with $\ell = m-1$.
We conclude:

\begin{theorem}[Hankel-total positivity of multivariate Fuss--Narayana polynomials of negative type]
   \label{thm.Qn-Pn-.hankelTP}
For each integer $m \ge 1$,
the sequences $(P^{(m)-}_n)_{n \ge 0}$ and $(Q^{(m)-}_n)_{n \ge 0}$
are coefficientwise Hankel-totally positive,
jointly in the indeterminates $x_0,\ldots,x_{m-1}$.
\end{theorem}

Once again we can pass to the limit $m \to\infty$ and define
the formal power series
\begin{eqnarray}
   Q^{(\infty)-}_{n,k}(\bfx)
   & = &
   {k+1 \over n+1}
   \sum_{\begin{scarray}
            j_0,j_1,\ldots \ge 0 \\
            \sum j_i = n-k
         \end{scarray}}
   \prod_{i=0}^\infty \binom{n+j_i}{j_i} \, x_i^{j_i}
      \label{def.Q-infty}  \\[2mm]
   P^{(\infty)-}_n(\bfx)
   & = &
   {1 \over n}
   \sum_{\begin{scarray}
            j_0,j_1,\ldots \ge 0 \\
            \sum j_i = n
         \end{scarray}}
   \!
   \binom{n+j_0}{j_0 -1} \, x_0^{j_0} \:
   \prod_{i=1}^\infty \binom{n+j_i-1}{j_i} \, x_i^{j_i}
   \qquad
      \label{def.P-infty}
\end{eqnarray}
(and of course $Q^{(\infty)-}_n \eqdef Q^{(\infty)-}_{n,0}$)
where $\bfx = (x_i)_{i \ge 0}$.
Here $Q^{(\infty)-}_{n,k}$ is a symmetric function in $\bfx$,
while $P^{(\infty)-}_n$ is a power series in $x_0$
whose coefficients are symmetric functions of $x_1,x_2,\ldots\;$.
The functional equation \reff{eq.functeqn.G} with $m=\infty$
can be rewritten as
\be
   G(t)  \;=\;  H(tG(t))
\ee
where $H(t)$ is the generating function \reff{def.H}
for the complete homogeneous symmetric functions.
Then Lagrange inversion gives the dual of
\reff{eq.Qinfty.symfun.all} and \reff{eq.Qinfty.symfun.h.all}:
\begin{subeqnarray}
   Q^{(\infty)-}_{n,k}(\bfx)  \:=\: [t^{n-k}] \, G(t)^{k+1}
   & = &
   {k+1 \over n+1} \, [t^{n-k}] \, H(t)^{n+1}
   \;=\;
   {k+1 \over n+1} \, [t^{n-k}] \, E(-t)^{-(n+1)}
       \nonumber \\[-2mm]  \\
   &  &
   \hspace*{-3.5cm} =\;
   \sum_{\lambda \vdash n-k} {k+1 \over n+1}
       \binom{n+1}{n+1-\ell(\lambda),\, m_1(\lambda),\, m_2(\lambda),\, \ldots}
       \; h_\lambda(\bfx)
  \slabel{eq.Q-infty.symfun_h} \\[2mm]
   &  &
   \hspace*{-3.5cm} =\;
   \sum_{\lambda \vdash n-k} {k+1 \over n+1} \,
       (-1)^{n - \ell(\lambda)} \,
       \binom{n+\ell(\lambda)}{n,\, m_1(\lambda),\, m_2(\lambda),\, \ldots}
       \; e_\lambda(\bfx)
   \;.
  \slabel{eq.Q-infty.symfun_e}
  \label{eq.Q-infty.symfun.all}
\end{subeqnarray}
That is, $Q^{(\infty)-}_{n,k} = \omega Q^{(\infty)}_{n,k}$
where $\omega$ is the involution of the ring of symmetric functions
defined by $\omega h_n = e_n$.
Also, the definition \reff{def.Q-infty} can trivially be rewritten as
\be
   Q^{(\infty)-}_{n,k}(\bfx)
   \;=\;
   {k+1 \over n+1}
   \sum_{\lambda \vdash n-k}
   \left[ \prod_i \binom{n+\lambda_i}{n} \right] \, m_\lambda(\bfx)
   \;.
  \label{eq.Q-infty.symfun_m}
\ee
We call $Q^{(\infty)-}_{n,k}(\bfx)$
the \textbfit{Fuss--Narayana symmetric functions of negative type}.

A further application of Lagrange inversion gives for $n \ge 1$
\begin{subeqnarray}
   & &
   \hspace*{-1.2cm}
   P^{(\infty)-}_n(\bfx)  \;=\; [t^n] \, \big( 1 - x_0 tG(t) \big)^{-1}
      \;=\; {1 \over n} \, [t^{n-1}] \, {x_0 \over (1- x_0 t)^2} \, H(t)^n
        \\[2mm]
   & &
   \;=\;
   {1 \over n} \, \sum_{j=1}^n j \, x_0^j \, [t^{n-j}] \, H(t)^n
        \\[2mm]
   & &
   \;=\;
   {1 \over n} \, \sum_{j=1}^n j \, x_0^j
   \!\!\!
   \sum_{\begin{scarray}
           r_1, r_2, \ldots \ge 0 \\
           \sum i r_i = n-j
         \end{scarray}}
   \!\!\!
   \binom{n}{n - \sum r_i ,\, r_1 ,\, r_2 ,\, \ldots}
      \, \prod_{i=1}^\infty h_i(\bfx)^{r_i}
      \qquad
              \\[2mm]
   & &
   \;=\;
   {1 \over n} \, \sum_{j=1}^n j \, x_0^j \,
      \sum_{\lambda \vdash n-j}
      \binom{n}{n - \ell(\lambda) ,\, m_1(\lambda) ,\, m_2(\lambda) ,\, \ldots}
       \; h_\lambda(\bfx)
      \qquad
  \label{eq.P-infty.symfun.all}
\end{subeqnarray}
as well as a similar formula in terms of $e_\lambda$.

Finally, we can again conclude:

\begin{theorem}[Hankel-total positivity of $\infty$-variate Fuss--Narayana series of negative type]
   \label{thm.Q-inftyP-infty.hankelTP}
The sequences $(P^{(\infty)-}_n)_{n \ge 0}$ and $(Q^{(\infty)-}_n)_{n \ge 0}$
are Hankel-totally positive with respect to the coefficientwise order
on the formal-power-series ring $\Z[[\bfx]]$.
\end{theorem}

\smallskip

Return now to the case of finite $m$,
and suppose that the multiset $\{x_0,\ldots,x_{m-1}\}$
consists of elements $\bfy = (y_1,\ldots,y_l)$
with multiplicities $\bfp = (p_1,\ldots,p_l)$
(where of course ${\sum p_i = m}$).
Then the functional equation \reff{eq.functeqn.G} becomes
\be
   G  \;=\;  \prod_{i=1}^l (1 - y_i t G)^{-p_i}
   \;,
 \label{eq.functeqn.G.bis}
\ee
and solving it by Lagrange inversion gives \cite[eq.~(3.4.5)]{Gessel_16}
\be
   G(t)^{k+1}
   \;=\;
   \sum_{n=k}^\infty \widetilde{Q}^{(m)-}_{n,k}(\bfy;\bfp) \: t^{n-k}
\ee
where
\begin{eqnarray}
   \widetilde{Q}^{(m)-}_{n,k}(\bfy;\bfp)
   & = &
   {k+1 \over n+1}
   \!\!
   \sum_{\begin{scarray}
            j_1,\ldots,j_l \ge 0 \\
            j_1+\ldots+j_l = n-k
         \end{scarray}}
   \!\!
   \prod_{i=1}^l \binom{p_i (n+1) + j_i-1}{j_i} \, y_i^{j_i}
   \;.
            \nonumber \\[-6mm]
 \label{def.Q-tilden.y1yl}
\end{eqnarray}
If we make the convention that $x_0 = y_1$,
then a further application of Lagrange inversion gives
\be
   f_0(t)
   \;=\;
   {1 \over 1 \,-\, y_1 t G(t)}
   \;=\;
   \sum_{n=0}^\infty \widetilde{P}^{(m)-}_n(\bfy;\bfp) \: t^n
\ee
where $\widetilde{P}^{(m)-}_0 = 1$ and
\begin{subeqnarray}
   \widetilde{P}^{(m)-}_n(\bfy;\bfp)
   & = &
   {1 \over n}
   \!\!
   \sum_{\begin{scarray}
            j_1,\ldots,j_l \ge 0 \\
            j_1+\ldots+j_l = n
         \end{scarray}}
   \!\!\!\!
   \binom{p_1 n + j_1}{j_1 -1} \, y_1^{j_1}
   \prod_{i=2}^l \binom{p_i n + j_i-1}{j_i} \, y_i^{j_i}
            \nonumber \\[-7mm]
 \slabel{def.P-tilden.y1yl}
\end{subeqnarray}
for $n \ge 1$.


\subsection{Combinatorial interpretations of the multivariate Fuss--Narayana polynomials}
   \label{subsec.periodic.combinatorial}

It is an immediate consequence of the definitions that
$P^{(m)}_n(\bfx)$ is the generating polynomial
for $m$-Dyck paths of length $(m+1)n$ in which
each $m$-fall to a height $i \bmod m+1$ gets weight $x_i$.
Moreover, \reff{eq.Qn.Pn+1} then implies that
$Q^{(m)}_n(\bfx)$ is the generating polynomial
for $m$-Dyck paths of length $(m+1)(n+1)$ in which
each $m$-fall to a height $i \bmod m+1$ gets weight $x_i$,
{\em except for}\/ the last step (an $m$-fall to height 0)
which gets weight 1 rather than $x_0$.
Or equivalently, from \reff{def.F.Qn} and the definition of $F$ we can see that
$Q^{(m)}_n(\bfx)$ is the generating polynomial for partial $m$-Dyck paths
from $(0,0)$ to $((m+1)n+m,m)$, with the usual weights.
Similar interpretations apply to the
multivariate Fuss--Narayana polynomials of negative type,
$P^{(m)-}_n(\bfx)$ and $Q^{(m)-}_n(\bfx)$.

In this subsection we would like to present some other, less trivial,
combinatorial interpretations of the multivariate Fuss--Narayana polynomials,
in terms of ordered trees (of various types), noncrossing partitions,
and parking functions.


\subsubsection[Fuss--Narayana polynomials of positive type: $(m+1)$-ary trees]{Fuss--Narayana polynomials of positive type: $\bm{(m+1)}$-ary trees}

\begin{figure}[!ht]
\begin{center}
\includegraphics[scale=1.8]{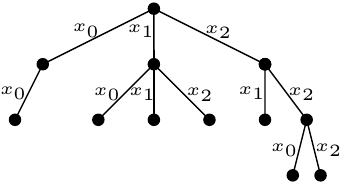}
\caption{\label{fig.positivearytree}
         A $3$-ary tree, with the edge weights shown.}
\end{center}
\end{figure}

Let us recall \cite[p.~295]{Stanley_86}
the recursive definition of an {\em $r$-ary tree}\/ ($r \ge 1$):
it is either empty or else consists of a root together with
an ordered list of $r$ subtrees, each of which is an $r$-ary tree
(which may be empty).  We draw an edge from each vertex to the root
of each of its nonempty subtrees;  an edge from a vertex to the root
of its $i$th subtree (let us number them $0 \le i \le r-1$)
will be called an {\em $i$-edge}\/ (see Figure~\ref{fig.positivearytree} for an example).

Let $f$ be the ordinary generating function for nonempty $r$-ary trees
in which each vertex gets a weight $t$ and each $i$-edge gets a weight $x_i$.
The recursive definition leads immediately to the functional equation
\be
   f  \;=\;  t \prod_{i=0}^{r-1} (1 + x_i f)
   \;,
\ee
which is precisely \reff{eq.functeqn.f}.  Hence:

\begin{proposition}
   \label{prop.fuss-narayana.trees}
The multivariate Fuss--Narayana polynomial $Q^{(m)}_n(\bfx)$
defined in \reff{def.Qn.x0xm} is the generating polynomial for
$(m+1)$-ary trees on $n+1$ vertices
in which each $i$-edge gets a weight $x_i$.

More generally,
the multivariate Fuss--Narayana polynomial $Q^{(m)}_{n,k}(\bfx)$
is the generating polynomial for ordered forests of
$(m+1)$-ary trees on $n+1$ total vertices with $k+1$ components
(hence $n-k$ edges)
in which each $i$-edge gets a weight $x_i$.
\end{proposition}

\noindent
This result (for trees) is due to Cigler \cite[Theorem~1]{Cigler_87},
who gave a combinatorial proof.

Since $P^{(m)}_n(\bfx) = x_0 \, Q^{(m)}_{n-1}(\bfx)$,
we also have:

\begin{proposition}
   \label{prop.fuss-narayana.trees.bis}
For $n \ge 1$,
the multivariate Fuss--Narayana polynomial $P^{(m)}_n(\bfx)$
defined in \reff{def.Pn.x0xm} is $x_0$ times the generating polynomial for
$(m+1)$-ary trees on $n$ vertices in which
each $i$-edge gets a weight $x_i$.
Equivalently, $P^{(m)}_n(\bfx)$ is the generating polynomial for
$(m+1)$-ary trees on $n+1$ vertices in which
the edge (if any) emanating from the root is a 0-edge
and in which each $i$-edge gets a weight $x_i$.
\end{proposition}

This result also has a direct combinatorial proof.
There is a bijection $\varphi$ between
$(m+1)$-ary trees on $n$ unlabeled vertices
and $m$-Dyck paths of length $(m+1)n$,
which can be defined recursively as follows\footnote{
   Generalizing what is done in \cite[p.~11]{Stanley_15} for the binary case.
}:
The empty tree $\emptyset$ maps onto the zero-step path;
if the tree $T$ consists of a root $r$ and subtrees $T_0,\ldots,T_m$,
then $\varphi(T) =
   \varphi(T_0) U \varphi(T_1) U \cdots \varphi(T_{m-1}) U \varphi(T_m) D$
where $U$ is a rise and $D$ is an $m$-fall.
Let us now say that the height $h(v)$ of a vertex $v \in T$
is the sum of the indices $i$ ($0 \le i \le m$)
on the edges connecting it to the root.
Then, corresponding to each vertex $v \in T$,
the path $\varphi(T)$ contains rises starting at heights
$h(v), h(v)+1, \ldots, h(v)+m-1$
and an $m$-fall starting at height $h(v)+m$.
Therefore, under \reff{eq.alpha.periodic} with $p=m+1$,
the $m$-Dyck path will get a weight $\prod_{v \in T} x_{h(v)}$.
This is the right weight for the root ($x_0$)
and for the children of the root ($x_i$ for child $i$)
but not for the remaining vertices.
However, the correct weights can be arranged by a simple bijection
that ``twists'' the tree:
for child $i$ of the root, we permute its subtrees cyclically by $i$;
we then do the same thing successively at lower levels of the tree.
After this bijection, the weights will be $x_0$ for the root
and $x_i$ for each vertex that is child $i$ of its parent,
as asserted in Proposition~\ref{prop.fuss-narayana.trees.bis}.


We can also pass to the limit $m \to\infty$.
Let us define recursively an {\em $\infty$-ary tree}\/:
it is either empty or else consists of a root together with
an ordered list of subtrees indexed by $\N$,
each of which is an $\infty$-ary tree (which may be empty)
{\em and only finitely many of which are nonempty}\/.
Please note that such a graph is necessarily finite
(as always, the recursion is carried out only finitely many times).
We define $i$-edges as before.
We then have:

\begin{proposition}
   \label{prop.fuss-narayana.trees.infty}
The $\infty$-variate Fuss--Narayana series $Q^{(\infty)}_n(\bfx)$
defined in \reff{def.Qinfty} is the generating formal power series for
$\infty$-ary trees on $n+1$ vertices in which
each $i$-edge gets a weight $x_i$.

More generally,
the $\infty$-variate Fuss--Narayana series $Q^{(\infty)}_{n,k}(\bfx)$
is the generating formal power series for ordered forests of
$\infty$-ary trees on $n+1$ total vertices with $k+1$ components
(hence $n-k$ edges)
in which each $i$-edge gets a weight $x_i$.
\end{proposition}

\begin{proposition}
   \label{prop.fuss-narayana.trees.infty.bis}
The $\infty$-variate Fuss--Narayana series $P^{(\infty)}_n(\bfx)$
defined in \reff{def.Pinfty} is $x_0$ times
the generating formal power series for
$\infty$-ary trees on $n$ vertices in which
each $i$-edge gets a weight $x_i$.
Equivalently, $P^{(\infty)}_n(\bfx)$ is the generating polynomial for
$\infty$-ary trees on $n+1$ vertices in which
the edge (if any) emanating from the root is a 0-edge
and in which each $i$-edge gets a weight $x_i$.
\end{proposition}

%

We can also express the foregoing results in terms of symmetric functions,
by viewing $r$-ary and $\infty$-ary trees
from a slightly different point of view.
Recall first
that an {\em ordered tree}\/ (also called {\em plane tree}\/)
is a rooted tree in which the children of each vertex
(or equivalently, the edges emanating outwards from that vertex)
are linearly ordered.
Then an $r$-ary (resp.\ $\infty$-ary) tree is simply an ordered tree
in which each edge carries a label $i \in \{0,\ldots,r-1\}$ (resp.\ $i \in \N$)
and the edges emanating outwards from each vertex consist, in order,
of zero or one edges labeled 0, then zero or one edges labeled 1,
and so forth.
An edge with label $i$ will be called an $i$-edge.
Since the choice of labels on the edges emanating outwards from a vertex $v$
can be made independently for each $v$,
Proposition~\ref{prop.fuss-narayana.trees.infty} can be rephrased as:

\begin{proposition}
   \label{prop.fuss-narayana.trees.infty.symfun}
The $\infty$-variate Fuss--Narayana series $Q^{(\infty)}_n(\bfx)$
defined in \reff{def.Qinfty} is the generating formal power series for
ordered trees on $n+1$ vertices in which
each vertex with $j$ children gets a weight $e_j(\bfx)$,
where $e_j$ is the elementary symmetric function.
\end{proposition}

\noindent
The other propositions have similar rephrasings.

But now we can use the well-known enumeration of
ordered forests of (unlabeled) ordered trees
with a given degree sequence \cite[pp.~30--36]{Stanley_99}.
Given a forest $F$ of rooted trees, we define the {\em type}\/ of $F$
to be the sequence $\bfr = (r_0,r_1,\ldots)$
where $r_i$ vertices of $F$ have out-degree $i$.
If $F$ has $n$~vertices and $k$ components,
then clearly $\sum\limits_{i \ge 0} r_i = n$
and $\sum\limits_{i \ge 0} i r_i = n-k$;
we abbreviate these conditions as $\bfr \to (n,k)$.
Then \cite[Theorem~5.3.10]{Stanley_99}
if~$\bfr$ is a sequence satisfying $\bfr \to (n,k)$,
the number of ordered forests of (unlabeled) ordered trees
on $n$ vertices and $k$ components with type $\bfr$ is
$\displaystyle {k \over n}  \binom{n}{r_0,r_1,\ldots}$.
Combining this with Proposition~\ref{prop.fuss-narayana.trees.infty.symfun},
we conclude that
\be
   Q^{(\infty)}_{n,k}(\bfx)
   \;=\;
   \sum_{\bfr \to (n+1,k+1)} {k+1 \over n+1} \binom{n+1}{r_0,r_1,\ldots}
      \, \prod_{j=0}^\infty e_j(\bfx)^{r_j}
   \;,
 \label{eq.Qinfty.symfun.all.bis}
\ee
recovering \reff{eq.Qinfty.symfun_e.a}.

\subsubsection[Fuss--Narayana polynomials of negative type: Multi-$m$-ary trees]{Fuss--Narayana polynomials of negative type: Multi-$\bm{m}$-ary trees}

\begin{figure}[!ht]
\begin{center}
\includegraphics[scale=1.8]{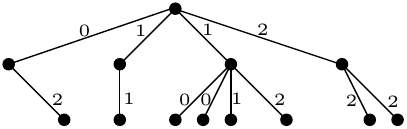}
\caption{\label{fig.multiarytree}
         A multi $3$-ary tree, with the edge labels shown.}
\end{center}
\end{figure}

We can also give a combinatorial interpretation for the
Fuss--Narayana polynomials of negative type.
Let us adopt the reinterpretation of $r$-ary and $\infty$-ary trees
as ordered trees with labeled edges;
and let us then consider the variant in which the number of edges
of each label emanating from a given vertex,
instead of being ``zero or one'', is ``zero or more'':
we call this a \textbfit{multi-$\bm{r}$-ary}
(resp.\ \textbfit{multi-$\bm{\infty}$-ary}) \textbfit{tree}
(see Figure~\ref{fig.multiarytree}).
[In particular, a multi-unary-tree ($r=1$) is simply an ordered tree.]

Let $g$ be the ordinary generating function
for nonempty multi-$r$-ary trees
in which each vertex gets a weight~$t$
and each $i$-edge gets a weight $x_i$.
The definition leads immediately to the functional equation
\be
   g  \;=\;  t \prod_{i=0}^{r-1} (1 - x_i g)^{-1}
   \;.
\ee
which is precisely \reff{eq.functeqn.g}.  Hence:

\begin{proposition}
   \label{prop.fuss-narayana.trees_extended}
The multivariate Fuss--Narayana polynomial of negative type
$Q^{(m)-}_n(\bfx)$ defined in \reff{def.Qn-.x0xm-1}
is the generating polynomial for multi-$m$-ary trees on $n+1$ vertices
in which each $i$-edge gets a weight $x_i$.

More generally,
the multivariate Fuss--Narayana polynomial of negative type
$Q^{(m)-}_{n,k}(\bfx)$
is the generating polynomial for ordered forests of multi-$m$-ary trees
on $n+1$ total vertices with $k+1$ components (hence $n-k$ edges)
in which each $i$-edge gets a weight $x_i$.
\end{proposition}

\begin{proposition}
   \label{prop.fuss-narayana.trees_extended.infty}
The $\infty$-variate Fuss--Narayana series of negative type
$Q^{(\infty)-}_n(\bfx)$ defined in \reff{def.Q-infty}
is the generating formal power series
for multi-$\infty$-ary trees on $n+1$ vertices
in which each $i$-edge gets a weight $x_i$.
\end{proposition}

Applying \reff{def.f0.Pn-} we then get:

\begin{proposition}
   \label{prop.fuss-narayana.trees_extended.P}
The multivariate Fuss--Narayana polynomial of negative type
$P^{(m)-}_n(\bfx)$ defined in \reff{def.Pn-.x0xm-1}
is the generating polynomial for multi-$m$-ary trees on $n+1$ vertices
in which all the edges emanating from the root are 0-edges
and in which each $i$-edge gets a weight $x_i$.
Equivalently, $P^{(m)-}_n(\bfx)$
is the generating polynomial for ordered forests of multi-$m$-ary trees
on $n$ vertices
in which each $i$-edge gets a weight $x_i$
and each component gets a weight $x_0$.
\end{proposition}

\begin{proposition}
   \label{prop.fuss-narayana.trees_extended.infty.P}
The $\infty$-variate Fuss--Narayana series of negative type
\linebreak
$P^{(\infty)-}_n(\bfx)$ defined in \reff{def.P-infty}
is the generating formal power series for multi-$\infty$-ary trees
on $n+1$ vertices
in which all the edges emanating from the root are 0-edges
and in which each $i$-edge gets a weight $x_i$.
Equivalently, $P^{(\infty)-}_n(\bfx)$
is the generating formal power series for ordered forests of
multi-$\infty$-ary trees on $n$ vertices
in which each $i$-edge gets a weight $x_i$
and each component gets a weight $x_0$.
\end{proposition}

Since the choice of labels on the edges emanating outwards from a vertex $v$
can be made independently for each $v$,
Propositions~\ref{prop.fuss-narayana.trees_extended.infty}
and \ref{prop.fuss-narayana.trees_extended.infty.P}
can be rephrased in terms of symmetric functions as:

\begin{proposition}
   \label{prop.fuss-narayana.trees_extended.infty.symfun}
The $\infty$-variate Fuss--Narayana series of negative type
$Q^{(\infty)-}_n(\bfx)$ defined in \reff{def.Q-infty}
is the generating formal power series
for ordered trees on $n+1$ vertices in which
each vertex with $j$ children gets a weight $h_j(\bfx)$,
where $h_j$ is the complete homogeneous symmetric function.

More generally,
the $\infty$-variate Fuss--Narayana series of negative type
$Q^{(\infty)-}_{n,k}(\bfx)$
is the generating formal power series
for ordered forests of ordered trees
on $n+1$ total vertices with $k+1$ components
in which each vertex with $j$ children gets a weight $h_j(\bfx)$.
\end{proposition}

\begin{proposition}
   \label{prop.fuss-narayana.trees_extended.infty.symfun.P}
The $\infty$-variate Fuss--Narayana series of negative type
$P^{(\infty)-}_n(\bfx)$ defined in \reff{def.P-infty}
is the generating formal power series
for ordered trees on $n+1$ vertices in which
each vertex with $j$ children gets a weight $h_j(\bfx)$,
except for the root, which gets a weight $x_0^k$.
\end{proposition}

By an argument similar to that leading to \reff{eq.Qinfty.symfun.all.bis},
we can conclude that
\be
   Q^{(\infty)-}_{n,k}(\bfx)
   \;=\;
   \sum_{\bfr \to (n+1,k+1)} {k+1 \over n+1} \binom{n+1}{r_0,r_1,\ldots}
      \, \prod_{j=0}^\infty h_j(\bfx)^{r_j}
   \;,
 \label{eq.Q-infty.symfun.all.bis}
\ee
recovering \reff{eq.Q-infty.symfun_h}.


\subsubsection{Noncrossing partitions}

The multivariate Fuss--Narayana polynomials and symmetric functions
also have interpretations in terms of noncrossing partitions.
Let us recall \cite{Simion_00a} that a partition $\pi = \{B_1,\ldots,B_k\}$
of the set $[n] \eqdef \{1,\ldots,n\}$ is called {\em noncrossing}\/
if there do not exist four elements $a < b < c < d$
such that $a,c \in B_i$, $b,d \in B_j$ and $i \ne j$.
We denote by $\NC_n$ the set of all noncrossing partitions of $[n]$;
its cardinality is the Catalan number $C_n$.

If $\pi$ is a partition of the set $[n]$,
we define the {\em type} of $\pi$
to be the list of the sizes of the blocks of $\pi$,
written in weakly decreasing order;
it is clearly a partition $\lambda(\pi)$ of the integer $n$.
Kreweras \cite[Th\'eor\`eme~4]{Kreweras_72} showed that,
for any $\lambda \vdash n$,
the number of noncrossing partitions of $[n]$ of type $\lambda$ is
\be
   {n (n-1) \cdots (n-\ell(\lambda)+2)
    \over
    m_1(\lambda)! \, m_2(\lambda)! \,\cdots\, m_n(\lambda)!
   }
   \;.
\ee
Comparing this with \reff{eq.Qinfty.symfun_e} and \reff{eq.Q-infty.symfun_h}
specialized to $k=0$, we see immediately that:

\begin{proposition}
   \label{prop.kreweras}
The multivariate Fuss--Narayana symmetric functions
have the interpretations
\begin{eqnarray}
   Q^{(\infty)}_n(\bfx)
   & = &
   \sum_{\pi \in \NC_n}  e_{\lambda(\pi)}(\bfx)
           \\[2mm]
   Q^{(\infty)-}_n(\bfx)
   & = &
   \sum_{\pi \in \NC_n}  h_{\lambda(\pi)}(\bfx)
\end{eqnarray}
In other words, $Q^{(\infty)}_n$ (resp.\ $Q^{(\infty)-}_n$)
is the generating polynomial for noncrossing partitions of $[n]$
in which each block of size $i$ gets a weight $e_i$ (resp.\ $h_i$).
\end{proposition}

\noindent
See Stanley \cite[Proposition~2.4]{Stanley_97}
for a second proof of Proposition~\ref{prop.kreweras},
based on a noncrossing analogue of the exponential formula
due to Speicher \cite{Speicher_94}.

The multivariate Fuss--Narayana polynomial $Q^{(m)}_n(\bfx)$
also has an interpretation as a generating polynomial
for multichains in the lattice of noncrossing partitions.
Recall that the set $\NC_n$ becomes a partially ordered set (poset)\footnote{
   See e.g.\ Stanley \cite[Chapter~3]{Stanley_86}
   for the basic definitions concerning posets.
}
when the partitions are ordered by refinement:
$\pi \le \pi'$ if every block of $\pi$ is contained in a block of $\pi'$.
This poset has a minimum element $\hat{0}$ (the partition into singletons)
and a maximum element $\hat{1}$ (the partition with a single block).
In fact, the poset $\NC_n$ is a lattice.
Moreover, it is graded of rank $n-1$:
the rank function is $\rk(\pi) = n - |\pi|$
where $|\pi|$ is the number of blocks in $\pi$.

Recall now that a {\em chain}\/ (resp.\ {\em multichain}\/)
of length $\ell$ in a poset $P$
is a sequence $\xi_0 < \xi_1 < \ldots < \xi_\ell$
(resp.\ $\xi_0 \le \xi_1 \le \ldots \le \xi_\ell$) in $P$.
Edelman \cite[Theorem~3.2]{Edelman_80} enumerated, by bijective arguments,
the chains $\xi_0 < \xi_1 < \ldots < \xi_\ell$ in $\NC_n$
where the elements of the chain have specified ranks
$r_0 < r_1 < \ldots < r_\ell$.  This enumeration is equivalent to:

\begin{proposition}
   \label{prop.edelman}
The multivariate Fuss--Narayana polynomial $Q^{(m)}_n(\bfx)$
defined in \reff{def.Qn.x0xm}
is the rank generating polynomial
for multichains of length $m-1$ in the lattice $\NC_{n+1}$:
\be
   Q^{(m)}_n(\bfx)
   \;=\;
   \!\!\!\!\!\!
   \sum_{\begin{scarray}
            \pi_0,\ldots,\pi_{m-1} \in \NC_{n+1} \\
            \pi_0 \le \ldots \le \pi_{m-1}
         \end{scarray}
        }
   \!\!\!\!\!
   x_0^{\rk(\pi_0)} x_1^{\rk(\pi_1) - \rk(\pi_0)}
     \,\cdots\,
   x_{m-1}^{\rk(\pi_{m-1}) - \rk(\pi_{m-2})}
   x_m^{n- \rk(\pi_{m-1})}
   \;.
 \label{eq.prop.edelman}
\ee
\end{proposition}

\noindent
Specializing to $\bfx = \bone$,
it follows that the number of multichains of length $m-1$ in $\NC_{n+1}$
is the Fuss--Catalan number $C_{n+1}^{(m+1)}$
\cite[Proposition~9.35]{Nica_06}.

\bigskip

{\bf Remarks.}
1.  As Stanley \cite{Stanley_97} observes,
it is curious that Proposition~\ref{prop.kreweras} refers to $\NC_n$,
while Proposition~\ref{prop.edelman} refers to $\NC_{n+1}$.

2.  Armstrong \cite{Armstrong_09} gives several combinatorial interpretations
of the Fuss--Narayana numbers and their generalizations to Coxeter groups.
\myendremark

\subsubsection{Parking functions}  \label{subsubsec.periodic.parking}

In a classic paper, Stanley \cite{Stanley_97}
has shown a connection between the Fuss--Narayana symmetric functions
(though he does not use this name), parking functions,
and noncrossing partitions.
Here we mention only his first result \cite[Proposition~2.2]{Stanley_97}
and refer to the original paper for many additional developments.
See also the very recent paper \cite{Stanley_18}
for some interesting generalizations.

A {\em parking function}\/ \cite{Yan_15} of length $n$
is a sequence $(a_1,\ldots,a_n)$ of positive integers
such that its increasing rearrangement $(b_1,\ldots,b_n)$
satisfies $b_i \le i$.\footnote{
   We follow Stanley \cite{Stanley_97} in using the convention in which
   $1 \le a_i \le n$;
   many other authors (e.g.\ \cite{Yan_15}) use the convention
   $0 \le a_i \le n-1$.
}
Let $\scrp_n$ denote the set of parking functions of length $n$.
The symmetric group $\Sym_n$ acts on $\scrp_n$ by permuting coordinates.
For each partition $\lambda$ of the integer $n$,
let $\tau_{\lambda,n}$ be the multiplicity of the irreducible character
of $\Sym_n$ indexed by $\lambda$ in the action of $\Sym_n$ on $\scrp_n$.
Then define the {\em parking-function symmetric function}\/ $\PF_n$ by
\be
   \PF_n  \;=\;  \sum_{\lambda \vdash n} \tau_{\lambda,n} \, s_\lambda
 \label{def.PFn}
\ee
where $s_\lambda$ are the Schur functions.
(This is the Frobenius characteristic of the action of $\Sym_n$ on $\scrp_n$:
 see \cite[section~7.18]{Stanley_99}.)
Then Stanley \cite{Stanley_97} shows (following Haiman \cite{Haiman_94}) that
\be
   \PF_n  \;=\;  {1 \over n+1} \, [t^n] \, H(t)^{n+1}
\ee
(see also \cite[Section~3.2]{Rattan_14}).
Comparing this with \reff{eq.Q-infty.symfun.all},
we see that $\PF_n = Q_n^{(\infty)-}$.
Applying the involution $\omega$ defined by $\omega h_n = e_n$,
we see also that
\be
   \omega \PF_n  \;=\;  {1 \over n+1} \, [t^n] \, E(t)^{n+1}
\ee
and hence by \reff{eq.Qinfty.symfun.all}
that $\omega\PF_n = Q_n^{(\infty)}$.
Stanley's parking-function symmetric functions are thus identical
to our Fuss--Narayana symmetric functions.

The equations \reff{eq.Qinfty.symfun_m}, 
\reff{eq.Q-infty.symfun_h} and \reff{eq.Q-infty.symfun_m}
appear in \cite[eqns.~(6), (5) and~(4)]{Stanley_97};
and \reff{eq.Q-infty.symfun_e} appears in \cite[eqn.~(1.1)]{Stanley_18}.
Stanley also gives formulae for $\PF_n$
in terms of the power-sum symmetric functions $p_\lambda$
and the Schur functions $s_\lambda$ \cite[eqns.~(2) and~(3)]{Stanley_97}.
See \cite{Stanley_97} for many further connections
between $\PF_n$, $\omega \PF_n$ and noncrossing partitions.

\subsection[Period $m\!+\!1$: Production-matrix method]{Period $\bm{m\!+\!1}$: Production-matrix method}
   \label{subsec.periodic.m+1.prodmat}

Consider again the case of period~$p = m+1$,
with $\alpha_i = x_{i+1 \bmod m+1}$.
Then Proposition~\ref{prop.prod.Sm.Jm.Tm}(b) tells us that
the production matrix for the lower-triangular matrix $\sfS^{(m)}(\balpha)$
of generalized $m$-Stieltjes--Rogers polynomials $S_{n,k}^{(m)}(\balpha)$
is the matrix $P^{(m)\mathrm{S}}(\balpha)$
defined by \reff{eq.prop.contraction},
which in the present case is
\begin{eqnarray}
   P^{(m)}(x_1,\ldots,x_m;x_0)
   & \eqdef &
   L(x_1,x_1,x_1,\ldots)
   \:
   L(x_2,x_2,x_2,\ldots)
   \:\cdots\:
   \hspace*{1cm}
       \nonumber \\
   & & \qquad
   L(x_m,x_m,x_m,\ldots)
   \:
   U^\star(x_0,x_0,x_0,\ldots)
   \;.
   \hspace*{1cm}
 \label{eq.prop.contraction.periodic.m+1}
\end{eqnarray}
The elements of this matrix are given by a simple explicit formula:

\begin{lemma}
   \label{lemma.periodic.m+1.prodmat}
The matrix $P^{(m)}(x_1,\ldots,x_m;x_0)$ defined in
\reff{eq.prop.contraction.periodic.m+1} has entries
\be
   P^{(m)}(x_1,\ldots,x_m;x_0)_{ij}
   \;=\;
   \begin{cases}
       e_{i-j+1}(x_0,x_1,\ldots,x_m)   & \textrm{if $j \ge 1$} \\[1mm]
       x_0 \, e_i(x_1,\ldots,x_m)      & \textrm{if $j = 0$}
   \end{cases}
 \label{eq.lemma.periodic.m+1.prodmat}
\ee
where $e_i$ are the elementary symmetric functions
(of course $e_i \eqdef 0$ for $i < 0$).
\end{lemma}

Let us in fact prove, for use in the next subsection,
a slightly more general formula in which
the first entry $x_0$ in the matrix $U^\star$
is replaced by a different value $y$:

\begin{lemma}
   \label{lemma.periodic.m+1.prodmat.bis}
The matrix
\begin{eqnarray}
   P^{(m)}(x_1,\ldots,x_m;y,x_0)
   & \eqdef &
   L(x_1,x_1,x_1,\ldots)
   \:
   L(x_2,x_2,x_2,\ldots)
   \:\cdots\:
   \hspace*{1cm}
       \nonumber \\
   & & \qquad
   L(x_m,x_m,x_m,\ldots)
   \:
   U^\star(y,x_0,x_0,\ldots)
   \hspace*{1cm}
 \label{eq.prop.contraction.periodic.m+1.bis}
\end{eqnarray}
has entries
\be
   P^{(m)}(x_1,\ldots,x_m;y,x_0)_{ij}
   \;=\;
   \begin{cases}
       e_{i-j+1}(x_0,x_1,\ldots,x_m)   & \textrm{if $j \ge 1$} \\[1mm]
       y \, e_i(x_1,\ldots,x_m)        & \textrm{if $j = 0$}
   \end{cases}
 \label{eq.lemma.periodic.m+1.prodmat.bis}
\ee
\end{lemma}

\proof
By induction on $m$.
When $m=0$ we have $P^{(0)}(\;;x_0) = U^\star(y,x_0,x_0,\ldots)$
and \reff{eq.lemma.periodic.m+1.prodmat.bis} holds.
When $m \ge 1$, we have by definition
\be
   P^{(m)}(x_1,\ldots,x_m;y,x_0) 
   \;=\;
   L(x_1,x_1,x_1,\ldots) \, P^{(m-1)}(x_2,\ldots,x_m;y,x_0)
   \;.
\ee
Then for all $i \ge 0$ and $j \ge 1$ we have by the inductive hypothesis
\begin{subeqnarray}
   P^{(m)}(x_1,\ldots,x_m;y,x_0)_{ij}
   & = &
   \sum_{\ell \ge 0} L(x_1,x_1,x_1,\ldots)_{i\ell} \,
                     P^{(m-1)}(x_2,\ldots,x_m;y,x_0)_{\ell j}
       \\[2mm]
   & = &
   x_1 \, e_{i-j}(x_0,x_2,\ldots,x_m)  \:+\: e_{i-j+1}(x_0,x_2,\ldots,x_m)
       \qquad\qquad \\[2mm]
   & = &
   e_{i-j+1}(x_0,x_1,\ldots,x_m)
   \;.
\end{subeqnarray}
Similarly, for all $i \ge 0$ we have
\begin{subeqnarray}
   P^{(m)}(x_1,\ldots,x_m;y,x_0)_{i0}
   & = &
   \sum_{\ell \ge 0} L(x_1,x_1,x_1,\ldots)_{i\ell} \,
                     P^{(m-1)}(x_2,\ldots,x_m;y,x_0)_{\ell 0}
       \qquad\qquad \\[2mm]
   & = &
   x_1 \, y \, e_{i-1}(x_2,\ldots,x_m)  \:+\: y \, e_i(x_2,\ldots,x_m)
       \qquad \\[2mm]
   & = &
   y \, e_i(x_1,\ldots,x_m)
   \;.
\end{subeqnarray}
\qed

Thus, the production matrix $P^{(m)}(x_1,\ldots,x_m;x_0)$
is a lower-Hessenberg matrix that is Toeplitz except for the zeroth column.
Let us give such matrices a name:
an \textbfit{AZ-matrix} is a lower-Hessenberg matrix of the form
\be
   \AZ(\ba,\bz)  \;=\;
   \begin{bmatrix}
      z_0 & a_0 &     &     &     &    \\
      z_1 & a_1 & a_0 &     &     &    \\
      z_2 & a_2 & a_1 & a_0 &     &    \\
      z_3 & a_3 & a_2 & a_1 & a_0 &    \\
      \vdots & \vdots  & \vdots & \vdots & \vdots    & \ddots
   \end{bmatrix}
 \label{def.AZ}
\ee
or in other words
\be
   \AZ(\ba,\bz)_{ij}
           \;=\;  \begin{cases}
                      z_i         &  \textrm{if $j=0$} \\
                      a_{i+1-j}   &  \textrm{if $1 \le j \le i+1$} \\
                      0           &  \textrm{if $j > i+1$}
                  \end{cases}
\ee
Here $\ba = (a_n)_{n \ge 0}$ and $\bz = (z_n)_{n \ge 0}$
are sequences in a commutative ring $R$;
they are called, respectively,
the {\em A-sequence}\/ and {\em Z-sequence}\/ of the matrix $\AZ(\ba,\bz)$.
The ordinary generating functions
$A(t) = \sum_{n=0}^\infty a_n t^n$ and $Z(t) = \sum_{n=0}^\infty z_n t^n$
are called the {\em A-series}\/ and {\em Z-series}\/.
We also write $\AZ(A,Z)$ as a synonym for $\AZ(\ba,\bz)$.

In our case we have
\begin{subeqnarray}
   a_i  & = &  e_i(x_0,x_1,\ldots,x_m)
        \;=\;  e_i(x_1,\ldots,x_m) \:+\: x_0 \, e_{i-1}(x_1,\ldots,x_m)
       \qquad \\[2mm]
   z_i  & = &  x_0 \, e_i(x_1,\ldots,x_m)
\end{subeqnarray}
or equivalently
\begin{subeqnarray}
   A(t)  & = &  E(t;x_0,x_1,\ldots,x_m)  \;=\;
                (1+ x_0 t) \, E(t;x_1,\ldots,x_m)   \qquad \\[2mm]
   Z(t)  & = &  x_0 \, E(t;x_1,\ldots,x_m)
                \;=\;  {x_0 \over 1 + x_0 t} \, A(t)
 \label{eq.periodm+1.AZ}
\end{subeqnarray}

AZ-matrices arise as the production matrices of Riordan arrays.
Let us review this connection briefly, as we shall make use of it.
Let $R$ be a commutative ring,
and let $f(t) = \sum_{n=0}^\infty f_n t^n$
and $g(t) = \sum_{n=1}^\infty g_n t^n$ be formal power series
with coefficients in $R$; note that $g$ has zero constant term.
Then the \textbfit{Riordan array} \cite{Shapiro_91,Sprugnoli_94,Barry_16}
associated to the pair $(f,g)$
is the infinite lower-triangular matrix
$\scrr(f,g) = (\scrr(f,g)_{nk})_{n,k \ge 0}$ defined by
\be
   \scrr(f,g)_{nk}
   \;=\;
   [t^n] \, f(t) g(t)^k
   \;.
 \label{def.riordan}
\ee
We call a Riordan array {\em normalized}\/ if $f_0 = g_1 = 1$,
or equivalently if all the diagonal elements $\scrr(f,g)_{nn} = f_0 g_1^n$
are equal to 1.
We then have the following well-known characterization
of Riordan arrays by their $A$- and $Z$-sequences
\cite{Deutsch_09,He_09,He_15,Barry_16}:

\begin{proposition}[Production matrices of Riordan arrays]
   \label{prop.riordan.production}
Let $L$ be a lower-triangular matrix (with entries in a commutative ring $R$)
with invertible diagonal entries,
and let $P = L^{-1} \Delta L$ be its production matrix.
Then $L$ is a Riordan array if and only~if $P$ is an AZ-matrix.

More precisely, $L = \scrr(f,g)$ if and only~if $P = \AZ(A,Z)$,
where the generating functions $\big( f(t), g(t) \big)$
and $\big( A(t), Z(t) \big)$ are connected by
\be
   g(t) \;=\; t \, A(g(t))  \;,\qquad
   f(t)  \;=\;  {f_0 \over 1 \,-\, t Z(g(t))}
   \;.
 \label{eq.prop.riordan.production}
\ee
\end{proposition}

\noindent
Proofs can be found in \cite{Deutsch_09,He_15,Barry_16,Sokal_totalpos}.

\bigskip

So, in the periodic case with period~$m+1$,
the matrix $\sfS^{(m)}(\balpha)$
of generalized $m$-Stieltjes--Rogers polynomials is a Riordan array
$\scrr(f,g)$.
We have $Z(t) = {x_0 \over 1+ x_0 t} A(t)$ and $f_0 = 1$,
hence $f(t) = 1 + x_0 g(t)$.
Using Lagrange inversion to solve \reff{eq.prop.riordan.production},
we find that
$P_{n,k}^{(m)}(\bfx) \eqdef S_{n,k}^{(m)}(\balpha)$
are given for $n>0$ by
\begin{subeqnarray}
   P_{n,k}^{(m)}(\bfx)
   & \eqdef &
   [t^n] \, f(t) \, g(t)^k
       \\[2mm]
   & = &
   [t^n] \, g(t)^k  \;+\;  x_0 \, [t^n] \, g(t)^{k+1}
       \\[2mm]
   & = &
   {k \over n} \, [t^{n-k}] \, A(t)^n
   \;+\;
   x_0 \, {k+1 \over n} \, [t^{n-k-1}] \, A(t)^{n+1}
       \\[2mm]
   & = &
   Q_{n-1,k-1}^{(m)}(\bfx)  \:+\: x_0 \, Q_{n-1,k}^{(m)}(\bfx)
\end{subeqnarray}
where $A(t) = E(t;x_0,x_1,\ldots,x_m)$
and $Q_{n,k}^{(m)}(\bfx)$ is given by \reff{def.Qn.x0xm}
[cf.\ also \reff{eq.Qinfty.symfun.all}].
When $k=0$ this agrees with \reff{def.Pn.x0xm}
[cf.\ also \reff{def.Pinfty}].
In conclusion, we have
\be
   P_{n,k}^{(m)}(\bfx)
   \;=\;
   \begin{cases}
         \delta_{k0}    & \textrm{if $n=0$}  \\[2mm]
         Q_{n-1,k-1}^{(m)}(\bfx)  \:+\: x_0 \, Q_{n-1,k}^{(m)}(\bfx)
                        & \textrm{if $n \ge 1$}
   \end{cases}
\ee

\subsection[Period $m\!+\!1$: Row-generating polynomials]{Period $\bm{m\!+\!1}$: Row-generating polynomials}
   \label{subsec.periodic.row-generating}

Having computed the lower-triangular matrix
$\sfS^{(m)} = (S_{n,k}^{(m)}(\balpha))_{n,k \ge 0}$
of generalized $m$-Stieltjes--Rogers polynomials
for the case of period~$m+1$,
let us now consider its row-generating polynomials
\be
   S_n^{(m)}(\balpha;\xi)  \;=\;  \sum_{k=0}^n S_{n,k}^{(m)}(\balpha) \, \xi^k
\ee
and more generally its row-generating matrix
\be
   S^{(m)}_{n,k}(\balpha;\xi)
   \;\eqdef\;
   \sum_{\ell=k}^n S^{(m)}_{n,\ell}(\balpha) \: \xi^{\ell-k}
 \label{def.Snk.rowgen.bis}
\ee
[cf.\ \reff{def.Sn.rowgen}/\reff{def.Snk.rowgen}],
where $\xi$ is an indeterminate.
Please recall that the definition \reff{def.Snk.rowgen.bis}
can be written in matrix form as
\be
   \sfS^{(m)}(\xi)  \;=\;  \sfS^{(m)} \, T_\xi
 \label{eq.Sm.Smxi}
\ee
where $T_\xi$ is the lower-triangular Toeplitz matrix of powers of $\xi$
[cf.\ \reff{def.Txi}].
We write $P_n^{(m)}(\bfx;\xi)$
for the row-generating polynomial $S_n^{(m)}(\balpha;\xi)$
specialized to the weights $\alpha_i = x_{i+1 \bmod m+1}$.

We have just shown that the production matrix for $\sfS^{(m)}$
for the case of period~$m+1$
is the matrix $P^{(m)}(x_1,\ldots,x_m;x_0)$ defined in
\reff{eq.prop.contraction.periodic.m+1}/\reff{eq.lemma.periodic.m+1.prodmat}.
It then follows immediately from \reff{eq.Sm.Smxi} 
and Lemma~\ref{lemma.production.AB}
that the production matrix for $\sfS^{(m)}(\xi)$
is $T_\xi^{-1} P^{(m)}(x_1,\ldots,x_m;x_0) \, T_\xi$.
It turns out that this latter matrix has a very simple form,
as a result of the following simple identities:

\begin{lemma}
   \label{lemma.Txi.L.U}
\hfill\break\vspace*{-7mm}
\begin{itemize}
   \item[(a)] Let $L(s_1,s_2,\ldots)$ be the lower-bidiagonal matrix
defined by \reff{def.L}.  Then
\be
   T_\xi^{-1} \, L(s_1,s_2,\ldots) \, T_\xi
   \;=\;
\Scale[0.95]{
   \begin{bmatrix}
      1              &                &               &     &      & \\[0.5mm]
      s_1            & 1              &               &     &      & \\[0.5mm]
      (s_2-s_1)\xi   & s_2            & 1             &      &     & \\[0.5mm]
      (s_3-s_2)\xi^2 & (s_3-s_2)\xi   & s_3           & 1    &     & \\[0.5mm]
      (s_4-s_3)\xi^3 & (s_4-s_3)\xi^2 & (s_4-s_3)\xi  & \quad s_4 \quad  & 1   & \\[0.5mm]
      \vdots         & \vdots         & \vdots        & \vdots & \ddots & \ddots
   \end{bmatrix}
}
\ee
or in other words
\be
   [T_\xi^{-1} \, L(s_1,s_2,\ldots) \, T_\xi]_{ij}
   \;=\;
   \begin{cases}
       1                          & \textrm{if $i=j$}    \\
       s_i                        & \textrm{if $i=j+1$}  \\
       (s_i - s_{i-1})\xi^{i-j-1} & \textrm{if $i > j+1$} \\
       0                          & \textrm{otherwise}
   \end{cases}
   \label{eq.lemma.Txi.L.U.a}
\ee
   \item[(b)] Let $U^\star(s_1,s_2,\ldots)$ be the upper-bidiagonal matrix
defined by \reff{def.Ustar}.  Then
\be
   T_\xi^{-1} \, U^\star(s_1,s_2,\ldots) \, T_\xi
   \;=\;
\Scale[0.95]{
   \begin{bmatrix}
      s_1 + \xi      & 1              &               &      &     & \\[0.5mm]
      (s_2-s_1)\xi   & s_2            & 1             &      &     & \\[0.5mm]
      (s_3-s_2)\xi^2 & (s_3-s_2)\xi   & s_3           & 1    &     & \\[0.5mm]
      (s_4-s_3)\xi^3 & (s_4-s_3)\xi^2 & (s_4-s_3)\xi  & \quad s_4 \quad  & 1   & \\[0.5mm]
      \vdots         & \vdots         & \vdots        & \vdots & \ddots & \ddots
   \end{bmatrix}
}
\ee
or in other words
\be
   [T_\xi^{-1} \, U^\star(s_1,s_2,\ldots) \, T_\xi]_{ij}
   \;=\;
   \begin{cases}
       1                          & \textrm{if $i=j-1$}    \\
       s_1 + \xi                  & \textrm{if $i=j=0$}    \\
       s_{i+1}                    & \textrm{if $i=j \ge 1$}  \\
       (s_{i+1} - s_i)\xi^{i-j}   & \textrm{if $i > j$}  \\
       0                          & \textrm{otherwise}
   \end{cases}
   \label{eq.lemma.Txi.L.U.b}
\ee
\end{itemize}
\end{lemma}

\proof
(a) We have $T_\xi^{-1} = L(-\xi,-\xi,\ldots)$ and hence
\be
   [T_\xi^{-1} \, L(s_1,s_2,\ldots)]_{ij}
   \;=\;
   \begin{cases}
       1                    & \textrm{if $i=j$}    \\
       s_i - \xi            & \textrm{if $i=j+1$}  \\
       -s_{i-1} \xi         & \textrm{if $i=j+2$} \\
       0                    & \textrm{otherwise}
   \end{cases}
\ee
Then
\begin{eqnarray}
   [T_\xi^{-1} \, L(s_1,s_2,\ldots) \, T_\xi]_{ik}
   & = &
   \xi^{i-k}  \, \mathrm{I}[k \le i]   \;+\;
   (s_i - \xi) \, \xi^{i-1-k} \, \mathrm{I}[k \le i-1]
           \qquad \nonumber \\
   & & \hspace{2.2cm}
   \;+\;
   (-s_{i-1} \xi) \, \xi^{i-2-k} \, \mathrm{I}[k \le i-2]
   \;,
   \qquad\qquad
\end{eqnarray}
where $\mathrm{I}[\hbox{\sl statement}]$ equals 1 if {\sl statement}
is true and 0 if it is false;
this coincides with \reff{eq.lemma.Txi.L.U.a}.

Part~(b) is proven similarly.
\qed

Specializing this lemma to $s_1 = s_2 = \ldots = s$, we have
\be
   T_\xi^{-1} \, L(s,s,s,\ldots) \, T_\xi
   \;=\;
   L(s,s,s,\ldots)
\ee
(i.e.\ $T_\xi$ and $L(s,s,s,\ldots)$ commute,
 because both are lower-triangular Toeplitz matrices)
and
\be
   T_\xi^{-1} \, U^\star(s,s,s,\ldots) \, T_\xi
   \;=\;
   U^\star(s+\xi,s,s,\ldots)
   \;.
\ee
Using the definition \reff{eq.prop.contraction.periodic.m+1}
of $P^{(m)}(x_1,\ldots,x_m;x_0)$, we see that
\begin{subeqnarray}
   & &
   \!\!\!\!
   T_\xi^{-1} \, P^{(m)}(x_1,\ldots,x_m;x_0) \, T_\xi
       \nonumber \\[2mm]
   & &
   \;=\;
   L(x_1,x_1,x_1,\ldots)
   \:
   L(x_2,x_2,x_2,\ldots)
   \:\cdots\:
   L(x_m,x_m,x_m,\ldots)
   \:
   U^\star(x_0 + \xi,x_0,x_0,\ldots)
       \nonumber \\ \slabel{eq.prodmat.periodic.rowgen.a} \\
   & & 
   \;=\;
   P^{(m)}(x_1,\ldots,x_m;x_0+\xi,x_0)
 \label{eq.prodmat.periodic.rowgen}
\end{subeqnarray}
as defined in \reff{eq.prop.contraction.periodic.m+1.bis}.
By Lemma~\ref{lemma.bidiagonal} this matrix
is coefficientwise totally positive,
jointly in the indeterminates $x_0,\ldots,x_m$ and $\xi$.
From Theorem~\ref{thm.iteration2bis} we therefore conclude:

\begin{theorem}[Hankel-total positivity of row-generating polynomials of multivariate Fuss--Narayana polynomials]
   \label{thm.Pn.hankelTP.rowgen}
For each integer $m \ge 1$,
the sequence $(P^{(m)}_n(\bfx;\xi))_{n \ge 0}$
is coefficientwise Hankel-totally positive,
jointly in the indeterminates $x_0,\ldots,x_m$ and $\xi$.
\end{theorem}

But we can say more:
by \reff{eq.prodmat.periodic.rowgen.a} we see that
the matrix ${P^{(m)}(x_1,\ldots,x_m;x_0+\xi,x_0)}$
is nothing other than the matrix $P^{(m)\mathrm{S}}(\balpha)$
defined by \reff{eq.prop.contraction}
where $\alpha_i = x_{i+1 \bmod m+1}$
{\em except that}\/ $\alpha_m = x_0 + \xi$.
Calling these new weights $\balpha'$,
we see that
\be
   P^{(m)}_n(\bfx;\xi)  \;=\;  S_n^{(m)}(\balpha')
   \;.
\ee
So the row-generating polynomials $P^{(m)}_n(\bfx;\xi)$
are themselves given by an $m$-branched S-fraction
that is a small modification of the one for 
$P^{(m)}_n(\bfx) = P^{(m)}_n(\bfx;0)$.

The coefficientwise total positivity of the production matrix
\reff{eq.prodmat.periodic.rowgen} also implies,
by Theorem~\ref{thm.iteration.homo},
that the lower-triangular row-generating matrix \reff{def.Snk.rowgen.bis}
is coefficientwise totally positive.
But this latter fact has a much simpler proof:
it is an immediate consequence of the matrix identity \reff{eq.Sm.Smxi}
and the coefficientwise total positivity of $\sfS^{(m)}$
(Theorem~\ref{thm.Stype.generalized.minors})
and $T_\xi$ (Lemma~\ref{lemma.toeplitz.power}).

\section{Weights eventually periodic of period $\bm{m\!+\!1}$ or \hbox{$\bm{m}$}}
   \label{sec.eventually_periodic}

The results of the previous section can be generalized to handle
some cases in which the weights $\balpha$
are eventually periodic of period~$m+1$ or $m$,
i.e.\ they consist of a finite sequence followed by a periodic sequence.
As a special case we will obtain the
{\em multivariate Aval polynomials}\/ \cite{Aval_08,Mlotkowski_12}.
We will use the recurrence method
from Sections~\ref{subsec.periodic.m+1.recurrence}
and \ref{subsec.periodic.m.recurrence}.

\subsection[Period $m\!+\!1$]{Period $\bm{m\!+\!1}$}
   \label{subsec.eventually_periodic.m+1.recurrence}

Fix an integer $m \ge 1$,
and let $\balpha$ be given by
\be
   \balpha  \;=\; (\alpha_i)_{i \ge m}  \;=\;
   y_1,\ldots,y_\ell,
   x_0,\ldots,x_m, x_0,\ldots,x_m, \ldots
 \label{eq.weights.eventually_periodic.m+1}
\ee
for some $\ell \ge 1$.
We again use the recurrence \reff{eq.mSRfk.1}.
The periodicity of the $\balpha$ starting at $\alpha_{m+\ell}$
implies the periodicity of the $f_k$ starting at $k=\ell$;
therefore the product $f_k \cdots f_{k+m}$
is independent of $k$ for $k \ge \ell$, let us call it $F$.
It is the same quantity $F$ that was calculated
in Section~\ref{subsec.periodic.m+1.recurrence}
and is given by \reff{def.F.Qn}/\reff{def.Qn.x0xm}.
Then $f_\ell(t) = 1 + x_0 t F(t)$ is the same quantity
that was given by \reff{def.f0.Pn}/\reff{def.Pn.x0xm};
and more generally $f_{\ell+i}(t) = 1 + x_i t F(t)$ for $0 \le i \le m$.

Now we can work backwards:  from \reff{eq.mSRfk.2} we have
\begin{subeqnarray}
   f_{\ell-1}(t)
   & = &
   {1 \over 1 \:-\: y_\ell t \, f_{\ell}(t) \,\cdots\, f_{\ell+m-1}(t)}
          \\[2mm]
   & = &
   {1 \over 1 \:-\: y_\ell t \, \displaystyle {F(t) \over f_{\ell+m}(t)}}
          \\[2mm]
   & = &
   {1 \over 1 \:-\: y_\ell t \, \displaystyle {F(t) \over 1 + x_m t F(t)}}
          \\[2mm]
   & = &
   {1 + x_m t F(t)  \over 1 + (x_m - y_\ell) t F(t)}
      \;.
 \slabel{eq.fl-1.d}
 \label{eq.fl-1}
\end{subeqnarray}
And then (if $\ell \ge 2$) we can write
\begin{subeqnarray}
   f_{\ell-2}(t)
   & = &
   {1 \over 1 \:-\: y_{\ell-1} t \, f_{\ell-1}(t) \,\cdots\, f_{\ell+m-2}(t)}
          \\[2mm]
   & = &
   {1 \over 1 \:-\: y_{\ell-1} t \, f_{\ell-1}(t) \, 
          \displaystyle {F(t) \over f_{\ell+m-1}(t) \, f_{\ell+m}(t)}}
          \\[2mm]
   & = &
   {1 \over 1 \:-\: y_{\ell-1} t \,
          \displaystyle {F(t) \over
           [1 + (x_m - y_\ell) t F(t)] \, [1 + x_{m-1} t F(t)]}
   }
\end{subeqnarray}
and so forth.
The equations get rather messy as $\ell$ grows,
so we limit ourselves to working out the case $\ell=1$.

Applying Lagrange inversion to \reff{eq.functeqn.F} and \reff{eq.fl-1.d}
with $\ell=1$, we obtain for $n \ge 1$
\be
   [t^n] \, f_0(t)
   \;=\;
   \!\!\!\!\!
   \sum_{\begin{scarray}
            j_\star,j_0,\ldots,j_m \ge 0 \\
            j_\star+j_0+\ldots+j_m = n
         \end{scarray}}
   \!\!\!\!
   {j_\star \over n} \, y_1 \, (y_1 - x_m)^{j_\star - 1}
   \prod_{i=0}^m  \binom{n}{j_i}  x_i^{j_i}
   \;.
\ee
We do not know whether these polynomials have any combinatorially
interesting specializations
(other than, of course, the periodic case $y_1 = x_m$).

\subsection[Period $m$]{Period $\bm{m}$}
   \label{subsec.eventually_periodic.m.recurrence}

Fix an integer $m \ge 1$,
and let $\balpha$ now be given by
\be
   \balpha  \;=\; (\alpha_i)_{i \ge m}  \;=\;
   y_1,\ldots,y_\ell,
   x_0,\ldots,x_{m-1}, x_0,\ldots,x_{m-1}, \ldots
 \label{eq.weights.eventually_periodic.m}
\ee
for some $\ell \ge 1$.
This time we use the recurrence \reff{eq.mSRfk.2}.
The periodicity of the $\balpha$ starting at $\alpha_{m+\ell}$
implies the periodicity of the $f_k$ starting at $k=\ell$;
therefore the product $f_{k+1} \cdots f_{k+m}$ is independent of $k$
for $k \ge \ell-1$
(note that this is now $\ell-1$, not~$\ell$), let us call it $G$.
It is the same quantity $G$ that was calculated
in Section~\ref{subsec.periodic.m.recurrence}
and is given by \reff{def.G.Qn-k}/\reff{def.Qn-.x0xm-1}.
Then $f_\ell(t) = 1/[1 - x_0 t G(t)]$ is the same quantity
that was given by \reff{def.f0.Pn-}/\reff{def.Pn-.x0xm-1};
and more generally $f_{\ell+i}(t) = 1/[1 - x_i t G(t)]$ for $0 \le i \le m-1$.

Now we can again work backwards:  from \reff{eq.mSRfk.2} we have
\begin{subeqnarray}
   f_{\ell-1}(t)
   & = &
   {1 \over 1 - y_\ell t  \, f_{\ell}(t) \cdots f_{\ell+m-1}(t)}
       \\[2mm]
   & = &
   {1 \over 1 - y_\ell t \, G(t)}
   \;,
 \slabel{eq.fl-1.periodm.b}
\end{subeqnarray}
which is even simpler than \reff{eq.fl-1}.
And then (if $\ell \ge 2$) we can write
\begin{subeqnarray}
   f_{\ell-2}(t)
   & = &
   {1 \over 1 \:-\: y_{\ell-1} t \, f_{\ell-1}(t) \,\cdots\, f_{\ell+m-2}(t)}
          \\[2mm]
   & = &
   {1 \over 1 \:-\: y_{\ell-1} t \, f_{\ell-1}(t) \, 
          \displaystyle {G(t) \over f_{\ell+m-1}(t)}}
          \\[2mm]
   & = &
   {1 \over 1 \:-\: y_{\ell-1} t \, G(t) \,
          \displaystyle {1 - x_{m-1} t G(t) \over 1 - y_\ell t G(t)}}
   \;.
\end{subeqnarray}
Once again the equations get rather messy as $\ell$ grows,
so we limit ourselves to working out the case $\ell=1$.

Before proceeding with the $\ell=1$ case,
it is convenient to rename
$y_1,x_0,\ldots,x_{m-1}$ as $x_0,x_1,\ldots,x_m$.
Then, applying Lagrange inversion to
\reff{eq.functeqn.G} and \reff{eq.fl-1.periodm.b},
we obtain for $n \ge 1$
\be
   [t^n] \, f_0(t)
   \;=\;
   \!\!\!\!
   \sum_{\begin{scarray}
            j_0,\ldots,j_m \ge 0 \\
            j_0+\ldots+j_m = n
         \end{scarray}}
   \!\!\!
   {j_0 \over n} \, x_0^{j_0}
   \prod_{i=1}^m  \binom{n+j_i-1}{j_i}  x_i^{j_i}
   \;.
 \label{def.aval}
\ee
We call the right-hand side of \reff{def.aval}
the \textbfit{multivariate Aval polynomial} of order~$m$,
written $A_n^{(m)}(x_0,x_1,\ldots,x_m)$,
because the coefficients in \reff{def.aval} appeared first
(to our knowledge) in Aval's paper
\cite[Remarks~2.5 and 3.2]{Aval_08},
along with a combinatorial interpretation in terms of $m$-Dyck paths
that is different from the one implicit in our $m$-S-fraction;
see also \cite{Mlotkowski_12}.
As an immediate consequence of Theorem~\ref{thm.Stype.minors} we have:

\begin{theorem}[Hankel-total positivity of multivariate Aval polynomials]
   \label{thm.aval.hankelTP}
For each integer $m \ge 1$,
the sequence $(A_n^{(m)})_{n \ge 0}$
of multivariate Aval polynomials of order~$m$
(with of course $A_0^{(m)} \eqdef 1$)
is coefficientwise Hankel-totally positive,
jointly in the indeterminates $x_0,\ldots,x_m$.
\end{theorem}

\section{Weights quasi-affine or factorized of period ${\bm{m\!+\!1}}$ or $\bm{m}$}
   \label{sec.quasi-affine}

The polynomials $P_n(x,y,u,v)$ defined by the S-fraction
\be
   \sum_{n=0}^\infty P_n(x,y,u,v) \: t^n
   \;=\;
   \cfrac{1}{1 - \cfrac{xt}{1 - \cfrac{yt}{1 - \cfrac{(x+u)t}{1- \cfrac{(y+v)t}{1 - \cdots}}}}}
 \label{eq.eulerian.fourvar.contfrac}
\ee
with coefficients
\begin{subeqnarray}
   \alpha_{2k-1}  & = &  x + (k-1) u \\
   \alpha_{2k}    & = &  y + (k-1) v
 \label{def.weights.eulerian.fourvar}
\end{subeqnarray}
contain many classical polynomials as specializations:
these include the Eulerian polynomials [see \reff{def.eulerian} below],
the Stirling cycle polynomials
[see \reff{def.stirlingcycle}/\reff{eq.stirlingcycle.Sfraction} below],
and the record-antirecord permutation polynomials
[see
 \reff{eq.euler.contfrac.BISBIS.2F0}/\reff{eq.thm.rF0}/\reff{eq.recantirec}
 below],
among many others.\footnote{
   See \cite{Sokal-Zeng_masterpoly} for further examples.
}
Since $P_n(x,y,u,v)$ is a homogeneous polynomial of degree~$n$,
it has three ``truly independent'' variables.
By Euler's continued fraction \reff{eq.nfact.contfrac},
we have $P_n(1,1,1,1) = n!$, so it is natural to expect that
$P_n(x,y,u,v)$ enumerates permutations of $[n]$
according to some natural trivariate statistic.
Very recently, Zeng and one of us \cite{Sokal-Zeng_masterpoly}
found two alternative versions of this trivariate statistic:
\begin{subeqnarray}
   P_n(x,y,u,v)
   & = &
   \sum_{\sigma \in \Sym_n}
      x^{\arec(\sigma)} y^{\erec(\sigma)}
         u^{n - \exc(\sigma) - \arec(\sigma)} v^{\exc(\sigma) - \erec(\sigma)}
     \slabel{eq.eulerian.fourvar.arec} \\[2mm]
   & = &
   \sum_{\sigma \in \Sym_n}
      x^{\cyc(\sigma)} y^{\erec(\sigma)}
         u^{n - \exc(\sigma) - \cyc(\sigma)} v^{\exc(\sigma) - \erec(\sigma)}
     \slabel{eq.eulerian.fourvar.cyc}
     \label{eq.eulerian.fourvar}
\end{subeqnarray}
where the sum runs over permutations of $[n]$,
and $\arec(\sigma), \erec(\sigma), \exc(\sigma), \cyc(\sigma)$
denote the number of antirecords, exclusive records, excedances
and cycles in $\sigma$.\footnote{
   Given a permutation $\sigma \in \Sym_n$,
   an index $i \in [n]$ is called a
   {\em record}\/ (or {\em left-to-right maximum}\/)
         if $\sigma(j) < \sigma(i)$ for all $j < i$
      [note in particular that the index 1 is always a record];
   an {\em antirecord}\/ (or {\em right-to-left minimum}\/)
         if $\sigma(j) > \sigma(i)$ for all $j > i$
      [note in particular that the index $n$ is always an antirecord];
   an {\em exclusive record}\/ if it is a record and not also
         an antirecord;
   and an {\em exclusive antirecord}\/ if it is an antirecord and not also
      a record.
   An index $i \in [n]$ is called an {\em excedance}\/
   if $\sigma(i) > i$.
}
We~say that the weights \reff{def.weights.eulerian.fourvar}
are \textbfit{quasi-affine} of period~2
because $\alpha_i$ is affine in~$i$
within each residue class of $i \bmod 2$.


Here we would like to consider the generalization of this setup
to $m$-branched S-fractions,
in which the coefficients $\balpha = (\alpha_i)_{i \ge m}$
are quasi-affine of period either $m+1$ or $m$.
For period $p$ we will write
\be
   \alpha_{m+j+pk}  \;=\;  x_j \,+\, k u_j
 \label{eq.alpha.quasi-affine}
\ee
where $\bfx = (x_0,\ldots,x_{p-1})$
and $\bfu = (u_0,\ldots,u_{p-1})$ are indeterminates.
We let $P_n^{(m,p)}(\bfx,\bfu)$
be the polynomials obtained by specializing
the $m$-Stieltjes--Rogers polynomials $S_n^{(m)}(\balpha)$
to the weights \reff{eq.alpha.quasi-affine}.
[When $\bfu = \bzero$ the weights $\balpha$ are periodic of period $p$,
and the polynomials $P_n^{(m,p)}(\bfx,\bzero)$
reduce to the polynomials $P_n^{(m,p)}(\bfx)$
defined at the beginning of Section~\ref{sec.periodic}.]
We are not yet able to obtain a complete combinatorial interpretation,
analogous to \reff{eq.eulerian.fourvar}, of the polynomials
$P_n^{(m,p)}(\bfx,\bfu)$,
but we are able to handle some interesting special cases.

In particular we will treat the case $\bfu = \bfx$, i.e.
\be
   \alpha_{m+j+pk}  \;=\;  (k+1) x_j
   \;,
 \label{eq.alpha.quasi-affine.u=x}
\ee
with period $p = m+1$ or $m$.
For $(m,p) = (1,2)$ these are the homogenized Eulerian polynomials
\be
   P_n^{(1,2)}(x_0,x_1,x_0,x_1)
   \;=\;
   P_n(x_0,x_1,x_0,x_1)
   \;=\;
   \sum_{k=0}^n \euler{n}{k} \, x_0^{n-k} x_1^k
   \;,
 \label{def.eulerian}
\ee
where the Eulerian number $\euler{n}{k}$
is the number of permutations of $[n]$ with $k$ excedances.\footnote{
   Here we are using the convention \cite{Graham_94} in which
   (for $n \ge 1$) $\euler{n}{k}$ is nonzero for $0 \le k \le n-1$.
   Many authors use the opposite convention in which
   $\euler{n}{k}$ is nonzero for $1 \le k \le n$;
   then $x_0$ and~$x_1$ have to be interchanged.
}${}^,$\footnote{
   The identity \reff{def.eulerian} --- that is, the S-fraction
   for the Eulerian polynomials --- was found by Stieltjes
   \cite[section~79]{Stieltjes_1894}.
   Stieltjes does not specifically mention the Eulerian polynomials,
   but he does state that the continued fraction
   is the formal Laplace transform of
   $(1-y) / (e^{t(y-1)} - y)$,
   which is well known to be the exponential generating function
   of the Eulerian polynomials.
   Stieltjes also refrains from showing the proof:
   ``Pour abr\'eger, je supprime toujours les artifices qu'il faut employer
     pour obtenir la transformation de l'int\'egrale d\'efinie
     en fraction continue'' (!).
   But a proof is sketched, albeit also without much explanation,
   in the book of Wall \cite[pp.~207--208]{Wall_48}.
   The J-fraction corresponding to
   the contraction of this S-fraction
   was proven, by combinatorial methods,
   by Flajolet \cite[Theorem~3B(ii) with a slight typographical error]{Flajolet_80}.
   Dumont \cite[Propositions~2 and 7]{Dumont_86}
   gave a direct combinatorial proof of the S-fraction,
   based on an interpretation of the Eulerian polynomials
   in terms of ``bipartite involutions of $[2n]$''
   and a bijection of these onto Dyck paths.
}
We therefore refer to \reff{eq.alpha.quasi-affine.u=x}
as the \textbfit{Eulerian-quasi-affine} weights,
and to the polynomials $P_n^{(m,m+1)}(\bfx,\bfx)$
as \textbfit{multivariate Eulerian polynomials}.
We will show that they contain the $m$th-order Eulerian polynomials
\cite{Barbero_15} as specializations.

This identification will allow us, among other things, to resolve a mystery
from \cite{Elvey-Price-Sokal_wardpoly}.
It was shown in \cite{Elvey-Price-Sokal_wardpoly}
that the reversed second-order Eulerian polynomials $\overline{E}_n^{[2]}(x)$
[defined in \reff{eq.reveuler} below]
are given by a T-fraction \reff{eq.f0.Tfrac}
with coefficients $\delta_i = (i-1)(x-1)$ and $\alpha_i = i$.
Since here $\delta_i$ is {\em not}\/ coefficientwise nonnegative
(or even pointwise nonnegative when $0 \le x < 1$),
the general theory of \cite{Sokal_totalpos}
--- that is, the $m=1$ case of Theorem~\ref{thm.Ttype.minors} ---
says nothing about the Hankel-total positivity of
the reversed second-order Eulerian polynomials.
And yet, it was found empirically that the sequence of
reversed second-order Eulerian polynomials
{\em is}\/ coefficientwise Hankel-totally positive:
this was tested through the $13 \times 13$ Hankel matrix.
Here we will show why (Corollary~\ref{cor.eulerianr} below):
in addition to the T-fraction found in \cite{Elvey-Price-Sokal_wardpoly},
the reversed second-order Eulerian polynomials
are {\em also}\/ given by a 2-branched S-fraction
with Eulerian-quasi-affine coefficients \reff{eq.alpha.quasi-affine.u=x},
namely $(x_0,x_1,x_2) = (1,1,x)$ or $(1,x,1)$.
The coefficientwise Hankel-total positivity
is then an immediate consequence of Theorem~\ref{thm.Stype.minors}.
In fact, we will show more generally that, for every $m \ge 1$,
the $m$th-order Eulerian polynomials have an $m$-branched S-fraction
with Eulerian-quasi-affine coefficients \reff{eq.alpha.quasi-affine.u=x}
in which one of $x_1,\ldots,x_m$ is $x$
and all the other~$x_i$ (including $x_0$) are 1.
This again implies the coefficientwise Hankel-total positivity.
In our opinion this story illustrates some of the power
of the branched continued fractions to resolve problems
that seemed mysterious (or simply intractable)
from the point of view of the classical continued fractions,
and indeed to provide multivariate generalizations.


In fact, we can handle some weights that are considerably more general
than \reff{eq.alpha.quasi-affine.u=x}, namely
\be
   \alpha_{m+j+pk}  \;=\;  (k+1) c_k x_j
 \label{eq.alpha.factorized}
\ee
where $\bfc = (c_k)_{k \ge 0}$ are indeterminates.
We call these weights \textbfit{factorized} of period~$p$,
and we write $\widehat{P}_n^{(m,p)}(\bfx,\bfc)$
for the corresponding polynomials.
Of course, they reduce to $P_n^{(m,p)}(\bfx,\bfx)$ when $\bfc = \bone$.
We call $\widehat{P}_n^{(m,m+1)}(\bfx,\bfc)$
the \textbfit{extended multivariate Eulerian polynomials}.

The plan of this section is as follows:
First we discuss briefly the special case $\bfx = \bfu = \bone$,
which leads to multifactorials
(Section~\ref{subsec.quasi-affine.multifactorials}).
Then we interpret the general Eulerian-quasi-affine case $\bfx = \bfu$
[i.e.\ \reff{eq.alpha.quasi-affine.u=x}] with period $m+1$,
and more generally the factorized case \reff{eq.alpha.factorized}
with period~$m+1$, in~terms of increasing $(m+1)$-ary trees;
and we also take the limit $m \to\infty$ to obtain the
Eulerian symmetric functions (Section~\ref{subsec.quasi-affine.m+1}).
By suitable modifications we can interpret
the Eulerian-quasi-affine case $\bfx = \bfu$ with period~$m$,
and more generally the factorized case with period~$m$,
in~terms of increasing multi-$m$-ary trees,
and take the limit $m \to\infty$ to obtain the
Eulerian symmetric functions of negative type;
however, some of these results are still conjectural
(Section~\ref{subsec.quasi-affine.m}).
In Section~\ref{subsec.quasi-affine.observation}
we make a brief observation on the analogy between the
periodic and quasi-affine cases.
And finally, we specialize the multivariate Eulerian polynomials
in order to relate them to $m$th-order Eulerian polynomials
(Section~\ref{subsec.quasi-affine.stirling}).


\subsection{Multifactorials}   \label{subsec.quasi-affine.multifactorials}

Let us denote by
\be
   F_n^{(r)}  \;\eqdef\;  \prod_{j=0}^{n-1} (1+jr)
 \label{def.multifactorial}
\ee
the $r$th-order multifactorials,
so that $F_n^{(0)} = 1$, $F_n^{(1)} = n!$, $F_n^{(2)} = (2n-1)!!$,
$F_n^{(3)} = (3n-2)!!!$ and so forth.
The multifactorials are specializations of the
homogenized Stirling cycle polynomials
\begin{subeqnarray}
   C_n(x,y)
   & \eqdef &
   \prod_{j=0}^{n-1} (x+jy)
        \\[2mm]
   & = &
   \sum_{k=0}^n
   \stirlingcycle{n}{k} \, x^k y^{n-k}
 \label{def.stirlingcycle}
\end{subeqnarray}
where $\stirlingcycle{n}{k}$ is the number of permutations of $[n]$
with $k$ cycles (or $k$ records, or $k$ antirecords):
we have
\be
   F_n^{(r)}  \;=\;  C_n(1,r)  \;=\;  r^n C_n(1/r,1)
   \;.
\ee
Already Euler \cite[section~26]{Euler_1760} \cite{Euler_1788}\footnote{
   See footnote~\ref{footnote_Euler_1760} in the Introduction
   for the history of \cite{Euler_1760}.
   The paper \cite{Euler_1788},
   which is E616 in Enestr\"om's \cite{Enestrom_13} catalogue,
   was apparently presented to the St.~Petersburg Academy in 1776,
   and published posthumously in 1788.
   See footnote~\ref{footnote_Euler} in Section~\ref{sec.hyper.rF0} below
   for discussion
   of Euler's method for proving \reff{eq.stirlingcycle.Sfraction}.
}
showed that the ordinary generating function
of the Stirling cycle polynomials has the S-fraction
\be
   \sum_{n=0}^\infty C_n(x,y) \: t^n
   \;=\;
   \cfrac{1}{1 - \cfrac{xt}{1 - \cfrac{yt}{1 - \cfrac{(x+y)t}{1- \cfrac{2yt}{1- \cdots}}}}}
  \label{eq.stirlingcycle.Sfraction}
\ee
with coefficients $\alpha_{2k-1} = x+(k-1)y$ and $\alpha_{2k} = ky$.
Therefore the multifactorials have a classical S-fraction
with coefficients $\alpha_{2k-1} = 1+(k-1)r$ and $\alpha_{2k} = kr$.
In~terms of \reff{eq.eulerian.fourvar} we have $C_n(x,y) = P_n(x,y,y,y)$.

Here we will show that the multifactorials are also given by
very simple $m$-branched S-fractions in which the weights
are Eulerian-quasi-affine of period $m+1$ or $m$:  namely,
$P_n^{(m,m+1)}(\bone,\bone) = F_n^{(m)}$ (Corollary~\ref{cor.quasi-affine})
and
$P_n^{(m,m)}(\bone,\bone) = F_n^{(m+1)}$
(Corollary~\ref{cor.quasi-affine.multi}).
These results will, in fact, be corollaries of more powerful results
in which we supply a combinatorial interpretation
for the multivariate polynomials
$P_n^{(m,m+1)}(\bfx,\bfx)$ and $P_n^{(m,m)}(\bfx,\bfx)$.


\bigskip

{\bf Remarks.}
1.  These different representations illustrate once again
the extreme non\-uniqueness of $m$-S-fractions with $m \ge 2$.
For instance, $(3n-2)!!!$ is given by at least three different
types of 3-S-fractions:
\begin{itemize}
   \item  $\balpha = 1,3,0,0,4,6,0,0,7,9,0,0,10,12,0,0,\ldots\,$
      (and uncountably many variants),
      coming from the classical S-fraction \reff{eq.stirlingcycle.Sfraction}
      via Proposition~\ref{prop.reduction};
   \item  $\balpha = 1,1,2,0,2,3,3,0,4,4,5,0,5,6,6,0,\ldots\,$
      (and uncountably many variants),
      coming from the 2-S-fraction $P_n^{(2,2)}(\bone,\bone) = F_n^{(3)}$
      via Proposition~\ref{prop.reduction}; and
   \item  $\balpha = 1,1,1,1,2,2,2,2,3,3,3,3,4,4,4,4,\ldots\,$,
      coming from the 3-S-fraction 
      \linebreak
      $P_n^{(3,4)}(\bone,\bone) = F_n^{(3)}$.
\end{itemize}

2.  For period $p=1$, the sequences $P_n^{(m,1)}(\bone,\bone)$
appear to have combinatorial interpretations as the number of
linear extensions of certain grid-like posets \cite{Pan_15}:
see \cite[A274644]{OEIS} for $(m,p) = (2,1)$ and $(3,2)$,
\cite[A274763]{OEIS} for $(m,p) = (3,1)$,
and \cite{Pan_15} for $(m,p) = (4,1)$.
We hope to investigate this connection in a future paper.

3.  For periods $p \ne 1,m,m+1$ we have been unable to find any
combinatorial interpretation of the polynomials $P_n^{(m,p)}(\bfx,\bfx)$
--- or even of the numbers $P_n^{(m,p)}(\bone,\bone)$ ---
other than the obvious one in terms of $m$-Dyck paths.
Indeed, the sequences $P_n^{(m,p)}(\bone,\bone)$ for $p \ne 1,m,m+1$
are not currently in \cite{OEIS}
(with the exception, noted above, of $(m,p) = (3,2)$).
\myendremark

\subsection[Period ${m+1}$: Increasing $(m+1)$-ary trees]{Period ${\bm{m+1}}$: Increasing $\bm{(m+1)}$-ary trees}
   \label{subsec.quasi-affine.m+1}

\begin{figure}[!ht]
\begin{center}
\includegraphics[scale=1]{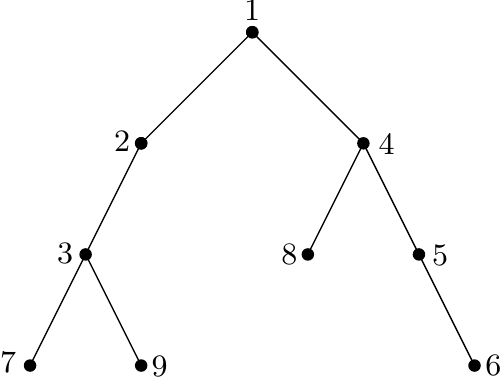}
\caption{\label{fig.increasingbinarytree} An increasing binary tree (the edges are left or right).}
\end{center}
\end{figure}

In Section~\ref{subsec.periodic.combinatorial}
we gave the recursive definition of an $r$-ary tree,
and we said what is meant by an $i$-edge in such a tree.
There the vertices were unlabeled; now we label them.
A {\em labeled $r$-ary tree}\/ on $n$ vertices
is an $r$-ary tree on $n$ vertices
in which the vertices are given distinct labels from the set $[n]$.
An \textbfit{increasing $\bm{r}$-ary tree} is a labeled $r$-ary tree in which
the label of each child is greater than the label of its parent;
otherwise put, the labels increase along every path downwards from the root
(see Figure~\ref{fig.increasingbinarytree} for an example).
The meaning of $i$-edge is unchanged, since it does not refer to the labels.

It is well known \cite[p.~24]{Stanley_86}
that the number of increasing binary trees on $n$ vertices is $n!$,
and more generally that the number of increasing $r$-ary trees
on $n$ vertices is the multifactorial $F_n^{(r-1)}$
\cite[p.~30, Example~1]{Bergeron_92}.
We would now like to generalize this enumeration to a multivariate polynomial.

\subsubsection{Multivariate Eulerian polynomials and Eulerian symmetric functions: Definitions and statement of results}

Fix an integer $m \ge 1$, and let $\bfx = (x_0,\ldots,x_m)$ be indeterminates.
Let $\scrq^{(m)}_n(\bfx)$ be the generating polynomial for
increasing $(m+1)$-ary trees on $n+1$ vertices
in~which each $i$-edge gets a weight $x_i$;
and more generally,
let $\scrq^{(m)}_{n,k}(\bfx)$ be the generating polynomial for
unordered forests of increasing $(m+1)$-ary trees
on $n+1$ total vertices with $k+1$ components
in~which each $i$-edge gets a weight $x_i$.
(If we were to consider instead ordered forests, it would just multiply
 the polynomial by $(k+1)!$, since the trees are labeled and hence
 distinguishable.)
And finally, define $\scrp^{(m)}_n(\bfx)$ by
\begin{subeqnarray}
   \scrp^{(m)}_0(\bfx)  & = &    1   \\
   \scrp^{(m)}_n(\bfx)  & = &    x_0 \, \scrq^{(m)}_{n-1}(\bfx)
        \qquad\hbox{for $n \ge 1$}
 \slabel{def.scrp.b}
 \label{def.scrp}
\end{subeqnarray}
$\scrp^{(m)}_n(\bfx)$ is therefore
the generating polynomial for
increasing $(m+1)$-ary trees on $n+1$ vertices in which
the only edge (if any) emanating from the root is a 0-edge
and in which each $i$-edge gets a weight $x_i$.
Both $\scrp^{(m)}_n(\bfx)$ and $\scrq^{(m)}_n(\bfx)$
are homogeneous polynomials of degree~$n$,
while $\scrq^{(m)}_{n,k}(\bfx)$ is homogeneous of degree~$n-k$;
we refer to them collectively as \textbfit{multivariate Eulerian polynomials}.
In particular, $\scrp^{(m)}_n(\bone) = \scrq^{(m)}_{n-1}(\bone) = F_n^{(m)}$
for $n \ge 1$.

We will prove:

\begin{theorem}[Branched continued fractions for the multivariate Eulerian polynomials]
 \label{thm.quasi-affine}
For every integer $m \ge 1$, we have:
\begin{itemize}
   \item[(a)] $\scrp^{(m)}_n(\bfx) = P_n^{(m,m+1)}(\bfx,\bfx)$,
which by definition equals
$S_n^{(m)}(\balpha)$ where the weights $\balpha$ are given by
the Eulerian-quasi-affine formula \reff{eq.alpha.quasi-affine.u=x}
with period $p=m+1$.
   \item[(b)] $\scrq^{(m)}_n(\bfx) =  J_n^{(m)}(\bbeta)$
where the weights $\bbeta$ are given by
\be
   \beta_i^{(\ell)}
   \;=\;
   {(i+1)! \over (i-\ell)!} \: e_{\ell+1}(x_0,\ldots,x_m)
 \label{eq.thm.quasi-affine}
\ee
(here $e_{\ell+1}$ is the elementary symmetric polynomial of degree $\ell+1$).
   \item[(c)] $\scrq^{(m)}_{n,k}(\bfx) =  J_{n,k}^{(m)}(\bbeta)$
where the weights $\bbeta$ are again given by \reff{eq.thm.quasi-affine}.
\end{itemize}
\end{theorem}

Since $\scrp^{(m)}_n(\bone) = F_n^{(m)}$, it follows that:

\begin{corollary}[Eulerian-quasi-affine $m$-S-fraction of period $m+1$ for multifactorials]
  \label{cor.quasi-affine}
For every integer $m \ge 1$, we have
$P_n^{(m,m+1)}(\bone,\bone) = F_n^{(m)}$.
\end{corollary}

Please note that the polynomials $\scrq^{(m)}_n(\bfx)$
and $\scrq^{(m)}_{n,k}(\bfx)$
are symmetric in $x_0,\ldots,x_m$,
since we can make a global permutation of the child-order indices $0,\ldots,m$.
It therefore follows from \reff{def.scrp.b}
and Theorem~\ref{thm.quasi-affine}(a)
that the polynomials $P_n^{(m,m+1)}(\bfx,\bfx)$
are, up to the prefactor $x_0$, also symmetric in $x_0,\ldots,x_m$ ---
a fact that seems far from obvious from the branched continued fraction
or from the combinatorial definition via $m$-Dyck paths.
The invariance of $P^{(m,m+1)}_n(\bfx,\bfx)$
under permutations of $x_1,\ldots,x_m$
also illustrates once again
the nonuniqueness of $m$-S-fractions with $m \ge 2$.

Combining Theorems~\ref{thm.quasi-affine}(a) and \ref{thm.Stype.minors},
and using \reff{def.scrp.b} to handle $\scrq^{(m)}_n$, we conclude:

\begin{corollary}[Hankel-total positivity of the multivariate Eulerian polynomials]
  \label{cor.quasi-affine.2}
For every integer ${m \ge 1}$,
the sequences $(\scrp^{(m)}_n(\bfx))_{n \ge 0}$
and $(\scrq^{(m)}_n(\bfx))_{n \ge 0}$
are coefficientwise Hankel-totally positive,
jointly in the indeterminates $x_0,\ldots,x_m$.
\end{corollary}

We will also show:

\begin{lemma}
   \label{lemma.Pm.beta.totalpos}
The production matrix $P^{(m)}(\bbeta)$ [cf.\ \reff{def.Pm}]
associated to the weights \reff{eq.thm.quasi-affine}
is coefficientwise totally positive,
jointly in the indeterminates $x_0,\ldots,x_m$.
\end{lemma}

Using Theorems~\ref{thm.Jtype.minors} and \ref{thm.quasi-affine}(b)
we obtain an alternate proof of the Hankel-total positivity of
$(\scrq^{(m)}_n(\bfx))_{n \ge 0}$;
and from Theorems~\ref{thm.Jtype.generalized.minors}
and \ref{thm.quasi-affine}(c) we conclude:

\begin{corollary}[Total positivity of lower-triangular matrix of multivariate Eulerian polynomials]
  \label{cor.quasi-affine.2a}
For every integer ${m \ge 1}$,
the lower-triangular matrix $(\scrq^{(m)}_{n,k}(\bfx))_{n,k \ge 0}$
is coefficientwise totally positive,
jointly in the indeterminates $x_0,\ldots,x_m$.
\end{corollary}

We also have explicit expressions for the multivariate Eulerian polynomials
$\scrp^{(m)}_n$ and $\scrq^{(m)}_n$,
and a recurrence for $\scrq^{(m)}_{n,k}$,
given by repeated application of a first-order linear differential operator.
We will prove:

\begin{proposition}[Differential expressions for the multivariate Eulerian polynomials]
   \label{prop.MVeulerian.differential}
For every integer ${m \ge 1}$, we have
\begin{eqnarray}
   \scrp^{(m)}_n(\bfx)  & = &  (\scrd_m \,+\, x_0)^n \: 1
          \label{eq.scrp.differential}  \\[2mm]
   \scrq^{(m)}_n(\bfx)  & = &
      \Bigl( \scrd_m \,+\, \sum_{i=0}^m x_i \Bigr)^{\! n} \: 1
          \label{eq.scrq.differential}  \\
   \scrq^{(m)}_{n,k}(\bfx)  & = &
      \Bigl( \scrd_m \,+\, (k+1) \sum_{i=0}^m x_i \Bigr)
      \scrq^{(m)}_{n-1,k}(\bfx)  \:+\:  \scrq^{(m)}_{n-1,k-1}(\bfx)
      \:+\: \delta_{n0} \delta_{k0}
      \qquad
          \label{eq.scrqnk.differential}
\end{eqnarray}
where
\be
   \scrd_m
   \;=\;
   \sum_{i=0}^m \biggl( \! x_i
                        \!\!\!\!
                            \sum_{\begin{scarray}
                                    0 \le j \le m \\
                                    j \ne i
                                  \end{scarray}}
                        \!\!\!\!\!\!
                            x_j \!
                \biggr) \, {\partial \over \partial x_i}
   \;.
 \label{def.dm}
\ee
\end{proposition}

{\bf Remarks.}
1.  For the classical ($m=1$) homogenized Eulerian polynomials
$A_n(x,y) = \sum_{k=0}^n \euler{n}{k} \, x^{n-k} y^k$,
the recurrence
$A_{n+1}(x,y) = [xy (\partial/\partial x + \partial/\partial y) + x] A_n(x,y)$
can easily be derived from the classic differential recurrence
for the univariate Eulerian polynomials $A_n(y) = A_n(1,y)$,
i.e.\ $A_{n+1}(y) = (1+ny) A_n(y) + y(1-y) A'_n(y)$
\cite[p.~10]{Petersen_15}.
The interesting thing here is that the $n$-dependent recurrence
for the univariate polynomials becomes an $n$-independent recurrence
for the homogeneous polynomials.
This observation --- which is perhaps not as well known as it should be ---
goes back to Carlitz \cite{Carlitz_73} (see also \cite{Visontai_13}).

2. The polynomials $\scrq^{(m)}_n(\bfx)$ for $m=2$
were introduced by Dumont \cite{Dumont_80},
who gave also the combinatorial interpretations
in terms of increasing ternary trees \cite[Proposition~1]{Dumont_80}
and Stirling permutations \cite[Proposition~3]{Dumont_80}
(for the latter, see Section~\ref{subsec.quasi-affine.stirling} below).
\myendremark

\bigskip

Whenever $m < m'$ we have
$\scrq^{(m)}_{n,k}(x_0,\ldots,x_m) = \scrq^{(m')}_{n,k}(x_0,\ldots,x_m,0,\ldots,0)$
and
$\scrp^{(m)}_n(x_0;x_1,\ldots,x_m) = \scrp^{(m')}_n(x_0;x_1,\ldots,x_m,0,\ldots,0)$:
this follows immediately from the combinatorial definition,
or alternatively
(for $\scrp^{(m)}_n$ and $\scrq^{(m)}_n$)
from the branched continued fraction
using Theorem~\ref{thm.quasi-affine}(a) and Proposition~\ref{prop.reduction}.
We can therefore pass to the limit ${m \to\infty}$:
let $\bfx = (x_i)_{i \ge 0}$ be indeterminates,
and define $\scrq^{(\infty)}_n(\bfx)$, $\scrq^{(\infty)}_{n,k}(\bfx)$
and $\scrp^{(\infty)}_n(\bfx)$ in terms of increasing $\infty$-ary trees
in the obvious way.
Here $\scrq^{(\infty)}_n(\bfx)$ and $\scrq^{(\infty)}_{n,k}(\bfx)$
are symmetric functions in the indeterminates $\bfx$,
while $\scrp^{(\infty)}_n(\bfx)$ is a symmetric function
up to the prefactor $x_0$;
they are homogeneous of degree~$n$ or $n-k$.
We refer to $\scrq^{(\infty)}_n(\bfx)$ and $\scrq^{(\infty)}_{n,k}(\bfx)$
as the \textbfit{Eulerian symmetric functions}.

Theorem~\ref{thm.quasi-affine}(b,c),
Corollaries~\ref{cor.quasi-affine.2} and \ref{cor.quasi-affine.2a},
and Proposition~\ref{prop.MVeulerian.differential}
extend immediately to $m=\infty$.
Let us state only the result on Hankel-total positivity:

\begin{corollary}[Hankel-total positivity of $\infty$-variate Eulerian series]
  \label{cor.quasi-affine.3}
The sequences $(\scrp^{(\infty)}_n(\bfx))_{n \ge 0}$
and $(\scrq^{(\infty)}_n(\bfx))_{n \ge 0}$
are Hankel-totally positive with respect to the coefficientwise order
on the formal-power-series ring $\Z[[\bfx]]$.
\end{corollary}

\smallskip

{\bf Remark.}
For $\scrq^{(\infty)}$ we can prove a stronger result:
namely, that the sequence $(\scrq^{(\infty)}_n(\bfx))_{n \ge 0}$
is Hankel-totally positive with respect to the Schur order
on the ring of symmetric functions,
which is stronger than the monomial (= coefficientwise) order.
See the Remark after the proof of
Lemma~\ref{lemma.Pm.beta.totalpos.factorized}
at the end of Section~\ref{subsubsec.quasi-affine.bijective}.
\myendremark

\smallskip

We can say more about the Eulerian symmetric functions $\scrq^{(\infty)}_n$
and $\scrq^{(\infty)}_{n,k}$
by using the interpretation of $\infty$-ary trees
as ordered trees in which each edge carries a label $i \in \N$
and the edges emanating outwards from each vertex consist, in order,
of zero or one edges labeled 0, then zero or one edges labeled 1, etc.\ 
(Section~\ref{subsec.periodic.combinatorial}).
Since the choice of labels on the edges emanating outwards from a vertex $v$
can be made independently for each $v$,
the definition of $\scrq^{(\infty)}_n$ can be rephrased as:
$\scrq^{(\infty)}_n(\bfx)$ is the generating formal power series for
increasing ordered trees on $n+1$ vertices in which
each vertex with $i$ children gets a weight $e_i(\bfx)$,
where $e_i$ is the elementary symmetric function;
and an analogous rephrasing holds for $\scrq^{(\infty)}_{n,k}$.
But now let us define
the {\em degree sequence}\/ $\lambda$ of a rooted tree $T$
to be the list of its vertices' out-degrees,
written in weakly decreasing order and with final zeroes removed.
Thus, if $T$ has $n+1$ vertices, $\lambda$ is a partition of $n$.
Then the definition of $\scrq^{(\infty)}_n$ can be further rephrased as:

\begin{proposition}[Eulerian symmetric functions in terms of elementary symmetric functions]
   \label{prop.eulerian.symfun}
For each partition $\lambda$ of $n$,
let $a_\lambda$ be the number of increasing ordered trees on $n+1$ vertices
that have degree sequence $\lambda$.  Then
\be
   \scrq^{(\infty)}_n(\bfx)
   \;=\;
   \sum_{\lambda \vdash n} a_\lambda \, e_\lambda(\bfx)
   \;.
 \label{eq.prop.eulerian.symfun}
\ee
\end{proposition}

For small $n$, the coefficients $a_\lambda$ can be computed
by explicit enumeration of ordered trees together with
the well-known result \cite[pp.~67 and 609, exercise~20]{Knuth_98}
(see also \cite[Lemma~2.1]{Gessel-Seo_06} for a refinement)
that the number of increasing labelings of a rooted tree on $n+1$ vertices
is $(n+1)!$ divided by the product of subtree sizes.
However, a more efficient method is to compute $\scrq^{(\infty)}_n$
by using the production matrix $P(\bbeta)$
[i.e.\ \reff{def.Pm} with $m=\infty$]
for the weights \reff{eq.thm.quasi-affine}.
The first few $\scrq^{(\infty)}_n$ are:
\begin{subeqnarray}
   \scrq^{(\infty)}_0  & = &   1  \\[1mm]
   \scrq^{(\infty)}_1  & = &   e_1  \\[1mm]
   \scrq^{(\infty)}_2  & = &   e_1^2 + 2e_2  \\[1mm]
   \scrq^{(\infty)}_3  & = &   e_1^3 + 8e_1 e_2 + 6e_3  \\[1mm]
   \scrq^{(\infty)}_4  & = &
       e_1^4 + 22 e_1^2 e_2 + 16 e_2^2 + 42 e_1 e_3 + 24e_4  \\[1mm]
   \scrq^{(\infty)}_5  & = &
       e_1^5 + 52 e_1^3 e_2 + 136 e_1 e_2^2 + 192 e_1^2 e_3 + 
           180 e_2 e_3 + 264 e_1 e_4 + 120 e_5
       \qquad \\[1mm]
   \scrq^{(\infty)}_6  & = &
   e_1^6 + 114 e_1^4 e_2 + 720 e_1^2 e_2^2 + 272 e_2^3 + 
 732 e_1^3 e_3 + 2304 e_1 e_2 e_3 
               \nonumber \\
     & & \quad\,
  +\, 540 e_3^2 + 1824 e_1^2 e_4 + 1248 e_2 e_4 + 1920 e_1 e_5 + 720 e_6
       \qquad
 \label{eq.euleriansymfun.0-6}
\end{subeqnarray}
Since the number of increasing ordered trees on $n+1$ vertices
is $(2n-1)!!$ \cite[p.~30, Corollary~1(iv)]{Bergeron_92},
we have $\sum_{\lambda \vdash n} a_\lambda = (2n-1)!!$;
and since the number of increasing $r$-ary trees
on $n+1$ vertices is the multifactorial $F_{n+1}^{(r-1)}$
\cite[p.~30, Example~1]{Bergeron_92},
we have
\be
   \sum_{\lambda \vdash n} a_\lambda \, \prod_i \binom{r}{\lambda_i}
   \;=\;
   F_{n+1}^{(r-1)}
   \;.
\ee
But these identities are insufficient to determine all the $a_\lambda$.

\begin{openproblem}
\rm
Find an explicit formula for the $a_\lambda$.
\end{openproblem}

We can, however, obtain an explicit expansion
for the Eulerian symmetric functions $\scrq^{(\infty)}_n$
and $\scrq^{(\infty)}_{n,k}$
in terms of the complete homogeneous symmetric functions~$h_\lambda$.
We use the method and notation of
Bergeron, Flajolet and Salvy \cite{Bergeron_92}.\footnote{
   See also \cite[Chapter~5]{Bergeron_98} for a more general context.
   Increasing ordered trees are treated in
   \cite[Example~5, pp.~364--365]{Bergeron_98}.
}
Let $Y_n = Y_n(\bphi)$ be the generating polynomial
for increasing ordered trees on $n$ vertices in which
each vertex with $i$ children gets a weight $\phi_i$, where $\phi_0 = 1$;
and more generally,
let $Y_{n,j} = Y_{n,j}(\bphi)$ be the generating polynomial
for unordered forests of increasing ordered trees on $n$ total vertices
with $j$ components in which
each vertex with $i$ children gets a weight $\phi_i$.
Here the weights $\bphi = (\phi_i)_{i \ge 1}$
are in the first instance indeterminates,
but they can later be specialized to values in any commutative ring
containing the rationals.
We will use the exponential generating functions
$Y(t) = \sum_{n=1}^\infty Y_n \, t^n/n!$
and $Y(t)^j/j! = \sum_{n=1}^\infty Y_{n,j} \, t^n/n!$
and the ordinary generating function
$\phi(w) = \sum_{k=0}^\infty\phi_k w^k$.
Then standard enumerative arguments \cite[Theorem~1]{Bergeron_92}
show that $Y(t)$ satisfies the ordinary differential equation
$Y'(t) = \phi(Y(t))$, leading to the implicit equation
\be
   t  \;=\; \int\limits_0^{Y(t)} {dw \over \phi(w)}
   \;.
 \label{eq.bergeron}
\ee
Introducing $\psi(w) = 1/\phi(w) = 1 + \sum_{i=1}^\infty \psi_i w^i$,
we then have
\be
   t   \;=\; Y(t) \, f(Y(t))
   \qquad\hbox{where}\qquad
   f(y) \;=\; 1 + \sum_{i=1}^\infty {\psi_i \over i+1} \, y^i
   \;.
 \label{eq.bergeron.fY}
\ee 
Solving $Y(t) = t/f(Y(t))$ by Lagrange inversion \reff{eq.lagrange.H}
with $H(u) = u^j/j!$ gives
\begin{subeqnarray}
   & &
   \hspace*{-1.2cm}
   Y_{n,j}
   \;=\;  n! \: [t^n] \, {Y(t)^j \over j!}
   \;=\; {(n-1)! \over (j-1)!} \: [y^{n-j}] \, f(y)^{-n}
        \\[2mm]
   & & \;
   =\;
   {(n-1)! \over (j-1)!}
    \!\! \sum_{\begin{scarray}
                   k_1, k_2, \ldots \ge 0 \\
                   \sum i k_i = n-j
                \end{scarray}}
          \!\!\!\!
          \binom{-n}{-n-\sum k_i,\, k_1,\, k_2,\, \ldots}
          \prod_{i=1}^\infty \Bigl( {\psi_i \over i+1} \Bigr) ^{\! k_i}
        \\[2mm]
   & & \;
   =\;
   {(n-1)! \over (j-1)!}
    \!\! \sum_{\begin{scarray}
                   k_1, k_2, \ldots \ge 0 \\
                   \sum i k_i = n-j
                \end{scarray}}
          \!\!
          (-1)^{\sum k_i} \,
          \binom{n+\sum k_i - 1}{n-1,\, k_1,\, k_2,\, \ldots}
          \prod_{i=1}^\infty \Bigl( {\psi_i \over i+1} \Bigr) ^{\! k_i}
   \,.
   \qquad
 \label{eq.bergeron.Yn}
\end{subeqnarray}
In our application we have $\phi_i = e_i$ and hence $\psi_i = (-1)^i h_i$.
This gives an explicit formula for the Eulerian symmetric functions:

\begin{proposition}[Eulerian symmetric functions in terms of complete homogeneous symmetric functions]
   \label{prop.eulerian.symfun.2}
We have
\begin{subeqnarray}
   \scrq^{(\infty)}_{n,k}(\bfx)
   \:=\:
   Y_{n+1,k+1}
   & \!=\! &
   {n! \over k!}
    \!\! \sum_{\begin{scarray}
                k_1, k_2, \ldots \ge 0 \\
                \sum i k_i = n-k
             \end{scarray}}
          \!\!
          \binom{n+\sum k_i}{n,\, k_1,\, k_2,\, \ldots}
          \prod_{i=1}^\infty \biggl( {(-1)^{i-1} h_i \over i+1} \biggr) ^{\! k_i}
    \qquad
       \\[2mm]
   & \!=\! &
   {1 \over k!}
   \sum_{\lambda \vdash n-k}
   \,
   (-1)^{n-k-\ell(\lambda)} \:
   {(n+\ell(\lambda))!  \over
    \displaystyle \prod\limits_i (i+1)^{m_i(\lambda)} \, m_i(\lambda)!
   }
   \; h_\lambda(\bfx)
      \;.
      \qquad\qquad
  \slabel{eq.prop.eulerian.symfun.2.b}
  \label{eq.prop.eulerian.symfun.2}
\end{subeqnarray}
\end{proposition}

\medskip

{\bf Remarks.}
1.  To see explicitly that the coefficient in
\reff{eq.prop.eulerian.symfun.2.b} is an integer, write it as
\be
   {\Bigl( k \,+\, \sum (i+1) m_i \Bigr)!
    \over
    k! \, \prod\limits_i (i+1)^{m_i} \, m_i!
   }
   \;=\;
   \binom{k \,+\, \sum (i+1) m_i}{k,\, 2m_1,\, 3m_2,\, \ldots}
   \prod\limits_i {((i+1)m_i)! \over (i+1)^{m_i} \, m_i!}
   \;.
\ee

2.  Computation of $\scrq_n^{(\infty)}$ for $n \le 10$
shows that some of the integers $a_\lambda$ in \reff{eq.prop.eulerian.symfun}
have large prime factors,
in contrast to the coefficients in \reff{eq.prop.eulerian.symfun.2},
which do not.
This suggests that there may not exist any simple formula
(at least not a multiplicative one) for the $a_\lambda$.

3.  Inspection of \reff{eq.euleriansymfun.0-6} shows that
our Eulerian symmetric functions $\scrq_n^{(\infty)}$ are {\em not}\/ equal
(or as far as we can tell, related)
to the ``Eulerian quasisymmetric functions''
(which are in fact symmetric!)\ defined by
Shareshian and Wachs \cite[Theorem~1.2]{Shareshian_10},
nor to their duals.

4.  Our multivariate Eulerian polynomials $\scrp_n^{(m)}$ and $\scrq_n^{(m)}$
are also not equal (or as far as we can tell, related)
to the multivariate Eulerian polynomials $A_n$, $\widetilde{A}_n$,
$C_n$ and $E_n$ defined by Haglund and Visontai \cite{Haglund_12}.
Indeed, in the latter polynomials the number of variables $\bfx$
is proportional to $n$, while in our polynomials it is a fixed number $m+1$.
\myendremark

\subsubsection{Extended multivariate Eulerian polynomials and Eulerian symmetric functions: Definitions and statement of results}
   \label{subsec.qusi-affine.extended}

But we can go much farther, and introduce additional indeterminates.
Given an increasing rooted tree $T$
(here it will be $(m+1)$-ary, but the notion is more general),
let us define the \textbfit{level} of a vertex $j \in T$
to be the number of children of the vertices $1,\ldots,j-1$ that are $> j$.
More generally, given an unordered forest $F$
consisting of $k+1$ increasing rooted trees,
let $r_j$ be the number of trees in $F$ that contain
at least one of the vertices $\{1,\ldots,j\}$;
then define the level of a vertex $j \in F$
to be the number of children of the vertices $1,\ldots,j-1$ that are $> j$,
plus $k+1 - r_j$.
(When the forest has only one tree,
 this definition reduces to the preceding one.)

Now let $\bfx = (x_0,\ldots,x_m)$ and $\bfc = (c_L)_{L \ge 0}$
be indeterminates,
and let $\scrq^{(m)}_n(\bfx,\bfc)$ be the generating polynomial for
increasing $(m+1)$-ary trees on $n+1$ vertices in which
each $i$-edge gets a weight $x_i$
and each vertex at level $L$ gets a weight $c_L$,
divided by $c_0$.
(Since vertex $n+1$ is at level 0,
 this definition can equivalently be rephrased as saying that
 each vertex at level $L$ {\em except for vertex $n+1$}\/
 gets a weight $c_L$, with no final division.
 In particular we see that $\scrq^{(m)}_n$ is a polynomial in $\bfc$.)
Similarly,
let $\scrq^{(m)}_{n,k}(\bfx,\bfc)$ be the generating polynomial for
unordered forests of increasing $(m+1)$-ary trees
on $n+1$ total vertices with $k+1$ components
in~which each $i$-edge gets a weight $x_i$
and each vertex at level $L$ gets a weight $c_L$,
divided by $c_0 c_1 \cdots c_k$.
(We will see later that $\scrq^{(m)}_{n,k}$ is indeed a polynomial in $\bfc$,
 among other things because it has an $m$-J-fraction with coefficients
 that are polynomials in $\bfc$.)
And finally, define $\scrp^{(m)}_n(\bfx,\bfc)$ by
\begin{subeqnarray}
   \scrp^{(m)}_0(\bfx,\bfc)  & = &    1   \\[1mm]
   \scrp^{(m)}_n(\bfx,\bfc)  & = &    c_0 x_0 \, \scrq^{(m)}_{n-1}(\bfx,\bfc)
        \qquad\hbox{for $n \ge 1$}
\end{subeqnarray}
Thus, $\scrp_n^{(m)}(\bfx,\bfc)$ is the same as $\scrq_n^{(m)}(\bfx,\bfc)$
but restricted to trees in which
the edge (if any) emanating from the root is a 0-edge.
Here $\scrp^{(m)}_n(\bfx,\bfc)$ and $\scrq^{(m)}_n(\bfx,\bfc)$
are homogeneous of degree~$n$ in $\bfx$
and also homogeneous of degree~$n$ in $\bfc$,
while $\scrq^{(m)}_{n,k}(\bfx)$ is homogeneous of degree~$n-k$ in $\bfx$
and in $\bfc$;
we refer to them collectively
as \textbfit{extended multivariate Eulerian polynomials}.
We will prove:


\begin{theorem}[Branched continued fractions for the extended multivariate Eulerian polynomials]
 \label{thm.factorized}
For every integer $m \ge 1$, we have:
\begin{itemize}
   \item[(a)] $\scrp^{(m)}_n(\bfx,\bfc) = \widehat{P}_n^{(m,m+1)}(\bfx,\bfc)$,
which by definition equals
$S_n^{(m)}(\balpha)$ where the weights $\balpha$ are given by
the factorized formula \reff{eq.alpha.factorized} with period $p=m+1$.
   \item[(b)] $\scrq^{(m)}_n(\bfx,\bfc) =  J_n^{(m)}(\bbeta)$
where the weights $\bbeta$ are given by
\be
   \beta_i^{(\ell)}
   \;=\;
   {(i+1)! \over (i-\ell)!} \:
   c_{i-\ell} c_{i-\ell+1} \cdots c_i
      \: e_{\ell+1}(x_0,\ldots,x_m)
   \;.
 \label{eq.thm.factorized}
\ee
   \item[(c)] $\scrq^{(m)}_{n,k}(\bfx,\bfc) =  J_{n,k}^{(m)}(\bbeta)$
where the weights $\bbeta$ are again given by \reff{eq.thm.factorized}.
\end{itemize}
\end{theorem}

As before, we get:

\begin{corollary}[Hankel-total positivity of the extended multivariate Eulerian polynomials]
  \label{cor.factorized.2}
For every integer ${m \ge 1}$,
the sequences $(\scrp^{(m)}_n(\bfx,\bfc))_{n \ge 0}$
and $(\scrq^{(m)}_n(\bfx,\bfc))_{n \ge 0}$
are coefficientwise Hankel-totally positive,
jointly in the indeterminates $\bfx$ and $\bfc$.
\end{corollary}

We will also show:

\begin{lemma}
   \label{lemma.Pm.beta.totalpos.factorized}
The production matrix $P^{(m)}(\bbeta)$ [cf.\ \reff{def.Pm}]
associated to the weights \reff{eq.thm.factorized}
is coefficientwise totally positive,
jointly in the indeterminates $\bfx$ and $\bfc$.
\end{lemma}

Once again this gives an alternate proof of the Hankel-total positivity of
$(\scrq^{(m)}_n(\bfx,\bfc))_{n \ge 0}$;
and by Theorem~\ref{thm.Jtype.generalized.minors} we get:

\begin{corollary}[Total positivity of lower-triangular matrix of extended multivariate Eulerian polynomials]
  \label{cor.factorized.2a}
For every integer ${m \ge 1}$,
the lower-triangular matrix $(\scrq^{(m)}_{n,k}(\bfx,\bfc))_{n,k \ge 0}$
is coefficientwise totally positive,
jointly in the indeterminates $\bfx$ and $\bfc$.
\end{corollary}

Of course, by specializing Theorem~\ref{thm.factorized},
Lemma~\ref{lemma.Pm.beta.totalpos.factorized}
and Corollaries~\ref{cor.factorized.2} and \ref{cor.factorized.2a}
to $\bfc = \bone$, we recover Theorem~\ref{thm.quasi-affine},
Lemma~\ref{lemma.Pm.beta.totalpos}
and Corollaries~\ref{cor.quasi-affine.2} and \ref{cor.quasi-affine.2a}.

Theorem~\ref{thm.factorized}(b,c) and
Corollaries~\ref{cor.factorized.2} and \ref{cor.factorized.2a}
also have obvious extensions to $m=\infty$, which we refrain from stating.
But it is worth introducing (with the obvious definition) the
\textbfit{extended Eulerian symmetric functions}
$\scrq^{(\infty)}_n(\bfx,\bfc)$ and $\scrq^{(\infty)}_{n,k}(\bfx,\bfc)$:
they are symmetric functions of $\bfx$
with coefficients in the polynomial ring $\Z[\bfc]$.
We then have the following generalization of
Proposition~\ref{prop.eulerian.symfun}:

\begin{proposition}[Extended Eulerian symmetric functions in terms of elementary symmetric functions]
   \label{prop.eulerian.symfun.factorized}
Let $\bfc = (c_L)_{L \ge 0}$ be indeterminates;
and for each partition $\lambda$ of~$n$,
let $a_\lambda(\bfc)$ be the sum,
over all increasing ordered trees on $n+1$ vertices
that have degree sequence $\lambda$,
of the product $\prod\limits_{j=1}^n c_{\lev(j)}$
where $\lev(j)$ is the level of vertex $j$
[note that $j=n+1$ is not included in this product].
Then
\be
   \scrq^{(\infty)}_n(\bfx,\bfc)
   \;=\;
   \sum_{\lambda \vdash n} a_\lambda(\bfc) \, e_\lambda(\bfx)
   \;.
\ee
\end{proposition}

\begin{openproblem}
\rm
Find an explicit formula for the $a_\lambda(\bfc)$,
and/or for the coefficients $b_\lambda(\bfc)$ in
$\scrq^{(\infty)}_n(\bfx,\bfc) =
   \sum\limits_{\lambda \vdash n} b_\lambda(\bfc) \, h_\lambda(\bfx)$.
\end{openproblem}

\bigskip

We are going to give two proofs of Theorems~\ref{thm.quasi-affine},
by completely different (but complementary) methods.
The first proof is by a bijection to labeled lattice paths:
this directly yields Theorem~\ref{thm.factorized}(b,c),
from which we deduce Theorem~\ref{thm.factorized}(a)
by applying odd contraction (Proposition~\ref{prop.contraction_odd});
then Theorem~\ref{thm.quasi-affine} is a special case.
The second proof is by the Euler--Gauss recurrence method:
this directly yields Theorem~\ref{thm.quasi-affine}(a),
from which we deduce Theorem~\ref{thm.quasi-affine}(b)
by applying odd contraction.
We have not yet succeeded in generalizing this second proof
to handle the case $\bfc \ne \bone$.
On the other hand, this second proof extends to handle the period-$m$ case
to be discussed in Section~\ref{subsec.quasi-affine.m},
for which we have been unable to complete the bijective proof.

\subsubsection{Proof by bijection to labeled lattice paths}
   \label{subsubsec.quasi-affine.bijective}

We are going to prove Theorem~\ref{thm.factorized}
by a bijection to labeled lattice paths.
This method for proving continued fractions
was invented by Flajolet \cite{Flajolet_80} and Viennot \cite{Viennot_83}
and used subsequently by many authors;
it can be summarized briefly as follows:\footnote{
   See also \cite{Sokal-Zeng_masterpoly} for a more detailed summary.
   What we here call a ``labeled Motzkin path''
   is (up~to small changes in notation)
   called a ``path diagramme'' by Flajolet \cite[p.~136]{Flajolet_80}
   and a ``history'' by Viennot \cite[p.~II-9]{Viennot_83}.
   Our set $\bfL$ of allowed labels is essentially what
   Flajolet \cite[p.~136]{Flajolet_80} calls a ``possibility function''.
}
For each $n \ge 0$, let $\scro_n$ be some class of combinatorial objects,
whose ordinary generating function (possibly with some weights)
we wish to expand as a classical continued fraction.
Now let $\scrm_n$ be the set of Motzkin paths of length~$n$:
we write them as $\omega = (\omega_0,\ldots,\omega_n)$
with $\omega_j = (j,h_j)$;
the steps are $s_j = h_j - h_{j-1} \in \{-1,0,1\}$.
For each $s \in \{-1,0,1\}$ and each $h \ge 0$,
let $L(s,h)$ be a finite (possibly empty) set of ``allowed labels''
for a step of type~$s$ starting at height~$h$;
we write $\bfL = (L(s,h))_{s \in \{-1,0,1\} ,\, h \ge 0}$.
Then an {\em $\bfL$-labeled Motzkin path}\/ of length~$n$
is a pair $(\omega,\xi)$ where $\omega \in \scrm_n$
and $\xi = (\xi_1,\ldots,\xi_n)$ is a set of labels satisfying
$\xi_j \in L(s_j,h_{j-1})$ for $1 \le j \le n$;
we write $\scrm^{\bfL}_n$ for the set of
$\bfL$-labeled Motzkin paths of length~$n$.
If we can find a bijection from $\scro_n$ to $\scrm^{\bfL}_n$,
then enumerating $\scro_n$ is equivalent to
enumerating Motzkin paths of length~$n$ with a weight $|L(s,h)|$
for each step of type~$s$ starting at height~$h$.
And by Flajolet's \cite{Flajolet_80} correspondence
of weighted Motzkin paths with J-fractions,
this means that the ordinary generating function
$\sum\limits_{n=0}^\infty |\scro_n| \, t^n$
is given by a J-fraction
with coefficients $\gamma_h = |L(0,h)|$
and $\beta_h = |L(1,h-1)| \: |L(-1,h)|$.
Moreover, if the elements of $\scro_n$ carry weights that,
under the bijection, correspond to a product of weights for each step
that moreover depend only on the step's type $s_j$,
its starting height $h_{j-1}$ and its label $\xi_j$,
then the weighted ordinary generating function
is also given by a J-fraction.
Similar considerations apply for mappings to other classes
of labeled lattice paths:
for instance, a mapping to labeled Dyck paths of length $2n$
gives rise to an S-fraction \reff{eq.f0.Sfrac},
and more generally a mapping to labeled $m$-Dyck paths of length $(m+1)n$
gives rise to an $m$-S-fraction \reff{eq.fk.mSfrac}.

\def\disco{\circle*{2}}
\newcommand{\mytree}{
  \setlength{\unitlength}{0.7mm}
  \begin{picture}(15,10)(-3,-4)
    \thicklines
    \put(0,-3){\disco}
    \put(0,-3){\line(3,5){5}}
    \put(5,5.7){\disco}
    \put(5,5.7){\line(3,-5){5}}
    \put(10,-3){\disco}
    \put(10,-3){\line(3,-5){5}}
    \put(15,-11.7){\disco}
  \end{picture}
  \vphantom{{\Huge $\int^\int_{A_{A_A}}$}}
}

Now, in Section~\ref{subsec.periodic.combinatorial}
we gave a combinatorial proof of
Proposition~\ref{prop.fuss-narayana.trees.bis}
by constructing a bijection from $(m+1)$-ary trees on $n$ vertices
to $m$-Dyck paths of length $(m+1)n$.
It would be natural to try to prove Theorem~\ref{thm.factorized}(a)
by constructing a bijection from
{\em increasing}\/ $(m+1)$-ary trees on $n$ vertices
to some set of suitably {\em labeled}\/ $m$-Dyck paths $(\omega,\xi)$;
and it would moreover be natural to define the path $\omega$
by re-using the same bijection that we have already devised
(using the tree without looking at its labels),
and then to define the labels $\xi$
by using in some way the labels on the tree.
But it turns out that this approach cannot work, not even for $m=1$.
Indeed, consider {\em any}\/ bijection from binary trees on $n$ vertices
to Dyck paths of length $2n$ (not necessarily the one constructed
in Section~\ref{subsec.periodic.combinatorial});
and consider the binary tree \mytree on 4~vertices
(the root is at the top).
This tree has 3 increasing labelings:
the root must be labeled 1, and then the left child of the root
can be labeled 2, 3 or 4.
But no matter how we map this tree to a Dyck path of length~8,
that path can reach height at most~4.
And the label sets $L(\pm 1, h)$ for $0 \le h \le 4$ must satisfy
$|L(1,h-1)| \: |L(-1,h)| = \alpha_h$
where $(\alpha_1,\alpha_2,\alpha_3,\alpha_4) = (1,1,2,2)$
from Euler's continued fraction \reff{eq.nfact.contfrac} for~$n!$.
But this means that each $|L(\pm 1, h)|$ must be 1 or 2,
so it is impossible for the Dyck path associated to this tree
to have 3 possible labelings.
So the approach of defining the Dyck path by looking only
at the underlying unlabeled tree cannot work.


The next-best thing would be to find a bijection from
increasing $(m+1)$-ary trees on $n$ vertices
to some set of suitably labeled $m$-Dyck paths $(\omega,\xi)$,
where not only $\xi$ but also $\omega$ would depend on the tree's labels.
We have not yet succeeded in doing this,
and we leave it as an open problem to find such a bijection
(which would moreover have to map the statistics in such a way
as to prove Theorem~\ref{thm.factorized}(a)
or at least Theorem~\ref{thm.quasi-affine}(a)).
Instead, we will prove Theorem~\ref{thm.factorized}
by constructing a bijection from
increasing $(m+1)$-ary trees on $n+1$~vertices
to a set of labeled reversed $m$-\L{}ukasiewicz paths of length $n$.
(A reversed $m$-\L{}ukasiewicz path is a path in $\Z \times \N$
 that uses steps $(1,r)$ with $-1 \le r \le m$.)
This will show that the corresponding ordinary generating function
is given by an $m$-J-fraction,
namely, the one asserted in Theorem~\ref{thm.factorized}(b).
We will then observe that the weights occurring in this $m$-J-fraction
are actually those arising from the
odd contraction \reff{eq.prop.contraction_odd} of an $m$-S-fraction,
namely, the one asserted in Theorem~\ref{thm.factorized}(a).

By generalizing the bijection from trees to forests
we will also prove Theorem~\ref{thm.factorized}(c).
But there is a slight twist:  the bijection works most naturally
for {\em ordered}\/ forests of increasing $(m+1)$-ary trees,
not unordered ones.
As mentioned earlier, the generating polynomial for ordered forests
--- let us call it $\widetilde{\scrq}^{(m)}_{n,k}(\bfx,\bfc)$ ---
is simply $(k+1)!$ times the generating polynomial
$\scrq^{(m)}_{n,k}(\bfx,\bfc)$ for unordered forests.
So we will show that $\widetilde{\scrq}^{(m)}_{n,k}(\bfx,\bfc)$
is $(k+1)!$ times the generalized $m$-Jacobi--Rogers polynomial
$J_{n,k}^{(m)}(\bbeta)$;
this will immediately prove the desired result
for $\scrq^{(m)}_{n,k}(\bfx,\bfc)$.

\medskip

\proofof{Theorem~\ref{thm.factorized}}
We will construct a bijection from
the set of ordered forests of increasing $(m+1)$-ary trees
on the vertex set $[n+1]$ with $k+1$ components
to the set of $\bfL$-labeled reversed partial $m$-\L{}ukasiewicz paths
from $(0,k)$ to $(n,0)$,
where the label sets $\bfL$ will be defined below.\footnote{
   When restricted to trees ($k=0$),
   our bijection generalizes the one
   described by Flajolet \cite[p.~143]{Flajolet_80}
   (following Fran\c{c}on \cite{Francon_78})
   for the case $m=1$ (increasing binary trees).
   Also, after completing this proof we discovered
   that essentially the same construction
   is contained (for general $m$, though only for trees)
   in a very interesting unpublished paper of Kuba and Varvak \cite{Kuba_09}.
   Indeed, Kuba and Varvak use somewhat more general weights than we do
   --- encoding the complete ``local type'' of children at each vertex ---
   because they are not concerned with total positivity.
}

\begin{figure}[!ht]
\begin{center}
(a) \qquad\qquad\qquad\quad
           \includegraphics[scale=1.5]{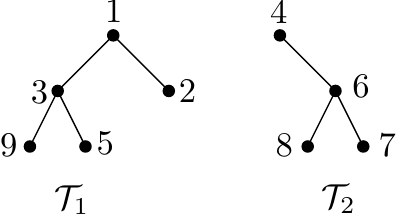}  \\[5mm]
(b) \qquad \includegraphics[scale=1.5]{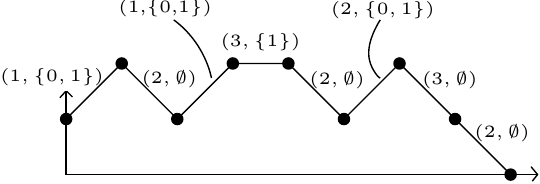}
\caption{
    (a) An ordered forest $(\mathcal T_1, \mathcal T_2)$
    of increasing binary trees on 9 total vertices.
    (b) Its corresponding labelled reversed partial 1-\L{}ukasiewicz path
    from $(0,1)$ to $(8,0)$.
}
\label{fig.bijection.forest.path}
\end{center}
\end{figure}

Given an ordered forest $F$ of increasing $(m+1)$-ary trees
on the vertex set $[n+1]$ with $k+1$ components,
we define a labeled reversed partial $m$-\L{}ukasiewicz path $(\omega,\xi)$
of length $n$ as follows (see Figure~\ref{fig.bijection.forest.path} for an example):

\medskip

{\bf Definition of the path $\bm{\omega}$.}
The path $\omega$ starts at height $h_0 = k$
and takes steps $s_1,\ldots,s_n$ with $s_i = \deg(i) - 1$,
where $\deg(i)$ is the number of children of vertex $i$
(we don't take any step $s_{n+1}$);
we have $-1 \le s_i \le m$ since $0 \le \deg(i) \le m+1$.\footnote{
   We remark that this definition of the path $\omega$
   is identical to what was done in Section~\ref{sec.mJR.trees}
   for unlabeled trees,
   except that here we use the given labeling of the vertices
   rather than the depth-first-search labeling.
   The proof of the interpretation of the height $h_j$
   given in Section~\ref{sec.mJR.trees} is simply a specialization
   of the one given in Lemma~\ref{lemma.hj} below.
   However, with the depth-first-search labeling
   we always finish processing one tree before moving on to the next,
   while that is {\em not}\/ in general the case here.
}
Recall that $r_j$ is the number of trees in $F$ that contain
at least one of the vertices $\{1,\ldots,j\}$.
We then claim:

\begin{lemma}
   \label{lemma.hj}
For $1 \le j \le n$,
the height $h_j = k+ \sum_{i=1}^j s_i$ has the following interpretations:
\begin{itemize}
   \item[(a)]  $h_j$ is the number of children of the vertices $\{1,\ldots,j\}$
      whose labels are $> j$, plus $k-r_j$.
   \item[(b)]  $h_{j-1}$ is the number of children of the vertices
      $\{1,\ldots,j-1\}$ whose labels are $> j$, plus $k+1-r_j$.
\end{itemize}
In particular, $h_j > k-r_j$ whenever $j$
is not the highest-numbered vertex of its tree,
and $h_j \ge k-r_j$ always;
and $h_{j-1}$ is the level of the vertex $j$ as we have defined it.
\end{lemma}

\proof
By induction on $j$.
For the base case $j=1$, the claims are clear
since $r_1 = 1$, $h_0 = k$ and $h_1 = k + \deg(1) - 1$.

For $j>1$, vertex $j$ is either the child of another node,
or the root of a tree.  We consider these two cases separately:

(i) Suppose that $j$ is the child of another node
(obviously numbered ${\le j-1}$).
By the inductive hypothesis~(a),
$h_{j-1}$ is the number of children of the vertices ${\{1,\ldots,j-1\}}$
whose labels are $\ge j$, plus $k-r_{j-1}$;
since one of these children is $j$,
it follows that
$h_{j-1} - 1$ is the number of children of the vertices $\{1,\ldots,j-1\}$
whose labels are $> j$, plus $k-r_{j-1}$.
Now vertex $j$ has $\deg(j)$ children, all of which have labels $> j$;
so $h_j = h_{j-1} + s_j = h_{j-1} - 1 + \deg(j)$
is the number of children of the vertices $\{1,\ldots,j\}$
whose labels are $> j$, plus $k-r_{j-1}$.
Since $r_j = r_{j-1}$, the preceding two sentences prove
claims (b) and (a), respectively.

(ii) Suppose that $j$ is a root.
By the inductive hypothesis~(a),
$h_{j-1}$ is the number of children of the vertices $\{1,\ldots,j-1\}$
whose labels are $\ge j$, plus $k-r_{j-1}$;
since $j$ is not one of these children,
$h_{j-1}$ is the number of children of the vertices $\{1,\ldots,j-1\}$
whose labels are $> j$, plus $k-r_{j-1}$.
Now all of the children of vertex $j$ have labels $> j$,
so $h_j = h_{j-1} +  \deg(j) - 1$
is the number of children of the vertices $\{1,\ldots,j\}$
whose labels are $> j$, plus $k-r_{j-1}-1$.
Since $r_j = r_{j-1} +1$, the preceding two sentences prove 
claims (b) and (a), respectively.
\qed

It follows from Lemma~\ref{lemma.hj}
that $h_1,\ldots,h_{n-1} \ge 0$ and $h_n = 0$.
So the path $\omega$ is indeed a reversed partial $m$-\L{}ukasiewicz path
from $(0,k)$ to $(n,0)$.

\medskip

{\bf Definition of the labels $\bm{\xi}$.}
The label $\xi_j$ will be an ordered pair $\xi_j = (\xi'_j,\xi''_j)$
where $\xi'_j$ is a positive integer
and $\xi''_j$ is a subset of $\{0,\ldots,m\}$;
more precisely, the label set $L(s,h)$
for a step $s$ starting at height $h$ will be
\be
   L(s,h)  \;=\;  [h+1] \,\times\, \binom{\{0,\ldots,m\}}{s+1}
\ee
where $\binom{S}{r}$ denotes the set of subsets of $S$ of cardinality $r$.
The label $\xi'_j$ is defined as 1 plus the number of vertices $> j$
that are either children of $\{1,\ldots,j-1\}$ or roots
and that precede $j$ in the depth-first-search order.\footnote{
   Here the depth-first-search order could be replaced by any chosen order
   on the vertices of $F$ that commutes with truncation.
   The key property we need is that
   the order on the truncated forest $F_{j-1}$ to be defined below
   is the restriction of the order on the full forest $F$.
}
The label $\xi''_j$ is the set of $i \in \{0,\ldots,m\}$
for which $j$ has an $i$-child.

\medskip

{\bf The inverse bijection.}
We claim that this mapping $F \mapsto (\omega,\xi)$ is a bijection
from the set of ordered forests of increasing $(m+1)$-ary trees
on the vertex set $[n+1]$ with $k+1$ components
to the set of $\bfL$-labeled reversed partial $m$-\L{}ukasiewicz paths
from $(0,k)$ to $(n,0)$.
To prove this, we explain the inverse mapping.

Given a labeled reversed partial $m$-\L{}ukasiewicz path $(\omega,\xi)$,
we build up the ordered forest $F$ vertex-by-vertex:
at stage $j$ we will have an ordered forest $F_j$
in which some of the vertices are labelled $1,\ldots,j$
and some others are unnumbered ``vacant slots''.
The starting forest $F_0$ has $k+1$ singleton components,
each of which is a vacant slot (these components are of course ordered).
We now ``read'' the path step-by-step, from $j=1$ through $j=n$.
When we read a step $s_j$ with labels $(\xi'_j,\xi''_j)$,
we insert a new vertex $j$ into one of the vacant slots of $F_{j-1}$,
namely, the $\xi'_j$th vacant slot
in the depth-first-search order of $F_{j-1}$.
We also create new vacant slots that are children of $j$:
namely, a vacant $i$-child for each $i \in \xi''_j$.
This defines $F_j$.
After stage $n$, the forest $F_n$ has only one vacant slot:
we insert vertex $n+1$ into this slot, thereby defining
the ordered forest $F$.


It is fairly clear that this insertion algorithm
defines a map $(\omega,\xi) \mapsto F$
that is indeed the inverse of the mapping $F \mapsto (\omega,\xi)$
defined previously.

\bigskip

Now recall that we want to enumerate ordered forests of
increasing $(m+1)$-ary trees on the vertex set $[n+1]$
with $k+1$ components in which
each $i$-edge gets a weight $x_i$
and each vertex at level $L$ gets a weight $c_L$
(divided by $c_0 \cdots c_k$).
We use the bijection to push these weights from the forests
to the labeled reversed partial $m$-\L{}ukasiewicz paths.
Given a forest $F$, each vertex $j \in [n+1]$
contributes a weight $c_{\lev(j)}$ [where $\lev(j)$ is its level]
and also a weight $x_i$ for each $i$-child.
Under the bijection, this vertex (if $j \le n$)
is mapped to a step $s_j = \deg(j) - 1$
from height $h_{j-1} = \lev(j)$ to height $h_j = h_{j-1} + s_j$.
Therefore, the weight in the labeled path $(\omega,\xi)$
corresponding to this vertex is 
$c_{h_{j-1}} \prod\limits_{i \in \xi''_j} x_i$,
and the weight of the labeled path $(\omega,\xi)$
[before division by $c_0 \cdots c_k$]
is the product of these weights over $1 \le j \le n$,
times the weight $c_0$ corresponding to the vertex $n+1$.

Now we sum over the labels $\xi$ to get the total weight
for each path $\omega$.
Summing over $\xi'_j$ gives a factor $h_{j-1} + 1$.
Summing over $\xi''_j$ gives a weight $e_r(x_0,\ldots,x_m)$
when vertex $j$ has $r$ ($= s_j + 1$) children.
Putting everything together,
the weight in the reversed partial $m$-\L{}ukasiewicz path
for a step~$s$ ($-1 \le s \le m$) starting from height~$h$ will be 
\be
   W(s,h) \;=\;  (h+1) \, c_h \, e_{s+1}(x_0,\ldots,x_m)
   \;.
 \label{eq.Wsh}
\ee

We now want to read this path backwards
(so that it becomes an ordinary partial $m$-\L{}ukasiewicz path
 in our definition);
then a step~$s$ starting at height~$h$
becomes a step~$s' = -s$ starting at height~$h' = h+s$.
Therefore, in the ordinary partial $m$-\L{}ukasiewicz path,
the weight will be
\be
   W'(s',h')  \;=\;  (h'+s'+1) \, c_{h'+s'} \, e_{1-s'}(x_0,\ldots,x_m)
   \;.
 \label{eq.Wsh.prime}
\ee

Finally, in the partial $m$-\L{}ukasiewicz path
we can associate to each $\ell$-fall from height $h'$
the corresponding rises from $h'-\ell \to h'-\ell+1 \to \ldots \to h'$;
and in addition there are rises $0 \to 1 \to \ldots \to k$
(exactly one of each type) that do not correspond to falls.
By our convention for $m$-Jacobi--Rogers polynomials
(which gives weight 1 to rises),
we must attribute the weights \reff{eq.Wsh.prime} for the rises
$h'-\ell \to h'-\ell+1 \to \ldots \to h'$
to the corresponding $\ell$-fall.
The resulting weight for an $\ell$-fall from height $i$ is thus
\begin{subeqnarray}
   \beta_i^{(\ell)}
   & = &
   W'(1, i-\ell) \, W'(1, i-\ell+1) \,\cdots\, W'(1,i-1) \:\times\: W(-\ell,i)
       \qquad \\[2mm]
   & = &
   (i-\ell+2) c_{i-\ell+1} \;
   (i-\ell+3) c_{i-\ell+2} \;\cdots\;
   (i+1) c_i
           \nonumber \\
   & & \hspace*{3cm}
   \;\times\;
   (i-\ell+1) c_{i-\ell} \: e_{\ell+1}(x_0,\ldots,x_m)
   \;,
\end{subeqnarray}
which is exactly \reff{eq.thm.factorized}.
Furthermore, the weight for the rises $0 \to 1 \to \ldots \to k$
that do not correspond to falls is
\be
   W'(1,0) \, W'(1,1) \,\cdots\, W'(1,k-1)
   \;=\;
   2 c_1 \: 3c_2 \:\cdots\: (k+1) c_k
   \;=\;
   (k+1)! \, c_1 \cdots c_k
   \;.
\ee
And finally, let us not forget the weight $c_0$ for the vertex $n+1$.
But now we recall that $\widetilde{Q}_{n,k}^{(m)}(\bfx,\bfc)$
was defined as the weighted sum over ordered forests
{\em divided by $c_0 \cdots c_k$}\/.
Therefore, we get
\be
   \widetilde{Q}_{n,k}^{(m)}(\bfx,\bfc)
   \;=\;
   (k+1)! \, J_{n,k}^{(m)}(\bbeta)
   \;.
\ee
Since
$\widetilde{Q}_{n,k}^{(m)}(\bfx,\bfc) = (k+1)! \, Q_{n,k}^{(m)}(\bfx,\bfc)$,
Theorem~\ref{thm.factorized}(c) follows.
And Theorem~\ref{thm.factorized}(b) is the special case $k=0$.


On the other hand, Lemma~\ref{lemma.contraction.factorized} below
with $b_k = (k+1) c_k$
shows that these $\bbeta$ are precisely those arising by
odd contraction \reff{eq.prop.contraction_odd}
from the $m$-S-fraction
specialized to the factorized weights \reff{eq.alpha.factorized} with $p=m+1$.
Then Proposition~\ref{prop.contraction_odd} tells us that
$S_n^{(m)}(\balpha) = \alpha_m J_{n-1}^{(m)}(\bbeta)
 = c_0 x_0 \scrq_{n-1}^{(m)}(\bfx,\bfc) = \scrp_n^{(m)}(\bfx,\bfc)$.
This proves Theorem~\ref{thm.factorized}(a).
\qed

\begin{lemma}
   \label{lemma.contraction.factorized}
Consider the production matrix $P^{(m)} = (p^{(m)}_{i,j})_{i,j \ge 0}$
defined by \reff{eq.prop.contraction_odd}
with $\alpha_{m+j+(m+1)k} = b_k x_j$ for $0 \le j \le m$ and $k \ge 0$:
\be
   P^{(m)}
   \;=\;
   U^\star(b_0 x_0, b_1 x_0, b_2 x_0, \ldots)
   \:
   L(b_0 x_1, b_1 x_1, b_2 x_1, \ldots)
   \:\cdots\:
   L(b_0 x_m, b_1 x_m, b_2 x_m, \ldots)
   \;.
\ee
Then
\be
   p^{(m)}_{i,j}
   \;=\;
      b_j b_{j+1} \cdots b_i
          \, e_{i-j+1}(x_0,\ldots,x_m)
\ee
where of course $e_{i-j+1}(x_0,\ldots,x_m) = 0$
whenever $j < i-m$ or $j > i+1$.
\end{lemma}

\proof
By induction on $m$.
If $m=0$, then
\be
   p_{i,j}^{(0)}
   \;=\;
   U^\star(b_0 x_0, b_1 x_0,\ldots)
   \;=\;
   \begin{cases}
        b_i x_0 \;=\; b_i e_1(x_0) &  \textrm{if $j=i$} \\
        1 \;=\; e_0(x_0)           &  \textrm{if $j=i+1$}  \\
        0                          &  \textrm{otherwise}
   \end{cases}
\ee
as claimed.
For $m >0$, we have by definition
$P^{(m)} = P^{(m-1)} \,  L(b_0 x_m, b_1 x_m,\ldots)$
and hence
\begin{subeqnarray}
   p^{(m)}_{i,j}
   & = &
   p^{(m-1)}_{i,j} \:+\:  b_j x_m \, p^{(m-1)}_{i,j+1}. 
         \\[2mm]
   & = &
   b_j b_{j+1} \cdots b_i \, e_{i-j+1}(x_0,\ldots,x_{m-1})
      \:+\:
   b_j x_m \, b_{j+1} \cdots b_i \, e_{i-j}(x_0,\ldots,x_{m-1})
         \nonumber \\ \\
   & = &
    b_j b_{j+1} \cdots b_i \, e_{i-j+1}(x_0,\ldots,x_m)
\end{subeqnarray}
where the second line used the inductive hypothesis.
%
%
%
%
%
\qed

Finally, let us prove Lemma~\ref{lemma.Pm.beta.totalpos.factorized}
(and hence also Lemma~\ref{lemma.Pm.beta.totalpos}):

\proofof{Lemma~\ref{lemma.Pm.beta.totalpos.factorized}}
The production matrix $P(\bbeta)$ associated to the weights
\reff{eq.thm.factorized} has matrix elements
\be
   p_{ij}  \;=\;  \beta_i^{(i-j)}
           \;=\; \widehat{c}_j \widehat{c}_{j+1} \cdots \widehat{c}_i
                 \: e_{i-j+1}(x_0,\ldots,x_m)
 \label{eq.Pbeta.ij}
\ee
where $\widehat{c}_k \eqdef (k+1) c_k$;
therefore it can be written as $P(\bbeta) = DTD'$
where $D = \diag(\widehat{c}_0,\, \widehat{c}_0 \widehat{c}_1,\,
                      \widehat{c}_0 \widehat{c}_1 \widehat{c}_2,\,\ldots)$,
$D' = \diag(1,\, 1/\widehat{c}_0,\, 1/(\widehat{c}_0 \widehat{c}_1),\,\ldots)$,
and $T$ is the unit-lower-Hessenberg Toeplitz matrix
$(e_{i-j+1}(\bfx))_{i,j \ge 0}$.
(Here we are working temporarily in the ring $\Q[\bfx,\bfc,\bfc^{-1}]$,
although in the end everything will lie in $\Z[\bfx,\bfc]$.)
By the Jacobi--Trudi formula \cite[Corollary~7.16.2]{Stanley_99},
all the minors of $T$ are skew Schur functions and hence monomial-positive
(in $\bfx$);
it follows that $P(\bbeta)$ is coefficientwise totally positive
(in $\bfx$ and $\bfc$).
\qed

{\bf Remark.}
Since the skew Schur functions are nonnegative linear combinations
of Schur functions \cite[Corollary~7.15.9]{Stanley_99},
this proof shows that in fact
$P(\bbeta)$ is totally positive with respect to the Schur order
on the ring of symmetric functions,
which is stronger than the coefficientwise (= monomial) order.\footnote{
      {\em A priori}\/ the foregoing proof shows only that the minors of
      $P(\bbeta)$ are linear combinations
      of Schur functions with nonnegative {\em rational}\/ coefficients
      (more precisely, coefficients in $\Q_+[\bfc]$).
      But it is an immediate consequence of \reff{eq.Pbeta.ij}
      that the minors of $P(\bbeta)$ are {\em integer}\/ linear combinations
      of the $e_\lambda$
      (more precisely, linear combinations with coefficients in $\Z[\bfc]$),
      and hence also of the $m_\lambda$;
      and the transition matrix from $m_\lambda$ to $s_\lambda$
      is integer-valued \cite[pp.~315--316]{Stanley_99}.
      So the minors of $P(\bbeta)$
      are linear combinations of Schur functions
      with nonnegative {\em integer}\/ coefficients
      (more precisely, coefficients in $\Z_+[\bfc]$).
}
Theorem~\ref{thm.Jtype.minors} then implies that
the sequence $(\scrq^{(\infty)}_n(\bfx))_{n \ge 0}$
is Hankel-totally positive with respect to the Schur order
(coefficientwise in $\bfc$);
and Theorem~\ref{thm.Jtype.generalized.minors} implies that
the lower-triangular matrix $(\scrq^{(\infty)}_{n,k}(\bfx,\bfc))_{n,k \ge 0}$
is totally positive with respect to the Schur order
(coefficientwise in $\bfc$).
\myendremark

%
%
%
%

\subsubsection{Proof by the Euler--Gauss recurrence method}
   \label{subsubsec.quasi-affine.euler-gauss}

We would now like to give a second proof of Theorem~\ref{thm.quasi-affine}(a,b),
based on the Euler--Gauss recurrence method for proving continued fractions,
generalized to $m$-S-fractions as in Proposition~\ref{prop.euler-gauss.mSR}.
This proof gives additional information not provided by the bijective proof:
namely, recurrences (and when $\bfx = \bone$, semi-explicit formulae)
for the quantities $g_k = f_0 f_1 \cdots f_k$.

Let us recall briefly the method:
if $(g_k(t))_{k \ge -1}$ are formal power series with constant term 1
(with coefficients in some commutative ring $R$)
satisfying a recurrence
\be
   g_k(t) - g_{k-1}(t)  \;=\; \alpha_{k+m} t \, g_{k+m}(t)
   \qquad\hbox{for } k \ge 0
 \label{eq.recurrence.gkm.0.bis}
\ee
for some coefficients $\balpha = (\alpha_i)_{i \ge m}$ in $R$,
then $g_0(t)/g_{-1}(t) = \sum_{n=0}^\infty S^{(m)}_n(\balpha) \, t^n$
and more generally
$g_k(t)/g_{-1}(t) = \sum_{n=0}^\infty S^{(m)}_{n|k}(\balpha) \, t^n$.

Here we will apply this method with the choice $g_{-1}(t) = 1$.
(By contrast, the application of this method to
 ratios of contiguous hypergeometric series,
 to be given in Sections~\ref{sec.hyper.rF0}--\ref{sec.hyper.rphis} below,
 will use in general $g_{-1}(t) \ne 1$.)
Here the weights $\balpha$ are the Eulerian-quasi-affine weights
\reff{eq.alpha.quasi-affine.u=x} with period $p=m+1$, i.e.
\be
   \alpha_k  \;=\;  \lfloor (k+1)/(m+1) \rfloor \, x_{k+1 \bmod m+1}
   \;.
 \label{eq.alphak.periodm+1}
\ee
We need to find series $(g_k(t))_{k \ge 0}$ with constant term 1
satisfying \reff{eq.recurrence.gkm.0.bis}.
Let us write $g_k(t) = \sum_{n=0}^\infty g_{k,n} \, t^n$
and define $g_{k,n} = 0$ for $n < 0$.
Then \reff{eq.recurrence.gkm.0.bis} can be written as
\be
   g_{k,n} \,-\, g_{k-1,n}  \;=\;  \alpha_{k+m} \, g_{k+m,n-1}
   \qquad\hbox{for } k,n \ge 0
   \;.
 \label{eq.recurrence.gkm.0.bis.gkn}
\ee
Here are the required $g_{k,n}$:

\begin{proposition}[Euler--Gauss recurrence for multivariate Eulerian polynomials]
   \label{prop.euler-gauss.quasi-affine}
Let $\bfx = (x_0,\ldots,x_m)$ be indeterminates;
we work in the ring $R = \Z[\bfx]$.
Set $g_{k,n} = \delta_{n0}$ for $k < 0$,
and then define $g_{k,n}$ for $k,n \ge 0$ by the recurrence
\be
   g_{k,n}
   \;=\;
   \Bigl( \scrd_m \,+\, \sum_{i=1}^m \alpha_{k+i} \Bigr) g_{k,n-1}
       \:+\:  g_{k-m,n}
 \label{def.gkn.quasi-affine}
\ee
where $\balpha$ are given by \reff{eq.alphak.periodm+1}
and $\scrd_m$ is given by \reff{def.dm}.
Then:
\begin{itemize}
   \item[(a)]  $g_{k,0} = 1$ for all $k \in \Z$.
\\[-5mm]
   \item[(b)]  $(g_{k,n})$ satisfies the recurrence
      \reff{eq.recurrence.gkm.0.bis.gkn}.
\\[-5mm]
   \item[(c)]  $(g_{k,n})$ also satisfies the recurrence
\be
   g_{k,n}
   \;=\;
   \Bigl( \scrd_m \,+\, \sum_{i=0}^m \alpha_{k+i} \Bigr) g_{k,n-1}
       \:+\:  g_{k-m-1,n}
   \qquad\hbox{for } k,n \ge 0
   \;.
 \label{def.gkn.quasi-affine.variant}
\ee
\\[-13mm]
   \item[(d)]  For $0 \le k \le m$, we have
$\displaystyle
   g_{k,n}
   \:=\:
   \Bigl( \scrd_m \,+\, \sum_{i=0}^k x_i \Bigr)^{\! n} \: 1
$.
\end{itemize}
Therefore $S^{(m)}_n(\balpha) = g_{0,n} = (\scrd_m + x_0)^n \, 1$,
and more generally $S^{(m)}_{n|k}(\balpha) = g_{k,n}$.
\end{proposition}

Here \reff{def.gkn.quasi-affine} is, for each $k \ge 0$,
a recurrence in $n$ for $g_{k,n}$, with $g_{k-m,n}$ being considered known;
it should therefore be interpreted as an outer recurrence on $k$
together with an inner recurrence on $n$.
Note also that we are using \reff{eq.alphak.periodm+1}
also to define $\alpha_i = 0$ for $-1 \le i \le m-1$.
This implies that \reff{eq.recurrence.gkm.0.bis.gkn}
holds also for $-(m+1) \le k \le -1$.

\medskip

\proofof{Proposition~\ref{prop.euler-gauss.quasi-affine}}
(a)  We see trivially using \reff{def.gkn.quasi-affine} that
$g_{k,0} = 1$ for all $k \in \Z$, i.e.\ $g_k(t)$ has constant term 1.

(b) We will now prove that \reff{eq.recurrence.gkm.0.bis.gkn} holds.
The proof will be by an outer induction on $k$ and an inner induction on $n$.
The base cases $k = -m,\ldots,-1$ hold trivially
because $g_{k,n} = g_{k-1,n} = \delta_{n0}$
and $\alpha_0,\ldots,\alpha_{m-1} = 0$.
Now suppose that \reff{eq.recurrence.gkm.0.bis.gkn} holds
for a given $k$ and all $n \ge 0$;
we want to prove it when $k$ is replaced by $k+m$, i.e.\ that
\be
   g_{k+m,n} \:-\: g_{k+m-1,n} \:-\: \alpha_{k+2m} \, g_{k+2m,n-1}
   \;=\;
   0
   \quad\hbox{for all $n \ge 0$}
   \;.
 \label{eq.proof.prop.quasi-affine.2}
\ee
We will prove \reff{eq.proof.prop.quasi-affine.2} by induction on $n$.
It holds for $n=0$ because $g_{k+m,0} = g_{k+m-1,0} = 1$ and $g_{k+2m,-1} = 0$.
If $n > 0$, then we use \reff{def.gkn.quasi-affine}
on each of the three terms on the left-hand side of
\reff{eq.proof.prop.quasi-affine.2}, giving
\begin{subeqnarray}
   g_{k+m,n}
   & = &
   \Bigl( \scrd_m \,+\, \sum_{i=1}^m \alpha_{k+m+i} \Bigr)  g_{k+m,n-1}
       \:+\:  g_{k,n}
       \\
   g_{k+m-1,n}
   & = &
   \Bigl( \scrd_m \,+\, \sum_{i=1}^m \alpha_{k+m-1+i} \Bigr)  g_{k+m-1,n-1}
       \:+\:  g_{k-1,n}
       \\
   \alpha_{k+2m} \, g_{k+2m,n-1}
   & = &
   \alpha_{k+2m} \,
   \Bigl( \scrd_m \,+\, \sum_{i=1}^m \alpha_{k+2m+i} \Bigr)  g_{k+2m,n-2}
       \:+\:  \alpha_{k+2m} \, g_{k+m,n-1}
       \nonumber \\[-3mm]
 \label{eq.proof.prop.quasi-affine.3}
\end{subeqnarray}
Using $\alpha_{j+m+1} = \alpha_j + x_{j+1 \bmod m+1}$,
and using the fact that $\scrd_m$ is a pure first-order differential operator
and hence satisfies the Leibniz rule, we have
\begin{eqnarray}
   \hbox{LHS of \reff{eq.proof.prop.quasi-affine.2}}
   & = &
   \Bigl( \scrd_m \,+\, \sum_{i=1}^m \alpha_{k+m-1+i} \Bigr)
      (g_{k+m,n-1} - g_{k+m-1,n-1} - \alpha_{k+2m} g_{k+2m,n-2})
          \nonumber \\
   & & \quad
   +\; g_{k+2m,n-2} \, (\scrd_m \, \alpha_{k+2m})
          \nonumber \\[2mm]
   & &  \quad
   +\; (g_{k,n} - g_{k-1,n} - \alpha_{k+m} g_{k+m,n-1})
         \nonumber \\[2mm]
   & & \quad
   +\; (\alpha_{k+2m} - \alpha_{k+m}) \, (g_{k+m,n-1} - g_{k+m,n-1})
          \nonumber \\[1mm]
   & & \quad
   -\; \alpha_{k+2m} \, \biggl( \sum_{i=1}^m x_{k+m+i \bmod m+1} \biggr)
                     \, g_{k+2m,n-2}
   \;.
 \label{eq.proof.prop.quasi-affine.4}
\end{eqnarray}
The first term vanishes by the hypothesis of the inner induction on $n$;
the third term vanishes by the hypothesis of the outer induction on $k$;
and the fourth term vanishes.
On the other hand, from \reff{eq.alphak.periodm+1}
and the definition of $\scrd_m$ we have
\begin{subeqnarray}
   \scrd^m \, \alpha_{k+2m}
   & = &
   \lfloor (k+2m+1)/(m+1) \rfloor \; \scrd_m \, x_{k-1 \bmod m+1}
       \\
   & = &
   \lfloor (k+2m+1)/(m+1) \rfloor \; x_{k-1 \bmod m+1} \:
      \biggl( \sum_{i=0}^m x_i \,-\, x_{k-1 \bmod m+1} \biggr)
       \nonumber \\[-2mm] \\[-2mm]
   & = &
   \alpha_{k+2m} \, \biggl( \sum_{i=0}^m x_i \,-\, x_{k-1 \bmod m+1} \biggr)
   \;,
 \label{eq.proof.prop.quasi-affine.5}
\end{subeqnarray}
so the second term in \reff{eq.proof.prop.quasi-affine.4}
cancels the fifth term.
This completes the proof of \reff{eq.recurrence.gkm.0.bis.gkn}.

(c)  Using \reff{eq.recurrence.gkm.0.bis.gkn} with $k \to k-m$
(which is valid as noted above since $k-m \ge -m \ge -(m+1)$),
we get $g_{k-m,n} = g_{k-m-1,n} + \alpha_k g_{k,n-1}$.
Substituting this on the right-hand side of \reff{def.gkn.quasi-affine}
gives \reff{def.gkn.quasi-affine.variant}.

(d) For $0 \le k \le m$, the term $g_{k-m-1,n}$ on the right-hand side of
\reff{def.gkn.quasi-affine.variant} is simply $\delta_{n0}$,
and $\sum_{i=0}^m \alpha_{k+i} = \sum_{i=0}^k x_i$.
So (d) is an immediate consequence.
\qed

\medskip

{\bf Remark.}  When $\bfx = \bone$, we can use
\be
   \left. \scrd_m \prod_{i=0}^m x_i^{n_i} \right|_{\bfx = \bone}
   \;=\;
   m \sum_{i=0}^m n_i
\ee
together with the fact that $g_{k,n-1}$ is homogeneous in $\bfx$
of degree $n-1$, to specialize \reff{def.gkn.quasi-affine.variant} as
\be
   g_{k,n}(\bone)
   \;=\;
   \big[ m(n-1) + k + 1 \bigr] g_{k,n-1}(\bone) \:+\:  g_{k-m-1,n}(\bone)
   \;.
 \label{def.gkn.bone}
\ee
It follows that for $0 \le k \le m$,
the $g_{k,n}(\bone)$ are shifted multifactorials:
\be
   g_{k,n}(\bone)  \;=\;  \prod_{j=0}^{n-1} (k+1+jm)  \;=\;  C_n(k+1,m)
   \;.
\ee
However, for $k > m$ the forms get more complicated:
$g_{k,n}(\bone)$ is a linear combination of $\lfloor (k+1)/(m+1) \rfloor + 1$
shifted multifactorials.
These can be computed by symbolic algebra for any given $m$ and $k$:
for instance, for $m=2$ the first few are
\begin{subeqnarray}
   g_{0,n}(\bone)  & = & (2n-1)!!  \;=\;  C_n(1,2)   \\[1mm]
   g_{1,n}(\bone)  & = & (2n)!!    \;=\;  C_n(2,2)   \\[1mm]
   g_{2,n}(\bone)  & = & (2n+1)!!  \;=\;  C_n(3,2)  \\[1mm]
   g_{3,n}(\bone)  & = & -(2n+1)!! + (2n+2)!!  \;=\;  -C_n(3,2) + 2 C_n(4,2)  \\[1mm]
   g_{4,n}(\bone)  & = & -(2n+2)!! + (2n+3)!!  \;=\;  -2C_n(4,2) + 3C_n(5,2) \\[1mm]
   g_{5,n}(\bone)  & = & -(2n+3)!! + {(2n+4)!! \over 2} \;=\; -3C_n(5,2) + 4C_n(6,2)  \qquad
\end{subeqnarray}
But we have not been able to find a general
expression for the coefficients in these sums.
\myendremark

\bigskip

We have thus shown that $P^{(m,m+1)}(\bfx,\bfx) = (\scrd_m + x_0)^n \, 1$.
To complete the proof of Theorem~\ref{thm.quasi-affine}(a)
we need to show that $\scrp^{(m)}_n(\bfx) = (\scrd_m + x_0)^n \, 1$,
as asserted in Proposition~\ref{prop.MVeulerian.differential}.

\proofof{Proposition~\ref{prop.MVeulerian.differential}}
We will treat the cases of $\scrp^{(m)}_n$ and $\scrq^{(m)}_n$ in parallel,
and write $A_n$ as a shorthand for either one.

\begin{figure}[!ht]
\begin{center}
\includegraphics[scale=1]{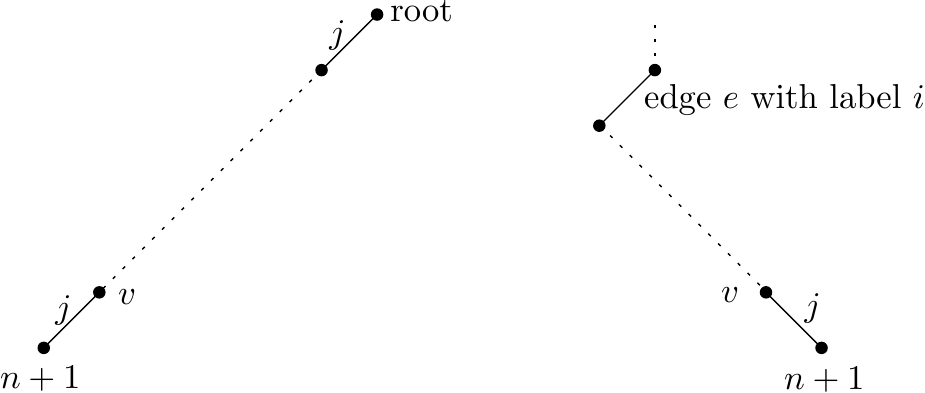} \\[3mm]
     {\em (i)} \hspace*{5cm} {\em (ii)}
\end{center}
\caption{The two possibilities for the vertex $n+1$.}
  \label{fig.positivetype}
\end{figure}

Consider vertex $n+1$ in an increasing ${(m+1)}$-ary tree on $n+1$ vertices:
it is necessarily a $j$-child for some $j \in \{0,\ldots,m\}$
of some vertex $v \in [n]$.
Follow the path from $v$ upwards towards the root
until the first $i$-edge with $i \ne j$ is reached (if one exists).
There are therefore two possibilities (see Figure~\ref{fig.positivetype} for an illustration):
\begin{itemize}
   \item[(i)]  No such edge exists: the path from $v$ to the root
        consists entirely of $j$-edges.
        (This includes the case in which $v$ is itself the root.)
        In the case of $\scrp^{(m)}_n$ this entails $j=0$.
   \item[(ii)]  Such an edge exists;  call it $e$.
\end{itemize}
And conversely, an increasing ${(m+1)}$-ary tree on $n+1$ vertices
is obtained from an increasing ${(m+1)}$-ary tree on $n$ vertices
by inserting vertex $n+1$ in one of these two ways:
either we choose a label $j$,
follow the path of $j$-edges (if any) downwards from the root
until there are no more $j$-edges,
and then attach vertex $n+1$ via a new $j$-edge;
or else we choose an edge $e$ (say that it is an $i$-edge),
choose a label $j \ne i$,
follow the path of $j$-edges (if any) downwards from the bottom end of $e$
until there are no more $j$-edges,
and then attach vertex $n+1$ via a new $j$-edge.

Possibility~(i) gives weight $x_0 A_{n-1}$ for $\scrp^{(m)}_n$,
or $\Bigl( \sum\limits_{j=0}^m x_j \Bigr) A_{n-1}$ for $\scrq^{(m)}_n$.
Possibility~(ii) gives weight $\scrd_m A_{n-1}$.
Adding these together gives
$A_n = (\scrd_m + x_0) A_{n-1}$ for $\scrp^{(m)}_n$, and
$A_n = {\Bigl( \scrd_m + \sum\limits_{j=0}^m x_j \Bigr) A_{n-1}}$
for $\scrq^{(m)}_n$.

$\scrq^{(m)}_{n,k}$ is handled similarly to $\scrq^{(m)}_n$,
with two small differences:
in possibility~(i) there are $k+1$ different roots
from which we could follow the path of $j$-edges downwards
to attach vertex $n+1$;
or $n+1$ could be made a new isolated vertex.
\qed

Finally, Theorem~\ref{thm.quasi-affine}(b)
follows from Theorem~\ref{thm.quasi-affine}(a) by odd contraction,
exactly as in Section~\ref{subsubsec.quasi-affine.bijective}.

\subsection[Period $m$: Increasing multi-$m$-ary trees]{Period $\bm{m}$: Increasing multi-$\bm{m}$-ary trees}
   \label{subsec.quasi-affine.m}

In Section~\ref{subsec.periodic.combinatorial}
we also defined the concept of a multi-$r$-ary tree.
We now define labeled multi-$r$-ary trees and increasing multi-$r$-ary trees
in the obvious way.

It is known \cite[p.~30, Corollary~1(iv)]{Bergeron_92}
that the number of increasing multi-unary trees on $n+1$ vertices
is $(2n-1)!!$;
more generally, by using \reff{eq.bergeron} with $\phi(w) = 1/(1-w)^r$,
it is easy to see that the number of increasing multi-$r$-ary trees
on $n+1$ vertices is the shifted multifactorial
$\prod\limits_{j=0}^{n-1} [r+j(r+1)]$.
We would now like to generalize this enumeration to a multivariate polynomial.

\subsubsection{Multivariate Eulerian polynomials and Eulerian symmetric functions of negative type: Definitions and statement of results}

Fix an integer $m \ge 1$,
and let $\bfx = (x_0,\ldots,x_{m-1})$ be indeterminates.
Let $\scrq^{(m)-}_n(\bfx)$ be the generating polynomial for
increasing multi-$m$-ary trees on $n+1$ vertices in which
each $i$-edge gets a weight $x_i$;
let $\scrq^{(m)-}_{n,k}(\bfx)$ be the generating polynomial for
unordered forests of increasing multi-$m$-ary trees
on $n+1$ total vertices with $k+1$ components
in~which each $i$-edge gets a weight $x_i$;
and let $\scrp^{(m)-}_n(\bfx)$ be the generating polynomial for
increasing multi-$m$-ary trees on $n+1$ vertices in which
all edges emanating from the root are 0-edges
and each $i$-edge gets a weight $x_i$.
Both $\scrp^{(m)-}_n(\bfx)$ and $\scrq^{(m)-}_n(\bfx)$
are homogeneous polynomials of degree~$n$,
while $\scrq^{(m)-}_{n,k}(\bfx)$ is homogeneous of degree~$n-k$;
we refer to them collectively as
\textbfit{multivariate Eulerian polynomials of negative type}.
When $\bfx = \bone$ we have $\scrp^{(m)-}_n(\bone) = F_n^{(m+1)}$
[see below after \reff{eq.Pinfty-.eulerian.explicit.e}]
and $\scrq^{(m)-}_n(\bone) = \prod\limits_{j=0}^{n-1} [m+j(m+1)]$.



We will prove:

\begin{theorem}[Branched continued fractions for the multivariate Eulerian polynomials of negative type]
 \label{thm.quasi-affine.multi}
For every integer $m \ge 1$, we have:
\begin{itemize}
   \item[(a)] $\scrp^{(m)-}_n(\bfx) = P_n^{(m,m)}(\bfx,\bfx)$,
which by definition equals
$S_n^{(m)}(\balpha)$ where the weights $\balpha$ are given by
the Eulerian-quasi-affine formula \reff{eq.alpha.quasi-affine.u=x}
with period $p=m$.
   \item[(b)] $\scrq^{(m)-}_n(\bfx) =  J_n^{(m)}(\bbeta')$
where the weights $\bbeta'$ are given by
\be
   \beta_i^{\prime(\ell)}
   \;=\;
   {(i+1)! \over (i-\ell)!} \: h_{\ell+1}(x_0,\ldots,x_m)
   \;.
 \label{eq.thm.quasi-affine.multi}
\ee
   \item[(c)] $\scrq^{(m)-}_{n,k}(\bfx) =  J_{n,k}^{(m)}(\bbeta')$
where the weights $\bbeta'$ are given by \reff{eq.thm.quasi-affine.multi}.
   \item[(d)]  $\scrp^{(m)-}_n(\bfx) =
    \displaystyle \sum\limits_{j=1}^n j! \, x_0^j \, J_{n-1,j-1}^{(m)}(\bbeta')$
    for $n \ge 1$, where $\bbeta'$ is given by \reff{eq.thm.quasi-affine.multi}.
\end{itemize} 
\end{theorem}

We will prove Theorem~\ref{thm.quasi-affine.multi}(b,c,d)
by a small modification of our bijective proof
of Theorems~\ref{thm.quasi-affine} and \ref{thm.factorized}
given in Section~\ref{subsubsec.quasi-affine.bijective}.
However, we are unable at present to use this method to prove part~(a):
the trouble is that here the relation
between $\scrp^{(m)-}_n(\bfx)$ and $\scrq^{(m)-}_n(\bfx)$
is nontrivial [see \reff{eq.Pinfty-.eulerian.explicit.Q} below],
in contrast to the relation \reff{def.scrp.b}
between $\scrp^{(m)}_n(\bfx)$ and $\scrq^{(m)}_n(\bfx)$,
which corresponds precisely to odd contraction.
So we will instead prove Theorem~\ref{thm.quasi-affine.multi}(a)
by a small modification of the Euler--Gauss recurrence proof
given in Section~\ref{subsubsec.quasi-affine.euler-gauss}.

For the special case $\bfx = \bone$, we have:

\begin{corollary}[Eulerian-quasi-affine $m$-S-fraction of period $m$ for multifactorials]
  \label{cor.quasi-affine.multi}
\hfill\break
For every integer $m \ge 1$, we have
$P_n^{(m,m)}(\bone,\bone) = F_n^{(m+1)}$.
\end{corollary}

Theorem~\ref{thm.quasi-affine.multi}(a) also implies
that the polynomials $P_n^{(m,m)}(\bfx,\bfx)$
are symmetric in $x_1,\ldots,x_{m-1}$ (though not also in $x_0$) ---
a fact that seems far from obvious from the branched continued fraction
or from the $m$-Dyck-path definition.

{}From Theorems~\ref{thm.quasi-affine.multi}(a)
and \ref{thm.Stype.minors} we deduce:

\begin{corollary}[Hankel-total positivity of $\scrp^{(m)-}$]
  \label{cor.quasi-affine.multi.2}
For every integer $m \ge 1$,
the sequence $(\scrp^{(m)-}_n)_{n \ge 0}$
is coefficientwise Hankel-totally positive,
jointly in the indeterminates $x_0,\ldots,x_{m-1}$.
\end{corollary}

We will also show (as a special case of the more general
Lemma~\ref{lemma.Pm.beta.totalpos.factorized.multi}):

\begin{lemma}
   \label{lemma.Pm.beta.totalpos.quasi-affine.multi}
The production matrix $P(\bbeta')$ [cf.\ \reff{def.Pm}]
associated to the weights \reff{eq.thm.quasi-affine.multi}
is coefficientwise totally positive,
jointly in the indeterminates $x_0,\ldots,x_{m-1}$.
\end{lemma}

Please note that in the present case the production matrix $P(\bbeta)$
is a full lower-Hessenberg matrix --- not $(m,1)$-banded
as in the case of positive type ---
since the complete homogeneous symmetric functions $h_n(x_0,\ldots,x_{m-1})$
are nonvanishing for all~$n$.

Theorems~\ref{thm.Jtype.minors} and \ref{thm.Jtype.generalized.minors}
then imply:

\begin{corollary}[Hankel-total positivity of $\scrq^{(m)-}$]
  \label{cor.quasi-affine.multi.2a}
For every integer $m \ge 1$,
the sequence $(\scrq^{(m)-}_n)_{n \ge 0}$
is coefficientwise Hankel-totally positive,
jointly in the indeterminates $x_0,\ldots,x_{m-1}$.
\end{corollary}

\begin{corollary}[Total positivity of lower-triangular matrix $\scrq^{(m)-}_{n,k}$]
  \label{cor.quasi-affine.multi.2b}
For every integer ${m \ge 1}$,
the lower-triangular matrix $(\scrq^{(m)-}_{n,k}(\bfx))_{n,k \ge 0}$
is coefficientwise totally positive,
jointly in the indeterminates $x_0,\ldots,x_{m-1}$.
\end{corollary}

We also have explicit expressions for the multivariate Eulerian polynomials
of negative type $\scrp^{(m)-}_n$ and $\scrq^{(m)-}_n$,
and a recurrence for $\scrq^{(m)-}_{n,k}$,
given by repeated application of a first-order linear differential operator.
We will prove:

\begin{proposition}[Differential expressions for the multivariate Eulerian polynomials of negative type]
   \label{prop.MVeulerian.differential.multi}
For every integer ${m \ge 1}$, we have
\begin{eqnarray}
   \scrp^{(m)-}_n(\bfx)  & = &  (\scrd^-_m \,+\, x_0)^n \: 1
          \label{eq.scrp.differential.multi}  \\[2mm]
   \scrq^{(m)-}_n(\bfx)  & = &
      \Bigl( \scrd^-_m \,+\, \sum_{i=0}^{m-1} x_i \Bigr)^{\! n} \: 1
          \label{eq.scrq.differential.multi}  \\
   \scrq^{(m)-}_{n,k}(\bfx)  & = &
      \Bigl( \scrd^-_m \,+\, (k+1) \sum_{i=0}^{m-1} x_i \Bigr)
      \scrq^{(m)-}_{n-1,k}(\bfx)  \:+\:  \scrq^{(m)-}_{n-1,k-1}(\bfx)
      \:+\: \delta_{n0} \delta_{k0}
      \qquad
          \label{eq.scrqnk.differential.multi}
\end{eqnarray}
where
\be
   \scrd^-_m
   \;=\;
   \sum_{i=0}^{m-1} \biggl( \! x_i^2  \:+\:
                            x_i \sum_{j=0}^{m-1} x_j \!
                    \biggr) \, {\partial \over \partial x_i}
   \;.
 \label{def.dm-}
\ee
\end{proposition}

{\bf Remarks.}
1.  $\scrd_m$ [cf.\ \reff{def.dm}] and $\scrd^-_m$ can be given a unified
definition as
\be
   \scrd^\pm_m
   \;=\;
   \sum_{i=0}^{p-1} \biggl( \sum_{j=0}^{p-1} x_j \,\mp\, x_i \!
                    \biggr) x_i \, {\partial \over \partial x_i}
 \label{def.dm+-}
\ee
where $p=m+1$ or $m$, respectively.

2.  For $m=1$ we have trivially $P_n^{(1)-}(x) = (2n-1)!! \, x^n$;
the S-fraction for this case is due to Euler \cite[section~29]{Euler_1760}.
For $m=2$, the triangular array of coefficients corresponding to
$P_n^{(2)-}(x,y)$ is not currently in \cite{OEIS}.
\myendremark

\bigskip

We can again pass to the limit $m \to \infty$ in the obvious way,
defining generating formal power series $\scrp^{(\infty)-}_n(\bfx)$,
$\scrq^{(\infty)-}_n(\bfx)$ and $\scrq^{(\infty)-}_{n,k}(\bfx)$
for increasing multi-$\infty$-ary trees.
Here $\scrq^{(\infty)-}_n(\bfx)$ and $\scrq^{(\infty)-}_{n,k}(\bfx)$
are symmetric functions in the indeterminates $\bfx$,
while $\scrp^{(\infty)-}_n(\bfx)$ is a symmetric function
of $(x_i)_{i \ge 1}$ with coefficients in $\Z[x_0]$;
they are homogeneous of degree~$n$ or $n-k$.
We refer to $\scrq^{(\infty)-}_n(\bfx)$ and $\scrq^{(\infty)-}_{n,k}(\bfx)$
as the \textbfit{Eulerian symmetric functions of negative type}.
Theorem~\ref{thm.quasi-affine.multi}(b,c,d),
Corollaries~\ref{cor.quasi-affine.multi.2}--\ref{cor.quasi-affine.multi.2b}
and Proposition~\ref{prop.MVeulerian.differential.multi}
have obvious counterparts at $m=\infty$, which we refrain from stating.


We can say more about the Eulerian symmetric functions
of negative type
by using the interpretation of multi-$\infty$-ary trees
as ordered trees in which each edge carries a label $i \in \N$
and the edges emanating outwards from each vertex consist, in order,
of zero or more edges labeled 0, then zero or more edges labeled 1, etc.\
(Section~\ref{subsec.periodic.combinatorial}).
Since the choice of labels on the edges emanating outwards from a vertex $v$
can be made independently for each $v$,
the definition of $\scrq^{(\infty)-}_n$ can be rephrased as:
$\scrq^{(\infty)-}_n(\bfx)$ is the generating formal power series for
increasing ordered trees on $n+1$ vertices in which
each vertex with $i$ children gets a weight $h_i(\bfx)$,
where $h_i$ is the complete homogeneous symmetric function;
and an analogous rephrasing holds for $\scrq^{(\infty)-}_{n,k}$.
It follows that the Eulerian symmetric functions of positive and negative type
are simply the duals of each other:
$\scrq^{(\infty)-}_n = \omega \scrq^{(\infty)}_n$
and $\scrq^{(\infty)-}_{n,k} = \omega \scrq^{(\infty)}_{n,k}$,
where $\omega$ is the involution of the ring of symmetric functions
defined by $\omega h_i = e_i$.
Therefore, the formulae \reff{eq.prop.eulerian.symfun},
\reff{eq.euleriansymfun.0-6} and \reff{eq.prop.eulerian.symfun.2}
for $\scrq^{(\infty)}_n$ and $\scrq^{(\infty)}_{n,k}$
imply dual formulae for $\scrq^{(\infty)-}_n$ and $\scrq^{(\infty)-}_{n,k}$,
in which $e$ and $h$ are interchanged.

We can also be explicit about $\scrp^{(\infty)-}_n$.
Starting from the definition of $\scrp^{(\infty)-}_n$
and classifying the trees by the sizes of the subtrees of the root,
we see immediately that
\be
   \scrp^{(\infty)-}_n(\bfx)
   \;=\;
   \sum_{\ell=0}^\infty x_0^\ell
          \!\!\!
      \sum_{\begin{scarray}
               k_1,\ldots,k_\ell \ge 1 \\
               \sum k_i = n
            \end{scarray}}
          \!\!\!\!
      \binom{n}{k_1,\,\ldots,\, k_\ell}
      \prod_{i=1}^\ell \scrq^{(\infty)-}_{k_i - 1}(\bfx)
   \;.
 \label{eq.Pinfty-.eulerian.explicit.Q}
\ee
This can be rewritten in terms of exponential generating functions as
\be
   {\bf P}(t)
   \;\eqdef\;
   \sum_{n=0}^\infty \scrp^{(\infty)-}_n(\bfx) \, {t^n \over n!}
   \;=\;
   {1 \over 1 \,-\, x_0 Y(t)}
 \label{eq.Pinfty-.eulerian.explicit.egf}
\ee
where $Y(t) \eqdef
       \sum\limits_{n=1}^\infty \scrq^{(\infty)-}_{n-1}(\bfx) \, t^n/n!$
satisfies \reff{eq.bergeron.fY} with $\psi_i = (-1)^i e_i$.
Solving \reff{eq.bergeron.fY}/\linebreak\reff{eq.Pinfty-.eulerian.explicit.egf}
by Lagrange inversion gives for $n \ge 1$
\begin{subeqnarray}
   & &
   \hspace*{-6mm}
   \scrp^{(\infty)-}_n(\bfx)
   \;=\;  n! \, [t^n] \, {\bf P}(t)
   \;=\; (n-1)! \, [y^{n-1}] \, {x_0 \over (1- x_0 y)^2} \, f(y)^{-n}
        \\[2mm]
   & &
   \hspace*{-3mm}
   \;=\;
   (n-1)! \sum_{j=1}^\infty j \, x_0^j \, [y^{n-j}] \, f(y)^{-n}
        \\[3mm]
   & &
   \hspace*{-3mm}
   \;=\;
   (n-1)! \sum_{j=1}^\infty j \, x_0^j
          \!\!\!\! \sum_{\begin{scarray}
                   k_1, k_2, \ldots \ge 0 \\
                   \sum i k_i = n-j
                \end{scarray}}
          \!\!\!\!
          \binom{-n}{-n-\sum k_i,\, k_1,\, k_2,\, \ldots}
          \prod_{i=1}^\infty \Bigl( {(-1)^i e_i \over i+1} \Bigr) ^{\! k_i}
      \qquad\qquad
        \\[3mm]
   & &
   \hspace*{-3mm}
   \;=\;
   (n-1)! \sum_{j=1}^\infty j \, x_0^j
          \!\!\!\! \sum_{\begin{scarray}
                   k_1, k_2, \ldots \ge 0 \\
                   \sum i k_i = n-j
                \end{scarray}}
          \!\!\!\!
          (-1)^{\sum k_i}
          \binom{n+\sum k_i - 1}{n-1,\, k_1,\, k_2,\, \ldots}
          \prod_{i=1}^\infty \Bigl( {(-1)^i e_i \over i+1} \Bigr) ^{\! k_i}
      \qquad\qquad
       \\[3mm]
   & &
   \hspace*{-3mm}
   \;=\;
   n! \, x_0^n
   \:+\:
   \sum_{j=1}^{n-1}  j \, x_0^j
   \sum_{\lambda \vdash n-j}
   \,
   (-1)^{n-j-\ell(\lambda)} \:
   {(n+\ell(\lambda)-1)!  \over
    \displaystyle \prod\limits_i (i+1)^{m_i(\lambda)} \, m_i(\lambda)!
   }
   \; e_\lambda(\bfx)
      \;.
 \slabel{eq.Pinfty-.eulerian.explicit.e.e}
 \label{eq.Pinfty-.eulerian.explicit.e}
\end{subeqnarray}
%
This equation agrees with Theorem~\ref{thm.quasi-affine.multi}(c,d)
combined with the dual of \reff{eq.prop.eulerian.symfun.2.b}.
In the special case where $\bfx$ consists of $m$ ones
(one of which is $x_0$) and the rest zeroes,
we have $\phi(w) = 1/(1-w)^m$, hence $Y(t) = 1 - [1 - (m+1)t]^{1/(m+1)}$
and ${\bf P}(t) = [1 - (m+1)t]^{-1/(m+1)}$,
so that $\scrp^{(m)-}_n(\bone) = F_n^{(m+1)}$.

\subsubsection{Extended multivariate Eulerian polynomials and Eulerian symmetric functions of negative type: Definitions and statement of results}

We can also introduce additional indeterminates $\bfc$
as in Section~\ref{subsec.qusi-affine.extended}.
Let $\bfx = (x_0,\ldots,x_{m-1})$ and $\bfc = (c_L)_{L \ge 0}$
be indeterminates;
let $\scrq^{(m)-}_n(\bfx,\bfc)$ be the generating polynomial for
increasing multi-$m$-ary trees on $n+1$ vertices in which
each $i$-edge gets a weight $x_i$
and each vertex at level $L$ gets a weight $c_L$,
divided by $c_0$;
let $\scrq^{(m)-}_{n,k}(\bfx,\bfc)$ be the generating polynomial for
unordered forests of increasing multi-$m$-ary trees
on $n+1$ total vertices with $k+1$ components
in~which each $i$-edge gets a weight $x_i$
and each vertex at level $L$ gets a weight $c_L$,
divided by $c_0 c_1 \cdots c_k$;
and let $\scrp^{(m)-}_n(\bfx,\bfc)$ be defined as $\scrq^{(m)-}_n(\bfx,\bfc)$
but restricted to trees in which
all edges emanating from the root are 0-edges.
Both $\scrp^{(m)-}_n(\bfx,\bfc)$ and $\scrq^{(m)-}_n(\bfx,\bfc)$
are homogeneous of degree~$n$ in $\bfx$ and in $\bfc$,
while $\scrq^{(m)-}_{n,k}(\bfx)$ is homogeneous of degree~$n-k$ in $\bfx$
and in $\bfc$;
we refer to them collectively as
\textbfit{extended multivariate Eulerian polynomials of negative type}.

By analogy with Theorems~\ref{thm.quasi-affine}, \ref{thm.factorized}
and \ref{thm.quasi-affine.multi}, we suspect that the following is true
(though we are unable at present to prove it):

\begin{conjecture}[Branched continued fractions for the extended multivariate Eulerian polynomials of negative type, part~(a)]
 \label{conj.factorized.multi.a}
For every integer $m \ge 1$, we have:
\begin{itemize}
   \item[(a)] $\scrp^{(m)-}_n(\bfx,\bfc) = \widehat{P}_n^{(m,m)}(\bfx,\bfc)$,
which by definition equals
$S_n^{(m)}(\balpha)$ where the weights $\balpha$ are given by
the factorized formula \reff{eq.alpha.factorized} with period $p=m$.
\end{itemize}
\end{conjecture}

We are, however, able to prove the remaining parts of the analogy:

\begin{theorem}[Branched continued fractions for the extended multivariate Eulerian polynomials of negative type, parts~(b) and (c)]
  \label{thm.factorized.multi.bc}
For every integer $m \ge 1$, we have:
\begin{itemize}
   \item[(b)] $\scrq^{(m)-}_n(\bfx,\bfc) =  J_n^{(m)}(\bbeta')$
where the weights $\bbeta'$ are given by
\be
   \beta_i^{\prime(\ell)}
   \;=\;
   {(i+1)! \over (i-\ell)!} \:
   c_{i-\ell} c_{i-\ell+1} \cdots c_i
      \, h_{\ell+1}(x_0,\ldots,x_m)
   \;.
 \label{eq.thm.factorized.multi}
\ee
   \item[(c)] $\scrq^{(m)-}_{n,k}(\bfx,\bfc) =  J_{n,k}^{(m)}(\bbeta')$
where the weights $\bbeta'$ are given by \reff{eq.thm.factorized.multi}.
   \item[(d)]  $\scrp^{(m)-}_n(\bfx) =
    \displaystyle \sum\limits_{j=1}^n j! \, c_0 c_1 \cdots c_{j-1} \, x_0^j \,
                                        J_{n-1,j-1}^{(m)}(\bbeta')$
    for $n \ge 1$, where $\bbeta'$ is given by \reff{eq.thm.factorized.multi}.
\end{itemize}
\end{theorem}

If Conjecture~\ref{conj.factorized.multi.a} holds, it implies:

\begin{conjecture}[Hankel-total positivity of $\scrp^{(m)-}_n(\bfx,\bfc)$]
  \label{conj.factorized.multi.2}
For every integer ${m \ge 1}$,
the sequence $(\scrp^{(m)-}_n(\bfx,\bfc))_{n \ge 0}$
is coefficientwise Hankel-totally positive,
jointly in the indeterminates $\bfx$ and $\bfc$.
\end{conjecture}

We will also show:

\begin{lemma}
   \label{lemma.Pm.beta.totalpos.factorized.multi}
The production matrix $P(\bbeta')$ [cf.\ \reff{def.Pm}]
associated to the weights \reff{eq.thm.factorized.multi}
is coefficientwise totally positive,
jointly in the indeterminates $\bfx$ and $\bfc$.
\end{lemma}

Using Theorems~\ref{thm.Jtype.minors} and \ref{thm.Jtype.generalized.minors}
we then conclude:

\begin{corollary}[Hankel-total positivity of $\scrq^{(m)-}(\bfx,\bfc)$]
  \label{cor.factorized.multi.2a}
For every integer $m \ge 1$,
the sequence $(\scrq^{(m)-}_n(\bfx,\bfc))_{n \ge 0}$
is coefficientwise Hankel-totally positive,
jointly in the indeterminates $\bfx$ and $\bfc$.
\end{corollary}

\begin{corollary}[Total positivity of lower-triangular matrix of $\scrq^{(m)-}(\bfx,\bfc)$]
  \label{cor.factorized.multi.2b}
For every integer ${m \ge 1}$,
the lower-triangular matrix $(\scrq^{(m)-}_{n,k}(\bfx,\bfc))_{n,k \ge 0}$
is coefficientwise totally positive,
jointly in the indeterminates $\bfx$ and $\bfc$.
\end{corollary}


Of course, by specializing
Conjecture~\ref{conj.factorized.multi.a},
Theorem~\ref{thm.factorized.multi.bc},
Conjecture~\ref{conj.factorized.multi.2},
Lemma~\ref{lemma.Pm.beta.totalpos.factorized.multi},
and Corollaries~\ref{cor.factorized.multi.2a} and
\ref{cor.factorized.multi.2b}
to $\bfc = \bone$,
we recover Theorem~\ref{thm.quasi-affine.multi},
Corollary~\ref{cor.quasi-affine.multi.2},
Lemma~\ref{lemma.Pm.beta.totalpos.quasi-affine.multi},
and Corollaries~\ref{cor.quasi-affine.multi.2a} and
\ref{cor.quasi-affine.multi.2b}.

These definitions, conjectures and results extend in an obvious way
to $m=\infty$.  In particular,
the definition of $\scrq^{(\infty)-}_n(\bfx,\bfc)$ can be rephrased as:
$\scrq^{(\infty)-}_n(\bfx,\bfc)$ is the generating formal power series for
increasing ordered trees on $n+1$ vertices in which
each vertex with $i$ children gets a weight $h_i(\bfx)$,
and each vertex at level $L$ gets a weight $c_L$, divided by $c_0$.
Similarly,
$\scrp^{(\infty)-}_n(\bfx,\bfc)$ is the generating formal power series for
increasing ordered trees on $n+1$ vertices in which
each vertex with $i$ children gets a weight $h_i(\bfx)$,
except for the root, which gets a weight $x_0^i$,
and each vertex at level $L$ gets a weight $c_L$, divided by $c_0$.

\subsubsection{Proof by bijection to labeled lattice paths}

The proof of Theorem~\ref{thm.factorized.multi.bc}
is almost identical to that of Theorem~\ref{thm.factorized},
with only some slight changes to replace
$m$-\L{}ukasiewicz paths by $\infty$-\L{}ukasiewicz paths,
sets by multisets,
and elementary symmetric functions
by complete homogeneous symmetric functions.
We therefore supply only a sketch.

\proofof{Theorem~\ref{thm.factorized.multi.bc}}
We will construct a bijection from
the set of ordered forests of increasing multi-$m$-ary trees
on the vertex set $[n+1]$ with $k+1$ components
to the set of $\bfL$-labeled reversed partial $\infty$-\L{}ukasiewicz paths
from $(0,k)$ to $(n,0)$,
where the label sets $\bfL$ will be defined below.

The definition of the steps $s_j$ is the same as in the proof of
Theorem~\ref{thm.factorized};
but since the number of children of a vertex is unbounded,
we have only $-1 \le s_j < \infty$.
The path $\omega$ is therefore a reversed \L{}ukasiewicz path
of length $n$.

The label $\xi_j$ will be an ordered pair $\xi_j = (\xi'_j,\xi''_j)$
where $\xi'_j$ is a positive integer
and $\xi''_j$ is a multiset on $\{0,\ldots,m-1\}$;
more precisely, the label set $L(s,h)$
for a step $s$ starting at height $h$ will be
\be
   L(s,h)  \;=\;  [h+1] \,\times\, \doublebinom{\{0,\ldots,m-1\}}{s+1}
\ee
where $\doublebinom{S}{r}$ denotes the set of $r$-element multisets on $S$
(that is, the set of multi-indices $(n_i)_{i \in S}$ of nonnegative integers
 such that $\sum n_i = r$).
As before, the label $\xi'_j$ will say into which of the
$h_{j-1}+1$ ``vacant slots''
(numbered in order of their creation) the node~$j$ is to be inserted\footnote{
   We trust that there will be no confusion between the heights $h_j$
   and the elementary symmetric functions $h_n(\bfx)$.
};
and the label $\xi''_j$ will say which children node~$j$ has.

%

The computation of the weights is the same as in
the proof of Theorem~\ref{thm.factorized};
but since we are here dealing with multi-$m$-ary trees
rather than $(m+1)$-ary trees,
the summation over $\xi''_j$ gives a
complete homogeneous symmetric function $h_{s+1}(x_0,\ldots,x_{m-1})$
instead of an elementary symmetric function $e_{s+1}(x_0,\ldots,x_m)$.
This gives the weights \reff{eq.thm.factorized.multi}
and proves Theorem~\ref{thm.factorized.multi.bc}(b,c).

To prove Theorem~\ref{thm.factorized.multi.bc}(d),
we observe that, for $n \ge 1$,
$\scrp_n^{(m)-}(\bfx,\bfc)$ is the same as $\scrq_n^{(m)-}(\bfx,\bfc)$
except that the edges emanating from the root are constrained to be 0-edges.
If the root has $j$ children, this means that the first step of the
reversed \L{}ukasiewicz path is $s_0 = j-1$,
and it gets a weight $c_0 x_0^j$ rather than the usual $c_0 h_j(\bfx)$.
Therefore, the {\em last}\/ step of the \L{}ukasiewicz path
goes from height $\ell = j-1$ to height 0.
We remove this last step, yielding a partial \L{}ukasiewicz path
from $(0,0)$ to $(n-1,j-1)$,
and then include explicitly the weight for this step:
instead of
$\beta_{j-1}^{\prime(j-1)} = \widehat{c}_0 \cdots \widehat{c}_{j-1} h_j(\bfx)$
where $\widehat{c} \eqdef (k+1) c_k$,
it is $\widehat{c}_0 \cdots \widehat{c}_{j-1} x_0^j$.
This proves Theorem~\ref{thm.factorized.multi.bc}(d).
\qed

\proofof{Lemma~\ref{lemma.Pm.beta.totalpos.factorized.multi}}
The proof is identical to that of
Lemma~\ref{lemma.Pm.beta.totalpos.factorized},
with the Jacobi--Trudi formula for $h_n$ \cite[Theorem~7.16.1]{Stanley_99}
replacing that for $e_n$.
\qed

The remark after the proof of Lemma~\ref{lemma.Pm.beta.totalpos.factorized}
concerning total positivity in the Schur order also holds here.

\subsubsection{Proof by the Euler--Gauss recurrence method}

We will now prove Theorem~\ref{thm.quasi-affine.multi}(a)
by the Euler--Gauss recurrence method;
we~use the same notation and preliminary definitions as in
Section~\ref{subsubsec.quasi-affine.euler-gauss}.
We need to prove
\be
   g_{k,n} \,-\, g_{k-1,n}  \;=\;  \alpha_{k+m} \, g_{k+m,n-1}
   \qquad\hbox{for } k,n \ge 0
   \;,
 \label{eq.recurrence.gkm.0.bis.gkn.multi}
\ee
where the weights $\balpha$ are the Eulerian-quasi-affine weights
\reff{eq.alpha.quasi-affine.u=x} with period $p=m$, i.e.
\be
   \alpha_k  \;=\;  \lfloor k/m \rfloor \, x_{k \bmod m}
   \;.
 \label{eq.alphak.periodm}
\ee

\begin{proposition}[Euler--Gauss recurrence for multivariate Eulerian polynomials of negative type]
   \label{prop.euler-gauss.quasi-affine.multi}
Let $\bfx = (x_0,\ldots,x_{m-1})$ be indeterminates;
we work in the ring $R = \Z[\bfx]$.
Set $g_{k,n} = \delta_{n0}$ for $k < 0$,
and then define $g_{k,n}$ for $k,n \ge 0$ by the recurrence
\be
   g_{k,n}
   \;=\;
   \Bigl( \scrd^-_m \,+\, \sum_{i=1}^m \alpha_{k+i} \Bigr) g_{k,n-1}
       \:+\:  g_{k-m,n}
 \label{def.gkn.quasi-affine.multi}
\ee
where $\balpha$ are given by \reff{eq.alphak.periodm}
and $\scrd^-_m$ is given by \reff{def.dm-}.
Then:
\begin{itemize}
   \item[(a)]  $g_{k,0} = 1$ for all $k \in \Z$.
\\[-5mm]
   \item[(b)]  $(g_{k,n})$ satisfies the recurrence
      \reff{eq.recurrence.gkm.0.bis.gkn.multi}.
\\[-5mm]
   \item[(c)]  $(g_{k,n})$ also satisfies the recurrence
\be
   g_{k,n}
   \;=\;
   \Bigl( \scrd^-_m \,+\, \sum_{i=0}^m \alpha_{k+i} \Bigr) g_{k,n-1}
       \:+\:  g_{k-m-1,n}
   \qquad\hbox{for } k,n \ge 0
   \;.
 \label{def.gkn.quasi-affine.multi.variant}
\ee
\\[-13mm]
   \item[(d)]  For $0 \le k \le m-1$, we have
$\displaystyle
   g_{k,n}
   \:=\:
   \Bigl( \scrd^-_m \,+\, \sum_{i=0}^k x_i \Bigr)^{\! n} \: 1
$,
while for $k=m$ we~have
$\displaystyle
   g_{m,n}
   \:=\:
   \Bigl( \scrd^-_m \,+\, \sum_{i=0}^{m-1} x_i \,+\, 2x_0 \Bigr)^{\! n} \: 1
$.
\end{itemize}
Therefore $S^{(m)}_n(\balpha) = g_{0,n} = (\scrd^-_m + x_0)^n \, 1$,
and more generally $S^{(m)}_{n|k}(\balpha) = g_{k,n}$.
\end{proposition}

Note that we are using \reff{eq.alphak.periodm}
also to define $\alpha_i = 0$ for $0 \le i \le m-1$.
This implies that \reff{eq.recurrence.gkm.0.bis.gkn.multi}
holds also for $-m \le k \le -1$.

\medskip

\proofof{Proposition~\ref{prop.euler-gauss.quasi-affine.multi}}
(a) holds trivially as in Proposition~\ref{prop.euler-gauss.quasi-affine}.

(b) We will now prove that \reff{eq.recurrence.gkm.0.bis.gkn.multi} holds.
The proof will be by an outer induction on $k$ and an inner induction on $n$.
The base cases $k = -m,\ldots,-1$ hold trivially.
Now suppose that \reff{eq.recurrence.gkm.0.bis.gkn.multi} holds
for a given $k$ and all $n \ge 0$;
we want to prove it when $k$ is replaced by $k+m$, i.e.\ that
\be
   g_{k+m,n} \:-\: g_{k+m-1,n} \:-\: \alpha_{k+2m} \, g_{k+2m,n-1}
   \;=\;
   0
   \quad\hbox{for all $n \ge 0$}
   \;.
 \label{eq.proof.prop.quasi-affine.multi.2}
\ee
We will prove \reff{eq.proof.prop.quasi-affine.multi.2} by induction on $n$.
It holds for $n=0$ because $g_{k+m,0} = g_{k+m-1,0} = 1$ and $g_{k+2m,-1} = 0$.
If $n > 0$, then we use \reff{def.gkn.quasi-affine.multi}
on each of the three terms on the left-hand side of
\reff{eq.proof.prop.quasi-affine.multi.2}, giving
\begin{subeqnarray}
   g_{k+m,n}
   & = &
   \Bigl( \scrd^-_m \,+\, \sum_{i=1}^m \alpha_{k+m+i} \Bigr)  g_{k+m,n-1}
       \:+\:  g_{k,n}
       \\
   g_{k+m-1,n}
   & = &
   \Bigl( \scrd^-_m \,+\, \sum_{i=1}^m \alpha_{k+m-1+i} \Bigr)  g_{k+m-1,n-1}
       \:+\:  g_{k-1,n}
       \\
   \alpha_{k+2m} \, g_{k+2m,n-1}
   & = &
   \alpha_{k+2m} \,
   \Bigl( \scrd^-_m \,+\, \sum_{i=1}^m \alpha_{k+2m+i} \Bigr)  g_{k+2m,n-2}
       \:+\:  \alpha_{k+2m} \, g_{k+m,n-1}
       \nonumber \\[-3mm]
 \label{eq.proof.prop.quasi-affine.multi.3}
\end{subeqnarray}
Using $\alpha_{j+m} = \alpha_j + x_{j \bmod m}$,
and using the fact that $\scrd^-_m$ is a pure first-order differential operator
and hence satisfies the Leibniz rule, we have
\begin{eqnarray}
   \hbox{LHS of \reff{eq.proof.prop.quasi-affine.multi.2}}
   & = &
   \Bigl( \scrd^-_m \,+\, \sum_{i=1}^m \alpha_{k+m-1+i} \Bigr)
      (g_{k+m,n-1} - g_{k+m-1,n-1} - \alpha_{k+2m} g_{k+2m,n-2})
          \nonumber \\
   & & \quad
   +\; g_{k+2m,n-2} \, (\scrd^-_m \, \alpha_{k+2m})
          \nonumber \\[2mm]
   & &  \quad
   +\; (g_{k,n} - g_{k-1,n} - \alpha_{k+m} g_{k+m,n-1})
         \nonumber \\[2mm]
   & & \quad
   +\; x_{k \bmod m} \, (g_{k+m,n-1} - g_{k+m,n-1})
          \nonumber \\[1mm]
   & & \quad
   -\; \alpha_{k+2m} \, \Bigl( x_{k \bmod m} \,+\, \sum_{i=0}^{m-1} x_i \Bigr)
                     \, g_{k+2m,n-2}
   \;.
 \label{eq.proof.prop.quasi-affine.multi.4}
\end{eqnarray}
The first term vanishes by the hypothesis of the inner induction on $n$;
the third term vanishes by the hypothesis of the outer induction on $k$;
and the fourth term vanishes.
On the other hand, from \reff{eq.alphak.periodm}
and the definition of $\scrd^-_m$ we have
\begin{subeqnarray}
   \scrd^-_m \, \alpha_{k+2m}
   & = &
   (\lfloor k/m \rfloor + 2) \: \scrd^-_m \, x_{k \bmod m}
       \\
   & = &
   (\lfloor k/m \rfloor + 2) \: x_{k \bmod m} \,
      \Bigl( x_{k \bmod m} \,+\, \sum_{i=0}^{m-1} x_i \Bigr)
       \\
   & = &
   \alpha_{k+2m} \, \Bigl( x_{k \bmod m} \,+\, \sum_{i=0}^{m-1} x_i \Bigr)
   \;,
 \label{eq.proof.prop.quasi-affine.multi.5}
\end{subeqnarray}
so the second term in \reff{eq.proof.prop.quasi-affine.multi.4}
cancels the fifth term.

(c)  Using \reff{eq.recurrence.gkm.0.bis.gkn.multi} with $k \to k-m$
(which is valid as noted above since $k-m \ge -m$),
we get $g_{k-m,n} = g_{k-m-1,n} + \alpha_k g_{k,n-1}$.
Substituting this on the right-hand side of \reff{def.gkn.quasi-affine.multi}
gives \reff{def.gkn.quasi-affine.multi.variant}.

(d) For $0 \le k \le m$, the term $g_{k-m-1,n}$ on the right-hand side of
\reff{def.gkn.quasi-affine.multi.variant} is simply $\delta_{n0}$.
Then $\sum_{i=0}^m \alpha_{k+i} = \sum_{i=0}^k x_i$ for $0 \le k \le m-1$,
or $\sum_{i=0}^{m-1} x_i \,+\, 2x_0$ for $k=m$.
So~(d) is an immediate consequence.
\qed

\medskip

{\bf Remark.}  When $\bfx = \bone$, we can use
\be
   \left. \scrd^-_m \prod_{i=0}^{m-1} x_i^{n_i} \right|_{\bfx = \bone}
   \;=\;
   (m+1) \sum_{i=0}^{m-1} n_i
\ee
together with the fact that $g_{k,n-1}$ is homogeneous in $\bfx$
of degree $n-1$, to specialize \reff{def.gkn.quasi-affine.multi} as
\be
   g_{k,n}(\bone)
   \;=\;
   \big[ (m+1)(n-1) + k +1 \bigr] g_{k,n-1}(\bone) \,+\, g_{k-m,n}(\bone)
   \;.
 \label{def.gkn.bone.multi}
\ee
It follows that for $0 \le k \le m-1$,
the $g_{k,n}(\bone)$ are shifted multifactorials:
\be
   g_{k,n}(\bone)  \;=\;  \prod_{j=0}^{n-1} [k+1+j(m+1)]  \;=\;  C_n(k+1,m+1)
   \;.
 \label{def.gkn.bone.multi.ans1}
\ee
Moreover, for $k=m$ we can specialize \reff{def.gkn.quasi-affine.multi.variant}
as
\be
   g_{m,n}(\bone)
   \;=\;
   \big[ (m+1)(n-1) + m+2 \bigr] g_{m,n-1}(\bone) \,+\, g_{-1,n}(\bone)
   \;,
\ee
to yield
\be
   g_{m,n}(\bone)  \;=\;  \prod_{j=0}^{n-1} [m+2+j(m+1)]  \;=\;  C_n(m+2,m+1)
 \label{def.gkn.bone.multi.ans2}
\ee
[note that we have skipped one step compared to \reff{def.gkn.bone.multi.ans1}].
However, for $k > m$ the forms get more complicated:
$g_{k,n}(\bone)$ is a linear combination of $\lfloor (k+1)/(m+1) \rfloor + 1$
shifted multifactorials.
These can be computed by symbolic algebra for any given $m$ and $k$:
for instance, for $m=2$ the first few are
\begin{subeqnarray}
   g_{0,n}(\bone)  & = & (3n-2)!!!  \;=\;  C_n(1,3)   \\[1mm]
   g_{1,n}(\bone)  & = & (3n-1)!!!  \;=\;  C_n(2,3)   \\[1mm]
   g_{2,n}(\bone)  & = & (3n+1)!!!  \;=\;  C_n(4,3)   \\[1mm]
   g_{3,n}(\bone)  & = & -(3n+1)!!! + (3n+2)!!!  \;=\; -C_n(4,3) + 2 C_n(5,3)  \\[1mm]
   g_{4,n}(\bone)  & = & - {(3n+2)!!! \over 2} + {(3n+4)!!! \over 2}  \;=\;  -C_n(5,3) + 2 C_n(7,3) \\[1mm]
   g_{5,n}(\bone)  & = & -(3n+4)!!! + {(3n+5)!!! \over 2} \;=\; -4C_n(7,3) + 5C_n(8,3)  \qquad
\end{subeqnarray}
But we have not been able to find a general
expression for the coefficients in these sums.
\myendremark

\bigskip

We have thus shown that $P^{(m,m)}(\bfx,\bfx) = (\scrd^-_m + x_0)^n \, 1$.
To complete the proof of Theorem~\ref{thm.quasi-affine.multi}(a)
we need to show that $\scrp^{(m)-}_n(\bfx) = (\scrd^-_m + x_0)^n \, 1$,
as asserted in Proposition~\ref{prop.MVeulerian.differential.multi}.

\proofof{Proposition~\ref{prop.MVeulerian.differential.multi}}
We will treat the cases of $\scrp^{(m)-}_n$ and $\scrq^{(m)-}_n$ in parallel,
and write $A_n$ as a shorthand for either one.

\begin{figure}[!ht]
\begin{center}
\includegraphics[scale=1]{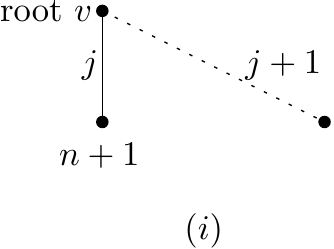}  \hspace*{1cm}
\includegraphics[scale=1]{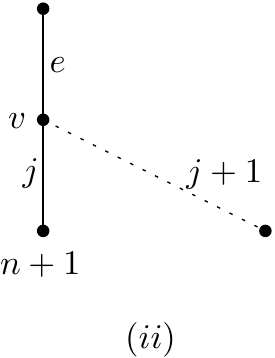} \hspace*{1cm}
\includegraphics[scale=1]{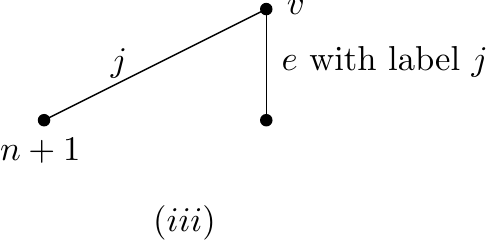}
\caption{
      The three possibilities for the vertex $n+1$.
}
\label{fig.negativetype}
\end{center}
\end{figure}

Consider vertex $n+1$ in an increasing multi-$m$-ary tree on $n+1$ vertices:
it is necessarily a $j$-child for some $j \in \{0,\ldots,m-1\}$
of some vertex $v \in [n]$.  There are three possibilities
(see Figure~\ref{fig.negativetype}):
\begin{itemize}
   \item[(i)]  $v$ is the root, and $n+1$ is its rightmost $j$-child.
        [In the case of $\scrp^{(m)-}_n$ this entails $j=0$.]
   \item[(ii)]  $v$ is a non-root vertex, and $n+1$ is its rightmost $j$-child.
        Let $e$ be the edge from $v$ to its parent.
   \item[(iii)]  $n+1$ is not the rightmost $j$-child of $v$;
        therefore $n+1$ has a $j$-sibling immediately to its right.
        Let $e$ be the edge from $v$ to this sibling.
\end{itemize}

\noindent
And conversely, an increasing multi-$m$-ary tree on $n+1$ vertices
is obtained from an increasing multi-$m$-ary tree on $n$ vertices
by inserting vertex $n+1$ in one of these three ways:
either we make $n+1$ the rightmost $j$-child of the root,
for some $j$ [$=0$ in the case of $\scrp^{(m)-}_n$];
or we choose an edge $e$ and insert $n+1$ in way~(ii) or way~(iii),
again for some $j$.

Possibility~(i) gives weight $x_0 A_{n-1}$ for $\scrp^{(m)-}_n$,
or $\Bigl( \sum\limits_{j=0}^{m-1} x_j \Bigr) A_{n-1}$ for $\scrq^{(m)-}_n$.

Possibility~(ii) gives weight
$\Bigl( \sum\limits_{j=0}^{m-1} x_j \Bigr)
 \Bigl( \sum\limits_{i=0}^{m-1} x_i \, {\partial \over \partial x_i} \Bigr)
 A_{n-1}$.

Possibility~(iii) gives weight
$\Bigl( \sum\limits_{j=0}^{m-1} x_j^2 \, {\partial \over \partial x_j} \Bigr)
 A_{n-1}$.

Adding these together gives $A_n = (\scrd^-_m + x_0) A_{n-1}$
for $\scrp^{(m)-}_n$, and
$A_n = {\Bigl( \scrd^-_m + \sum\limits_{j=0}^{m-1} x_j \Bigr)} A_{n-1}$
for $\scrq^{(m)-}_n$.

$\scrq^{(m)-}_{n,k}$ is handled similarly to $\scrq^{(m)-}_n$,
with two small differences:
in possibility~(i) there are $k+1$ different roots
where vertex $n+1$ could be attached;
or it could be made a new isolated vertex.
\qed

It is curious that the combinatorial arguments in the proofs of
Propositions~\ref{prop.MVeulerian.differential} and
\ref{prop.MVeulerian.differential.multi} are so different.
It would be of interest to find a unified understanding of these two results,
and ideally a generalization.

\subsection{A final observation}  \label{subsec.quasi-affine.observation}

Please note that there is a perfect analogy
between the periodic case (of period ${m+1}$ or $m$)
studied in Section~\ref{sec.periodic}
and the Eulerian-quasi-affine case (again of period~${m+1}$ or $m$)
studied here.
The Fuss--Narayana polynomials and series of positive type
($P_n^{(m)}$, $Q_n^{(m)}$, $Q_{n,k}^{(m)}$,
 $P_n^{(\infty)}$, $Q_n^{(\infty)}$, $Q_{n,k}^{(\infty)}$),
which correspond to period~${m+1}$,
and those of negative type
($P_n^{(m)-}$, $Q_n^{(m)-}$, $Q_{n,k}^{(m)-}$,
 $P_n^{(\infty)-}$, $Q_n^{(\infty)-}$, $Q_{n,k}^{(\infty)-}$),
which correspond to period~$m$,
are the generating polynomials or series
for {\em unlabeled}\/ ordered trees
(or ordered forests of unlabeled ordered trees)
of suitable types:
$(m+1)$-ary, $\infty$-ary, multi-$m$-ary or multi-$\infty$-ary, respectively,
and with ($P$) or without ($Q$) constraints on the labels of the edges
emanating from the root
(see Propositions~\ref{prop.fuss-narayana.trees}--%
\ref{prop.fuss-narayana.trees_extended.infty.symfun.P}).
Similarly, the
Eulerian polynomials and series of positive type
($\scrp_n^{(m)}$, $\scrq_n^{(m)}$, $\scrq_{n,k}^{(m)}$,
 $\scrp_n^{(\infty)}$, $\scrq_n^{(\infty)}$, $\scrq_{n,k}^{(\infty)}$),
which correspond to period~$m+1$,
and those of negative type
($\scrp_n^{(m)-}$, $\scrq_n^{(m)-}$, $\scrq_{n,k}^{(m)-}$,
 $\scrp_n^{(\infty)-}$, $\scrq_n^{(\infty)-}$, $\scrq_{n,k}^{(\infty)-}$),
which correspond to period~$m$,
are the generating polynomials or series
for {\em increasingly labeled}\/ ordered trees
(or unordered forests of increasingly labeled ordered trees)
of exactly the same types, with exactly the same weights.
In~both cases the simplest objects are the symmetric functions
$Q_{n,k}^{(\infty)} = \omega Q_{n,k}^{(\infty)-}$
and $\scrq_{n,k}^{(\infty)} = \omega \scrq_{n,k}^{(\infty)-}$,
in which all vertices of the tree are treated in an
identical and symmetric manner.
Then the $m$-Stieltjes--Rogers polynomials
$P_n^{(m)}, P_n^{(m)-}, \scrp_n^{(m)}, \scrp_n^{(m)-}$
and their $m \to \infty$ limits
$P_n^{(\infty)}, P_n^{(\infty)-}, \scrp_n^{(\infty)}, \scrp_n^{(\infty)-}$
are variants in which the edges emanating from the root
are constrained to be 0-edges;
as a consequence, these quantities are symmetric functions of $(x_i)_{i \ge 1}$
with a non-symmetric dependence on $x_0$.

\subsection[Stirling permutations and $r$th-order Eulerian polynomials]{Stirling permutations and $\bm{r}$th-order Eulerian polynomials}
   \label{subsec.quasi-affine.stirling}

A word $\bfw = w_1 \cdots w_L$ on a totally ordered alphabet $\mathbb{A}$
is called a \textbfit{Stirling word} if $i < j < k$ and $w_i = w_k$ imply
$w_j \ge w_i$:
that is, between any two occurrences of any letter $a$,
only letters that are larger than or equal to $a$ are allowed.
(Equivalently, between any two successive occurrences of the letter $a$,
only letters that are larger than $a$ are allowed.)
Now let $\bfr = (r_1,\ldots,r_n)$ be a finite sequence of nonnegative integers,
and define the multiset $M_\bfr = \{1^{r_1}, 2^{r_2}, \ldots, n^{r_n}\}$
consisting of $r_i$ copies of the letter $i$;
it has total cardinality $|\bfr| = \sum r_i$.
A \textbfit{permutation} of $M_\bfr$ is a word $w_1 \cdots w_{|\bfr|}$
containing $r_i$ copies of the letter $i$, for each $i \in [n]$;
it is called a \textbfit{Stirling permutation} of $M_\bfr$
if it is also a Stirling word.
The number of Stirling permutations of $M_\bfr$ is
\cite{Janson_11,Dzhumadildaev_14}
\be
   (r_1 + 1) (r_1 + r_2 + 1) \cdots (r_1 + \ldots + r_{n-1} + 1)
   \;,
 \label{eq.enum.stirling}
\ee
which is easy to see by induction since the $r_n$ copies of $n$
have to occur as a block,
and this block can be inserted in any of
$r_1 + \ldots + r_{n-1} + 1$ positions
in a Stirling permutation of $\{1^{r_1}, 2^{r_2}, \ldots, (n-1)^{r_{n-1}}\}$.

Stirling permutations were introduced by Gessel and Stanley \cite{Gessel_78}
for the case $r_1 = \ldots = r_n = 2$;
this was generalized to $r_1 = \ldots = r_n = r$
[which we denote by the shorthand $\bfr = (r^n)$]
by Gessel \cite{Gessel_78a} and Park \cite{Park_94a,Park_94b},
and to general multisets $M_\bfr$ by
Brenti \cite{Brenti_89,Brenti_98} and
others \cite{Janson_11,Dzhumadildaev_14}.
We refer to Stirling permutations of $M_{(r^n)}$
as \textbfit{$\bm{r}$-Stirling permutations} of order~$n$.
It follows from \reff{eq.enum.stirling}
that the number of $r$-Stirling permutations of order~$n$ is
the multifactorial $F_n^{(r)}$.

For any word $\bfw = w_1 \cdots w_L$
on a totally ordered alphabet $\mathbb{A}$,
a pair $(i,i+1)$ with $1 \le i \le L-1$
is called a \textbfit{descent} (resp.\ \textbfit{ascent}, \textbfit{plateau})
if $w_i > w_{i+1}$ (resp.\ $w_i < w_{i+1}$, $w_i = w_{i+1}$).
We call the index $i$ (resp.\ $i+1$) the beginning (resp.\ end)
of the descent, ascent or plateau.
We write $\des(\bfw)$ for the number of descents in the word $\bfw$.
If $\bfw$ is a Stirling word, then only the last occurrence of each letter $j$
can be the beginning of a descent,
and this can happen only for letters $j$
that are not the smallest letter of $\mathbb{A}$;
so $0 \le \des(\bfw) \le |\mathbb{A}| -1$ for $\bfw \ne \emptyset$
(of course $\des(\emptyset) = 0$).
We write $\euler{\bfr}{k}$ for the number of Stirling permutations
of $M_\bfr$ with $k$ descents,
and call them the {\em multiset Eulerian numbers}\/;
we then define the {\em $\bfr$-Eulerian polynomial}\/
$A_\bfr(x) = \sum\limits_{k=0}^n \euler{\bfr}{k} x^k$.
When $\bfr = (r^n)$ we write $\euler{\bfr}{k} = \euler{n}{k}^{\! (r)}$
and call them the \textbfit{$\bm{r^{th}}$-order Eulerian numbers};
for $r=1,2$ we also use the notations $\euler{n}{k}$
and $\eulersecond{n}{k}$, respectively.\footnote{
   Here we follow the convention of Graham {\em et al.}\/ \cite{Graham_94}
   that (when $n \ge 1$)
   $\euler{n}{k}^{\! (r)}$ is nonzero for $0 \le k \le n-1$.
   Many other authors \cite{Gessel_78,Park_94a,Janson_11,Dzhumadildaev_14}
   use the convention that the final index $L$
   is always the beginning of a descent
   (e.g.~by imposing the boundary condition $w_{L+1} = 0$),
   so that (when $n \ge 1$)
   $\euler{n}{k}^{\! (r)}$ is nonzero for $1 \le k \le n$.
}
We define the \textbfit{$\bm{r^{th}}$-order Eulerian polynomials}
\be
   A_n^{(r)}(x) 
   \;=\;
   \sum\limits_{k=0}^n \euler{n}{k}^{\!\! (r)} x^k
   \;.
 \label{def.eulerianr}
\ee
The $r$th-order Eulerian numbers satisfy the recurrence
\be
   \euler{n}{k}^{\! (r)}   \;=\;
      (k+1) \euler{n-1}{k}^{\!\! (r)}   \,+\,
      \big[ rn-(r-1)-k \big] \euler{n-1}{k-1}^{\!\! (r)}
   \quad\hbox{for } n \ge 1
 \label{eq.eulerian.recurrence}
\ee
with initial condition $\euler{0}{k}^{\! (r)} = \delta_{0k}$.
It follows easily from \reff{eq.eulerian.recurrence} that
$\euler{n}{n-1}^{\! (r)} = F_n^{(r-1)}$ and
$\sum\limits_{k=0}^n \euler{n}{k}^{\! (r)} = F_n^{(r)}$.

A bijection from $r$-Stirling permutations of order~$n$
to increasing $(r+1)$-ary trees with $n$ vertices
was found by Gessel \cite{Gessel_78a} (see \cite{Park_94a})
and independently by
Janson {\em et al.}\/ \cite{Janson_11} (see also \cite{Kuba_09}),
and can be described as follows:
Decompose the Stirling permutation as
$\bfw = \bfw_0 1 \bfw_1 1 \cdots \bfw_{r-1} 1 \bfw_r$.
Put vertex 1 at the root, and for each nonempty subword $\bfw_i$,
create an $i$-child of the root provisionally labeled $\bfw_i$;
then repeat this process recursively on each subword,
using the smallest letter of the subword in place of 1;
at the end, each node of the tree will have a label $j \in [n]$.
See \cite{Janson_11} for a description of the inverse bijection.

We can use this bijection to map statistics between
Stirling permutations and trees.
It is not difficult to prove recursively that:

\begin{lemma}
   \label{lemma.stirling.descents}
Let $\bfw$ be an $r$-Stirling permutation of order~$n$,
and let $T$ be the corresponding increasing $(r+1)$-ary tree
on the vertex set $[n]$.  Then:
\begin{itemize}
   \item[(i)] For $1 \le j \le n$ and $1 \le i \le r-1$,
       the following are equivalent:
       \begin{itemize}
           \item[(a)]  The vertex $j \in T$ has an $i$-child.
           \item[(b)]  In the word $\bfw$, between the $i$th and $(i+1)$st
               occurrences of the letter $j$ there is a nonempty subword.
               [Of course, all entries in this subword are $> j$.]
           \item[(c)]  In the word $\bfw$, the $i$th occurrence of the
               letter $j$ is the beginning of an ascent.
           \item[(d)]  In the word $\bfw$, the $(i+1)$st occurrence of the
               letter $j$ is the end of a descent.
       \end{itemize}
   \item[(ii)] For $1 \le j \le n$ and $i=0$, the following are equivalent:
       \begin{itemize}
           \item[(a)]  The vertex $j \in T$ has a $0$-child.
           \item[(d)]  In the word $\bfw$, the first occurrence of the
               letter $j$ is the end of a descent.
       \end{itemize}
   \item[(iii)] For $1 \le j \le n$ and $i=r$, the following are equivalent:
       \begin{itemize}
           \item[(a)]  The vertex $j \in T$ has an $r$-child.
           \item[(c)]  In the word $\bfw$, the last (i.e.\ $r$th)
               occurrence of the letter $j$ is the beginning of an ascent.
       \end{itemize}
\end{itemize}
\end{lemma}

It follows from Lemma~\ref{lemma.stirling.descents}
that the multivariate Eulerian polynomial $\scrq_{n-1}^{(r)}(\bfx)$
is the generating polynomial for $r$-Stirling permutations of order~$n$
in which we give a weight $x_i$ ($0 \le i \le r-1$)
whenever the $(i+1)$st occurrence of a letter is the end of a descent,
and a weight $x_r$ whenever the last (i.e.\ $r$th) occurrence of a letter
is the beginning of an ascent.
Therefore, taking $x_0 = \ldots = x_{r-1} = x$ and $x_r = 1$,
we see that $\scrq_{n-1}^{(r)}(x,\ldots,x,1)$
gives a weight $x$ for all descents
and hence is the $r$th-order Eulerian polynomial $A_n^{(r)}(x)$.
But $\scrq_{n-1}^{(r)}(\bfx)$ is symmetric in its arguments,
so we can equivalently take $x_0 = 1$ and $x_1 = \ldots = x_r = x$,
and write
\be
   A_n^{(r)}(x)
   \;=\;
   \scrq_{n-1}^{(r)}(x,\ldots,x,1)
   \;=\;
   \scrq_{n-1}^{(r)}(1,x,\ldots,x)
   \;=\;
   \scrp_n^{(r)}(1,x,\ldots,x)
   \;.
 \label{eq.scrp.eulerianr}
\ee
By Theorem~\ref{thm.quasi-affine}(a) and Corollary~\ref{cor.quasi-affine.2}
we can conclude:

\begin{corollary}[$r$th-order Eulerian polynomials $A_n^{(r)}$]
   \label{cor.eulerianr.0} 
For any integer $r \ge 1$,
the $r$th-order Eulerian polynomial $A_n^{(r)}(x)$
defined in \reff{def.eulerianr} satisfies
\begin{subeqnarray}
   A_n^{(r)}(x)
   \;=\;
   \scrq_{n-1}^{(r)}(1,x,\ldots,x)
   & = &
   \scrp_n^{(r)}(1,x,\ldots,x)
         \\[2mm]
   & = &
   P_n^{(r,r+1)}(\bfx,\bfx) \; \hbox{where $\bfx = (1,x,\ldots,x)$} \,.
      \qquad\quad
 \label{eq.cor.eulerianr.0}
\end{subeqnarray}
That is, it equals the $r$-Stieltjes--Rogers polynomial
$S_n^{(r)}(\balpha)$ where the weights $\balpha$ are given by
the Eulerian-quasi-affine formula \reff{eq.alpha.quasi-affine.u=x}
with period $p=r+1$ and $\bfx = (1,x,\ldots,x)$.

Therefore the sequence $(A_n^{(r)}(x))_{n \ge 0}$
is coefficientwise Hankel-totally positive in $x$.
\end{corollary}

A slightly different definition of the $r$th-order Eulerian polynomials
was used in \cite{Elvey-Price-Sokal_wardpoly}:
\begin{subeqnarray}
   E_0^{[r]}(x)  & = &  1  \\
   E_n^{[r]}(x)  & = &  x \,A_n^{(r)}(x)
       \;=\;  \sum\limits_{k=0}^{n-1} \euler{n}{k}^{\!\! (r)} x^{k+1}
           \quad\hbox{for $n \ge 1$}
 \label{eq.reveuler0}
\end{subeqnarray}
Then the reversed $r$th-order Eulerian polynomials
are defined by $\overline{E}_n^{[r]}(x) = x^n \, E_n^{[r]}(1/x)$, i.e.
\begin{subeqnarray}
   \overline{E}_0^{[r]}(x)  & = &  1  \\[2mm]
   \overline{E}_n^{[r]}(x)  & = &  x^{n-1} A_n^{(r)}(1/x)
       \;=\;  \sum\limits_{k=0}^{n-1} \euler{n}{k}^{\!\! (r)} x^{n-1-k}
           \quad\hbox{for $n \ge 1$}
 \slabel{eq.reveuler.b}
 \label{eq.reveuler}
\end{subeqnarray}
and satisfy $\overline{E}_n^{[r]}(0) = F_n^{(r-1)}$
and $\overline{E}_n^{[r]}(1) = F_n^{(r)}$.
Then \reff{eq.reveuler.b} and \reff{eq.scrp.eulerianr} give
\be
   \overline{E}_n^{[r]}(x)  \;=\;  x^{n-1} A_n^{(r)}(1/x)
      \;=\; \scrq_{n-1}^{(r)}(x,1,\ldots,1)
   \;.
\ee
But since $\scrq_{n-1}^{(r)}(\bfx)$ is symmetric in its arguments,
we can equivalently take one of $x_1,\ldots,x_r$ (say, $x_r$) to be $x$
and all the other~$x_i$ (including $x_0$) to be 1;
then
\be
   \overline{E}_n^{[r]}(x)
   \;=\;
   \scrq_{n-1}^{(r)}(1,\ldots,1,x)
   \;=\;
   \scrp_n^{(r)}(1,\ldots,1,x)
   \;.
 \label{eq.scrp.eulerianr.bis}
\ee
By Theorem~\ref{thm.quasi-affine}(a) and Corollary~\ref{cor.quasi-affine.2}
we can conclude:

\begin{corollary}[Reversed $r$th-order Eulerian polynomials $\overline{E}_n^{[r]}$]
   \label{cor.eulerianr} 
For any integer $r \ge 1$,
the reversed $r$th-order Eulerian polynomial $\overline{E}_n^{[r]}(x)$
defined in \reff{eq.reveuler0}/\reff{eq.reveuler}
satisfies
\begin{subeqnarray}
   \overline{E}_n^{[r]}(x)
   \;=\;
   \scrq_{n-1}^{(r)}(1,\ldots,1,x)
   & = &
   \scrp_n^{(r)}(1,\ldots,1,x)
         \\[2mm]
   & = &
   P_n^{(r,r+1)}(\bfx,\bfx) \; \hbox{where $\bfx = (1,\ldots,1,x)$} \,.
      \qquad\quad
 \label{eq.cor.eulerianr}
\end{subeqnarray}
That is, it equals the $r$-Stieltjes--Rogers polynomial
$S_n^{(r)}(\balpha)$ where the weights $\balpha$ are given by
the Eulerian-quasi-affine formula \reff{eq.alpha.quasi-affine.u=x}
with period $p=r+1$ and $\bfx = (1,\ldots,1,x)$.

Therefore the sequence $(\overline{E}_n^{[r]}(x))_{n \ge 0}$
is coefficientwise Hankel-totally positive in $x$.
\end{corollary}

Of course, this univariate result
--- which even for $r=2$ was a mystery for us
\cite{Elvey-Price-Sokal_wardpoly} until recently ---
is merely a very special case
of the results for multivariate Eulerian polynomials.

\section{Ratios of contiguous hypergeometric series I: $\bm{\FHyper{m+1}{0}}$}
   \label{sec.hyper.rF0}

As mentioned in the introduction,
Euler \cite[section~21]{Euler_1760} found in 1746 the continued fraction
\be
   \sum_{n=0}^\infty n! \: t^n
   \;=\;
   \cfrac{1}{1 - \cfrac{1t}{1 - \cfrac{1t}{1 - \cfrac{2t}{1- \cfrac{2t}{1- \cdots}}}}}
 \label{eq.nfact.contfrac.bis}
\ee
with coefficients $\alpha_{2k-1} = \alpha_{2k} = k$.
In fact, in the same paper Euler \cite[section~26]{Euler_1760}
found the more general continued fraction
\be
   \sum_{n=0}^\infty a(a+1)(a+2) \cdots (a+n-1) \: t^n
   \;=\;
   \cfrac{1}{1 - \cfrac{at}{1 - \cfrac{1t}{1 - \cfrac{(a+1)t}{1- \cfrac{2t}{1- \cdots}}}}}
 \label{eq.xupperfact.contfrac}
\ee
with coefficients $\alpha_{2k-1} = a+k-1$ and $\alpha_{2k} = k$.
And this is, in turn, the $b=1$ special case of the beautiful continued
fraction for ratios of contiguous hypergeometric series $\FHyper{2}{0}$
\cite[section~92]{Wall_48}:
\be
   {\FHYPERbottomzero{2}{a,b}{t}
    \over
    \FHYPERbottomzero{2}{a,b-1}{t}
   }
   \;\:=\;\:
   \cfrac{1}{1 - \cfrac{at}{1 - \cfrac{bt}{1 - \cfrac{(a+1)t}{1- \cfrac{(b+1)t}{1 - \cdots}}}}}
 \label{eq.euler.contfrac.BISBIS.2F0}
\ee
with coefficients $\alpha_{2k-1} = a+k-1$ and $\alpha_{2k} = b+k-1$.
(We do not know who was the first to discover
 \reff{eq.euler.contfrac.BISBIS.2F0},
 which is a limiting case of Gauss' \cite{Gauss_1813}
 continued fraction for ratios of contiguous $\tfo$.)

Here we will show that all these continued fractions
have ``higher'' generalizations:
namely, for every integer $m \ge 1$,
the ratio of contiguous hypergeometric series $\FHyper{m+1}{0}$
has a nice $m$-branched continued fraction.

\begin{theorem}[$m$-branched continued fraction for ratios of contiguous $\FHyper{m+1}{0}$]
   \label{thm.rF0}
\nopagebreak\hfill\break\nopagebreak
Fix an integer $m \ge 1$,
and define the polynomials $P_n^{(m)}(a_1,\ldots,a_m;a_{m+1})$ by
\be
   \sum_{n=0}^\infty P_n^{(m)}(a_1,\ldots,a_m;a_{m+1}) \: t^n
   \;\:=\;\:
   {\FHYPERbottomzero{m+1}{a_1,\ldots,a_{m+1}}{t}
    \over
    \FHYPERbottomzero{m+1}{a_1,\ldots,a_m,a_{m+1}-1}{t}
   }
   \;\,.
 \label{eq.thm.rF0}
\ee
Then $P_n^{(m)}(a_1,\ldots,a_m;a_{m+1}) = S_n^{(m)}(\balpha)$
where $\balpha = (\alpha_i)_{i \ge m}$ is given by
\be
   \balpha
   \;=\;
   a_1 \cdots a_m, \, 
   a_2 \cdots a_{m+1}, \, 
   a_3 \cdots a_{m+1} (a_1 + 1), \, 
   a_4 \cdots a_{m+1} (a_1 + 1)(a_2 + 1), \, 
   \ldots
   \;\,,
 \label{eq.thm.rF0.alphas}
\ee
which can be seen as products of $m$ successive pre-alphas:
\be
   \balphapre 
   \;=\;
   a_1,\ldots,a_{m+1}, a_1+1,\ldots,a_{m+1}+1, a_1+2,\ldots,a_{m+1}+2, \ldots
   \;\,.
 \label{eq.thm.rF0.prealphas}
\ee
\end{theorem}

Please note that here the ``pre-alphas'' \reff{eq.thm.rF0.prealphas}
are quasi-affine of period~$m+1$.
So this is a kind of ``higher'' generalization of the situation
studied in the preceding section:
no longer are the $\balpha$ themselves quasi-affine (when $m > 1$);
rather, they are products of $m$ successive pre-alphas that are quasi-affine.

Please note also that the polynomials $P_n^{(m)}(a_1,\ldots,a_m;a_{m+1})$
are symmetric in $a_1,\ldots,a_m$ (while $a_{m+1}$ plays a distinguished role).
But the $\balpha$ defined in \reff{eq.thm.rF0.alphas}
are not symmetric in $a_1,\ldots,a_m$ (when $m > 1$).
This illustrates once again the nonuniqueness of $m$-S-fractions when $m \ge 2$.

The proof of Theorem~\ref{thm.rF0} will be based on
the Euler--Gauss recurrence method
as generalized to $m > 1$ in Proposition~\ref{prop.euler-gauss.mSR}.\footnote{
   As mentioned earlier, this method (for $m=1$) was used implicitly
   by Euler \cite[section~21]{Euler_1760}
   for proving \reff{eq.nfact.contfrac.bis}.
   In \cite[section~26]{Euler_1760}, Euler stated that
   the same method can be applied to the more general series
   \reff{eq.xupperfact.contfrac},
   which reduces to \reff{eq.nfact.contfrac.bis} when $a=1$;
   but he did not provide the details, and he instead
   proved \reff{eq.xupperfact.contfrac} by an alternative method.
   Three decades later, however, Euler \cite{Euler_1788}
   returned to his original method and presented the details
   of the derivation of \reff{eq.xupperfact.contfrac}.
 \label{footnote_Euler}
}
Namely, fix an integer $m \ge 1$,
and let $(g_k(t))_{k \ge -1}$ be a sequence of formal power series
(with coefficients in some commutative ring $R$)
with constant term 1,
which satisfies a linear three-term recurrence of the form
\be
   g_k(t) - g_{k-1}(t)  \;=\; \alpha_{k+m} t \, g_{k+m}(t)
   \qquad\hbox{for } k \ge 0
 \label{eq.recurrence.gkm}
\ee
for some coefficients $\balpha = (\alpha_i)_{i \ge m}$ in $R$.
Then $g_0(t)/g_{-1}(t) = \sum_{n=0}^\infty S^{(m)}_n(\balpha) \, t^n$
where $S^{(m)}_n(\balpha)$ is the $m$-Stieltjes--Rogers polynomial
evaluated at the specified values $\balpha$.

Now it is easy to show that a suitably defined sequence of
contiguous hypergeometric series $\FHyper{m+1}{0}$
satisfies a recurrence of the form \reff{eq.recurrence.gkm}.
The key fact is the following:

\begin{lemma}[Three-term contiguous relation for $\FHyper{r}{0}$]
   \label{lemma.recurrence.rF0}
Fix an integer $r \ge 1$.  Then the hypergeometric series $\FHyper{r}{0}$
satisfies
\begin{eqnarray}
   & & \hspace*{-1cm}
   \FHYPERbottomzero{r}{a_1,\ldots,a_{i-1},a_i+1, a_{i+1},\ldots,a_r}{t}
   \:-\:
   \FHYPERbottomzero{r}{a_1,\ldots,a_r}{t}
        \hspace*{1cm}
        \nonumber \\[4mm]
   & & \hspace*{2cm}
   \;=\;
   a_1 \,\cdots \not\mathrel{\,a}_i \cdots\, a_r \, t \;\:
   \FHYPERbottomzero{r}{a_1 +1,\ldots,a_r +1}{t}
 \label{eq.lemma.recurrence.rF0}
\end{eqnarray}
where $\not\mathrel{\,a}_i$ indicates that $a_i$ is omitted from the product.
\end{lemma}

\noindent
The identity \reff{eq.lemma.recurrence.rF0} is easily proven
by comparing coefficients of $t^n$ on both sides;
it~is a special case of a more general identity for
hypergeometric series $\FHyper{r}{s}$ [see \reff{eq.second.1} below].

Now define
\begin{subeqnarray}
   g_{-1}(t)
   & = &
   \FHYPERbottomzero{m+1}{a_1,\ldots,a_m,a_{m+1}-1}{t}
           \\[2mm]
   g_0(t)
   & = &
   \FHYPERbottomzero{m+1}{a_1,\ldots,a_{m+1}}{t}
\end{subeqnarray}
and define $g_1,g_2,\ldots$ by successively incrementing
$a_1,\ldots,a_{m+1}$ by 1, continuing cyclically;  thus
\be
   g_k(t)
   \;=\;
   \FHYPERbottomzero{m+1}{a_1 + \lceil {k \over m+1} \rceil,\,
                          a_2 + \lceil {k-1 \over m+1} \rceil,\,
                          \ldots,\,
                          a_{m+1} + \lceil {k-m \over m+1} \rceil
                         }{t}
\ee
for all $k \ge -1$.
It is easy to see, using \reff{eq.lemma.recurrence.rF0},
that the sequence $(g_k(t))_{k \ge -1}$
satisfies the recurrence \reff{eq.recurrence.gkm}
with the $\balpha$ stated in \reff{eq.thm.rF0.alphas}.
This completes the proof of Theorem~\ref{thm.rF0}.

When $a_{m+1} = 1$, Theorem~\ref{thm.rF0} simplifies because
the denominator hypergeometric series reduces to 1,
and we get:

\begin{corollary}[$m$-branched continued fraction for products of Stirling cycle polynomials]
   \label{cor.rF0}
Define the homogenized Stirling cycle polynomials by
\be
   C_n(x,y)  \;\eqdef\;  \prod_{j=0}^{n-1} (x+jy)
\ee
[cf.\ \reff{def.stirlingcycle}].
Then the polynomials
\be
   P_n^{(m)}(a_1,\ldots,a_m;1)
   \;=\;
   \prod_{i=1}^m C_n(a_i,1)
 \label{eq.cor.rF0.2}
\ee
equal $S_n^{(m)}(\balpha)$ where $\balpha = (\alpha_i)_{i \ge m}$ is given by
\reff{eq.thm.rF0.alphas} with $a_{m+1} = 1$.
More generally, the polynomials
\be
   \prod_{i=1}^m C_n(x_i,y_i)
   \;=\;
   \biggl( \prod\limits_{i=1}^m y_i \biggr)^{\! n} \,
      P_n^{(m)}(x_1/y_1,\ldots,x_m/y_m;1)
 \label{eq.cor.rF0.3}
\ee
equal $S_n^{(m)}(\balpha)$ where $\balpha = (\alpha_i)_{i \ge m}$ is given by
\reff{eq.thm.rF0.alphas} evaluated at $a_i = x_i/y_i$ and $a_{m+1} = 1$
and multiplied by $\prod\limits_{i=1}^m y_i$, i.e.
\be
   \balpha
   \;=\;
   x_1 \cdots x_m, \, 
   x_2 \cdots x_m y_1 ,\, 
   x_3 \cdots x_m y_2 (x_1 + y_1) ,\,
   x_4 \cdots x_m y_3 (x_1 + y_1) (x_2 + y_2) ,\,
   \ldots
   \;\,.
\ee
\end{corollary}

\bigskip

{\bf Examples.}
1.  Take $a_1 = \ldots = a_m = 1$ in Corollary~\ref{cor.rF0};
then $P_n^{(m)} = (n!)^m$.
The pre-alphas are 
\be
   \balphapre
   \;=\;
   \underbrace{1,\ldots,1}_{\text{$m+1$ times}},
   \underbrace{2,\ldots,2}_{\text{$m+1$ times}},
   \underbrace{3,\ldots,3}_{\text{$m+1$ times}},
   \ldots
\ee
and the $\balpha$ are products of $m$ successive pre-alphas.
For $m=2$ this is shown in \reff{eq.mCF.nfactsq}.

2.  Take $a_1 = \ldots = a_m = \smhalf$ in Corollary~\ref{cor.rF0};
then $P_n^{(m)} = (2n-1)!!^m / 2^{mn}$.
We can get the sequence $(2n-1)!!^m$ by multiplying all the $\balpha$ by $2^m$,
or equivalently by multiplying all the $\balphapre$ by $2$.
The pre-alphas are thus
\be
   \balphapre
   \;=\;
   \underbrace{1,\ldots,1}_{\text{$m$ times}}, 2,
   \underbrace{3,\ldots,3}_{\text{$m$ times}}, 4,
   \underbrace{5,\ldots,5}_{\text{$m$ times}}, 6,
   \ldots
\ee
and the $\balpha$ are products of $m$ successive pre-alphas.
Similar methods handle $(3n-2)!!!$, $(3n-1)!!!$ and so forth.

3.  Take $a_j = j/m$ for $1 \le j \le m$ in Corollary~\ref{cor.rF0};
then $P_n^{(m)} = (mn)! / (m^m)^n$.
We can get the sequence $(mn)!$ by multiplying all the $\balpha$ by $m^m$,
or equivalently by multiplying all the $\balphapre$ by $m$.
The pre-alphas are thus
\be
   \balphapre
   \;=\;
   1,2,3,\ldots,m,m, m+1,m+2,m+3,\ldots,2m,2m,
   \ldots
\ee
and the $\balpha$ are products of $m$ successive pre-alphas.
For $m=2$ this is shown in \reff{eq.mCF.2nfact}.
\myendremark

\medskip

Combining Theorem~\ref{thm.rF0} with Theorem~\ref{thm.iteration2bis},
we conclude:

\begin{corollary}[Hankel-total positivity for ratios of contiguous $\FHyper{m+1}{0}$]
   \label{cor.rF0.TP}
The sequence of polynomials $(P_n^{(m)}(a_1,\ldots,a_m;a_{m+1}))_{n \ge 0}$
defined in \reff{eq.thm.rF0}
is coefficientwise Hankel-totally positive,
jointly in all the indeterminates $a_1,\ldots,a_{m+1}$.

In particular, the product \reff{eq.cor.rF0.3}
of homogenized Stirling cycle polynomials
is coefficientwise Hankel-totally positive,
jointly in all the indeterminates $x_1,\ldots,x_m$ and $y_1,\ldots,y_m$.
\end{corollary}

{\bf Remark.}
This result was conjectured a few years ago by one of us
\cite{Sokal_unpub_2015}, based on computations of Hankel minors
up to $8 \times 8$ ($m=2$) and $7 \times 7$ ($m=3$);
but at that time he had no idea how to prove it.
Even the fact that $P_n^{(m)}(a_1,\ldots,a_m;a_{m+1})$
has nonnegative coefficients is not completely trivial;
but this was proven by Gessel \cite{Gessel_private}
using \reff{eq.lemma.recurrence.rF0}.
\myendremark

\bigskip

It follows from \reff{def.stirlingcycle}/\reff{eq.cor.rF0.2}
that the polynomials $P_n^{(m)}(a_1,\ldots,a_m;1)$
have a simple combinatorial interpretation in terms of
$m$ independent permutations $\sigma_1,\ldots,\sigma_m$ of $[n]$:
\begin{subeqnarray}
   P_n^{(m)}(a_1,\ldots,a_m;1)
   & = &
   \sum_{\sigma_1,\ldots,\sigma_m \in \Sym_n}
   a_1^{\cyc(\sigma_1)} \,\cdots\, a_m^{\cyc(\sigma_m)}
           \\[4mm]
   & = &
   \sum_{\sigma_1,\ldots,\sigma_m \in \Sym_n}
   a_1^{\rec(\sigma_1)} \,\cdots\, a_m^{\rec(\sigma_m)}
\end{subeqnarray}
where $\cyc(\sigma)$ [resp.\ $\rec(\sigma)$]
denotes the numbers of cycles (resp.\ records) in the permutation $\sigma$.
It is natural to guess that the more general polynomial
$P_n^{(m)}(a_1,\ldots,a_m;a_{m+1})$ has an interpretation
\begin{subeqnarray}
   P_n^{(m)}(a_1,\ldots,a_m;a_{m+1})
   & = &
   \sum_{\sigma_1,\ldots,\sigma_m \in \Sym_n}
   \!\!
   a_1^{\cyc(\sigma_1)} \,\cdots\, a_m^{\cyc(\sigma_m)}
        \, a_{m+1}^{\mysteryone(\sigma_1,\ldots,\sigma_m)}
        \qquad
           \\[4mm]
   & = &
   \sum_{\sigma_1,\ldots,\sigma_m \in \Sym_n}
   \!\!
   a_1^{\rec(\sigma_1)} \,\cdots\, a_m^{\rec(\sigma_m)}
        \, a_{m+1}^{\mysterytwo(\sigma_1,\ldots,\sigma_m)}
        \qquad
\end{subeqnarray}
for some ``mystery statistics'' $\mysteryone$ and $\mysterytwo$;
but for $m > 1$ we have been unable to determine what these statistics are.
For $m=1$ the answer is known:
\begin{subeqnarray}
   P_n^{(1)}(a;b)
   & = &
   \sum_{\sigma \in \Sym_n}
   a^{\cyc(\sigma)} \, b^{\earec(\sigma)}
        \qquad
         \slabel{eq.recantirec.a}  \\[4mm]
   & = &
   \sum_{\sigma \in \Sym_n}
   a^{\rec(\sigma)} \, b^{\earec(\sigma)}
        \qquad
         \slabel{eq.recantirec.b}
         \label{eq.recantirec}
\end{subeqnarray}
where $\earec(\sigma)$ is the number of exclusive antirecords in $\sigma$
(that is, antirecords that are not also records).
Here \reff{eq.recantirec.b} was found three decades ago
by Dumont and Kreweras \cite{Dumont_88},
while \reff{eq.recantirec.a} is a specialization of
the result \reff{eq.eulerian.fourvar.cyc}
found recently by Zeng and one of us \cite{Sokal-Zeng_masterpoly}.
So we suspect that for $m > 1$ the ``mystery statistics''
have something to do with exclusive antirecords;
but we have been unable to find any viable candidate statistics.

\begin{openproblem}
\rm
Find the ``mystery statistics'' $\mysteryone$ and $\mysterytwo$.
\end{openproblem}

\section{Ratios of contiguous hypergeometric series II: $\bm{\FHyper{r}{s}}$}
   \label{sec.hyper.rFs}

By a slight generalization of the method employed in the previous section,
we can obtain a branched continued fraction for
the ratio of contiguous hypergeometric series $\FHyper{r}{s}$
for arbitrary integers $r,s$.
We will in fact consider {\em three}\/ different types of ratios of
contiguous hypergeometric series:
\begin{eqnarray}
   {\FHYPER{r}{s\,}{a_1,\ldots,a_r}{b_1,\ldots,b_s}{t}
    \over
    \FHYPER{r}{s\,}{a_1,\ldots,a_{r-1},a_r-1}{b_1,\ldots,b_{s-1},b_s-1}{t}
   }
   & \eqdef &
   \sum_{n=0}^\infty R_n^{(r,s)}(\bfa,\bfb)
                     \: t^n
   \quad\hbox{for $r,s \ge 1$} \quad
       \label{def.Rrs}  \\[4mm]
   {\FHYPER{r}{s\,}{a_1,\ldots,a_r}{b_1,\ldots,b_s}{t}
    \over
    \FHYPER{r}{s\,}{a_1,\ldots,a_{r-1},a_r-1}{b_1,\ldots,b_s}{t}
   }
   & \eqdef &
   \sum_{n=0}^\infty U_n^{(r,s)}(\bfa,\bfb)
                     \: t^n
   \quad\hbox{for $r \ge 1$, $s \ge 0$} \quad
       \label{def.Urs}  \\[4mm]
   {\FHYPER{r}{s\,}{a_1,\ldots,a_r}{b_1,\ldots,b_s}{t}
    \over
    \FHYPER{r}{s\,}{a_1,\ldots,a_r}{b_1,\ldots,b_{s-1},b_s-1}{t}
   }
   & \eqdef &
   \sum_{n=0}^\infty V_n^{(r,s)}(\bfa,\bfb)
                     \: t^n
   \quad\hbox{for $r \ge 0$, $s \ge 1$} \quad
      \label{def.Vrs}
\end{eqnarray}
We refer to \reff{def.Rrs}--\reff{def.Vrs}, respectively,
as the {\em first}\/, {\em second}\/ and {\em third}\/
ratios of contiguous hypergeometric series.
We will give branched continued fractions for all three:
in the first and second cases, the branching order is $m = \max(r-1,s)$,
while in the third case it is $m = \max(r,s)$.\footnote{
   Let us remark that the trivial case
   of the second continued fraction for $\ofz$
   can indeed be viewed as yielding a 0-branched continued fraction,
   since the ratio of contiguous series is $1/(1-t)$.
   However, we will henceforth exclude from consideration
   this degenerate case.
}
When $(r,s) = (2,1)$, the first continued fraction reduces to
Gauss' \cite{Gauss_1813} continued fraction for ratios of contiguous $\tfo$.
The second continued fraction may possibly be new,
even in the classical case of $\tfo$
(see the remarks below).

The proofs of these continued fractions will be based, once again,
on verifying the recurrence \reff{eq.recurrence.gkm}
for a suitably defined sequence $(g_k(t))_{k \ge 1}$
of hypergeometric series $\FHyper{r}{s}$.
And the verifications of this recurrence will be based, once again,
on three-term contiguous relations for $\FHyper{r}{s}$.
The needed contiguous relations are not new
\cite{Krattenthaler_contiguous,Krattenthaler_HYP};
but since they do not seem to be well known,
we devote a short subsection to
recalling them.
We then turn to the statement and proofs of the branched continued fractions.
We begin with the first continued fraction for $\FHyper{m+1}{m}$,
for which we need to use only one type of three-term contiguous relation.
Then we treat the first continued fraction
for $\FHyper{r}{s}$ for general $r,s$,
which requires combining two different
three-term contiguous relations in a suitable order,
or alternatively taking limits from $\FHyper{m+1}{m}$.
Then we consider analogously the second and third continued fractions.

\subsection[Three-term contiguous relations for general $\FHyper{r}{s}$]{Three-term contiguous relations for general $\bm{\FHyper{r}{s}}$}
   \label{subsec.hyper.rFs.contig}

The hypergeometric series $\FHyper{r}{s}$ is defined by
\be
   \FHYPER{r}{s}{a_1,\ldots,a_r}{b_1,\ldots,b_s}{t}
   \;=\;
   \sum_{n=0}^\infty
   {a_1^{\overline{n}} \,\cdots\, a_r^{\overline{n}}
    \over
    b_1^{\overline{n}} \,\cdots\, b_s^{\overline{n}}
   }
   \: {t^n \over n!}
   \;,
 \label{def.pFq}
\ee
where we have used the notation $a^{\overline{n}} = a(a+1) \cdots (a+n-1)$.
We consider \reff{def.pFq} as belonging to
the formal-power-series ring $R[[t]]$,
where $R$ is the ring $\Q(\bfb)[\bfa]$
of polynomials in the indeterminates $\bfa = (a_1,\ldots,a_r)$
whose coefficients are rational functions
in the indeterminates $\bfb = (b_1,\ldots,b_s)$.

\begin{lemma}[Three-term contiguous relations for $\FHyper{r}{s}$]
   \label{lemma.second}
For any indices $i,j$ we have:
\begin{eqnarray}
   & & \hspace*{-1cm}
   \FHYPER{r}{s\,}{a_1,\ldots,a_{i-1},a_i+1, a_{i+1},\ldots,a_r}
                {b_1,\ldots,b_s}{t}
   \:-\:
   \FHYPER{r}{s\,}{a_1,\ldots,a_r}{b_1,\ldots,b_s}{t}
        \hspace*{1cm}
        \nonumber \\[2mm]
   & & \hspace*{2cm}
   \;=\;
   {a_1 \,\cdots \not\mathrel{\,a}_i \cdots\, a_r
    \over
    b_1 \,\cdots\, b_s
   }
   \; t \;\:
   \FHYPER{r}{s\,}{a_1 +1,\ldots,a_r +1}{b_1 +1,\ldots,b_s +1}{t}
       \label{eq.second.1}  \\[4mm]
   & & \hspace*{-1cm}
   \FHYPER{r}{s\,}{a_1,\ldots,a_r}
                {b_1,\ldots,b_{i-1},b_i+1, b_{i+1},\ldots,b_s}{t}
   \:-\:
   \FHYPER{r}{s\,}{a_1,\ldots,a_r}{b_1,\ldots,b_s}{t}
        \hspace*{1cm}
        \nonumber \\[2mm]
   & & 
   \;=\;
   - \:
   {a_1 \,\cdots \, a_r
    \over
    (b_i +1) \: b_1 \,\cdots\, b_s
   }
   \; t \;\:
   \FHYPER{r}{s\,}{a_1 +1,\ldots,a_r +1}
             {b_1+1,\ldots,b_{i-1}+1,b_i+2, b_{i+1}+1,\ldots,b_s+1}{t}
       \label{eq.second.2}  \nonumber \\
%
\end{eqnarray}
\begin{eqnarray}
   & & \hspace*{-1cm}
   \FHYPER{r}{s\,}{a_1,\ldots,a_{i-1},a_i+1, a_{i+1},\ldots,a_r}
                {b_1,\ldots,b_{j-1},b_j+1, b_{j+1},\ldots,b_s}{t}
   \:-\:
   \FHYPER{r}{s\,}{a_1,\ldots,a_r}{b_1,\ldots,b_s}{t}
        \hspace*{1cm}
        \nonumber \\[2mm]
   & & \hspace*{-8mm}
   \;=\;
   {(b_j - a_i) \: a_1 \,\cdots \not\mathrel{\,a}_i \cdots\, a_r
    \over
    (b_j +1) \: b_1 \,\cdots\, b_s
   }
   \; t \;\:
   \FHYPER{r}{s\,}{a_1 +1,\ldots,a_r +1}
             {b_1+1,\ldots,b_{j-1}+1,b_j+2, b_{j+1}+1,\ldots,b_s+1}{t}
       \label{eq.second.4} \nonumber \\
\end{eqnarray}
where $\not\mathrel{\,a}_i$ indicates that $a_i$ is omitted from the product.
%
\end{lemma}

These formulae are easily proven by extracting the coefficient of $t^n$
on both sides.
%
%
%
%
%
%
See \cite{Petreolle-Sokal_contiguous} for further discussion
of these contiguous relations (and three other similar ones) and their history.

\subsection[Branched continued fraction for first ratio of contiguous $\FHyper{m+1}{m}$]{Branched continued fraction for first ratio of contiguous $\bm{\FHyper{m+1}{m}}$}
   \label{subsec.hyp.first.m+1Fm}

In this subsection we will obtain an $m$-branched continued fraction
for the first ratio \reff{def.Rrs} of contiguous $\FHyper{m+1}{m}$,
by using the contiguous relation \reff{eq.second.4}.
This generalizes Gauss' \cite{Gauss_1813} use of the $\tfo$ case
of \reff{eq.second.4} to derive a classical S-fraction
for the first ratio of contiguous $\tfo$ \cite[Chapter~XVIII]{Wall_48}.
Since the proof for $\FHyper{m+1}{m}$ parallels that of Theorem~\ref{thm.rF0}
but the weights $\balpha$ turn out to be significantly more complicated,
we shall give the proof first and let that motivate the definition
of the weights.

As before, we want to satisfy the recurrence \reff{eq.recurrence.gkm}
with suitably chosen coefficients $\balpha$.
Following the pattern of the proof of Theorem~\ref{thm.rF0},
we define
\begin{subeqnarray}
   g_{-1}(t)
   & = &
    \FHYPER{m+1}{m}{a_1,\ldots,a_m,a_{m+1}-1}{b_1,\ldots,b_{m-1},b_m-1}{t}
      \slabel{eq.hyp.first.g-1} \\[2mm]
   g_0(t)
   & = &
    \FHYPER{m+1}{m}{a_1,\ldots,a_{m+1}}{b_1,\ldots,b_m}{t}
\end{subeqnarray}
and then define $g_1,g_2,\ldots$ by incrementing first $a_1$ and $b_1$ by 1,
then $a_2$ and $b_2$ by 1, etc., continuing cyclically.
Please note that this cyclicity is
of period~$m+1$ for the numerator parameters
but period~$m$ for the denominator parameters
(strange, perhaps, but that is the way it is;
 when $m=1$ it coincides with the even-odd alternation
 used by Gauss \cite{Gauss_1813}).
Thus
\be
   g_k(t)
   \;=\;
   \FHYPER{m+1}{m}{a_1 + \lceil {k \over m+1} \rceil,\,
                   a_2 + \lceil {k-1 \over m+1} \rceil,\,
                   \ldots,\,
                   a_{m+1} + \lceil {k-m \over m+1} \rceil
                  }
                  {b_1 + \lceil {k \over m} \rceil,\, 
                   b_2 + \lceil {k-1 \over m} \rceil,\,
                   \ldots,\,
                   b_m + \lceil {k-(m-1) \over m} \rceil
                  }{t}
 \label{eq.gk.m+1Fm.0}
\ee
for all $k \ge -1$.
Defining
\begin{subeqnarray}
   a_{i,k}  & = &  a_i \,+\, \Bigl\lceil {k+1-i \over m+1} \Bigr\rceil
       \slabel{def.aik}  \\[2mm]
   b_{i,k}  & = &  b_i \,+\, \Bigl\lceil {k+1-i \over m} \Bigr\rceil
       \slabel{def.bik}
       \label{def.abik}
\end{subeqnarray}
we can write this simply as
\be
   g_k(t)
   \;=\;
   \FHYPER{m+1}{m}{a_{1,k} ,\, a_{2,k} ,\, \ldots,\, a_{m+1,k}}
                  {b_{1,k} ,\, b_{2,k} ,\, \ldots,\, b_{m,k}}
                  {t}
   \;.
 \label{eq.gk.m+1Fm}
\ee
Now, at stage $k$ the ``active'' variables
in the contiguous relation \reff{eq.second.4}
will be $a'_k \eqdef a_{[(k-1) \bmod (m+1)]+1,k}$
and $b'_k \eqdef b_{[(k-1) \bmod m]+1,k}$.
We then see that the recurrence \reff{eq.recurrence.gkm} is satisfied with
\be
   \alpha_{m+k}
   \;=\;
   { \big( b'_k - a'_k \big)
     \!\!\!
     \prod\limits_{\begin{scarray}
                       i \not\equiv k \bmod m+1
                   \end{scarray}}
     \!\!\!   a_{i,k}
     \over
     \big( b'_k - 1 \bigr) \;
     \prod\limits_{i=1}^m b_{i,k}
   }
   \;.
   \qquad
 \label{eq.alphas.m+1Fm}
\ee
Hence:

\begin{theorem}[$m$-branched continued fraction for first ratio of contiguous $\FHyper{m+1}{m}$]
   \label{thm.m+1Fm}
\hfill\break
Fix an integer $m \ge 1$,
and define $R_n^{(m+1,m)}(a_1,\ldots,a_m;a_{m+1};b_1,\ldots,b_{m-1};b_m)$ by
\be
   \sum_{n=0}^\infty R_n^{(m+1,m)}(\bfa,\bfb)
                     \: t^n
   \;\:=\;\:
   {\FHYPER{m+1}{m}{a_1,\ldots,a_{m+1}}{b_1,\ldots,b_m}{t}
    \over
    \FHYPER{m+1}{m}{a_1,\ldots,a_m,a_{m+1}-1}{b_1,\ldots,b_{m-1},b_m-1}{t}
   }
   \;\,.
 \label{eq.thm.m+1Fm}
\ee
Then $R_n^{(m+1,m)}(\bfa,\bfb) = S_n^{(m)}(\balpha)$
where the $\balpha$ are given by \reff{eq.alphas.m+1Fm}.
\end{theorem}

\subsection[Branched continued fraction for first ratio of contiguous $\FHyper{r}{s}$ for general $r,s$]{Branched continued fraction for first ratio of contiguous $\bm{\FHyper{r}{s}$} for general $\bm{r,s}$}
   \label{subsec.hyp.first.rFs}

We now consider arbitrary integers $r,s \ge 1$
and derive an $m$-branched continued fraction
for the first ratio \reff{def.Rrs} of contiguous $\FHyper{r}{s}$,
where $m = \max(r-1,s)$.
Since we have already treated $\FHyper{m+1}{m}$,
we can assume that $r \ne s+1$.
There are thus two cases:
\begin{itemize}
   \item[(a)]  $r > s+1$, hence $r=m+1$: so we are treating
$\FHyper{m+1}{s}$ with $1 \le s < m$.
   \item[(b)]  $r < s+1$, hence $s=m$: so we are treating
$\FHyper{r}{m}$ with $1 \le r < m+1$.
\end{itemize}
For each case we will give two proofs:
the first is based on taking limits in the $\FHyper{m+1}{m}$ result;
the~second is a direct proof of the recurrence \reff{eq.recurrence.gkm}
by using the identities
\reff{eq.second.1} and \reff{eq.second.4} [when $r > s+1$]
or \reff{eq.second.2} and \reff{eq.second.4} [when $r < s+1$]
in a suitable order.

\medskip

{\bf Case $\bm{\FHyper{m+1}{s}}$ with $\bm{1 \le s < m}$.}
In Theorem~\ref{thm.m+1Fm} we replace~$t$ by $b_1 \cdots b_{m-s} t$
and let $b_1, \ldots, b_{m-s} \to \infty$.
On the right-hand side of \reff{eq.thm.m+1Fm}
this gives precisely the desired ratio of $\FHyper{m+1}{s}$
(after a relabeling $b_i \to b_{i-(m-s)}$).
For the coefficient $\alpha_{m+k}$, this multiplies \reff{eq.alphas.m+1Fm}
by $b_1 \cdots b_{m-s}$ and then sends $b_1, \ldots, b_{m-s} \to \infty$.
There are two cases, depending on whether
$k \equiv 1,\ldots,m-s \bmod m$ or $k \equiv m-s+1,\ldots,m \bmod m$.

(i)  If $k \equiv 1,\ldots,m-s \bmod m$,
then the factor $(b'_k - a'_k)/(b'_k - 1)$ tends to 1,
and the terms $1 \le i \le m-s$ disappear from the denominator product.
We thus have (before relabeling)
\be
   \alpha_{m+k}
   \;=\;
   { \prod\limits_{\begin{scarray}
                       i \not\equiv k \bmod m+1
                   \end{scarray}}
     \!\!\!   a_{i,k}
     \over
     \prod\limits_{i=m-s+1}^m b_{i,k}
   }
   \;.
   \qquad
 \label{eq.alphas.rFs.case1.i}
\ee

(ii) If $k \equiv m-s+1,\ldots,m \bmod m$,
then the factor $(b'_k - a'_k)/(b'_k - 1)$ remains as~is,
and the terms $1 \le i \le m-s$ again disappear from the denominator product.
We thus have (before relabeling)
\be
   \alpha_{m+k}
   \;=\;
   { \big( b'_k - a'_k \big)
     \!\!\!
     \prod\limits_{\begin{scarray}
                       i \not\equiv k \bmod m+1
                   \end{scarray}}
     \!\!\!   a_{i,k}
     \over
     \big( b'_k - 1 \bigr) \;
     \prod\limits_{i=m-s+1}^m b_{i,k}
   }
   \;.
   \qquad
 \label{eq.alphas.rFs.case1.ii}
\ee

We now relabel $b_i \to b_{i-(m-s)}$:
to handle this, we replace \reff{def.bik} by
\be
   \widehat{b}_{i,k}
   \;=\;
   b_i \,+\, \Bigl\lceil {k+1-i-(m-s) \over m} \Bigr\rceil
 \label{def.bhat.ik}
\ee
and define
$\widehat{b}'_k \eqdef \widehat{b}_{[(k-1)-(m-s) \bmod m]+1,k}$.
With $a'_k \eqdef a_{[(k-1) \bmod (m+1)]+1,k}$ as before,
we then get
\be
   \alpha_{m+k}
   \;=\;
   \begin{cases}
      \displaystyle
      { \prod\limits_{\begin{scarray}
                          i \not\equiv k \bmod m+1
                      \end{scarray}}
        \!\!\!   a_{i,k}
        \over
        \prod\limits_{i=1}^s \widehat{b}_{i,k}
      }
      & \textrm{if $k \equiv 1,\ldots,m-s \bmod m$}   \\[12mm]
      \displaystyle
      { \big( \,\widehat{b}'_k - a'_k \big)
        \!\!\!
        \prod\limits_{\begin{scarray}
                          i \not\equiv k \bmod m+1
                      \end{scarray}}
        \!\!\!   a_{i,k}
        \over
        \big( \,\widehat{b}'_k - 1 \bigr) \;
        \prod\limits_{i=1}^s \widehat{b}_{i,k}
      }
      & \textrm{if $k \equiv m-s+1,\ldots,m \bmod m$}
   \end{cases}
 \label{def.alpha.m+1Fs}
\ee
Of course, this formula holds also when $s=m$
since it is then equivalent to \reff{eq.alphas.m+1Fm}.

We can alternatively get this same result by working directly
with the recurrences.  It suffices to define
\be
   g_k(t)
   \;=\;
   \FHYPER{m+1}{s\,}{a_{1,k} ,\, a_{2,k} ,\, \ldots,\, a_{m+1,k}}
                  {\widehat{b}_{1,k} ,\, \widehat{b}_{2,k} ,\, \ldots,\, \widehat{b}_{s,k}}
                  {t}
   \;,
 \label{eq.gk.m+1Fs}
\ee
i.e.\ the same as \reff{eq.gk.m+1Fm}
but with $\widehat{b}_{i,k}$ replacing $b_{i,k}$.
Concretely, this means that we increment the denominator variables
$b_1,\ldots,b_s$ in order and then pause for $m-s$ steps
before continuing cyclically;
the underlying cyclicity (in the ``time'' variable~$k$)
is thus of period $m$, as seen in \reff{def.bhat.ik}.
Of course, we increment the numerator variables cyclically of period $m+1$
as usual.
Then the recurrence \reff{eq.recurrence.gkm}
with the weights \reff{def.alpha.m+1Fs}
can be verified by using \reff{eq.second.1}
when $k \equiv 1,\ldots,m-s \bmod m$,
and \reff{eq.second.4} when  $k \equiv m-s+1,\ldots,m \bmod m$.
As this verification is purely computational,
we leave it to the reader.\footnote{
   We remark that the principal difficulty in this second proof
   is guessing the correct $(g_k)$ and $(\alpha_i)$.
   It is therefore convenient that we were able to determine the $(\alpha_i)$
   by the limiting method.
}

\bigskip

{\bf Case $\bm{\FHyper{r}{m}}$ with $\bm{1 \le r < m+1}$.}
In Theorem~\ref{thm.m+1Fm} we replace~$t$ by $t/(a_1 \cdots a_{m+1-r})$
and let $a_1, \ldots, a_{m+1-r} \to \infty$.
On the right-hand side of \reff{eq.thm.m+1Fm}
this gives precisely the desired ratio of $\FHyper{r}{m}$
(after a relabeling $a_i \to a_{i-(m+1-r)}$).
For the coefficient $\alpha_{m+k}$, it divides \reff{eq.alphas.m+1Fm}
by $a_1 \cdots a_{m+1-r}$ and then sends $a_1, \ldots, a_{m+1-r} \to \infty$.
There are two cases, depending on whether
$k \equiv 1,\ldots,m+1-r \bmod m+1$ or $k \equiv m+2-r,\ldots,m+1 \bmod m+1$.

(i)  If $k \equiv 1,\ldots,m+1-r \bmod m+1$,
then the factor $b'_k - a'_k$ becomes $-1$,
and the terms $1 \le i \le m+1-r$ disappear from the numerator product.
We thus have (before relabeling)
\be
   \alpha_{m+k}
   \;=\;
   - \:
   { \prod\limits_{i=m+2-r}^{m+1} \!   a_{i,k}
     \over
     \big( b'_k - 1 \big) \; \prod\limits_{i=1}^m b_{i,k}
   }
   \;.
   \qquad
 \label{eq.alphas.rFs.case2.i}
\ee

(ii) If $k \equiv m+2-r,\ldots,m+1 \bmod m+1$,
then the factor $b'_k - a'_k$ remains as~is,
and the terms $1 \le i \le m+1-r$ again disappear from the numerator product.
We thus have (before relabeling)
\be
   \alpha_{m+k}
   \;=\;
   { \big( b'_k - a'_k \big)
     \,
     \prod\limits_{i=m+2-r}^{m+1} \!  a_{i,k}
     \over
     \big( b'_k - 1 \big) \; \prod\limits_{i=1}^m b_{i,k}
   }
   \;.
   \qquad
 \label{eq.alphas.rFs.case2.ii}
\ee

We now relabel $a_i \to a_{i-(m+1-r)}$:
to handle this, we replace \reff{def.aik} by
\be
   \widehat{a}_{i,k}
   \;=\;
   a_i \,+\, \Bigl\lceil {k+1-i-(m+1-r) \over m+1} \Bigr\rceil
 \label{def.ahat.ik}
\ee
and define
$\widehat{a}'_k \eqdef \widehat{a}_{[(k-1)-(m+1-r) \bmod m+1]+1,k}$.
We then get
\be
   \alpha_{m+k}
   \;=\;
   \begin{cases}
      \displaystyle
      - \:
      { \prod\limits_{i=1}^{r} \!   \widehat{a}_{i,k}
        \over
        \big( b'_k - 1 \big) \; \prod\limits_{i=1}^m b_{i,k}
      }
         & \textrm{if $k \equiv 1,\ldots,m+1-r \bmod m+1$}   \\[12mm]
      \displaystyle
      { \big( b'_k - \widehat{a}'_k \big)
        \!\!
        \prod\limits_{i \not\equiv k-(m+1-r) \bmod m+1} \!\!\!  \widehat{a}_{i,k}
        \over
        \big( b'_k - 1 \big) \; \prod\limits_{i=1}^m b_{i,k}
      }
      & \textrm{if $k \equiv m+2-r,\ldots,m+1 \bmod m+1$}
   \end{cases}
 \label{def.alpha.rFm}
\ee
Of course, this formula holds also when $r=m+1$
since it is then equivalent to \reff{eq.alphas.m+1Fm}.

We can alternatively get this same result by working directly
with the recurrences.  It suffices to define
\be
   g_k(t)
   \;=\;
   \FHYPER{r}{m\,}{\widehat{a}_{1,k} ,\, \widehat{a}_{2,k} ,\, \ldots,\, \widehat{a}_{r,k}}
                  {b_{1,k} ,\, b_{2,k} ,\, \ldots,\, b_{m,k}}
                  {t}
   \;,
 \label{eq.gk.rFm}
\ee
i.e.\ the same as \reff{eq.gk.m+1Fm}
but with $\widehat{a}_{i,k}$ replacing $a_{i,k}$.
Concretely, this means that we increment the numerator variables
$a_1,\ldots,a_r$ in order and then pause for $m+1-r$ steps
before continuing cyclically;
the underlying cyclicity (in the ``time'' variable~$k$)
is thus of period $m+1$, as seen in \reff{def.ahat.ik}.
Of course, we increment the denominator variables cyclically of period $m$
as usual.
Then the recurrence \reff{eq.recurrence.gkm}
with the weights \reff{def.alpha.rFm}
can be verified by using \reff{eq.second.2}
when $k \equiv 1,\ldots,m+1-r \bmod m+1$,
and \reff{eq.second.4} when $k \equiv m+2-r,\ldots,m+1 \bmod m+1$.
We again leave this verification to the reader.

In summary:

\begin{theorem}[$m$-branched continued fraction for first ratio of contiguous $\FHyper{r}{s}$]
   \label{thm.rFs}
\hfill\break
Fix integers $r,s \ge 1$ and let $m = \max(r-1,s)$.
Define $R_n^{(r,s)}(a_1,\ldots,a_{r-1};a_r;$ $b_1,\ldots,b_{s-1};b_s)$ by
\be
   \sum_{n=0}^\infty R_n^{(r,s)}(\bfa,\bfb)
                     \: t^n
   \;\:=\;\:
   {\FHYPER{r}{s\,}{a_1,\ldots,a_r}{b_1,\ldots,b_s}{t}
    \over
    \FHYPER{r}{s\,}{a_1,\ldots,a_{r-1},a_r-1}{b_1,\ldots,b_{s-1},b_s-1}{t}
   }
   \;\,.
 \label{eq.thm.rFs}
\ee
Then $R_n^{(r,s)}(\bfa,\bfb) = S_n^{(m)}(\balpha)$
where the $\balpha$ are given by \reff{def.alpha.m+1Fs} when $r \ge s+1$
and by \reff{def.alpha.rFm} when $r \le s+1$.
\end{theorem}

For $\FHyper{m+1}{s}$ with $1 \le s \le m$,
we can deduce from this continued fraction
a simple {\em sufficient}\/ condition for Hankel-total positivity:

\begin{corollary}
   \label{cor.rFs}
Fix integers $1 \le s \le m$.
Let $b_1,\ldots,b_s$ be \emph{real numbers} satisfying
$b_1,\ldots,b_{s-1} > 0$ and $b_s > 1$,
and define $B = \min\limits_{1 \le i \le s} b_i$.
Now let $R$ be a partially ordered commutative ring
containing the real numbers (with their usual ordering),
and let $a_1,\ldots,a_{m+1}$ be elements of $R$
satisfying $0 \le a_i \le B$.
Then the sequence $(R_n^{(m+1,s)}(\bfa,\bfb))_{n \ge 0}$
defined by \reff{eq.thm.rFs} is Hankel-totally positive in $R$.

In particular, if $R = \R$, then the sequence
$(R_n^{(m+1,s)}(\bfa,\bfb))_{n \ge 0}$
is a Stieltjes moment sequence.
\end{corollary}

\proof
By Theorem~\ref{thm.Stype.minors}
it suffices to verify that the weights \reff{def.alpha.m+1Fs}
are nonnegative.
We manifestly have $a_{i,k} \ge a_i \ge 0$
and $\widehat{b}_{i,k} \ge b_i > 0$,
so the cases with $k \equiv 1,\ldots,m-s \bmod m$ are handled.
Also, $\widehat{b}'_k - 1 = b_s -1 > 0$ when $k=0$,
while $\widehat{b}'_k - 1 \ge b_{[(k-1) \bmod m] + 1} > 0$
when $k \ge m-s+1$ since the ceiling in \reff{def.bik} is at least 1.
Finally,
$\widehat{b}'_k - a'_k \ge
 b_{[(k-1)-(m-s) \bmod m] + 1} - a_{[(k-1) \bmod m+1] + 1} \ge 0$
since the ceiling in \reff{def.bik}
is at least as large as the ceiling in \reff{def.aik}
[the numerators in the ceilings are identical,
 but \reff{def.aik} has a denominator $m+1$ in place of $m$].
So the cases with $k \equiv m-s+1,\ldots,m \bmod m$ are also handled.
\qed

For $s=0$ there is a corresponding result, namely Corollary~\ref{cor.rF0.TP},
which implies Hankel-total positivity when
$R = \R[\bfa]$ with the coefficientwise order.
But when $s \ge 1$ and $R$ is a polynomial ring over the reals
with the coefficientwise order,
the upper bound $a_i \le B$ in the hypothesis of Corollary~\ref{cor.rFs}
forces the $a_i$ to be constants,
so that we are effectively working in $R = \R$.
We have nevertheless chosen to state the corollary in greater generality
because there may exist some interesting applications
to other partially ordered commutative rings
(for instance, polynomial rings with a pointwise order
 \cite{Brumfiel_79,Lam_84,Prestel_01,Marshall_08}).

\subsection[Branched continued fraction for second ratio of contiguous $\FHyper{m+1}{m}$]{Branched continued fraction for second ratio of contiguous $\bm{\FHyper{m+1}{m}$}}
   \label{subsec.hyp.second.m+1Fm}

We now turn to the $m$-branched continued fraction
for the second ratio \reff{def.Urs} of contiguous $\FHyper{m+1}{m}$.
As before, we want to satisfy the recurrence \reff{eq.recurrence.gkm}
with suitably chosen coefficients $\balpha$.
We define $g_k(t)$ for $k \ge 0$ by exactly the same formula
\reff{eq.gk.m+1Fm.0}/\reff{eq.gk.m+1Fm}
that was used in Section~\ref{subsec.hyp.first.m+1Fm} for the first ratio.
However --- and here is the slight twist --- we define $g_{-1}(t)$ by
\be
   g_{-1}(t)
   \;=\;
   \FHYPER{m+1}{m}{a_1,\ldots,a_m,a_{m+1}-1}{b_1,\ldots,b_{m-1},b_m}{t}
   \;,
\ee
which is different from
\reff{eq.hyp.first.g-1}/\reff{eq.gk.m+1Fm.0}/\reff{eq.gk.m+1Fm}
because it has $b_m$ in the denominator instead of $b_m - 1$.
We are therefore incrementing the variables $a_i$ and $b_i$
cyclically as before {\em except for}\/ a special treatment
of the step from $g_{-1}$ to $g_0$.

If $k \ge 1$, then the recurrence \reff{eq.recurrence.gkm}
holds with exactly the same coefficients $(\alpha_i)_{i \ge m+1}$
that were defined in \reff{eq.alphas.m+1Fm},
because the $g_j$ for $j \ge 0$ are exactly the same
as those defined in \reff{eq.gk.m+1Fm.0}/\reff{eq.gk.m+1Fm}.
On the other hand, for $k=0$ we prove directly,
using the contiguous relation \reff{eq.second.1}, that
\be
   g_0(t) \,-\, g_1(t)
   \;=\;
   {a_1 \cdots a_m \over b_1 \cdots b_m} \, t \; g_m(t)
   \;.
\ee
We therefore have:

\begin{theorem}[$m$-branched continued fraction for second ratio of contiguous $\FHyper{m+1}{m}$]
   \label{thm.m+1Fm.second}
Fix an integer $m \ge 1$,
and define $U_n^{(m+1,m)}(a_1,\ldots,a_m;a_{m+1};b_1,\ldots,b_{m-1};b_m)$ by
\be
   \sum_{n=0}^\infty U_n^{(m+1,m)}(\bfa,\bfb)
                     \: t^n
   \;\:=\;\:
   {\FHYPER{m+1}{m}{a_1,\ldots,a_{m+1}}{b_1,\ldots,b_m}{t}
    \over
    \FHYPER{m+1}{m}{a_1,\ldots,a_m,a_{m+1}-1}{b_1,\ldots,b_m}{t}
   }
   \;\,.
 \label{eq.thm.m+1Fm.second}
\ee
Then $R_n^{(m+1,m)}(\bfa,\bfb) = S_n^{(m)}(\balpha)$
where $\balpha$ is given by
\be
   \alpha_{m+k}
   \;=\;
   \begin{cases}
       \displaystyle {a_1 \cdots a_m \over b_1 \cdots b_m}
           & \textrm{if $k = 0$} \\[8mm]
       \displaystyle
       { \big( b'_k - a'_k \big)
         \!\!\!
         \prod\limits_{\begin{scarray}
                           i \not\equiv k \bmod m+1
                       \end{scarray}}
         \!\!\!   a_{i,k}
         \over
         \big( b'_k - 1 \bigr) \;
         \prod\limits_{i=1}^m b_{i,k}
       }
           & \textrm{if $k \ge 1$}
   \end{cases}
 \label{eq.alphas.m+1Fm.second}
\ee
\end{theorem}

\noindent
Here the coefficients $\balpha$ are identical to those in
Theorem~\ref{thm.m+1Fm} except for the first coefficient $\alpha_m$.

When $m=1$, Theorem~\ref{thm.m+1Fm.second} gives a classical S-fraction
\be
   {\FHYPER{2}{1}{a,b}{c}{t}
    \over
    \FHYPER{2}{1}{a,b-1}{c}{t}
   }
   \;=\;
   \cfrac{1}{1 - \cfrac{\displaystyle{a \over c} \, t}
                       {1 - \cfrac{\displaystyle{b(c-a) \over c(c+1)} \, t}
                                  { 1 - \cfrac{\displaystyle{(a+1)(c-b+1) \over (c+1)(c+2)}\, t}
                                              {1 - \ldots}}}}
\ee
with coefficients
\begin{subeqnarray}
   \alpha_1       & = &  {a \over c}  \\[2mm]
   \alpha_{2k-1}  & = &  {(a+k-1)(c-b+k-1) \over (c+2k-3)(c+2k-2)}
              \quad\hbox{for $k \ge 2$}  \\[2mm]
   \alpha_{2k}    & = &  {(b+k-1)(c-a+k-1) \over (c+2k-2)(c+2k-1)}
\end{subeqnarray}
We would be very surprised if this S-fraction were new,
but we have been unable to find it in the literature.
We would be grateful to any readers who could supply a relevant reference.
Frank \cite[eq.~(2.3)]{Frank_56} gives a continued fraction
for ${\tfo(a,b;c;t)}/{\tfo(a,b-1;c;t)}$,
but it is a T-fraction, not an S-fraction.

\subsection[Branched continued fraction for second ratio of contiguous $\FHyper{r}{s}$ for general $r,s$]{Branched continued fraction for second ratio of contiguous $\bm{\FHyper{r}{s}}$ for general $\bm{r,s}$}
   \label{subsec.hyp.second.rFs}

We now consider arbitrary integers $r \ge 1$ and $s \ge 0$
and derive an $m$-branched continued fraction
for the second ratio \reff{def.Urs} of contiguous $\FHyper{r}{s}$,
where $m = {\max(r-1,s)}$.
The case $(r,s) = (1,0)$ is trivial, so we exclude it henceforth;
then $m \ge 1$.
Since we have already treated $\FHyper{m+1}{m}$,
we can assume that $r \ne s+1$.
There are thus two cases:
\begin{itemize}
   \item[(a)]  $r > s+1$, hence $r=m+1$: so we are treating
$\FHyper{m+1}{s}$ with $0 \le s < m$.
   \item[(b)]  $r < s+1$, hence $s=m$: so we are treating
$\FHyper{r}{m}$ with $1 \le r < m+1$.
\end{itemize}
For each case there are again two proofs:
a proof based on taking limits in the $\FHyper{m+1}{m}$ result,
and a direct proof of the recurrence \reff{eq.recurrence.gkm}
by using the identities \reff{eq.second.1} and \reff{eq.second.4}
in a suitable order.
We will be brief since the proofs are similar to those
in Section~\ref{subsec.hyp.first.rFs}.

\medskip

{\bf Case $\bm{\FHyper{m+1}{s}}$ with $\bm{0 \le s < m}$.}
In Theorem~\ref{thm.m+1Fm.second} we replace~$t$ by $b_1 \cdots b_{m-s} t$
and let $b_1, \ldots, b_{m-s} \to \infty$.
On the right-hand side of \reff{eq.thm.m+1Fm.second}
this gives precisely the desired ratio of $\FHyper{m+1}{s}$
(after a relabeling $b_i \to b_{i-(m-s)}$).
For the coefficient $\alpha_{m+k}$,
this multiplies \reff{eq.alphas.m+1Fm.second}
by $b_1 \cdots b_{m-s}$ and then sends $b_1, \ldots, b_{m-s} \to \infty$.
For $k=0$, the resulting coefficient $\alpha_m$
(before relabeling)
is simply $(a_1 \cdots a_m)/(b_{m-s+1} \cdots b_m)$.
For $k > 0$ the coefficients \reff{eq.alphas.m+1Fm.second}
of the second continued fraction
are identical to the coefficients \reff{eq.alphas.m+1Fm}
of the first continued fraction,
so our analysis in Section~\ref{subsec.hyp.first.rFs}
holds without change.
After relabeling $b_i \to b_{i-(m-s)}$, the result is
\be
   \alpha_{m+k}
   \;=\;
   \begin{cases}
      \displaystyle {a_1 \cdots a_m \over b_1 \cdots b_s}
      & \textrm{if $k = 0$}   \\[6mm]
      \displaystyle
      { \prod\limits_{\begin{scarray}
                          i \not\equiv k \bmod m+1
                      \end{scarray}}
        \!\!\!   a_{i,k}
        \over
        \prod\limits_{i=1}^s \widehat{b}_{i,k}
      }
      & \textrm{if $k \equiv 1,\ldots,m-s \bmod m$}   \\[12mm]
      \displaystyle
      { \big( \,\widehat{b}'_k - a'_k \big)
        \!\!\!
        \prod\limits_{\begin{scarray}
                          i \not\equiv k \bmod m+1
                      \end{scarray}}
        \!\!\!   a_{i,k}
        \over
        \big( \,\widehat{b}'_k - 1 \bigr) \;
        \prod\limits_{i=1}^s \widehat{b}_{i,k}
      }
      & \textrm{if $k \equiv m-s+1,\ldots,m \bmod m$ and $k \neq 0$}
   \end{cases}
 \label{def.alpha.m+1Fs.second}
\ee
These coefficients are identical to those in \reff{def.alpha.m+1Fs}
except that $\alpha_m$ is different.
Of course, the result \reff{def.alpha.m+1Fs.second} holds also when $s=m$
since it is then equivalent to \reff{eq.alphas.m+1Fm};
for $s=0$ it coincides with \reff{eq.thm.rF0.alphas}.

We can alternatively get this same result by working directly
with the recurrences.  It suffices to define $g_k(t)$ for $k \ge 0$
by the same formula \reff{eq.gk.m+1Fs} as was used for the
first continued fraction.  However, we now define $g_{-1}(t)$ by
\be
   g_{-1}(t)
   \;=\;
   \FHYPER{m+1}{s}{a_1,\ldots,a_m,a_{m+1}-1}{b_1,\ldots,b_s}{t}
   \;.
\ee
If $k \ge 1$, then the recurrence \reff{eq.recurrence.gkm}
holds with exactly the same coefficients $(\alpha_i)_{i \ge m+1}$
that were defined in \reff{def.alpha.m+1Fs}
because the $g_j$ for $j \ge 0$ are exactly the same
as those defined in \reff{eq.gk.m+1Fs}.
On the other hand, for $k=0$ we prove directly,
using the contiguous relation \reff{eq.second.1}, that
\be
   g_0(t) \,-\, g_1(t)
   \;=\;
   {a_1 \cdots a_m \over b_1 \cdots b_s} \, t \; g_m(t)
   \;,
\ee
which gives the correct coefficient $\alpha_m$.

\bigskip

{\bf Case $\bm{\FHyper{r}{m}}$ with $\bm{1 \le r < m+1}$.}
In Theorem~\ref{thm.m+1Fm.second} we replace~$t$ by $t/(a_1 \cdots a_{m+1-r})$
and let $a_1, \ldots, a_{m+1-r} \to \infty$.
On the right-hand side of \reff{eq.thm.m+1Fm.second}
this gives precisely the desired ratio of $\FHyper{r}{m}$
(after a relabeling $a_i \to a_{i-(m+1-r)}$).
For the coefficient $\alpha_{m+k}$, it divides \reff{eq.alphas.m+1Fm.second}
by $a_1 \cdots a_{m+1-r}$ and then sends $a_1, \ldots, a_{m+1-r} \to \infty$.
For $k=0$, the resulting coefficient $\alpha_m$
(before relabeling)
is simply $(a_{m+2-r} \cdots a_{m+1})/(b_1 \cdots b_m)$.
For $k > 0$ the coefficients \reff{eq.alphas.m+1Fm.second}
of the second continued fraction
are identical to the coefficients \reff{eq.alphas.m+1Fm}
of the first continued fraction,
so our analysis in Section~\ref{subsec.hyp.first.rFs}
holds without change.
After relabeling $b_i \to b_{i-(m-s)}$, the result is
\be
   \alpha_{m+k}
   \;=\;
   \begin{cases}
      \displaystyle {a_1 \cdots a_r \over b_1 \cdots b_m}
         & \textrm{if $k=0$}  \\[6mm]
      \displaystyle
      - \:
      { \prod\limits_{i=1}^{r} \!   \widehat{a}_{i,k}
        \over
        \big( b'_k - 1 \big) \; \prod\limits_{i=1}^m b_{i,k}
      }
         & \textrm{if $k \equiv 1,\ldots,m+1-r \bmod m+1$}   \\[12mm]
      \displaystyle
      { \big( b'_k - \widehat{a}'_k \big)
        \!\!
        \prod\limits_{i \not\equiv k-(m+1-r) \bmod m+1} \!\!\!  \widehat{a}_{i,k}
        \over
        \big( b'_k - 1 \big) \; \prod\limits_{i=1}^m b_{i,k}
      }
      & \parbox{7cm}{\textrm{if $k \equiv m+2-r,\ldots,m+1 \bmod m+1$ and $k \ne 0$}}
   \end{cases}
 \label{def.alpha.rFm.second}
\ee
These coefficients are identical to those in \reff{def.alpha.rFm}
except that $\alpha_m$ is different.
Of course, the result \reff{def.alpha.rFm.second} holds also when $r=m+1$
since it is then equivalent to \reff{eq.alphas.m+1Fm}.

We can alternatively get this same result by working directly
with the recurrences.  It suffices to define $g_k(t)$ for $k \ge 0$
by the same formula \reff{eq.gk.rFm} as was used for the
first continued fraction.  However, we now define $g_{-1}(t)$ by
\be
   g_{-1}(t)
   \;=\;
   \FHYPER{r}{m\,}{a_1,\ldots,a_{r-1},a_r-1}{b_1,\ldots,b_m}{t}
   \;.
\ee
If $k \ge 1$, then the recurrence \reff{eq.recurrence.gkm}
holds with exactly the same coefficients $(\alpha_i)_{i \ge m+1}$
that were defined in \reff{def.alpha.rFm}
because the $g_j$ for $j \ge 0$ are exactly the same
as those treated in Section~\ref{subsec.hyp.first.rFs}.
On the other hand, for $k=0$ we prove directly,
using the contiguous relation \reff{eq.second.1}, that
\be
   g_0(t) \,-\, g_1(t)
   \;=\;
   {a_1 \cdots a_r \over b_1 \cdots b_m} \, t \; g_m(t)
   \;,
\ee
which gives the correct coefficient $\alpha_m$.

In summary:

\begin{theorem}[$m$-branched continued fraction for second ratio of contiguous $\FHyper{r}{s}$]
   \label{thm.rFs.second}
Define $U_n^{(r,s)}(\bfa,\bfb)$ by \reff{def.Urs}.  Then:
\begin{itemize}
   \item[(a)]  For $0 \le s \le m$,
      we have $U_n^{(m+1,s)}(\bfa,\bfb) = S_n^{(m)}(\balpha)$
      where the $\balpha$ are given by \reff{def.alpha.m+1Fs.second}.
   \item[(b)]  For $1 \le r \le m+1$,
      we have $U_n^{(r,m)}(\bfa,\bfb) = S_n^{(m)}(\balpha)$
      where the $\balpha$ are given by \reff{def.alpha.rFm.second}.
\end{itemize}
\end{theorem}

When $r=m=1$, Theorem~\ref{thm.rFs.second}(b) gives a classical S-fraction
for ${\ofo(a;b;t)}/{\ofo(a-1;b;t)}$.
Once again we would be surprised if this S-fraction were new,
but we have been unable to find it in the literature.

For $\FHyper{m+1}{s}$ with $0 \le s \le m$,
we can deduce from Theorem~\ref{thm.rFs.second}(a)
a simple {\em sufficient}\/ condition for Hankel-total positivity:

\begin{corollary}
   \label{cor.rFs.second}
Fix integers $0 \le s \le m$.
Let $b_1,\ldots,b_s$ be \emph{real numbers} satisfying
$b_1,\ldots,b_s > 0$,
and define $B = \min\limits_{1 \le i \le s} b_i$
(or $B = +\infty$ in~case $s=0$).
Now let $R$ be a partially ordered commutative ring
containing the real numbers (with their usual ordering),
and let $a_1,\ldots,a_{m+1}$ be elements of $R$
satisfying $0 \le a_i \le B$.
Then the sequence $(U_n^{(m+1,s)}(\bfa,\bfb))_{n \ge 0}$
defined by \reff{def.Urs} is Hankel-totally positive in $R$.

In particular, if $R = \R$, then the sequence
$(U_n^{(m+1,s)}(\bfa,\bfb))_{n \ge 0}$
is a Stieltjes moment sequence.
\end{corollary}

\proof
Essentially the same as the proof of Corollary~\ref{cor.rFs},
except that here the restriction $b_s > 1$ is not needed:
in Corollary~\ref{cor.rFs} it was employed to handle $\alpha_m$ (i.e.\ $k=0$),
but now $\alpha_m = (a_1 \cdots a_r)/(b_1 \cdots b_m)$.
\qed

For $s=0$ this coincides with Corollary~\ref{cor.rF0.TP}.
For $s \ge 1$, however, the upper bound $a_i \le B$
imposes the same limitations as were discussed after
Corollary~\ref{cor.rFs} when $R$ is a polynomial ring
with the coefficientwise order.

\subsection[Branched continued fraction for third ratio of contiguous $\FHyper{0}{m}$]{Branched continued fraction for third ratio of contiguous $\bm{\FHyper{0}{m}}$}

For the third ratio \reff{def.Vrs} of contiguous hypergeometric series,
the simplest case is $\FHyper{0}{m}$
because it needs only the contiguous relation \reff{eq.second.2}.
In order to obtain a continued fraction with nonnegative coefficients
$\balpha$, it is convenient to replace $t$ by $-t$ in the generating function.

\begin{theorem}[$m$-branched continued fraction for third ratio of contiguous $\FHyper{0}{m}$]
   \label{thm.0Fm}
\hfill\break
Fix an integer $m \ge 1$, and define $V_n^{(0,m)}(b_1,\ldots,b_{m-1};b_m)$ by
\be
   \sum_{n=0}^\infty (-1)^n \, V_n^{(0,m)}(b_1,\ldots,b_{m-1};b_m) \: t^n
   \;\:=\;\:
   {\FHYPERtopzero{m}{b_1,\ldots,b_m}{-t}
    \over
    \FHYPERtopzero{m}{b_1,\ldots,b_{m-1},b_m-1}{-t}
   }
   \;\,.
 \label{eq.thm.0Fm}
\ee
Then $(-1)^n \, V_n^{(0,m)}(\bfb) = S_n^{(m)}(\balpha)$
where $\balpha = (\alpha_i)_{i \ge m}$ is given by
\be
   \balpha
   \;=\;
   {1 \over (b_m-1) b_1 \cdots b_m}, \,
   {1 \over b_1 \cdots b_m (b_1+1)}, \,
   {1 \over b_2 \cdots b_m (b_1+1)(b_2+1)}, \,
   \ldots
   \;\,,
 \label{eq.thm.0Fm.alphas}
\ee
which can be seen as products of $m+1$ successive pre-alphas:
\be
   \balphapre 
   \;=\;
   {1 \over b_m -1}, {1 \over b_1}, {1 \over b_2}, \ldots, {1 \over b_m},
                     {1 \over b_1+1}, {1 \over b_2+1}, \ldots, {1 \over b_m+1},
                     {1 \over b_1+2}, \ldots
   \;\,.
 \label{eq.thm.0Fm.prealphas}
\ee
\end{theorem}

\proof
Define
\begin{subeqnarray}
   g_{-1}(t)
   & = &
   \FHYPERtopzero{m}{b_1,\ldots,b_{m-1},b_m-1}{-t}
           \\[2mm]
   g_0(t)
   & = &
   \FHYPERtopzero{m}{b_1,\ldots,b_m}{-t}
\end{subeqnarray}
and define $g_1,g_2,\ldots$ by successively incrementing
$b_1,\ldots,b_m$ by 1, continuing cyclically;  thus
\be
   g_k(t)
   \;=\;
   \FHYPERtopzero{m}{b_1 + \lceil {k \over m} \rceil,\,
                     b_2 + \lceil {k-1 \over m} \rceil,\,
                     \ldots,\,
                     b_m + \lceil {k-(m-1) \over m} \rceil
                    }{-t}
\ee
for all $k \ge -1$.
It is easy to see, using \reff{eq.second.2},
that the sequence $(g_k(t))_{k \ge -1}$
satisfies the recurrence \reff{eq.recurrence.gkm}
with the weights \reff{eq.thm.0Fm.alphas}.
\qed

When $m=1$ this reduces to the classical S-fraction
for ratios of contiguous $\zfo$ \cite[eqn.~(91.4)]{Wall_48},
which includes Lambert's \cite{Lambert_1768} continued fraction for tangent 
as the $b = 3/2$ special case.

Because of the denominators in \reff{eq.thm.0Fm.alphas},
we cannot use Theorem~\ref{thm.0Fm} to obtain {\em coefficientwise}\/
Hankel-total positivity,
but we can anyway deduce:

\begin{corollary}
   \label{cor.0Fm}
Let $b_1,\ldots,b_m$ be \emph{real numbers}
satisfying $b_1,\ldots,b_{m-1} > 0$ and ${b_m > 1}$.
Then the sequence
$\big( (-1)^n \, V_n^{(0,m)}(b_1,\ldots,b_{m-1};b_m) \big)_{n \ge 0}$
is a Stieltjes moment sequence.
\end{corollary}

\bigskip

Let us remark, finally, that a variant of Corollary~\ref{cor.0Fm}
can be proven without using continued fractions.
Namely, let us define the ratio of contiguous hypergeometric functions
in which {\em all}\/ of the parameters are incremented by 1:
\be
   {\FHYPER{r}{s\,}{a_1+1,\ldots,a_r+1}{b_1+1,\ldots,b_s+1}{t}
    \over
    \FHYPER{r}{s\,}{a_1,\ldots,a_r}{b_1,\ldots,b_s}{t}
   }
   \;\eqdef\;
   \sum_{n=0}^\infty W_n^{(r,s)}(\bfa,\bfb) \: t^n
   \quad\hbox{for $r,s \ge 0$} \,. \quad
      \label{def.Wrs}
\ee
It is easy to see that
\be
   {\FHYPER{r}{s\,}{a_1+1,\ldots,a_r+1}{b_1+1,\ldots,b_s+1}{t}
    \over
    \FHYPER{r}{s\,}{a_1,\ldots,a_r}{b_1,\ldots,b_s}{t}
   }
   \;=\;
   {b_1 \,\cdots\, b_s  \over  a_1 \,\cdots\, a_r} \;  {d \over dt} \,
       \log \: \FHYPER{r}{s\,}{a_1,\ldots,a_r}{b_1,\ldots,b_s}{t}
   \;.
\ee
Moreover, when $r=0$ it is known \cite{Hurwitz_1890,Hille_29,Ki_00}
that for all $b_1,\ldots,b_s > 0$, the function
$\FHYPERtopzero{s}{b_1,\ldots,b_s}{\,\cdot\,}$
is a real entire function of order $1/(s+1)$
that belongs to the Laguerre--P\'olya class $LP^+$
(that is, it is a limit, uniformly on compact subsets of $\C$,
of polynomials that have only negative real zeros).
And it is easy to show that if $f \in LP^+$
and we define real numbers $d_n$ by
\be
   {f'(t) \over f(t)}
   \;=\;
   \sum_{n=0}^\infty (-1)^n \, d_n \: t^n
   \;,
 \label{eq.prop.logderiv}
\ee
then $(d_n)_{n \ge 0}$ is a Stieltjes moment sequence.
This proves:

\begin{proposition}
   \label{prop.cor.0Fm.bis}
Let $b_1,\ldots,b_m$ be real numbers satisfying $b_1,\ldots,b_m > 0$.
Then the sequence
$\big( (-1)^n \, W_n^{(0,m)}(b_1,\ldots,b_m) \big)_{n \ge 0}$
is a Stieltjes moment sequence.
\end{proposition}

When $m=1$ this coincides with the statement of Corollary~\ref{cor.0Fm}
and therefore provides an alternate proof of it.
Similar reasoning was used in \cite[proof of Lemma~2]{Lassalle_12}.

\subsection[Branched continued fraction for third ratio of contiguous $\FHyper{m}{m}$]{Branched continued fraction for third ratio of contiguous $\bm{\FHyper{m}{m}}$}
   \label{subsec.hyp.third.mFm}

For the third ratio \reff{def.Vrs} of contiguous hypergeometric series
$\FHyper{r}{s}$,
the natural starting case is $\FHyper{m}{m}$
({\em not}\/ $\FHyper{m+1}{m}$);
we will obtain an $m$-branched continued fraction for it.
As before, we want to satisfy the recurrence \reff{eq.recurrence.gkm}
with suitably chosen coefficients $\balpha$.
We begin with
\begin{subeqnarray}
   g_{-1}(t)
   & = &
    \FHYPER{m}{m}{a_1,\ldots,a_m}{b_1,\ldots,b_m -1}{t}
           \\[2mm]
   g_0(t)
   & = &
    \FHYPER{m}{m}{a_1,\ldots,a_m}{b_1,\ldots,b_m}{t}
\end{subeqnarray}
We then define $g_1,g_2,\ldots$ by incrementing first $a_1$ and $b_1$ by 1,
continuing cyclically with period $m$ for the $b_i$,
and with period $m+1$ for the $a_i$, taking a pause after incrementing $a_m$.
This can be written for all $k \ge -1$ as
\be
   g_k(t)
   \;=\;
   \FHYPER{m}{m}{a_{1,k} ,\, a_{2,k} ,\, \ldots,\, a_{m,k}}
                  {b_{1,k} ,\, b_{2,k} ,\, \ldots,\, b_{m,k}}
                  {t}
   \;.
 \label{eq.gk.mFm.third}
\ee
where
\begin{subeqnarray}
   a_{i,k}  & = &  a_i \,+\, \Bigl\lceil {k+1-i \over m+1} \Bigr\rceil
       \slabel{def.aik.bis} \\[2mm]
   b_{i,k}  & = &  b_i \,+\, \Bigl\lceil {k+1-i \over m} \Bigr\rceil
       \slabel{def.bik.bis}
\end{subeqnarray}
Please note that this definition is identical to \reff{def.abik},
despite the different interpretation:
here there is no variable $a_{m+1}$, but the denominator $m+1$
in \reff{def.aik.bis} implements the pause.

We will now use two different contiguous relations depending on
the value of the ``time'' $k$:
\begin{itemize}
   \item[(i)]  If $k \equiv 0 \bmod m+1$,
     then at stage $k$ we will use the contiguous relation \reff{eq.second.2}
     with the ``active'' variable $b'_k \eqdef b_{[(k-1) \bmod m]+1,k}$.
     (This corresponds to the ``pause'' in the $a_i$.)
   \item[(ii)]  If $k \not\equiv 0 \bmod m+1$,
     then at stage $k$ we will use the contiguous relation \reff{eq.second.4}
     with the ``active'' variables
     $a'_k \eqdef a_{[(k-1) \bmod (m+1)]+1,k}$
     and $b'_k \eqdef b_{[(k-1) \bmod m]+1,k}$.
\end{itemize}

We then see that the recurrence \reff{eq.recurrence.gkm} is satisfied with
\be
   \alpha_{m+k}
   \;=\;
   \begin{cases}
      \displaystyle
      - \:
      { \prod\limits_{i=1}^m   a_{i,k}
        \over
        \big( b'_k -1 \big) \: \prod\limits_{i=1}^m  b_{i,k}
      }
      & \textrm{if $k \equiv 0 \bmod m+1$}   \\[12mm]
      \displaystyle
      { \big( b'_k - a'_k \big)
        \!\!\!
        \prod\limits_{\begin{scarray}
                          i \not\equiv k \bmod m+1
                      \end{scarray}}
        \!\!\!   a_{i,k}
        \over
        \big( b'_k - 1 \bigr) \;
        \prod\limits_{i=1}^m  b_{i,k}
      }
      & \textrm{if $k \not\equiv 0 \bmod m+1$}
   \end{cases}
 \label{def.alpha.mFm.third}
\ee
We therefore have:

\begin{theorem}[$m$-branched continued fraction for third ratio of contiguous $\FHyper{m}{m}$]
   \label{thm.mFm.third}
Fix an integer $m \ge 1$,
and define $V_n^{(m,m)}(a_1,\ldots,a_m;b_1,\ldots,b_{m-1};b_m)$ by
\be
   \sum_{n=0}^\infty V_n^{(m,m)}(\bfa,\bfb)
                     \: t^n
   \;\:=\;\:
   {\FHYPER{m}{m\,}{a_1,\ldots,a_m}{b_1,\ldots,b_m}{t}
    \over
    \FHYPER{m}{m\,}{a_1,\ldots,a_m}{b_1,\ldots,b_{m-1},b_m-1}{t}
   }
   \;\,.
 \label{eq.thm.mFm.third}
\ee
Then $V_n^{(m,m)}(\bfa,\bfb) = S_n^{(m)}(\balpha)$
where $\balpha$ is given by \reff{def.alpha.mFm.third}.
\end{theorem}

{\bf Remark.}
For the classical case of $\ofo$, this third continued fraction
--- for $\ofo(a;b;t)/\ofo(a;b-1;t)$ with weights \reff{def.alpha.rFm.third}
--- is known, but apparently not very well known.
It can be found in the books of
Perron \cite[3rd ed., vol.~2, p.~124]{Perron}
and Khovanskii \cite[p.~140]{Khovanskii_63}
and in the Wikpedia article on Gauss' continued fraction \cite{wikipedia},
but not (as~far as we can tell) in any of the other standard books
on continued fractions \cite{Wall_48,Jones_80,Lorentzen_92,Cuyt_08}.
However, as Perron observes, it is equivalent to the first
continued fraction for $\ofo(a;b;t)/\ofo(a-1;b-1;t)$
by virtue of the identity $\ofo(a;b;t) = {e^t \, \ofo(b-a;b;-t)}$.
\myendremark

\subsection[Branched continued fraction for third ratio of contiguous $\FHyper{r}{s}$ for general $r,s$]{Branched continued fraction for third ratio of contiguous $\bm{\FHyper{r}{s}}$ for general $\bm{r,s}$}
   \label{subsec.hyp.third.rFs}

We now consider arbitrary integers $r \ge 0$ and $s \ge 1$
and derive an $m$-branched continued fraction
for the third ratio \reff{def.Vrs} of contiguous $\FHyper{r}{s}$,
where $m = {\max(r,s)}$.
Since we have already treated $\FHyper{m}{m}$,
we can assume that $r \ne s$.
There are thus two cases:
\begin{itemize}
   \item[(a)]  $r > s$, hence $r=m$: so we are treating
$\FHyper{m}{s}$ with $1 \le s < m$.
   \item[(b)]  $r < s$, hence $s=m$: so we are treating
$\FHyper{r}{m}$ with $0 \le r < m$.
\end{itemize}
For brevity we will give only the
proof based on taking limits in the $\FHyper{m}{m}$ result;
however, there is also a proof using directly the contiguous relations.

\medskip

{\bf Case $\bm{\FHyper{m}{s}}$ with $\bm{1 \le s < m}$.}
In Theorem~\ref{thm.mFm.third} we replace~$t$ by $b_1 \cdots b_{m-s} t$
and let $b_1, \ldots, b_{m-s} \to \infty$.
On the right-hand side of \reff{eq.thm.mFm.third}
this gives precisely the desired ratio of $\FHyper{m}{s}$
(after a relabeling $b_i \to b_{i-(m-s)}$).
For the coefficient $\alpha_{m+k}$,
this multiplies \reff{def.alpha.mFm.third}
by $b_1 \cdots b_{m-s}$ and then sends $b_1, \ldots, b_{m-s} \to \infty$.
Since the coefficients \reff{eq.thm.mFm.third} involve $k \bmod m+1$
but the limits involve $k \bmod m$,
we now need to distinguish (alas!)\ four cases:
\begin{itemize}
   \item[(i)]  If $k \equiv 0 \bmod m+1$ and $k \equiv 1,\ldots,m-s \bmod m$,
then (before relabeling) we have
\be
   \alpha_{m+k}
   \;=\;
      - \:
      { \prod\limits_{i=1}^m   a_{i,k}
        \over
        \big( b'_k -1 \big) \: \prod\limits_{i=1}^m  b_{i,k}
      }
      \,\times\, b_1 \cdots b_{m-s}
   \;\to\; 0
\ee
because $b'_k - 1 \to \infty$
(the other $b_{i,k}$ with $1 \le i \le m-s$ cancel against
 the multiplier $b_1 \cdots b_{m-s}$).
   \item[(ii)]  If $k \equiv 0 \bmod m+1$ and $k \not\equiv 1,\ldots,m-s \bmod m$,
then (before relabeling) we have
\be
   \alpha_{m+k}
   \;=\;
      - \:
      { \prod\limits_{i=1}^m   a_{i,k}
        \over
        \big( b'_k -1 \big) \: \prod\limits_{i=1}^m  b_{i,k}
      }
      \,\times\, b_1 \cdots b_{m-s}
   \;\to\;
      - \:
      { \prod\limits_{i=1}^m   a_{i,k}
        \over
        \big( b'_k -1 \big) \: \prod\limits_{i=m-s+1}^m  b_{i,k}
      }
    \;.
\ee
   \item[(iii)]  If $k \not\equiv 0 \bmod m+1$ and $k \equiv 1,\ldots,m-s \bmod m$,
then (before relabeling) we have
\be
   \alpha_{m+k}
   \;=\;
      { \big( b'_k - a'_k \big)
        \!\!\!
        \prod\limits_{\begin{scarray}
                          i \not\equiv k \bmod m+1
                      \end{scarray}}
        \!\!\!   a_{i,k}
        \over
        \big( b'_k - 1 \bigr) \;
        \prod\limits_{i=1}^m  b_{i,k}
      }
      \,\times\, b_1 \cdots b_{m-s}
   \;\to\;
   { \prod\limits_{\begin{scarray}
                          i \not\equiv k \bmod m+1
                      \end{scarray}}
        \!\!\!   a_{i,k}
        \over
        \prod\limits_{i=m-s+1}^m  b_{i,k}
   }
   \;.
\ee
   \item[(iv)]  If $k \not\equiv 0 \bmod m+1$ and $k \not\equiv 1,\ldots,m-s \bmod m$,
then (before relabeling) we have
\be
   \alpha_{m+k}
   \;=\;
      { \big( b'_k - a'_k \big)
        \!\!\!
        \prod\limits_{\begin{scarray}
                          i \not\equiv k \bmod m+1
                      \end{scarray}}
        \!\!\!   a_{i,k}
        \over
        \big( b'_k - 1 \bigr) \;
        \prod\limits_{i=1}^m  b_{i,k}
      }
      \,\times\, b_1 \cdots b_{m-s}
   \;\to\;
      { \big( b'_k - a'_k \big)
        \!\!\!
        \prod\limits_{\begin{scarray}
                          i \not\equiv k \bmod m+1
                      \end{scarray}}
        \!\!\!   a_{i,k}
        \over
        \big( b'_k - 1 \bigr) \;
        \prod\limits_{i=m-s+1}^m  b_{i,k}
      }
   \;.
\ee
\end{itemize}

We now relabel $b_i \to b_{i-(m-s)}$:
to handle this, we define
\be
   \widehat{b}_{i,k}
   \;=\;
   b_i \,+\, \Bigl\lceil {k+1-i-(m-s) \over m} \Bigr\rceil
\ee
as before [cf.\ \reff{def.bhat.ik}]
and $\widehat{b}'_k \eqdef \widehat{b}_{[(k-1)-(m-s) \bmod m]+1,k}$.
We then get:
\be
   \alpha_{m+k}
   \;=\;
   \begin{cases}
       0   & \parbox{5cm}{\textrm{if $k \equiv 0 \bmod m+1$ and\\ \hphantom{if} $k \equiv 1,\ldots,m-s \bmod m$}}
            \\[6mm]
      \displaystyle
      - \:
      { \prod\limits_{i=1}^m   a_{i,k}
        \over
        \big(\, \widehat{b}'_k -1 \big) \: \prod\limits_{i=1}^s  \widehat{b}_{i,k}
      }
         & \parbox{5cm}{\textrm{if $k \equiv 0 \bmod m+1$ and\\ \hphantom{if} $k \not\equiv 1,\ldots,m-s \bmod m$}}
            \\[12mm]
      \displaystyle
   { \prod\limits_{\begin{scarray}
                          i \not\equiv k \bmod m+1
                      \end{scarray}}
        \!\!\!   a_{i,k}
        \over
        \prod\limits_{i=1}^s  \widehat{b}_{i,k}
   }
          & \parbox{5cm}{\textrm{if $k \not\equiv 0 \bmod m+1$ and\\ \hphantom{if} $k \equiv 1,\ldots,m-s \bmod m$}}
            \\[12mm]
      \displaystyle
      { \big(\, \widehat{b}'_k - a'_k \big)
        \!\!\!
        \prod\limits_{\begin{scarray}
                          i \not\equiv k \bmod m+1
                      \end{scarray}}
        \!\!\!   a_{i,k}
        \over
        \big(\, \widehat{b}'_k - 1 \bigr) \;
        \prod\limits_{i=1}^s  \widehat{b}_{i,k}
      }
          & \parbox{5cm}{\textrm{if $k \not\equiv 0 \bmod m+1$ and\\ \hphantom{if} $k \not\equiv 1,\ldots,m-s \bmod m$}}
    \end{cases}
 \label{def.alpha.mFs.third}
\ee

\medskip

{\bf Case $\bm{\FHyper{r}{m}}$ with $\bm{0 \le r < m}$.}
In Theorem~\ref{thm.mFm.third} we replace~$t$ by $t/(a_1 \cdots a_{m+1-r})$
and let $a_1, \ldots, a_{m+1-r} \to \infty$.
On the right-hand side of \reff{eq.thm.mFm.third}
this gives precisely the desired ratio of $\FHyper{r}{m}$
(after a relabeling $a_i \to a_{i-(m+1-r)}$).
For the coefficient $\alpha_{m+k}$,
this divides \reff{def.alpha.mFm.third}
by $a_1 \cdots a_{m+1-r}$ and then sends 
$a_1, \ldots, a_{m+1-r} \to \infty$.
There are three cases, depending on whether
$k \equiv 0 \bmod m+1$, 
$k \equiv 1,\ldots,m+1-r \bmod m+1$,
or $k \equiv m+2-r,\ldots,m \bmod m+1$.

(i) If $k \equiv 0 \bmod m+1$,
then (before relabeling) we have
\be
   \alpha_{m+k}
   \;=\;
   - \:
      { \prod\limits_{i=1}^m   a_{i,k}
        \over
        \big( b'_k -1 \big) \: \prod\limits_{i=1}^m  b_{i,k}
      }
      \,\times\, {1 \over a_1 \cdots a_{m+1-r}}
   \;\to\;
   - \:
      { \prod\limits_{i=m+2-r}^m   a_{i,k}
        \over
        \big( b'_k -1 \big) \: \prod\limits_{i=1}^m  b_{i,k}
      }
   \;.
 \label{eq.rFm.third.alphas.i}
\ee

(ii)  If $k \equiv 1,\ldots,m+1-r \bmod m+1$,
then (before relabeling) we have
\be
   \alpha_{m+k}
   \;=\;
   { \big( b'_k - a'_k \big)
        \!\!\!
        \prod\limits_{\begin{scarray}
                          i \not\equiv k \bmod m+1
                      \end{scarray}}
        \!\!\!   a_{i,k}
        \over
        \big( b'_k - 1 \bigr) \;
        \prod\limits_{i=1}^m  b_{i,k}
      }
      \,\times\, {1 \over a_1 \cdots a_{m+1-r}}
   \;\to\;
   - \:
      { \prod\limits_{i=m+2-r}^m   a_{i,k}
        \over
        \big( b'_k -1 \big) \: \prod\limits_{i=1}^m  b_{i,k}
      }
   \;,
 \label{eq.rFm.third.alphas.ii}
\ee
which is identical to \reff{eq.rFm.third.alphas.i}.

(iii) If $k \equiv m+2-r,\ldots,m+1 \bmod m+1$,
then (before relabeling) we have
\be
   \alpha_{m+k}
   \;=\;
   { \big( b'_k - a'_k \big)
        \!\!\!
        \prod\limits_{\begin{scarray}
                          i \not\equiv k \bmod m+1
                      \end{scarray}}
        \!\!\!   a_{i,k}
        \over
        \big( b'_k - 1 \bigr) \;
        \prod\limits_{i=1}^m  b_{i,k}
      }
      \,\times\, {1 \over a_1 \cdots a_{m+1-r}}
   \;\to\;
      {\big(b'_k - a'_k \big) \: \prod\limits_{i=m+2-r}^m   a_{i,k}
        \over
        \big( b'_k -1 \big) \: \prod\limits_{i=1}^m  b_{i,k}
      }
    \;.
\ee

We now relabel $a_i \to a_{i-(m+1-r)}$:
to handle this, we replace \reff{def.aik} by
\be
   \widehat{a}_{i,k}
   \;=\;
   a_i \,+\, \Bigl\lceil {k+1-i-(m+1-r) \over m+1} \Bigr\rceil
 \label{def.ahat.ik.bis}
\ee
and define
$\widehat{a}'_k \eqdef \widehat{a}_{[(k-1)-(m+1-r) \bmod m+1]+1,k}$.
We then get
\be
   \alpha_{m+k}
   \;=\;
   \begin{cases}
      \displaystyle
      - \:
      { \prod\limits_{i=1}^{r} \!   \widehat{a}_{i,k}
        \over
        \big( b'_k - 1 \big) \; \prod\limits_{i=1}^m b_{i,k}
      }
         & \textrm{if $k \equiv 0,\ldots,m+1-r \bmod m+1$}   \\[12mm]
      \displaystyle
      { \big( b'_k - \widehat{a}'_k \big)
        \!\!
        \prod\limits_{i \not\equiv k-(m+1-r) \bmod m+1} \!\!\!  \widehat{a}_{i,k}
        \over
        \big( b'_k - 1 \big) \; \prod\limits_{i=1}^m b_{i,k}
      }
      & \textrm{if $k \equiv m+2-r,\ldots,m \bmod m+1$}
   \end{cases}
 \label{def.alpha.rFm.third}
\ee
When $r=0$ this is simply \reff{eq.thm.0Fm.alphas}
[with a sign change because in Theorem~\ref{thm.0Fm}
 we used $-t$ in place of $t$].

In summary:

\begin{theorem}[$m$-branched continued fraction for third ratio of contiguous $\FHyper{r}{s}$]
   \label{thm.rFs.third}
Define $V_n^{(r,s)}(\bfa,\bfb)$ by \reff{def.Vrs}.  Then:
\begin{itemize}
   \item[(a)]  For $1 \le s \le m$,
      we have $V_n^{(m,s)}(\bfa,\bfb) = S_n^{(m)}(\balpha)$
      where the $\balpha$ are given by \reff{def.alpha.mFs.third}.
   \item[(b)]  For $0 \le r \le m$,
      we have $V_n^{(r,m)}(\bfa,\bfb) = S_n^{(m)}(\balpha)$
      where the $\balpha$ are given by \reff{def.alpha.rFm.third}.
\end{itemize}
\end{theorem}

Because of the minus signs in the coefficients
\reff{def.alpha.mFs.third} and \reff{def.alpha.rFm.third},
we are unable to give any result for Hankel-total positivity,
except for the case $\FHyper{0}{m}$
that was already treated in Corollary~\ref{cor.0Fm}.

\section{Ratios of contiguous hypergeometric series III: $\bm{\phiHyper{r}{s}}$}
   \label{sec.hyper.rphis}

The branched continued fractions for ratios of
hypergeometric series $\FHyper{r}{s}$ derived in the preceding section
can be straightforwardly generalized to
branched continued fractions for ratios of
basic hypergeometric series $\phiHyper{r}{s}$.

The basic hypergeometric series $\phiHyper{r}{s}$ is defined by \cite{Gasper_04}
\be
   \phiHYPER{r}{s}{a_1,\ldots,a_r}{b_1,\ldots,b_s}{q}{t}
   \;=\;
   \sum_{n=0}^\infty
   {(a_1;q)_n \, (a_2;q)_n \,\cdots\, (a_r;q)_n
    \over
    (b_1;q)_n \, (b_2;q)_n \,\cdots\, (b_s;q)_n \, (q;q)_n
   }
   \:
   \Bigl(\! (-1)^n q^{n(n-1)/2} \!\Bigr) ^{\! s+1-r}
   \:
   t^n
 \label{def.rphis}
\ee
where
\be
   (a;q)_n  \;\eqdef\;  \prod_{j=0}^{n-1} (1-aq^j)
   \;.
 \label{def.aq}
\ee
We consider \reff{def.rphis} as belonging to
the formal-power-series ring $R[[t]]$,
where $R$ is the ring $\Q(\bfb,q)[\bfa]$
of polynomials in the indeterminates $\bfa = (a_1,\ldots,a_r)$
whose coefficients are rational functions
in the indeterminates $\bfb = (b_1,\ldots,b_s)$ and $q$.\footnote{
   {\bf Warning:}  The older books of Bailey \cite{Bailey_35} and
   Slater \cite{Slater_66} define $\phiHyper{r}{s}$ {\em without}\/
   the factor $[(-1)^n q^{n(n-1)/2}]^{s+1-r}$.
   This should be borne in mind when interpreting the formulae
   from these books whenever $r \neq s+1$.
}

The proofs of the branched continued fractions for basic hypergeometric series
are completely analogous to those given in the
previous two sections for ordinary hypergeometric series,
but using the appropriately modified contiguous relations
\cite{Krattenthaler_contiguous,Krattenthaler_HYP}.
For brevity we will show only the principal case,
namely, the $m$-branched continued fraction for the first ratio
of $\phiHyper{m+1}{m}$.
This generalizes Heine's \cite{Heine_1847} continued fraction
for ratios of contiguous $\phiHyper{2}{1}$
\cite[p.~395]{Cuyt_08}.
The contiguous relation that we will need to use is
\cite[eq.~(3.3)]{Krattenthaler_contiguous}
\begin{eqnarray}
   & & \hspace*{-1cm}
   \phiHYPER{r}{s\,}{a_1,\ldots,a_{i-1},qa_i, a_{i+1},\ldots,a_r}
                {b_1,\ldots,b_{j-1},qb_j, b_{j+1},\ldots,b_s}{q}{t}
   \:-\:
   \phiHYPER{r}{s\,}{a_1,\ldots,a_r}{b_1,\ldots,b_s}{q}{t}
        \hspace*{1cm}
        \nonumber \\[2mm]
   & & 
   \;=\;
   (-1)^{s+1-r} \:
   {(a_i - b_j) \: (1-a_1) \,\cdots \overline{1-a_i} \cdots\, (1-a_r)
    \over
    (1 - qb_j) \: (1-b_1) \,\cdots\, (1-b_s)
   }
   \; t \;\:
           \nonumber \\
   & &  \hspace{2cm}
   \times \;
   \phiHYPER{r}{s\,}{qa_1,\ldots,qa_r}
             {qb_1,\ldots,qb_{j-1},q^2 b_j, qb_{j+1},\ldots,qb_s}{q}{q^{s+1-r} t}
       \label{eq.phi.contig.3}
\end{eqnarray}
where $\overline{\hbox{\sl stuff}}$ indicates that $\hbox{\sl stuff}$
is omitted from the product.

\subsection[Branched continued fraction for ratio of contiguous $\phiHyper{m+1}{m}$]{Branched continued fraction for ratio of contiguous $\bm{\phiHyper{m+1}{m}}$}
   \label{subsec.hyper.m+1Fm.contig}

The proof for $\phiHyper{m+1}{m}$
is closely analogous to the one given in Section~\ref{subsec.hyp.first.m+1Fm}
for $\FHyper{m+1}{m}$.
As before, we want to satisfy the recurrence \reff{eq.recurrence.gkm}
with suitably chosen coefficients $\balpha$.
Following the pattern of the proof of Theorem~\ref{thm.m+1Fm},
we define
\begin{subeqnarray}
   g_{-1}(t)
   & = &
    \phiHYPER{m+1}{m}{a_1,\ldots,a_m,a_{m+1}/q}{b_1,\ldots,b_{m-1},b_m/q}{q}{t}
         \\[2mm]
   g_0(t)
   & = &
    \phiHYPER{m+1}{m}{a_1,\ldots,a_{m+1}}{b_1,\ldots,b_m}{q}{t}
\end{subeqnarray}
and then define $g_1,g_2,\ldots$ by multiplying first $a_1$ and $b_1$ by $q$,
then multiplying $a_2$ and $b_2$ by $q$, etc., continuing cyclically.
As before, this cyclicity is of period~$m+1$ for the numerator parameters
but period~$m$ for the denominator parameters.
Thus
\be
   g_k(t)
   \;=\;
   \phiHYPER{m+1}{m}{a_{1,k} ,\, a_{2,k} ,\, \ldots,\, a_{m+1,k}}
                  {b_{1,k} ,\, b_{2,k} ,\, \ldots,\, b_{m,k}}
                  {q}{t}
   \;.
 \label{eq.gk.m+1phim}
\ee
for all $k \ge -1$,
where
\begin{subeqnarray}
   a_{i,k}  & = &  q^{\lceil (k+1-i)/(m+1) \rceil} \: a_i
       \slabel{def.aik.phi}  \\[2mm]
   b_{i,k}  & = &  q^{\lceil (k+1-i)/m \rceil} \: b_i
       \slabel{def.bik.phi}
       \label{def.abik.phi}
\end{subeqnarray}
Now, at stage $k$ the ``active'' variables
in the contiguous relation \reff{eq.phi.contig.3}
will be $a'_k \eqdef a_{[(k-1) \bmod (m+1)]+1,k}$
and $b'_k \eqdef b_{[(k-1) \bmod m]+1,k}$.
We then see that the recurrence \reff{eq.recurrence.gkm} is satisfied with
\be
   \alpha_{m+k}
   \;=\;
   { \big( a'_k - b'_k \big)
     \!\!\!
     \prod\limits_{\begin{scarray}
                       i \not\equiv k \bmod m+1
                   \end{scarray}}
     \!\!\!   (1 - a_{i,k})
     \over
     \big( 1 - b'_k/q \bigr) \;
     \prod\limits_{i=1}^m (1 - b_{i,k})
   }
   \;.
   \qquad
 \label{eq.alphas.m+1phim}
\ee
Note that the factors $(-1)^{s+1-r}$ and $q^{s+1-r}$
in \reff{eq.phi.contig.3} disappear here because we are in the case $r=s+1$.

\begin{theorem}[$m$-branched continued fraction for ratio of contiguous $\phiHyper{m+1}{m}$]
   \label{thm.m+1phim}
\hfill\break
Fix an integer $m \ge 1$,
and define $\Phi_n^{(m+1,m)}(a_1,\ldots,a_m;a_{m+1};b_1,\ldots,b_{m-1};b_m;q)$ by
\be
   \sum_{n=0}^\infty \Phi_n^{(m+1,m)}(\bfa,\bfb,q)
                     \: t^n
   \;\:=\;\:
   {\phiHYPER{m+1}{m}{a_1,\ldots,a_{m+1}}{b_1,\ldots,b_m}{q}{t}
    \over
    \phiHYPER{m+1}{m}{a_1,\ldots,a_m,a_{m+1}/q}{b_1,\ldots,b_{m-1},b_m/q}{q}{t}
   }
   \;\,.
 \label{eq.thm.m+1phim}
\ee
Then $\Phi_n^{(m+1,m)}(\bfa,\bfb,q) = S_n^{(m)}(\balpha)$
where the $\balpha$ are given by \reff{eq.alphas.m+1phim}.
\end{theorem}

%
%
%


\section{Some final remarks}  \label{sec.final}

\subsection{A conjecture on higher-order Genocchi numbers and Gandhi polynomials}

The {\em Genocchi numbers}\/ $G_2,G_4,G_6,\ldots$ \cite[A110501]{OEIS}
are positive integers defined by the exponential generating function
\be
   t \, \tan(t/2)
   \;=\;
   \sum_{n=1}^\infty G_{2n} \, {t^{2n} \over (2n)!}
   \;.
\ee
Their ordinary generating function has a classical S-fraction expansion
\be
   \sum_{n=0}^\infty G_{2n+2} \, t^n
   \;=\;
   \cfrac{1}{1 - \cfrac{1 \cdot 1 t}{1 - \cfrac{1 \cdot 2t}{1 - \cfrac{2 \cdot 2t}{1- \cfrac{2 \cdot 3 t}{1-\cdots}}}}}
\ee
where $\balpha = (\alpha_i)_{i \ge 1}$
are products of successive pairs of the pre-alphas
$\balpha^{\rm pre} = 1,1,2,2,$ $3,3,\ldots\:$.

More generally, define the {\em Gandhi polynomials}\/ $\scrg_n(y)$
for $n \ge 1$ by the recurrence
\be
   \scrg_n(y)  \;=\;  (y+1)^2 \, \scrg_{n-1}(y+1) \:-\: y^2 \, \scrg_{n-1}(y)
   \quad\hbox{for $n \ge 2$}
 \label{eq.recurrence.gandhi}
\ee
with initial condition $\scrg_1(y) = 1$.
They satisfy $\scrg_n(0) = G_{2n}$ and $\scrg_n(1) = G_{2n+2}$,
and a slightly modified ordinary generating function
has a classical S-fraction expansion\footnote{
   See also \cite{Dumont_95b} for a two-variable generalization.
}
\be
   1 \,+\, y \, \sum_{n=1}^\infty \scrg_n(y) \, t^n
   \;=\;
   \cfrac{1}{1 - \cfrac{y \cdot 1 \, t}{1 - \cfrac{1 \cdot (y+1)\, t}{1 - \cfrac{(y+1) \cdot 2\, t}{1- \cfrac{2 \cdot (y+2) \, t}{1-\cdots}}}}}
   \label{eq.gandhi.Sfrac}
\ee
where $\balpha = (\alpha_i)_{i \ge 1}$
are products of successive pairs of the pre-alphas
$\balpha^{\rm pre} = y,1,{y+1},2,$ ${y+2},3,\ldots\:$.

Following Han \cite{Han_93} and Domaratzki \cite{Domaratzki_04},
we now generalize the recurrence \reff{eq.recurrence.gandhi} to higher order.
Define, for each integer $m \ge 1$, the
{\em $m$th-order Gandhi polynomials}\/ $\scrg^{[m]}_n(y)$ for $n \ge 1$ by
the recurrence
\be
   \scrg^{[m]}_n(y)
   \;=\;
   (y+1)^{m+1} \, \scrg^{[m]}_{n-1}(y+1) \:-\: y^{m+1} \, \scrg^{[m]}_{n-1}(y)
   \quad\hbox{for $n \ge 2$}
 \label{eq.recurrence.gandhi.m}
\ee
with initial condition $\scrg^{[m]}_1(y) = 1$.
We then define the {\em $m$th-order Genocchi numbers}\/
$G_{2n+2}^{[m]} = \scrg^{[m]}_n(1)$;
they also satisfy $G_{2n}^{[m]} = \scrg^{[m]}_n(0)$.
We conjecture an $m$-branched continued fraction that generalizes
\reff{eq.gandhi.Sfrac}:

\begin{conjecture}
   \label{conj.gandhi}
We have $y^m \, \scrg^{[m]}_n(y) = S_n^{(m)}(\balpha)$
where $\balpha = (\alpha_i)_{i \ge m}$
are products of successive $(m+1)$-tuples of the pre-alphas
\be
   \balpha^{\rm pre} 
   \;=\;
      \underbrace{y,\ldots,y}_{\text{$m$ times}}, 1,
      \underbrace{y+1,\ldots,y+1}_{\text{$m$ times}}, 2,
      \underbrace{y+2,\ldots,y+2}_{\text{$m$ times}}, 3, \ldots
   \;.
\ee
\end{conjecture}

\noindent
In particular, when $y=1$ we conjecture an $m$-branched S-fraction
with integer coefficients for the $m$th-order Genocchi numbers.

The $\balpha$ occurring in Conjecture~\ref{conj.gandhi}
are very reminiscent of those occurring in Corollary~\ref{cor.rF0}
for the products of Stirling cycle polynomials,
and more generally those occurring in Theorem~\ref{thm.rF0}
for the ratios of hypergeometric series $\FHyper{m+1}{0}$.
In both cases the pre-alphas $\balpha^{\rm pre}$ are the same
(always quasi-affine of period $m+1$),
but here we take products of successive $(m+1)$-tuples
rather than successive $m$-tuples.

\subsection{Limitations of our method}

Nearly all of the forgeoing examples of Hankel-totally positive sequences
have been constructed using the $m$-Stieltjes--Rogers polynomials,
because this is the easiest method:  it suffices to prove the
$m$-S-fraction and verify that the coefficients $\balpha$ are nonnegative.
But, just as for $m=1$  \cite{Sokal_flajolet,Sokal_totalpos},
there are some limitations to this method:
not every sequence $\ba = (a_n)_{n \ge 0}$ with $a_0 = 1$
in a commutative ring $R$ can be expressed as $a_n = S_n^{(m)}(\balpha)$
with coefficients $\alpha_i \in R$,
even if we allow $m$ to be taken arbitrarily large.
Indeed, there is a very simple necessary condition:
since every $m$-Dyck path of nonzero length must end with an $m$-fall,
each polynomial $S_n^{(m)}(\balpha)$ with $n \ge 1$
is divisible by $\alpha_m = S_1^{(m)}(\balpha)$;
therefore, for a sequence $\ba = (a_n)_{n \ge 0}$ in a commutative ring $R$
to be expressible as $a_n = S_n^{(m)}(\balpha)$ with $\alpha_i \in R$,
a necessary condition is that $a_n$ is divisible in $R$ by $a_1$
for all $n \ge 1$.  Already this shows that many interesting
combinatorial and number-theoretic sequences,
such as the Ap\'ery numbers \cite[A005259]{OEIS}
\begin{subeqnarray}
   A_n
   & = &
   \sum_{k=0}^n \binom{n}{k}^{\! 2} \binom{n+k}{k}^{\! 2}
        \\[2mm]
   (A_n)_{n \ge 0}
   & = &
   1, 5, 73, 1445, 33001, 819005, 21460825, 584307365,
     \ldots
   \;,
 \label{def.apery}
\end{subeqnarray}
cannot be written as an $m$-S-fraction with {\em integer}\/ coefficients
$\balpha$.  (Of course, they can be written as $m$-S-fractions with
{\em rational}\/ coefficients.  It is an open question whether
any of these $m$-S-fractions have tractable closed-form expressions
and/or understandable combinatorial interpretations.)
Similarly, many interesting sequences of combinatorial polynomials,
such as the Narayana polynomials of type B,
\begin{subeqnarray}
   N_n^{\rm B}(x)
   & = &
   \sum_{k=0}^n \binom{n}{k}^{\! 2} \, x^k
        \\[2mm]
   \big( N_n^{\rm B}(x) \big)_{n \ge 0}
   & = &
   1,\, 1+x,\, 1+4x+x^2 ,\, \ldots \;,
 \label{def.narayanaB}
\end{subeqnarray}
cannot be written as an $m$-S-fraction with coefficients
in the polynomial ring $\R[x]$.

This limitation does not apply to the $m$-Thron--Rogers polynomials
$T_n^{(m)}(\balpha,\bdelta)$, so it is possible that some of these
sequences $\ba$ might be expressible as $a_n = T_n^{(m)}(\balpha,\bdelta)$
with coefficients $\alpha_i, \delta_i \in R$.
But the high degree of nonuniqueness of $m$-T-fractions
makes it very difficult to find tractable expansions of this kind;
even when $m=1$, only a few nontrivial examples are known
\cite{Elvey-Price-Sokal_wardpoly}.
We consider it to be an important open problem
to find T-fractions or $m$-T-fractions
for interesting combinatorial sequences.

Yet another possible approach is to use $m$-J-fractions.
As shown in Theorem~\ref{thm.Jtype.minors},
the sequence of $m$-Jacobi--Rogers polynomials is Hankel-totally positive
whenever the associated production matrix is totally positive.
When $m=1$ (i.e.\ for the ordinary Jacobi--Rogers polynomials),
this is a very useful method \cite{Sokal_totalpos},
because there are convenient criteria for proving the total positivity
of a tridiagonal matrix:
the contiguous-principal-minors criterion
\cite[p.~98]{Pinkus_10} \cite{Sokal_totalpos}
and the comparison theorem \cite{Sokal_totalpos,Zhu_18}.
But for $m > 1$, we have available fewer methods for proving
the total positivity of an $(m,1)$-banded matrix.
One sufficient condition, of course, is to write the matrix as a product of
lower-bidiagonal and upper-bidiagonal matrices with nonnegative entries;
but by Proposition~\ref{prop.contraction}
this is tantamount to saying that the $m$-Jacobi--Rogers polynomials
arise as the contraction of $m$-Stieltjes--Rogers polynomials
(at least if we use the order $L_1 L_2 \cdots L_m U^\star$),
so we get nothing new.
But there is a slight generalization of this approach that could be useful:
instead of writing $P = L_1 L_2 \cdots L_m U^\star$,
we write $P = L_1 L_2 \cdots L_{m-1} T$
where $T$ is tridiagonal and totally positive.
If the total positivity of $T$ has been established by some method
other than factorizing it as $L U^\star$,
this method could yield Hankel-totally positive sequences of
$m$-Jacobi--Rogers polynomials that do not simply arise as the
contraction of $m$-Stieltjes--Rogers polynomials.
%

\section*{Acknowledgments}

We are extremely grateful to J\'er\'emie Bouttier for drawing our attention
to \cite{Albenque_12},
and to Christian Krattenthaler for drawing our attention
to \cite{Krattenthaler_contiguous,Krattenthaler_HYP}.
We also wish to thank the organizers of the
80th S\'eminaire Lotharingien de Combinatoire (Lyon, March~2018)
for the opportunity to present a preliminary version of this work.
Finally, we thank an anonymous referee for a very careful reading
of the manuscript and for helpful comments.

This research was supported in part by
the U.K.~Engineering and Physical Sciences Research Council grant EP/N025636/1,
by a fellowship from the China Scholarship Council,
and by the National Natural Science Foundation of China grant No.~11971206.

\end{document}